\documentclass[12pt]{article}
\usepackage[a4paper,margin=1in]{geometry}

\usepackage{amsmath,amssymb,amsthm,mathtools,bm,mathrsfs}
\usepackage{booktabs,graphicx,float,enumitem,microtype,array}
\usepackage[round,colon,authoryear]{natbib}
\usepackage[colorlinks=true,linkcolor=blue,citecolor=blue,urlcolor=blue]{hyperref}

\newcommand{\dd}{\mathop{}\!\mathrm{d}}

\newcommand{\dist}{\operatorname{dist}}

\newtheorem{theorem}{Theorem}[section]
\newtheorem{proposition}[theorem]{Proposition}

\newtheorem{lemma}[theorem]{Lemma}
\newtheorem{corollary}[theorem]{Corollary}

\theoremstyle{definition}

\newtheorem{example}[theorem]{Example}
\theoremstyle{remark}
\newtheorem{remark}[theorem]{Remark}
\numberwithin{equation}{section}

\begin{document}

\title{Geometry, Local Minimax Risk, and the Limits of Sample Average Approximation in Linear Stochastic Optimization}
\date{}
\author{Hanbo Yang and Ming Yuan\\Columbia University}
\maketitle

\begin{abstract}
We study how the local geometry of a linear stochastic optimization problem determines its statistical difficulty, measured by expected decision regret. The optimal local rate is governed by how fast the support function of the feasible set grows away from its linearization at the true mean; depending on the instance, it can be $n^{-1/2}$, $n^{-1}$, $n^{-2}$, or a non-power rate, and sample average approximation (SAA) attains these rates in a Gaussian mean model. At smooth instances, SAA also attains the classical sharp local minimax constant. We then allow the data distribution to vary in a Kullback--Leibler (KL) ball of radius proportional to $1/n$ around a reference distribution $P$. The exact local minimax limit of this benchmark, which may depend on $P$, coincides with that of a Gaussian model whose local parameter set is the covariance ellipsoid of $P$. A simple coordinatewise truncated SAA rule, which does not use $P$, attains the minimax rate, and at smooth instances its constant matches the minimax constant as the neighborhood radius grows. Ordinary SAA, in contrast, can have strictly suboptimal worst-case expected regret when the local rate is fast relative to the tails of $P$, because rare tail events then dominate the expectation. The results delineate when SAA is sharply optimal, rate-optimal, or statistically suboptimal, and quantify the local price of distributional robustness.


\noindent\textbf{MSC2020 subject classifications:} Primary 62C20, 90C15; secondary 62F12, 52A20.
\end{abstract}

\section{Introduction}
\label{sec:introduction}

Sample average approximation (SAA) replaces an unknown expectation by its empirical counterpart and then solves the resulting optimization problem. It is model-free, computationally transparent, and supported by a mature theory of consistency, asymptotic distributions, and finite-sample guarantees; see, for example, \citet{Shapiro1993} and the monograph treatment in \citet{ShapiroDentchevaRuszczynski2021}. A more demanding question is whether SAA is \emph{statistically optimal for the particular optimization instance and uncertainty model at hand}. Standard worst-case rates do not answer this question. The attainable regret can vary sharply with the local geometry of the feasible region, and the optimality of the empirical mean can additionally depend on how distributional perturbations interact with the tails of the data.

We study these issues in the canonical linear stochastic optimization problem
\begin{equation}
\max_{\bm w\in\Omega}\ \bm w^\top\bm\mu,
\label{eq:intro_problem}
\end{equation}
where \(\Omega\subset\mathbb R^d\) is compact and convex and \(\bm\mu\) is unknown. Its optimal value is the support function
\[
f_*(\bm z)=\sup_{\bm w\in\Omega}\bm w^\top\bm z,
\]
and we write \(\bm w_*(\bm z)\) for a selected maximizer. Given $n$ independent observations, ${\bm x}_1,\ldots, {\bm x}_n$, with mean \(\bm\mu\), ordinary SAA is the plug-in policy \(\bm w_*(\bar{\bm x})\) where $\bar{\bm x}$ is the sample mean. We evaluate a data-dependent policy through decision regret and compare it with \emph{all} measurable sample-based policies.

Our analysis considers two local statistical experiments. We first study Gaussian observations with covariance \(\Sigma\) and local mean perturbations \(\bm\mu+\bm t/\sqrt n\). We then allow the entire data-generating distribution to vary in a Kullback--Leibler neighborhood of a reference distribution \(P\),
\[
\mathcal Q_{n,M}^{\mathrm{KL}}(P)
=\left\{Q\ll P:\mathbb E_Q\|\bm x\|_2<\infty,\
D_{\mathrm{KL}}(Q\|P)\le \frac{M^2}{2n}\right\}.
\]
The second experiment is locally distributionally robust: it permits changes in the mean, covariance, and tails, subject only to a shrinking entropy budget. In both experiments, the central deterministic object is the local support-function remainder
\[
r_{\bm\mu}(\bm h)
=f_*(\bm\mu+\bm h)-f_*(\bm\mu)
-\bm w_*(\bm\mu)^\top\bm h.
\]
Statistical uncertainty determines the scale and directions in which \(\bm\mu\) can be learned; the local growth of \(r_{\bm\mu}\) converts that uncertainty into regret.

This viewpoint yields both an optimality result and a limitation of SAA. In the Gaussian experiment, SAA adapts to the local geometry and attains the local minimax rate: the regret is of order \(n^{-1}\) at ordinary smooth instances and \(n^{-1/2}\) at nonsmooth instances, and degenerate curvature can produce faster rates such as \(n^{-2}\) or non-power rates. On KL neighborhoods of radius \(1/n\), the exact local minimax benchmark is governed by the same support-function geometry, and a coordinatewise truncated version of SAA attains the minimax rate. Ordinary SAA, however, need not attain even the rate: rare but KL-admissible tail perturbations can dominate its expected regret when the local geometric rate is faster than the available tail control for the empirical mean. Thus there is no universal answer to whether SAA is optimal. It can be sharply optimal, merely rate-optimal, or strictly suboptimal, depending jointly on local geometry and distributional tails.

\subsection{Main contributions}
\label{subsec:intro_contributions}

More specifically, our results can be summarized as follows.
\paragraph{General geometry and exact Gaussian local risk.}
If
\[
r_{\bm\mu}(\varepsilon\bm t)
=a(\varepsilon)\psi(\bm t)+o(a(\varepsilon))
\]
uniformly on bounded sets, the Gaussian local minimax regret, normalized by
$a(n^{-1/2})$, converges to the Bayes Jensen-gap functional
$V_{M,\Sigma}(\psi)$ in \eqref{eq:defineV}. When $\psi$ is nonaffine on the local parameter ball, this constant is positive and SAA attains the rate.
This covers nonsmooth, quadratic, degenerate, and non-power geometry in a common decision-theoretic representation.
One fixed strongly convex feasible set exhibits rates $n^{-1/2}$, $n^{-1}$, and $n^{-2}$ at different centers; a separate compact convex set gives $1/(\sqrt n\log n)$.
Finite-sample upper and lower bounds in terms of a symmetric local optimization modulus dispense with a local expansion altogether.

\paragraph{Exact local minimax risk on KL neighborhoods.}
For the distributional neighborhood \(\mathcal Q_{n,M}^{\mathrm{KL}}(P)\), we establish the exact local minimax limit
\[
\frac{\mathcal R_{n,M}^{\mathrm{KL}}(P)}{a(n^{-1/2})}
\longrightarrow
V_{M,I_d}\!\left(\psi\circ\Sigma^{1/2}\right),
\]
where $V_{M,I_d}$ is given by \eqref{eq:defineV} and $\Sigma$ is the covariance matrix of $\bm x$ under $P$, under a local exponential-moment condition; the matching lower bound requires only a finite second moment and nonsingular covariance. The rule attaining the exact constant uses knowledge of $P$ and $M$; the truncated mean introduced below does not, and attains the rate. The same limit arises in a matched Gaussian experiment with local alternatives $\bm\mu+\Sigma^{1/2}\bm u/\sqrt n$, so the KL neighborhood is asymptotically summarized by the covariance ellipsoid of the reference distribution while the loss is summarized by the local geometry \(\psi\). Whenever the limiting local modulus is nonaffine, a coordinatewise truncated empirical mean, followed by optimization, is rate-minimax. At smooth instances, its large-local-ball constant reduces to
\(\operatorname{tr}(\Sigma \nabla^2f_*(\bm\mu))/2\), which is also the local minimax constant. In this sense, local distributional robustness need not change first-order efficiency, although it changes the class of distributions against which an estimator must perform uniformly.

\paragraph{When ordinary SAA is suboptimal.}
We compare ordinary SAA with the truncated plug-in policy uniformly over the KL neighborhood. A geometry-weighted tail criterion supplies sufficient conditions for rate-optimality and for equality of the leading maximal regrets of ordinary and truncated SAA. For power-type geometry \(a(\varepsilon)=\varepsilon^\alpha\) and a reference distribution satisfying an exponential-power moment condition \(\mathbb E_P\exp(c\|\bm x\|_2^\zeta)<\infty\), ordinary SAA has maximal regret bounded by
\[
O\!\left(n^{-\alpha/2}+n^{-\zeta}\right).
\]
It is rate-optimal when \(\alpha\le 2\zeta\), and it has the same leading maximal-regret constant as truncated SAA when \(\alpha<2\zeta\). Matching examples establish rate sharpness over the stated class. In particular, for Laplace tails \((\zeta=1)\) and quartic local geometry \((\alpha=4)\), the local minimax and truncated-SAA risks are of order \(n^{-2}\), whereas ordinary SAA has worst-case regret of order \(n^{-1}\). The failure is not visible from the Gaussian weak limit of the sample mean; it arises from lack of uniform integrability under rare distributional perturbations.

\paragraph{Smooth efficiency as a special case.}
At a $C^2$ instance, $a(\varepsilon)=\varepsilon^2$ and
$\psi(\bm t)=\bm t^\top H_{\bm\mu}\bm t/2$, where
$H_{\bm\mu}=\nabla^2f_*(\bm\mu)$. The expanding-ball Gaussian minimax constant is
$\operatorname{tr}(\Sigma H_{\bm\mu})/2$, attained by SAA without knowledge of $\Sigma$.
Under the boundary regularity conditions in Section~\ref{sec:geometry}, this equals
\[
\frac{1}{2\|\bm\mu\|_2}
\sum_{i=1}^{d-1}\frac{\bm u_i^\top\Sigma\bm u_i}{\kappa_i(\bm\mu)},
\]
where $\bm u_i$ and $\kappa_i(\bm\mu)$ are the principal tangent directions and curvatures.
Thus tangent noise is weighted by inverse curvature. Truncated SAA has the same smooth efficiency constant on KL neighborhoods after the local radius tends to infinity.

\subsection{Relation to existing literature}
\label{subsec:related_work}

\paragraph{Classical SAA asymptotics.}
The statistical theory of SAA establishes consistency, asymptotic distributions of optimal values and solutions, and sample-size guarantees under increasingly general stochastic-programming models; representative references include \citet{Shapiro1993,Shapiro2008,ShapiroDentchevaRuszczynski2021}. Our focus differs in two respects. First, the target is regret rather than estimation error of the optimizer or optimal value. Second, we compare SAA with the local minimax risk over all measurable policies. This asks an efficiency question rather than only deriving the sampling distribution of SAA. The support-function structure of \eqref{eq:intro_problem} makes it possible to identify sharp constants, characterize nonquadratic rates, and distinguish weak convergence from convergence of expected regret under distributional perturbations.

\paragraph{Local minimax and instance-optimal stochastic optimization.}
The closest conceptual antecedent is \citet{DuchiRuan2021}, who develop a H\'ajek--Le Cam-type local minimax theory for stochastic convex optimization and relate statistical difficulty to tilt stability, together with asymptotically optimal online algorithms. Our linear setting permits a different level of geometric resolution. We obtain an exact Bayes--minimax representation in terms of the support-function remainder, sharp constants expressed through principal curvatures, finite-sample modulus bounds, and a theory of degenerate and non-power rates. We also extend the local analysis from a Gaussian mean experiment to population-centered KL neighborhoods and identify when the empirical-mean implementation ceases to attain the benchmark. The smooth trace constant serves as the classical benchmark; the additional conclusions are the exact fixed-radius Jensen-gap limit for general local support growth, the corresponding full-KL minimax limit, and a separation of regret rates caused by rare tails. 

A recent and complementary contribution is \citet{JiangPananjady2026}, who study finite-sample instance optimality for unconstrained smooth strongly convex stochastic optimization with additive and multiplicative oracle noise. They show that SAA and averaged stochastic approximation can be substantially suboptimal at realistic sample sizes and construct an instance-optimal variance-reduced method. Their loss, oracle model, and finite-sample objective differ from ours. Together, the two perspectives emphasize that the optimality of SAA is experiment-, loss-, and instance-dependent.

\paragraph{Distributionally robust optimization and shrinking divergence neighborhoods.}
KL and more general divergence neighborhoods have been used both to quantify model sensitivity and to construct statistically calibrated robust procedures. \citet{Lam2016} derives infinitesimal expansions for worst-case performance over shrinking KL neighborhoods. \citet{Lam2019} and \citet{DuchiGlynnNamkoong2021} connect empirical divergence-based DRO with empirical likelihood, asymptotically calibrated confidence bounds, and variance regularization, while \citet{HeLam2025} develop higher-order expansions and Bartlett corrections for divergence-based DRO. From a large-deviations perspective, \citet{VanParysEsfahaniKuhn2021} show that relative-entropy DRO is optimal for controlling out-of-sample disappointment. Wasserstein and other data-driven ambiguity sets, which are centered at the empirical distribution, are studied in \citet{MohajerinEsfahaniKuhn2018,BertsimasGuptaKallus2018,GaoKleywegt2023}. Complementary work asks whether robust or optimistic reformulations can improve the out-of-sample behavior of SAA in particular model classes; see \citet{AndersonPhilpott2022} and \citet{GotohKimLim2023}.

Our KL neighborhood is centered at a fixed reference distribution rather than the empirical distribution, and our target is different. We characterize the minimax decision regret over all measurable policies, not the value of a particular robust optimization formulation or the coverage of an uncertainty set. The radius \(1/n\) generates a local asymptotic experiment in which covariance and support geometry determine the minimax benchmark. The comparison between ordinary and truncated SAA then exposes an additional phenomenon: two procedures can have the same Gaussian weak limit but different maximal expected regret because KL-admissible rare events affect uniform integrability. This is the mechanism behind the suboptimality of ordinary SAA at sufficiently fast geometric rates.

\paragraph{Predict-then-optimize and decision-focused learning.}
For linear downstream optimization, the smart predict-then-optimize framework of \citet{ElmachtoubGrigas2022} evaluates predictions through induced decision regret rather than prediction error. Generalization analyses such as \citet{ElBalghitiEtAl2023} highlight the role of margins around parameter values that lead to nonunique decisions and introduce geometric strength conditions on the feasible set. Our ambiguity set plays a related geometric role, but our question is different: we fix a local statistical experiment and characterize its minimax regret, including exact constants and the effect of distributional tails, rather than bounding generalization error over prediction classes.

\paragraph{Nonsmooth functional and norm estimation.}
A separate literature studies minimax estimation of nonsmooth functionals and norms in Gaussian models. The connection is direct: for \(p^\star=p/(p-1)\) and \(\Omega=B_{p^\star}^d=\{\bm w:\|\bm w\|_{p^\star}\le1\}\), we have \(f_*(\bm\mu)=\|\bm\mu\|_p\), so estimating an \(\ell_p\)-norm is estimating the \emph{optimal value} of the same stochastic optimization problem. Relevant work includes \citet{LepskiNemirovskiSpokoiny1999} on \(L_r\)-norms in nonparametric Gaussian models, \citet{CaiLow2011} on the \(\ell_1\) functional in the Gaussian sequence model, and the extensions to general \(L_r\)-norms \citep{HanJiaoMukherjeeWeissman2020} and to \(\ell_\gamma\)-functionals over sparse normal means \citep{CollierCommingesTsybakov2020}.

Our target is different: we evaluate a data-dependent decision through regret rather than estimating \(f_*(\bm\mu)\) or \(\bm w_*(\bm\mu)\). Regret removes the affine first-order component of the value functional, so the local remainder becomes the governing object, and value estimation and decision regret can have different local orders. At a differentiable nonzero point, \(\|\bm\mu\|_p\) is a regular scalar functional with root-\(n\) estimation error, whereas decision regret is typically second order. At degenerate points the contrast is larger: for \(p>2\) and \(\bm\mu=r\bm e_1\), the value remains first-order regular in the radial direction, while the transverse decision regret has order \(n^{-p/2}\). These faster-order instances are also where tail stabilization can become necessary under KL robustness; see Examples~\ref{ex:lq_value_decision} and~\ref{ex:div_matching_tail_geometry}.
\paragraph{Empirical optimization, model structure, and ETO/IEO.}
A closely related line of work led by Lam asks when empirical optimization can be statistically improved. \citet{Lam2021} proves, under suitable smoothness with respect to the data-generating distribution, that the first-order optimality-gap behavior of alternative data-driven formulations is second-order stochastically dominated by empirical optimization. Our smooth, nondegenerate results are consistent with this efficiency principle: ordinary SAA is sharply optimal in the quadratic regime of the Gaussian experiment. Our KL counterexamples delineate a different regime. When the leading quadratic curvature vanishes and a faster regret rate is attainable, rare tail perturbations can make the raw empirical mean suboptimal even though a truncated empirical optimizer remains rate-minimax. 

Other work identifies mechanisms that improve on SAA through additional structure. \citet{LamWang2026} show that first-order improvements can arise when side information is correctly specified and geometrically effective. The comparison of SAA with estimate-then-optimize (ETO) and integrated-estimation-optimization (IEO) is similar in spirit: \citet{ElmachtoubLamZhangZhao2025} show that under correct parametric specification ETO asymptotically dominates IEO in the sense of stochastic dominance of regret and, when the comparison applies, can also dominate SAA, whereas under misspecification the ordering can reverse. \citet{ElmachtoubLamLanZhang2025} isolate the impact of misspecification, and \citet{LanEtAl2025} study local misspecification through contiguity and derive an explicit bias--variance tradeoff among SAA, ETO, and IEO. In these works the gain comes from model structure that reduces effective estimation variance. The mechanism in our KL results is different: tail stabilization can improve maximal expected regret without altering the Gaussian local information. Consequently, no unconditional ranking among data-driven procedures should be expected; which procedure is preferable depends on whether model structure or tail behavior is the binding constraint.

\subsection{Organization}
Section~\ref{sec:problem} contrasts the global Gaussian-width benchmark with a local quadratic-growth bound. Section~\ref{sect:fast_rate} develops the sharp Gaussian local minimax theory, including the Bayes representation, smooth efficiency bound, general local scales, and a finite-sample optimization modulus. Section~\ref{sec:geometry} translates these statistical results into boundary curvature, ambiguity, and contact geometry, and illustrates the resulting range of rates. Section~\ref{sect:divergence_neighborhoods} first introduces a matched ellipsoidal Gaussian benchmark, then establishes the exact local minimax risk over shrinking KL neighborhoods and quantifies the price of robustness. Section~\ref{sec:saa_limits} compares ordinary and truncated SAA and gives sharp conditions for optimality or suboptimality. Section~\ref{sec:numerical_experiments} illustrates the Gaussian rates and the geometry--tail transition numerically. Section~\ref{sec:discussion} concludes. All proofs, intermediate technical results, and detailed calculations are collected in the Online Supplement. Table~\ref{tab:benchmark-guide} in Section~\ref{sec:saa_limits} summarizes the benchmarks and procedure-specific risks used throughout.

\section{Global Complexity and Local Geometric Fast Rates}
\label{sec:problem}

\subsection{Problem setup}

Let $\Omega\subset\mathbb R^d$ be a compact convex set containing $\bm 0$, and define its support function by
\begin{equation}
f_*(\bm z):=\max_{\bm w\in\Omega}\bm w^\top\bm z,
\qquad \bm z\in\mathbb R^d.
\label{eq:origin_problem}
\end{equation}
For each $\bm z$, let
\[
\mathcal W_*(\bm z)
:=\arg\max_{\bm w\in\Omega}\bm w^\top\bm z,
\]
The measurable maximum theorem gives a Borel measurable selection; fix one and denote it by $\bm w_*(\bm z)\in\mathcal W_*(\bm z)$ \citep[Theorem~18.19]{AliprantisBorder2006}. The ambiguity set is
\begin{equation}
A:=\{\bm z\in\mathbb R^d:\mathcal W_*(\bm z)
\text{ is not a singleton}\}.
\label{eq:ambiguity_set}
\end{equation}
For a true objective vector $\bm\mu$ and a feasible decision $\bm w\in\Omega$, define the regret
\begin{equation}
R_{\bm\mu}(\bm w):=f_*(\bm\mu)-\bm w^\top\bm\mu.
\label{eq:pointwise_regret}
\end{equation}
We suppress the subscript $\bm\mu$ when the underlying instance is clear. Unless stated otherwise, $d$, $\Omega$, and $\bm\mu$ are fixed as $n\to\infty$, and
\[
R_\Omega:=\sup_{\bm w\in\Omega}\|\bm w\|_2<\infty.
\]

We observe
\begin{equation}
\bm x_1,\ldots,\bm x_n
\overset{\mathrm{i.i.d.}}{\sim}\mathcal N(\bm\mu,\Sigma),
\qquad
\Sigma\in\mathcal C
:=\left\{\Gamma=\Gamma^\top:
\underline\lambda I_d\preceq\Gamma\preceq\overline\lambda I_d\right\},
\label{eq:cov_model}
\end{equation}
where $0<\underline\lambda\le\overline\lambda<\infty$ are fixed constants. The covariance matrix $\Sigma$ is a fixed element of $\mathcal C$ that enters only through the benchmarks and the constants below: the SAA policy defined next does not use $\Sigma$ or the constants $\underline\lambda,\overline\lambda$. Writing
\[
\bar{\bm x}:=\frac1n\sum_{i=1}^n\bm x_i,
\]
sample average approximation (SAA) is the plug-in policy
\[
\widehat{\bm w}_{\mathrm{SAA}}=\bm w_*(\bar{\bm x}).
\]
For a compact set $D\subset\mathbb R^d$, let
\begin{equation}
\mathcal G(D):=
\mathbb E\left[\sup_{\bm w\in D}\bm w^\top\bm G\right],
\qquad \bm G\sim\mathcal N(\bm 0,I_d),
\label{eq:gaussian_width}
\end{equation}
be its Gaussian width.

\subsection{A global benchmark governed by Gaussian width}

We first consider performance uniformly over all $\bm\mu\in\mathbb R^d$. For a fixed covariance matrix $\Sigma$, define the global minimax regret
\begin{equation}
\mathfrak R_n^\Sigma(\Omega)
:=\inf_{\widehat{\bm w}_n}
\sup_{\bm\mu\in\mathbb R^d}
\mathbb E_{\bm\mu,\Sigma}
R_{\bm\mu}(\widehat{\bm w}_n),
\label{eq:global_minimax_fixed_cov}
\end{equation}
where the infimum is over measurable maps
$\widehat{\bm w}_n:(\mathbb R^d)^n\to\Omega$.

\begin{proposition}[Global minimax regret]
\label{prop:global_gaussian_width}
For every positive definite $\Sigma$ and every $n\ge1$,
\begin{equation}
\left(1-\frac1{\sqrt2}\right)
\frac{\mathcal G(\Sigma^{1/2}\Omega)}{\sqrt n}
\le \mathfrak R_n^\Sigma(\Omega)
\le \sup_{\bm\mu\in\mathbb R^d}
\mathbb E_{\bm\mu,\Sigma}
R_{\bm\mu}(\bm w_*(\bar{\bm x}))
\le 2\frac{\mathcal G(\Sigma^{1/2}\Omega)}{\sqrt n}.
\label{eq:global_width_sandwich}
\end{equation}
\end{proposition}

The global benchmark depends on the geometry of the entire feasible set through $\mathcal G(\Sigma^{1/2}\Omega)$. Moreover, $\bm w_*(\bar{\bm x})$ attains the global minimax rate without using $\Sigma$. For a sequence of problems, global minimax consistency is possible precisely, up to universal constants, when $\mathcal G(\Sigma^{1/2}\Omega)=o(\sqrt n)$. Under the fixed-dimensional asymptotics used in the remainder of the paper, this condition holds automatically; the role of Proposition~\ref{prop:global_gaussian_width} is to provide a global
baseline against which the instance-specific local rates can be compared.

\subsection{A local fast rate from quadratic growth}

The global rate does not use the location of $\bm\mu$. Faster rates become possible when the optimizer is locally stable. Let $\|\cdot\|_\Omega$ and $\|\cdot\|_K$ be a pair of dual norms, and define
\begin{equation}
\dist(\bm z,A):=\inf_{\bm a\in A}\|\bm z-\bm a\|_K.
\label{eq:distance_to_ambiguity}
\end{equation}
A natural choice is available when $\Omega=-\Omega$ and
$\bm0\in\operatorname{int}\Omega$: then
\begin{equation}
\|\bm w\|_\Omega:=\inf\{a>0:\bm w\in a\Omega\},
\qquad
K:=\{\bm z:f_*(\bm z)\le1\},
\qquad
\|\bm z\|_K=f_*(\bm z).
\label{eq:polar_norm_convention}
\end{equation}

Assume that, for some $c_\Omega>0$,
\begin{equation}
\bm z^\top\bigl(\bm w_*(\bm z)-\bm w\bigr)
\ge c_\Omega \dist(\bm z,A)
\|\bm w-\bm w_*(\bm z)\|_\Omega^2,
\qquad \bm z\in\mathbb R^d,\quad \bm w\in\Omega.
\label{eq:fast_geometric_condition}
\end{equation}
This is a quadratic-growth condition for the exposed support point. It quantifies how rapidly the objective deteriorates as the decision moves away from the optimizer, with the strength of the inequality vanishing as $\bm z$ approaches the ambiguity set.

\begin{proposition}[A covariance-dependent local fast-rate bound]
\label{prop:cov_geometric_upper}
Suppose \eqref{eq:fast_geometric_condition} holds at
$\bm z=\bm\mu$ and $\dist(\bm\mu,A)>0$. Then
\begin{equation}
\mathbb E_{\bm\mu,\Sigma}
R_{\bm\mu}(\bm w_*(\bar{\bm x}))
\le
\frac{\mathbb E\|\Sigma^{1/2}\bm G\|_K^2}
{c_\Omega n \dist(\bm\mu,A)},
\qquad \bm G\sim\mathcal N(\bm0,I_d).
\label{eq:cov_geometric_moment_bound}
\end{equation}
If, in addition, the polar-norm convention \eqref{eq:polar_norm_convention} is in force, then
\begin{equation}
\mathcal G(\Sigma^{1/2}\Omega)^2
\le \mathbb E\|\Sigma^{1/2}\bm G\|_K^2
\le\left(1+\frac\pi2\right)
\mathcal G(\Sigma^{1/2}\Omega)^2,
\label{eq:cov_width_moment_comparison}
\end{equation}
and consequently
\begin{equation}
\mathbb E_{\bm\mu,\Sigma}
R_{\bm\mu}(\bm w_*(\bar{\bm x}))
\le \frac{1+\pi/2}{c_\Omega}
\frac{\mathcal G(\Sigma^{1/2}\Omega)^2}
{n \dist(\bm\mu,A)}
\le \frac{(1+\pi/2)\overline\lambda}{c_\Omega}
\frac{\mathcal G(\Omega)^2}
{n \dist(\bm\mu,A)}.
\label{eq:fast_ddopt_bound}
\end{equation}
\end{proposition}

\begin{remark}[Midpoint strong convexity]
\label{rem:strong-convexity-sufficient}
If $\Omega=-\Omega$, $\bm0\in\operatorname{int}\Omega$, and $\Omega$ is $c_\Omega$-strongly convex in the midpoint sense relative to $\|\cdot\|_\Omega$ (so that $(\bm u+\bm v)/2+\bm s\in\Omega$ whenever $\|\bm s\|_\Omega\le\frac{c_\Omega}{2}\|\bm u-\bm v\|_\Omega^2$), then $A=\{\bm0\}$ and \eqref{eq:fast_geometric_condition} holds with $\dist(\bm z,A)=\|\bm z\|_K$.
\end{remark}

Proposition~\ref{prop:global_gaussian_width} shows that worst-case global regret is of order $n^{-1/2}$ and is determined by the width of the entire feasible set. Proposition~\ref{prop:cov_geometric_upper} shows how local stability can square the estimation error and improve the pointwise rate to $n^{-1}$. Its coefficient still combines a global quantity, $\mathcal G(\Omega)$, with the local margin $\dist(\bm\mu,A)$, so it is a benchmark rather than a sharp characterization. The next section replaces this mixed bound by a local minimax theory based on the full support-function remainder.

\section{Local Minimax Risk under Gaussian Sampling}
\label{sect:fast_rate}

\subsection{Local experiments and minimax risks}
In the original statistical problem, a policy has the form
\begin{equation}
    \hat{\bm{w}}_n=\hat{\bm{w}}_n(\bm{x}_1,\ldots,\bm{x}_n),
    \qquad \hat{\bm{w}}_n:(\mathbb{R}^d)^n\to\Omega.
    \label{eq:fast_sample_only_policy}
\end{equation}
Fix $\bm{\mu}\in\mathbb{R}^d$ and $M>0$ to evaluate local performance around $\bm{\mu}$. For the lower-bound argument, define
\begin{equation}
    R_{n,M}^{\Sigma}(\bm{\mu})
    =\inf_{\hat{\bm{w}}}\sup_{\bm{t}\in B_2(M)}
       \mathbb{E}_{\bm\mu_{n,\bm{t}},\Sigma}
       \left[f_*(\bm{\mu}_{n,\bm{t}})
       -\hat{\bm{w}}^\top\bm{\mu}_{n,\bm{t}}\right].
    \label{eq:fast_local_minimax_definition}
\end{equation}
where $\bm{\mu}_{n,\bm{t}}=\bm{\mu}+\bm{t}/\sqrt{n}$ and $B_2(M)=\{\bm{t}\in\mathbb{R}^d:\|\bm{t}\|_2\le M\}$. This is a pointwise local asymptotic minimax benchmark at the fixed reference center $\bm\mu$; we refer to optimality in this sense as \emph{instance optimality}. As is standard in local asymptotic theory, a rule in $R_{n,M}^{\Sigma}(\bm\mu)$ may be indexed by the benchmark center $\bm\mu$ and the fixed covariance $\Sigma$, but not by the unknown local perturbation $\bm t$. This convention gives the benchmark procedure more information than ordinary SAA, which is the sample-only rule $\bm w_*(\bar{\bm x})$ that uses neither the center nor $\Sigma$. In expressions evaluated under $\mathbb P_{\bm\mu_{n,\bm t},\Sigma}$, the notation $R(\bm w)$ means $f_*(\bm\mu_{n,\bm t})-\bm w^\top\bm\mu_{n,\bm t}$. For every rule of the form \eqref{eq:fast_sample_only_policy},
\begin{equation}
    \sup_{\bm{t}\in B_2(M)}
       \mathbb{E}_{\bm\mu_{n,\bm{t}},\Sigma}R(\hat{\bm{w}}_n)
    \ge R_{n,M}^{\Sigma}(\bm{\mu}).
    \label{eq:fast_sample_only_comparison}
\end{equation}
Define the local remainder by
\begin{equation}
    r_{\bm{\mu}}(\bm{u})
    =f_*(\bm{\mu}+\bm{u})-f_*(\bm{\mu})
       -\bm{w}_*(\bm{\mu})^\top\bm{u}.
    \label{eq:fast_local_remainder}
\end{equation}
This is nonnegative and convex, even when $f_*$ is not differentiable.

At an ambiguous instance, the displayed remainder depends on the selected element of $\mathcal W_*(\bm\mu)$, but all statistical quantities below do not. If $\bm w_0,\bm w_1\in\mathcal W_*(\bm\mu)$ and $r^{(j)}_{\bm\mu}$ is defined using $\bm w_j$, then
\[
    r^{(1)}_{\bm\mu}(\bm h)
    =r^{(0)}_{\bm\mu}(\bm h)
      -(\bm w_1-\bm w_0)^\top\bm h.
\]
Thus the two versions differ only by a linear functional. Such a term cancels from the Bayes Jensen gaps defining $R_{n,M}^{\Sigma}$ and $V_{M,\Sigma}$ and from the symmetric modulus $D_{\bm\mu}$. The minimax risks, rates, and constants are therefore invariant to the measurable selection.

The next subsection gives an exact decision-theoretic representation of $R_{n,M}^{\Sigma}(\bm\mu)$.

\subsection{A Bayes representation of local regret}
The key simplification is that the optimization problem enters the statistical decision problem only through the convex remainder $r_{\bm\mu}$. We now present an identity that makes this reduction exact.

Let $\mathcal{P}(B_2(M))$ be the set of probability measures supported on $B_2(M)$. For each $\pi\in\mathcal{P}(B_2(M))$, consider the following hierarchical model 
\begin{equation}
	\bm{T}\sim\pi,\qquad 
	\bm{x}_1,\ldots,\bm{x}_n\mid\bm{T}
	\overset{\mathrm{i.i.d.}}{\sim}\
	\mathcal{N}(\bm{\mu}_{n,\bm{T}},\Sigma).
	\label{eq:fast_hierarchical_prior}
\end{equation}
Thus $\pi$ is a prior on the local perturbation $\bm{T}$ and induces a prior on the true mean $\bm{\mu}_{n,\bm{T}}$. Under this model, denote $\bm{Z}=\sqrt{n}(\bar{\bm{x}}-\bm{\mu}_{n,\bm{T}})$, $\bm{Y}=\sqrt{n}(\bar{\bm{x}}-\bm{\mu})=\bm{T}+\bm{Z}$, and $\bm{m}_{\pi,\Sigma}(\bm{Y})=\mathbb{E}_{\pi,\Sigma}[\bm{T}\mid\bm{Y}]$. It is clear that $\bm{Z}\sim\mathcal{N}(\bm{0},\Sigma)$ and is independent of $\bm{T}$.

\begin{lemma}[Local Bayes--minimax identity]
\label{lem:fast_bayes_minimax}
Under the hierarchical model \eqref{eq:fast_hierarchical_prior},
\begin{equation}
\begin{aligned}
    R_{n,M}^{\Sigma}(\bm{\mu})
    &=\sup_{\pi\in\mathcal{P}(B_2(M))}
       \inf_{\hat{\bm{w}}}
       \mathbb{E}_{\pi,\Sigma}\left[
          f_*(\bm{\mu}_{n,\bm{T}})
          -\hat{\bm{w}}(\bm{x}_1,\ldots,\bm{x}_n)^\top
             \bm{\mu}_{n,\bm{T}}
       \right]\\
    &=\sup_{\pi\in\mathcal{P}(B_2(M))}
       \mathbb{E}_{\pi,\Sigma}\left[
       r_{\bm{\mu}}\left(\frac{\bm{T}}{\sqrt{n}}\right)
       -r_{\bm{\mu}}\left(
          \frac{\bm{m}_{\pi,\Sigma}(\bm{Y})}{\sqrt{n}}\right)
       \right].
\end{aligned}
    \label{eq:fast_exact_jensen_gap}
\end{equation}
For each fixed $\pi$, the inner infimum is attained by the Bayes policy
\begin{equation}
    \hat{\bm{w}}_\pi^{\mathrm{B}}(\bm{x}_1,\ldots,\bm{x}_n)
    =\bm{w}_*\left(
       \bm{\mu}+\frac{\bm{m}_{\pi,\Sigma}(\bm{Y})}{\sqrt{n}}
    \right).
    \label{eq:fast_bayes_policy}
\end{equation}
\end{lemma}

The identity is a Jensen-gap representation. The prior determines how difficult it is to infer the local perturbation $\bm T$ from the Gaussian shift $\bm Y=\bm T+\bm Z$, while $r_{\bm\mu}$ converts the remaining posterior uncertainty into regret. In particular, the rate is not fixed by the $n^{-1/2}$ estimation scale alone; it depends on how rapidly the remainder grows at the reference instance.

\subsection{Smooth instances and the sharp quadratic constant}
At a smooth instance the remainder is locally quadratic. The Hessian of the support function therefore acts as the loss matrix in the Gaussian shift experiment.
\begin{theorem}[Local asymptotic minimax regret]
\label{thm:fast_local_minimax}
Suppose that $f_*$ is twice continuously differentiable on a neighborhood of the true value $\bm{\mu}\ne\bm{0}$, and write $H_{\bm{\mu}}=\nabla^2f_*(\bm{\mu})$. Keep $d$, $\Omega$, and $\bm{\mu}$ fixed as $n\to\infty$. For every fixed $M>0$ and $\Sigma\in\mathcal C$,
\begin{equation}
    \lim_{n\to\infty}n R_{n,M}^{\Sigma}(\bm{\mu})
       =\frac12\mathcal{B}_{M,\Sigma}(H_{\bm{\mu}}),
    \label{eq:fast_fixed_neighborhood_limit}
\end{equation}
where $\mathcal{B}_{M,\Sigma}(H_{\bm{\mu}})
       :=\sup_{\pi\in\mathcal{P}(B_2(M))}
       \mathbb{E}_{\pi,\Sigma}\left[
       (\bm{T}-\bm{m}_{\pi,\Sigma}(\bm{Y}))^\top
       H_{\bm{\mu}}(\bm{T}-\bm{m}_{\pi,\Sigma}(\bm{Y}))\right]$.
Moreover,
\begin{equation}
    \lim_{M\to\infty}\lim_{n\to\infty}
       n R_{n,M}^{\Sigma}(\bm{\mu})
       =\frac12\operatorname{tr}(\Sigma H_{\bm{\mu}}).
    \label{eq:fast_sharp_local_minimax}
\end{equation}
The SAA policy $\bm{w}_*(\bar{\bm{x}})$ attains this local asymptotic minimax bound in that, for every fixed $M$,
\begin{equation}
    \lim_{n\to\infty}\sup_{\Sigma\in\mathcal C}\sup_{\bm{t}\in B_2(M)}
       \left|n\mathbb{E}_{\bm\mu_{n,\bm{t}},\Sigma}
       R(\bm{w}_*(\bar{\bm{x}}))
       -\frac12\operatorname{tr}(\Sigma H_{\bm{\mu}})\right|=0.
    \label{eq:fast_uniform_plugin_limit}
\end{equation}
\end{theorem}

For a fixed $M$, the minimax constant $\mathcal B_{M,\Sigma}(H_{\bm\mu})/2$ can be smaller than the plug-in constant $\operatorname{tr}(\Sigma H_{\bm\mu})/2$. Sharp optimality of SAA is therefore understood in the standard local asymptotic minimax sense \citep{Hajek1972,LeCam1972,vanderVaart1998}: first take $n\to\infty$ on a fixed local ball and then let the local radius $M\to\infty$. At a fixed radius, procedures that exploit the prior information $\|\bm t\|_2\le M$ in the spirit of Stein-type shrinkage \citep{Stein1956,JamesStein1961} can have a strictly smaller constant, and shrinkage-type corrections to SAA are known to matter in high-dimensional, small-sample regimes \citep{SmithWinkler2006,GuptaRusmevichientong2021}; such gains vanish as $M\to\infty$ and do not affect the fixed-dimensional rate or the sharp constant established here. The theorem also shows that SAA attains the sharp constant without using $\Sigma$, uniformly over $\Sigma\in\mathcal C$.

\subsection{General local growth and arbitrary rates}
Quadratic behavior as described in the previous subsection is only one possibility. More generally, suppose that for some $a:(0,\varepsilon_0)\to(0,\infty)$ with $a(\varepsilon)\to0$ as $\varepsilon\downarrow0$, and some finite continuous function $\psi:\mathbb{R}^d\to\mathbb{R}$, for every $\bm{t}\in\mathbb{R}^d$,
\begin{equation}
	\lim_{\varepsilon\downarrow0}
	\frac{r_{\bm{\mu}}(\varepsilon\bm{t})}{a(\varepsilon)}
	=\psi(\bm{t}).
	\label{eq:fast_general_local_expansion0}
\end{equation}
$\psi$ is automatically convex if \eqref{eq:fast_general_local_expansion0} holds true. Since each $r_{\bm\mu}(\varepsilon\,\cdot)/a(\varepsilon)$ is convex and $\psi$ is finite and continuous, pointwise convergence in \eqref{eq:fast_general_local_expansion0} is automatically uniform on bounded sets.
Let
\begin{equation}
	\label{eq:defineV}
	V_{M,\Sigma}(\psi):=\sup_{\pi\in\mathcal{P}(B_2(M))}
	\mathbb{E}_{\pi,\Sigma}\left[
	\psi(\bm{T})-\psi(\bm{m}_{\pi,\Sigma}(\bm{Y}))\right].
\end{equation}
Note that $V_{M,\Sigma}(\psi)$ is finite and nonnegative. It is strictly positive if and only if $\psi$ is not affine on $B_2(M)$. 

\begin{theorem}[Local minimax regret under a general scaling]
\label{thm:fast_general_local_rate}
Suppose that \eqref{eq:fast_general_local_expansion0} holds. Then
\begin{equation}
    \lim_{n\to\infty}
       \frac{R_{n,M}^{\Sigma}(\bm{\mu})}{a(n^{-1/2})}=V_{M,\Sigma}(\psi),
    \label{eq:fast_general_local_minimax}
\end{equation}
Moreover, whenever $V_{M,\Sigma}(\psi)>0$, the SAA policy
$\bm{w}_*(\bar{\bm{x}})$ attains the local minimax rate: there are
constants $0<c\le C<\infty$ such that, for all sufficiently large $n$,
\begin{equation}
    c a(n^{-1/2})
       \le R_{n,M}^{\Sigma}(\bm{\mu})\le\sup_{\bm{t}\in B_2(M)}
          \mathbb{E}_{\bm\mu_{n,\bm{t}},\Sigma}
             R(\bm{w}_*(\bar{\bm{x}}))
       \le C a(n^{-1/2}).
    \label{eq:fast_general_plugin_rate}
\end{equation}
The constants may depend on the fixed $d,\Omega,\bm{\mu},M,\underline\lambda,\overline\lambda$ and can be chosen uniformly over $\Sigma\in\mathcal C$.
\end{theorem}

Quantitatively, the proof yields
\begin{equation}
    \left|\frac{R_{n,M}^{\Sigma}(\bm{\mu})}{a(n^{-1/2})}-V_{M,\Sigma}(\psi)\right|
    \le 2\sup_{\bm{t}\in B_2(M)}
       \left|\frac{r_{\bm{\mu}}(\bm{t}/\sqrt{n})}{a(n^{-1/2})}
          -\psi(\bm{t})\right|.
    \label{eq:fast_general_local_error}
\end{equation}
When $\psi$ is not affine on $B_2(M)$, $R_{n,M}^{\Sigma}(\bm{\mu})\asymp a(n^{-1/2})$; otherwise $R_{n,M}^{\Sigma}(\bm{\mu})=o(a(n^{-1/2}))$.

The theorem gives a direct scale transformation:
\[
\text{estimation uncertainty at scale }n^{-1/2}
\quad\longmapsto\quad
\text{decision regret at scale }a(n^{-1/2}).
\]
It also separates exact minimax characterization from algorithmic attainment: $V_{M,\Sigma}(\psi)$ is the benchmark, while ordinary SAA is guaranteed to attain its rate whenever the limiting convex modulus is nonaffine.

\subsection{A finite-sample local optimization modulus}
The asymptotic expansion is not needed to compare the statistical difficulty of two instances. Define, for any compact convex $\Omega$, any $\bm{\mu}\in\mathbb{R}^d$,
and $\varepsilon,M>0$, the symmetric local modulus
\begin{equation}
	\begin{aligned}
		D_{\bm{\mu}}(\varepsilon,M)
		&:=\sup_{\bm{t}\in B_2(M)}
		\left\{\frac{f_*(\bm{\mu}+\varepsilon\bm{t})
			+f_*(\bm{\mu}-\varepsilon\bm{t})}{2}
		-f_*(\bm{\mu})\right\}\\
		&=\frac12\sup_{\bm{t}\in B_2(M)}
		\{r_{\bm{\mu}}(\varepsilon\bm{t})
		+r_{\bm{\mu}}(-\varepsilon\bm{t})\}.
	\end{aligned}
	\label{eq:fast_symmetric_difference_definition}
\end{equation}
Then $D_{\bm{\mu}}(\varepsilon,M)\ge0$.

\begin{proposition}[A finite-sample bound without a local expansion]
\label{prop:fast_symmetric_difference}
For every integer $n\ge1$,
\begin{equation}
    2\Phi(-M/\sqrt{\lambda_{\min}(\Sigma)})D_{\bm{\mu}}(n^{-1/2},M)
       \le R_{n,M}^{\Sigma}(\bm{\mu})
       \le 2D_{\bm{\mu}}(n^{-1/2},M),
    \label{eq:fast_symmetric_difference_bounds}
\end{equation}
where $\Phi$ is the standard normal distribution function. 
\end{proposition}

Note that Proposition \ref{prop:fast_symmetric_difference} does not assume differentiability of $f_*$ or \eqref{eq:fast_general_local_expansion0}. It shows that, up to constants depending only on the fixed local radius and spectral bounds,
\[
R_{n,M}^{\Sigma}(\bm\mu)\asymp D_{\bm\mu}(n^{-1/2},M).
\]
This is the most assumption-light form of the geometry-to-regret principle: inject statistical uncertainty at scale $n^{-1/2}$ and evaluate the symmetric growth of the support function at that scale.

\section{Geometry-to-Regret Correspondence}
\label{sec:geometry}

\subsection{Nonsmooth faces and first-order rates}
When the exposed face $\mathcal W_*(\bm\mu)$ is not a singleton, the first-order directional derivative of the support function is already nonlinear. This yields the slow local regime. More specifically, for any compact convex $\Omega$, let
\begin{equation}
		\psi_1(\bm{t})
		:=\sup_{\bm{w}\in \mathcal W_*(\bm{\mu})}
		(\bm{w}-\bm{w}_*(\bm{\mu}))^\top\bm{t}.
	\label{eq:fast_exposed_face_remainder}
\end{equation}

\begin{corollary}[First-order limit for a compact convex set]
\label{cor:fast_first_order_limit}
Condition \eqref{eq:fast_general_local_expansion0} always holds with $a(\varepsilon)=\varepsilon$ and $\psi=\psi_1$. For every $M>0$, $V_{M,\Sigma}(\psi_1)>0$ if and only if $\mathcal W_*(\bm{\mu})$ contains more than one point. 
\end{corollary}

Together with Theorem \ref{thm:fast_general_local_rate}, this immediately implies that
\begin{equation}
	\lim_{n\to\infty}\sqrt{n}R_{n,M}^{\Sigma}(\bm{\mu})
	=V_{M,\Sigma}(\psi_1)>0,
	\label{eq:fast_first_order_minimax}
\end{equation}
when $\mathcal W_*(\bm{\mu})$ is not a singleton. Otherwise $\psi_1\equiv0$, and faster rates can be achieved, i.e.,
\begin{equation*}
	R_{n,M}^{\Sigma}(\bm{\mu})=o(n^{-1/2}).
\end{equation*}

Hence decision ambiguity is precisely the source of a positive $n^{-1/2}$ local constant. A unique optimizer eliminates the first-order term, but it does not determine the next rate: the remainder may be quadratic, higher order, or non-power.

\subsection{Distance to ambiguity and normalized local shape}
Positive homogeneity separates radial distance from directional shape. Suppose that $\Omega$ is not a singleton, so that $\bm0\in A$ and $\dist(\bm\mu,A)\le\|\bm\mu\|_K<\infty$. Recall that $\|\cdot\|_K$ is a fixed norm used to define $\dist(\cdot,A)$. Since $A$ is a cone that need not be closed, $\dist(\bm\mu,A)>0$ requires more than differentiability of $f_*$ at $\bm\mu$; it holds whenever $\bm\mu\neq\bm0$ and $f_*$ is $C^2$ on a neighborhood of $\bm\mu$, because then $A$ does not meet that neighborhood. For such $\bm\mu$, let
\begin{equation}
	\bm{\nu}_{\bm{\mu}}:=\frac{\bm{\mu}}{\dist(\bm{\mu},A)},\qquad
	H_{\bm{\nu}_{\bm{\mu}}}
	:=\nabla^2 f_*(\bm{\nu}_{\bm{\mu}}).
	\label{eq:fast_distance_normalization}
\end{equation}
The rescaling ensures that $\dist(\bm{\nu}_{\bm{\mu}},A)=1$ and $H_{\bm{\mu}}
={\dist(\bm{\mu},A)^{-1}}H_{\bm{\nu}_{\bm{\mu}}}$ so that
\begin{equation}
	\lim_{n\to\infty}nR_{n,M}^{\Sigma}(\bm{\mu})
	=\frac{\mathcal{B}_{M,\Sigma}(H_{\bm{\nu}_{\bm{\mu}}})}
	{2\dist(\bm{\mu},A)}.
	\label{eq:fast_distance_fixed_M}
\end{equation}
This implies that
\begin{corollary}[The sharp constant and the distance to the ambiguity set]
\label{cor:fast_distance_smooth}
For any $\bm\mu\neq \bm0$ such that $f_\ast$ is $C^2$ on a neighborhood of $\bm\mu$,
\begin{equation}
    \lim_{M\to\infty}\lim_{n\to\infty}nR_{n,M}^{\Sigma}(\bm{\mu})
       =\frac{\operatorname{tr}(\Sigma H_{\bm{\nu}_{\bm{\mu}}})}
               {2\dist(\bm{\mu},A)}=\frac{[\mathcal{G}(\Sigma^{1/2}\Omega)]^2}{\dist(\bm{\mu},A)}
          \chi_{\Omega,\Sigma}(\bm{\mu}),
    \label{eq:fast_distance_sharp_constant}
\end{equation}
where $\chi_{\Omega,\Sigma}(\bm{\mu}):=\operatorname{tr}(\Sigma H_{{\bm{\nu}_{\bm{\mu}}}})/(2{[\mathcal{G}(\Sigma^{1/2}\Omega)]^2}).$
\end{corollary}

The distance $\dist(\bm\mu,A)$ therefore controls only one component of difficulty. The normalized Hessian $H_{\bm\nu_{\bm\mu}}$ retains the directional geometry, and two instances at the same distance from ambiguity can have different rates when one has degenerate curvature.

At a smooth point with $\operatorname{tr}(\Sigma H_{\bm{\mu}})>0$, the power $n^{-1}$ is locally optimal. The fast rate in \eqref{eq:fast_ddopt_bound} is sharp up to constant factors precisely when
\begin{equation}
\operatorname{tr}(\Sigma\nabla^2f_*(\bm{\mu}))
\asymp\frac{[\mathcal{G}(\Sigma^{1/2}\Omega)]^2}{\dist(\bm{\mu},A)},
\label{eq:fast_curvature_comparability}
\end{equation}
or equivalently, $\chi_{\Omega,\Sigma}(\bm{\mu})\asymp1$, with the strong-convexity constant in \eqref{eq:fast_ddopt_bound} controlled. 

For example, for $\Omega=B_2^d$ with $d\ge2$,
\begin{equation*}
f_*(\bm{\mu})=\|\bm{\mu}\|_2,\qquad
\operatorname{tr}(\Sigma\nabla^2f_*(\bm{\mu}))
=\frac{\operatorname{tr}(\Sigma)-\bm\mu^\top\Sigma\bm\mu/\|\bm\mu\|_2^2}{\|\bm{\mu}\|_2},\qquad
\dist(\bm{\mu},A)=\|\bm{\mu}\|_2,
\end{equation*}
and $[\mathcal{G}(\Sigma^{1/2}B_2^d)]^2\asymp d$ under the fixed spectral bounds. Hence $\chi_{\Omega,\Sigma}(\bm\mu)\asymp1$ and the fast-rate bound
matches the local minimax scale.
In contrast, consider the strongly convex ellipse
\begin{equation*}
\Omega_L=\{\bm{w}\in\mathbb{R}^2:w_1^2/L^2+w_2^2\le1\},
\qquad L\ge1,\qquad \bm{\mu}=\bm{e}_1.
\end{equation*}
Its support function is $f_*(\bm{\mu})=(L^2\mu_1^2+\mu_2^2)^{1/2}$,
so $A=\{\bm{0}\}$ and
\begin{equation*}
\dist(\bm{\mu},A)=L,\qquad
[\mathcal{G}(\Sigma^{1/2}\Omega_L)]^2\asymp L^2+1,\qquad
\operatorname{tr}(\Sigma\nabla^2f_*(\bm{\mu}))=\Sigma_{22}/L.
\end{equation*}
As a linear image of the Euclidean ball, this ellipse satisfies \eqref{eq:fast_geometric_condition} with $c_{\Omega_L}=1/2$ in its own norm, uniformly in $L$. Under the fixed spectral bounds, $\Sigma_{22}\asymp1$, so that
\[
\chi_{\Omega_L,\Sigma}(\bm e_1)\asymp\frac{1}{L^2}\longrightarrow0
\qquad(L\to\infty).
\]
Thus \eqref{eq:fast_curvature_comparability} fails. The bound \eqref{eq:fast_ddopt_bound} is of order $L/n$, whereas the sharp local constant is $\operatorname{tr}(\Sigma H_{\bm\mu})/(2n)=\Sigma_{22}/(2Ln)$; the bound overshoots by a factor of order $L^2$, although the quadratic-growth constant does not deteriorate. A bound that combines the global width $\mathcal G(\Omega)$ with the local margin $\dist(\bm\mu,A)$ therefore cannot be sharp in general.

\subsection{Smooth boundaries and inverse curvature}
The identity between the tangential Hessian of a support function and the inverse Weingarten map is classical; see \citet[Section~2]{HowardHug2007}. The next result gives its local form and applies it to the statistical constant in Theorem \ref{thm:fast_local_minimax}.

Let $d\ge2$, let $\Omega$ be compact and convex with nonempty interior, and fix $\bm{\mu}\ne\bm{0}$ with a unique maximizer $\bm{w}_*(\bm{\mu})$. Denote $\bm{n}_{\bm{\mu}}:=\bm{\mu}/\|\bm{\mu}\|_2$ as the outward unit normal to the supporting hyperplane selected by $\bm\mu$, $\mathcal{T}_{\bm{\mu}}:=\bm{n}_{\bm{\mu}}^\perp$ as the linear subspace parallel to that hyperplane, and $U_{\bm{\mu}}$ as an orthonormal basis matrix of $\mathcal{T}_{\bm{\mu}}$, with $U_{\bm{\mu}}^\top U_{\bm{\mu}}=I_{d-1}$ and $\operatorname{range}(U_{\bm{\mu}})=\mathcal{T}_{\bm{\mu}}$. When $\partial\Omega$ is a $C^2$ hypersurface near $\bm w_*(\bm\mu)$, $\mathcal T_{\bm\mu}$ is its tangent space at $\bm w_*(\bm\mu)$; in this case, let $\bm N$ be the outward unit normal map on the boundary neighborhood and define
\begin{equation}
	\mathscr{S}_{\bm{\mu}}
	:=U_{\bm{\mu}}^\top
	D\bm{N}(\bm{w}_*(\bm{\mu}))U_{\bm{\mu}},\qquad
	\kappa_i(\bm{\mu}):=\lambda_i(\mathscr{S}_{\bm{\mu}}),
	\quad i=1,\ldots,d-1.
	\label{eq:fast_boundary_curvature_matrix}
\end{equation}
Here, for any differentiable vector-valued function $F:\mathbb{R}^m\to\mathbb{R}^k$, the derivative $DF(\bm{x})$ is defined as a linear map with
\begin{equation*}
	F(\bm{x}+\bm{h})=F(\bm{x})+DF(\bm{x})[\bm{h}]+o(\|\bm{h}\|_2).
\end{equation*}
For a $C^2$ boundary, $\bm N$ is $C^1$, so its differential exists; it is taken along boundary curves and acts on tangent vectors.
We use the sign convention under which a Euclidean sphere of radius $R$ has $\kappa_i=1/R$. 

\begin{corollary}[Boundary principal curvatures and the sharp constant]
\label{cor:fast_principal_curvatures}
Assume that $\partial\Omega$ is a $C^2$ hypersurface on a neighborhood of $\bm{w}_*(\bm{\mu})$ and that $\mathscr{S}_{\bm{\mu}}\succ0$, equivalently $\kappa_i(\bm{\mu})>0$ for every $i$. Then $f_*$ is $C^2$ on a neighborhood of $\bm{\mu}$, and
\begin{equation}
    H_{\bm{\mu}}
       =\frac{1}{\|\bm{\mu}\|_2}
          U_{\bm{\mu}}\mathscr{S}_{\bm{\mu}}^{-1}U_{\bm{\mu}}^\top.
    \label{eq:fast_hessian_principal_curvatures}
\end{equation}
In particular, $\operatorname{rank}(H_{\bm{\mu}})=d-1$ and
$\ker(H_{\bm{\mu}})=\operatorname{span}\{\bm{\mu}\}$.
Choose the columns $\bm u_1,\ldots,\bm u_{d-1}$ of
$U_{\bm\mu}$ along orthonormal principal directions. The sharp local
asymptotic minimax constant is
\begin{equation}
\begin{aligned}
    \lim_{M\to\infty}\lim_{n\to\infty}nR_{n,M}^{\Sigma}(\bm\mu)
       &=\frac1{2\|\bm\mu\|_2}
          \sum_{i=1}^{d-1}\frac{\bm u_i^\top\Sigma\bm u_i}{\kappa_i(\bm\mu)}\\
       &=\frac{\dist(\bm n_{\bm\mu},A)}{2\dist(\bm\mu,A)}
          \sum_{i=1}^{d-1}\frac{\bm u_i^\top\Sigma\bm u_i}{\kappa_i(\bm\mu)}.
\end{aligned}
    \label{eq:fast_principal_curvature_minimax}
\end{equation}
The second equality uses any fixed norm $\|\cdot\|_K$ defining $\dist(\cdot,A)$.
\end{corollary}

The formula has a direct interpretation. Noise normal to the boundary does not change the exposed decision to first order. Tangential noise does, and its contribution is divided by the corresponding principal curvature. Flatter boundary directions are therefore statistically harder.

\begin{remark}[Quadratic growth implies positive curvature]
\label{rem:growth_implies_curvature}
Let $d\ge2$, let $\Omega$ have nonempty interior, and fix $\bm\mu\ne\bm0$ with $\dist(\bm\mu,A)>0$ such that $\partial\Omega$ is a $C^2$ hypersurface near $\bm w_*(\bm\mu)$. If the quadratic-growth condition \eqref{eq:fast_geometric_condition} holds at $\bm z=\bm\mu$, then $\mathscr{S}_{\bm\mu}\succ0$, so Corollary~\ref{cor:fast_principal_curvatures} applies without a separate positive-curvature assumption (Lemma~\ref{lem:growth_implies_curvature} in the Online Supplement).
\end{remark}

\subsection{Boundary contact and rate exponents}
The curvature formula extends to nonsmooth or degenerate boundaries through the order of contact with the supporting hyperplane. First, we introduce the following local parametrization of the boundary $\partial\Omega$. 

\begin{proposition}[Local boundary parametrization and growth conditions]
\label{prop:boundary-local-parametrization}
Suppose that $\Omega$ has a unique supporting hyperplane at $\bm w_*(\bm\mu)$. Then there exist $\rho>0$ and a continuous convex function $\varrho_{\bm\mu}$, defined on a neighborhood of $\{\bm z:\|\bm z\|_2\le\rho\}$, such that $\varrho_{\bm\mu}(\bm0)=0$, $0\le\varrho_{\bm\mu}(\bm z)<\rho$, and
\begin{equation}
\bm w_*(\bm\mu)+U_{\bm\mu}\bm z-\bm n_{\bm\mu}b\in\Omega
\quad\Longleftrightarrow\quad \varrho_{\bm\mu}(\bm z)\le b,
\qquad \|\bm z\|_2\le\rho,\quad |b|\le\rho.
\label{eq:boundary-local-parametrization}
\end{equation}
Here $\varrho_{\bm\mu}(\bm z)$ is differentiable at $\bm0$, with $\nabla\varrho_{\bm\mu}(\bm0)=\bm0$ and $\varrho_{\bm\mu}(\bm z)=o(\|\bm z\|_2)$ as $\bm z\to\bm0$. 
The following conditions give lower and upper bounds for $\varrho_{\bm\mu}$.
\begin{enumerate}[label=(\roman*),itemsep=2pt]
\item \emph{Lower bound.} Suppose that, for some $\beta>0$ and $\underline\kappa>0$,
\begin{equation*}
\bm n_{\bm\mu}^\top(\bm w_*(\bm\mu)-\bm w)
\ge\frac{\underline\kappa}{1+\beta}\|U_{\bm\mu}^\top(\bm w-\bm w_*(\bm\mu))\|_2^{1+\beta},
\qquad \bm w\in\partial\Omega\text{ near }\bm w_*(\bm\mu).
\end{equation*}
Then, after decreasing $\rho$ if necessary,
\begin{equation*}
\varrho_{\bm\mu}(\bm z)\ge\frac{\underline\kappa}{1+\beta}\|\bm z\|_2^{1+\beta},
\qquad \|\bm z\|_2\le\rho.
\end{equation*}
\item \emph{Upper bound.} If $\partial\Omega$ is a $C^{1,\beta}$ hypersurface near $\bm w_*(\bm\mu)$, where $0<\beta\le1$, then $\varrho_{\bm\mu}\in C^{1,\beta}$ locally and, after decreasing $\rho$ if necessary, there exists $\bar\kappa>0$ such that
\begin{equation*}
0\le\varrho_{\bm\mu}(\bm z)\le\frac{\bar\kappa}{1+\beta}\|\bm z\|_2^{1+\beta},
\qquad \|\bm z\|_2\le\rho.
\end{equation*}
\end{enumerate}
\end{proposition}

With this local parametrization, two-sided power growth of the boundary profile $\varrho_{\bm{\mu}}$ determines the order of the local minimax risk and the maximal expected regret of SAA.

\begin{theorem}[Two-sided boundary contact determines the minimax rate]
\label{thm:geometry-contact-rate}
Suppose that $\Omega$ admits the local representation \eqref{eq:boundary-local-parametrization} of Proposition~\ref{prop:boundary-local-parametrization}, with a continuous convex $\varrho_{\bm\mu}$. Assume that, for some $\beta>0$ and constants $0<\underline\kappa\le\bar\kappa<\infty$,
\[
\frac{\underline\kappa}{1+\beta}\|\bm z\|_2^{1+\beta}
\le \varrho_{\bm\mu}(\bm z)\le
\frac{\bar\kappa}{1+\beta}\|\bm z\|_2^{1+\beta},
\qquad \|\bm z\|_2\le\rho.
\] Then, for every fixed $M>0$,
\[
R_{n,M}^{\Sigma}(\bm\mu)
\asymp
\sup_{\|\bm t\|_2\le M}
\mathbb E_{\bm\mu_{n,\bm t},\Sigma}R(\bm w_*(\bar{\bm x}))
\asymp
\|\bm\mu\|_2^{-1/\beta}n^{-(1+\beta)/(2\beta)},
\]
with implied constants that may be chosen uniformly over $\Sigma\in\mathcal C$.
If, in addition, $\varrho_{\bm\mu}(\delta\bm z)/\delta^{1+\beta}$ converges locally uniformly to a homogeneous convex function $\varrho_0$, then the local expansion in Theorem~\ref{thm:fast_general_local_rate} holds with
\[
a(\varepsilon)=\varepsilon^{(1+\beta)/\beta},
\qquad
\psi(\bm t)=\|\bm\mu\|_2^{-1/\beta}\varrho_0^*(U_{\bm\mu}^\top\bm t),
\]
where $\varrho_0^*(\bm y):=\sup_{\bm z\in\mathbb R^{d-1}}\{\bm y^{\top}\bm z-\varrho_0(\bm z)\}$ is the convex conjugate.
\end{theorem}

A more detailed statement is given in Theorem~\ref{thm:fast_boundary_growth_lower} of the Online Supplement. 

The exponent provides a useful dictionary. Quadratic boundary growth corresponds to $\beta=1$ and the $n^{-1}$ regime. When $0<\beta<1$, the boundary is sharper than quadratic and the support remainder has order greater than two, yielding rates faster than $n^{-1}$ such as $n^{-2}$. When $\beta>1$, the boundary is flatter and the rate lies between $n^{-1/2}$ and $n^{-1}$.

\begin{remark}[Sufficient conditions for quadratic boundary growth]
\label{rem:fast_holder_boundary_growth}
Under the geometric setup above, either Euclidean strong convexity together with a local $C^{1,1}$ boundary, or a local $C^2$ boundary with $\mathscr S_{\bm\mu}\succ0$, gives the two-sided boundary growth condition with $\beta=1$; precise statements and proofs are given in Lemma~\ref{lem:quadratic_boundary_growth} of the Online Supplement. Theorem~\ref{thm:geometry-contact-rate} then gives an $n^{-1}$ local minimax rate, attained by SAA. In the $C^2$ case, the quadratic profile also recovers the inverse-curvature formula in Corollary~\ref{cor:fast_principal_curvatures}. This complements Remark~\ref{rem:growth_implies_curvature}: positive curvature yields local quadratic contact, while the earlier remark derives positive curvature from the quadratic-growth condition \eqref{eq:fast_geometric_condition}.
\end{remark}

\subsection{Examples}
The examples below show why neither global strong convexity nor distance to the ambiguity set alone determines the local rate. Detailed calculations are given in Section~E of the Online Supplement.

\begin{example}[Rectangular constraints]
\label{ex:fast_rectangle}
Let $\Omega=[0,1]^d$, so $f_*(\bm\mu)=\sum_{j=1}^d(\mu_j)_+$. Denote $\mathcal S_0=\{j\in[d]:\mu_j=0\}$ and $s=\|\bm\mu\|_0=d-|\mathcal S_0|$. 

For $j\in[d]$, if $\mu_j=0$, then the optimal $w_{*,j}$ can take any value in $[0,1]$, whereas if $\mu_j\neq0$, $w_{*,j}=\mathbb I\{\mu_j>0\}$ is unique. Hence the ambiguity set is $A=\bigcup_{j=1}^d\{\bm \mu:\mu_j=0\}$, and $\bm\mu\in A$ exactly when $d>s$.

For fixed $M>0$ and sufficiently small $\varepsilon>0$, if $\mu_j\neq0$,  $\mu_j+\varepsilon t_j$ has the same sign as $\mu_j$, uniformly over $\bm{t}\in B_2(M)$. Hence $r_{\bm\mu}(\varepsilon\bm t)=\varepsilon\psi(\bm t)$, where
\begin{equation*}
\psi(\bm t)=\sum_{j\in\mathcal S_0}\left[(t_j)_+-w_{*,j}(\bm\mu)t_j\right].
\end{equation*}
If $s<d$, $\bm\mu$ is in the ambiguity set, and Corollary~\ref{cor:fast_first_order_limit} gives $V_{M,\Sigma}(\psi)>0$ because $\psi$ is nonaffine on $B_2(M)$. By Theorem~\ref{thm:fast_general_local_rate}, $R_{n,M}^{\Sigma}(\bm\mu)\asymp n^{-1/2}$, and SAA attains this local minimax rate.

If $s=d$, the entire local ball $B_2(M)$ shares the same optimal policy for sufficiently large $n$, so $R_{n,M}^{\Sigma}(\bm\mu)=0$. The maximal expected regret of SAA on $B_2(M)$ decays exponentially. Calculations and details are provided in Section~\ref{app:rectangle-details} of the Online Supplement.
\end{example}

\begin{example}[One strongly convex set, three local rates]
\label{ex:fast_center_dependent_rates}
For $\Omega=B_{4/3}^2:=\{\bm{w}\in\mathbb{R}^2:\|\bm{w}\|_{4/3}\le 1\}$, $f_*(\bm\mu)=\|\bm\mu\|_4$, and strict convexity gives $A=\{\bm0\}$. Fix $M>0$ and compare the local rates at three centers.

For $\bm{\mu}=\bm0$, any $\bm{w}\in\Omega$ is optimal. Corollary~\ref{cor:fast_first_order_limit} and Theorem~\ref{thm:fast_general_local_rate} therefore give $R_{n,M}^{\Sigma}\left(\bm{0}\right)\asymp n^{-1/2}$.

For $\bm{\mu}=2^{-1/4}(1,1)^\top$, direct calculation gives
\begin{equation*}
    H_{\bm{\mu}}
	=\frac{3}{2\sqrt{2}}
	\begin{pmatrix}1&-1\\-1&1\end{pmatrix}.
\end{equation*}
Hence by Theorem \ref{thm:fast_local_minimax}, $R_{n,M}^{\Sigma}\left(2^{-1/4}(1,1)^\top\right)\asymp n^{-1}$, where the sharp constant when $M\to\infty$ is
\begin{equation*}
	\lim_{M\to\infty}\lim_{n\to\infty}
	nR_{n,M}^{\Sigma}\bigl(2^{-1/4}(1,1)^\top\bigr)
	=\frac12\operatorname{tr}(\Sigma H_{\bm{\mu}})
	=\frac{3}{4\sqrt{2}}(\Sigma_{11}+\Sigma_{22}-2\Sigma_{12}).
\end{equation*}
The same $n^{-1}$ rate applies to any $\bm{\mu}=(\mu_1,\mu_2)^\top$ with $\mu_1\mu_2\neq0$.

For $\bm{\mu}=\bm{e}_1:=(1,0)^\top$, note that $H_{\bm{e}_1}=\bm{0}$. Uniformly for bounded $\bm{t}$,
\begin{equation*}
	f_*(\bm{e}_1+\varepsilon\bm{t})
	=1+\varepsilon t_1+\frac{\varepsilon^4t_2^4}{4}
	+o(\varepsilon^4).
\end{equation*}
Theorem \ref{thm:fast_general_local_rate} therefore applies with
$a(\varepsilon)=\varepsilon^4$ and $\psi(\bm{t})=t_2^4/4$, yielding $R_{n,M}^{\Sigma}(\bm{e}_1)\asymp n^{-2}$.
More generally, at any nonzero point on either coordinate axis, the same $n^{-2}$ rate holds. Under the polar-norm convention \eqref{eq:polar_norm_convention}, $\dist(\bm\mu,A)=\|\bm\mu\|_4$, so the two displayed nonzero centers both have distance $1$ from $A$, yet have different local minimax rates. Theorem~\ref{thm:fast_general_local_rate} shows that SAA attains all three rates; at the smooth centers with nonzero Hessian, Theorem~\ref{thm:fast_local_minimax} gives sharp attainment in the expanding-neighborhood limit. More details are given in Section~\ref{app:l4/3-value-decision} of the Online Supplement.
\end{example}

\begin{example}[Optimal value versus optimal decision on $\ell_{p^\star}$ balls]
\label{ex:lq_value_decision}
Let $d\ge2$, $1<p<\infty$, and let $p^\star=p/(p-1)$. For
$\Omega=B_{p^\star}^d=\{\bm w\in\mathbb{R}^d:\|\bm w\|_{p^\star}\le1\}$,
\begin{equation*}
    f_*(\bm\mu)=\|\bm\mu\|_p,
    \qquad
    w_{*,j}(\bm\mu)
      =\frac{\operatorname{sgn}(\mu_j)|\mu_j|^{p-1}}
             {\|\bm\mu\|_p^{p-1}}
      \quad(\bm\mu\ne\bm0).
\end{equation*}
For $\bm\mu=\mu_1\bm e_1$, $\mu_1>0$, the unique optimizer is $\bm w_*(\bm\mu)=\bm e_1$ and, uniformly for bounded $\bm t$,
\begin{equation*}
    f_*(\bm\mu+\varepsilon\bm t)
    =\mu_1+\varepsilon t_1
      +\frac{\varepsilon^p}{p\mu_1^{p-1}}
         \sum_{j=2}^d|t_j|^p+o(\varepsilon^p).
\end{equation*}
Estimating the optimal value still involves the first-order radial term $\varepsilon t_1$, whereas decision regret subtracts that affine term and gives  $a(\varepsilon)=\varepsilon^p$ and $\psi(\bm t)=(p\mu_1^{p-1})^{-1}\sum_{j=2}^d|t_j|^p$. Since $d\ge2$, $\psi(\bm t)$ is nonaffine on every $B_2(M)$, $M>0$. Theorem~\ref{thm:fast_general_local_rate} yields, for fixed $M>0$,
\begin{equation*}
    R_{n,M}^{\Sigma}(\bm\mu)
    \asymp\sup_{\bm t\in B_2(M)}
       \mathbb E_{\bm\mu_{n,\bm t},\Sigma}R(\bm w_*(\bar{\bm x}))
    \asymp n^{-p/2}.
\end{equation*}
The case $p=4$ recovers the axis-point rate in Example~\ref{ex:fast_center_dependent_rates}.

Thus, in the present fixed-dimensional Gaussian experiment, the same problem can have a first-order $n^{-1/2}$ value-estimation scale (under expected absolute loss), but much faster optimal-decision regret. Detailed calculations are given in Section~\ref{app:lq-value-decision} of the Online Supplement.
\end{example}

\begin{example}[A non-power rate]
\label{ex:fast_nonpower_rate}
Choose $0<\delta<e^{-2}$ and define, for $|s|\le\delta$,
\begin{equation*}
	g(0)=0,\qquad
	g(s)=\frac{|s|}{\log(1/|s|)}\quad(s\ne0).
\end{equation*}
Then $g$ is even, convex, and $C^1$, with $g'(0)=0$. Define \begin{equation*}
	\Omega=\operatorname{conv}\left(
	\{\bm{0}\}\cup
	\left\{\bigl(1+g(s)-sg'(s),\,g'(s)\bigr):|s|\le\delta\right\}
	\right).
\end{equation*}
For $\bm{\mu}=\bm{e}_1$, there exists a unique optimizer $\bm w_*(\bm e_1)=\bm e_1$. For bounded $\bm t$ and sufficiently small $\varepsilon>0$, direct calculation gives
\begin{equation*}
    r_{\bm{e}_1}(\varepsilon\bm{t})=(1+\varepsilon t_1)
		g\left(\frac{\varepsilon t_2}{1+\varepsilon t_1}\right).
\end{equation*}
For every fixed $M>0$,
\begin{equation*}
\lim_{\varepsilon\downarrow0}\sup_{\bm t\in B_2(M)}
\left|\frac{r_{\bm e_1}(\varepsilon\bm t)}
{\varepsilon/\log(1/\varepsilon)}-|t_2|\right|=0,
\end{equation*}
Hence $a(\varepsilon)=\varepsilon/\log(1/\varepsilon)$, and $\psi(\bm{t})=|t_2|$ is not affine. Theorem~\ref{thm:fast_general_local_rate} gives
\begin{equation*}
R_{n,M}^{\Sigma}(\bm e_1)\asymp\frac1{\sqrt n\log n},
\end{equation*}
attained in rate by SAA. This rate is $o(n^{-1/2})$ but slower than $n^{-\alpha/2}$ for every $\alpha>1$. Thus uniqueness and differentiability do not force a power expansion, and a theory restricted to integer-order Taylor terms would miss legitimate local regimes. More details are given in Section~\ref{app:non_power_local_rate} of the Online Supplement.
\end{example}

\begin{table}[H]
\centering
\caption{Representative local geometries and nominal minimax rates.}
\label{tab:geometry-rate-summary}
\renewcommand{\arraystretch}{1.18}
\begin{tabular}{@{}>{\raggedright\arraybackslash}p{0.29\textwidth}>{\raggedright\arraybackslash}p{0.23\textwidth}>{\centering\arraybackslash}p{0.19\textwidth}>{\raggedright\arraybackslash}p{0.20\textwidth}@{}}
\toprule
Local geometry at the exposed decision & Local scale $a(\varepsilon)$ & Minimax regret & Representative example\\
\midrule
Nonunique exposed face & $\varepsilon$ & $n^{-1/2}$ & Rectangle at a zero coordinate\\
Positive quadratic curvature & $\varepsilon^2$ & $n^{-1}$ & Euclidean ball, generic smooth point\\
General $p$th-order contact & $\varepsilon^p$ & $n^{-p/2}$ & $\ell_{p^\star}$ ball at an axis point\\
Non-power contact & $\varepsilon/\log(1/\varepsilon)$ & $1/(\sqrt n\log n)$ & Example~\ref{ex:fast_nonpower_rate}\\
\bottomrule
\end{tabular}
\end{table}

\section{Local Minimax Risk on KL Neighborhoods}
\label{sect:divergence_neighborhoods}

\subsection{The robust local experiment}

For the following mean parametrization, assume that, for some $\eta>0$,
\begin{equation}
\begin{aligned}
    &\mathbb E_P\bm x=\bm\mu,\qquad\operatorname{Cov}_P(\bm x)=\Sigma\in\mathcal C,\\
    &\mathbb E_P e^{\bm\theta^\top(\bm x-\bm\mu)}<\infty,
       \qquad \|\bm\theta\|_2<\eta.
\end{aligned}
    \label{eq:ng_reference_distribution}
\end{equation}
Define
\begin{equation}
    \Lambda_P(\bm\theta)
       :=\log\mathbb E_P e^{\bm\theta^\top(\bm x-\bm\mu)},
       \qquad \|\bm\theta\|_2<\eta.
    \label{eq:ng_cumulant_definition}
\end{equation}
The dimension $d$ and the reference distribution $P$ are fixed.

Throughout, $\bm\mu=\mathbb E_P\bm x$ and $\Sigma=\operatorname{Cov}_P(\bm x)\succ0$, and $P$ satisfies \eqref{eq:ng_reference_distribution} except where a result is explicitly stated under a finite second moment only. For any probability measure $Q$ with $\mathbb E_Q\|\bm x\|_2<\infty$, write
\begin{equation*}
    \bm\mu_Q:=\mathbb E_Q\bm x,\qquad
    R_Q(\bm w):=f_*(\bm\mu_Q)-\bm w^\top\bm\mu_Q,\qquad \bm w\in\Omega.
\end{equation*}
The observations $\bm{x}_1,\ldots,\bm{x}_n$ are i.i.d.\ with law $Q$. The covariance assumption $\Sigma\in\mathcal{C}$ is only imposed on the reference law $P$.

\subsection{A matched Gaussian benchmark}
To isolate the statistical price of allowing full distributional perturbations, it is useful to compare the KL neighborhood with a Gaussian experiment having the same local mean geometry. For fixed $\Sigma\succ0$, define
\begin{equation}
\begin{aligned}
R_{n,M}^{\mathrm{ell},\Sigma}(\bm\mu)
:=\inf_{\widehat{\bm w}_n}\sup_{\bm u\in B_2(M)}
\mathbb E_{\bm\mu+\Sigma^{1/2}\bm u/\sqrt n,\Sigma}
\left[f_*\left(\bm\mu+\frac{\Sigma^{1/2}\bm u}{\sqrt n}\right)
-\widehat{\bm w}_n^\top
\left(\bm\mu+\frac{\Sigma^{1/2}\bm u}{\sqrt n}\right)\right].
\end{aligned}
\label{eq:matched_gaussian_benchmark}
\end{equation}
The local parameter set is the same covariance ellipsoid that will emerge from the KL constraint.

\begin{corollary}[Matched ellipsoidal Gaussian experiment]
\label{cor:matched_gaussian_benchmark}
Under the local expansion \eqref{eq:fast_general_local_expansion0}, for every fixed $M>0$,
\begin{equation}
    \lim_{n\to\infty}
    \frac{R_{n,M}^{\mathrm{ell},\Sigma}(\bm\mu)}{a(n^{-1/2})}
    =V_{M,I_d}(\psi\circ\Sigma^{1/2}).
    \label{eq:matched_gaussian_limit}
\end{equation}
If $\psi$ is nonaffine on $\Sigma^{1/2}B_2(M)$, ordinary SAA is rate-optimal for this Gaussian benchmark. If $f_*$ is $C^2$ near $\bm\mu$, then
\[
    \lim_{n\to\infty}nR_{n,M}^{\mathrm{ell},\Sigma}(\bm\mu)
    =\frac12\mathcal B_{M,I_d}
       (\Sigma^{1/2}H_{\bm\mu}\Sigma^{1/2}),
\]
and the expanding-neighborhood constant is $\operatorname{tr}(\Sigma H_{\bm\mu})/2$.
\end{corollary}

\subsection{KL neighborhoods and a truncated mean}

For $M>0$, define
\begin{equation}
\begin{aligned}
    \mathcal Q_{n,M}^{\mathrm{KL}}(P)
       &:=\left\{Q\ll P:\mathbb E_Q\|\bm x\|_2<\infty,\ D_{\mathrm{KL}}(Q\|P)\le\frac{M^2}{2n}\right\},\\
    \mathcal R_{n,M}^{\mathrm{KL}}(P)
       &:=\inf_{\hat{\bm w}_n}\sup_{Q\in\mathcal Q_{n,M}^{\mathrm{KL}}(P)}
               \mathbb E_{Q^{\otimes n}}R_Q(\hat{\bm w}_n).
\end{aligned}
    \label{eq:div_kl_benchmark}
\end{equation}
The infimum is over measurable maps $(\mathbb R^d)^n\to\Omega$. As in the preceding local benchmarks, a rule chosen for a fixed $(P,n,M)$ may depend on that reference model, but not on the unknown $Q$. The common truncated plug-in policy defined below uses neither $P$ nor $\bm\mu$ nor $\Sigma$.

Define the coordinatewise truncated observations and their mean by
\begin{equation}
\begin{aligned}
\tau_n&=\log^2(n+1),\qquad
x_{i,j}^{(n)}=\operatorname{sgn}(x_{i,j})\min\{|x_{i,j}|,\tau_n\},\\
\widetilde{\bm\mu}_n&=\frac1n\sum_{i=1}^n\bm x_i^{(n)},
\qquad
\bm Z_{n,Q}=\sqrt n(\widetilde{\bm\mu}_n-\bm\mu_Q).
\end{aligned}
\label{eq:div_truncated_mean}
\end{equation}
The truncation threshold $\tau_n$ is deterministic and does not use $P$, $\bm\mu$, or $\Sigma$. Uniform mean-ellipsoid bounds, moment bounds, and a Gaussian approximation for $\bm Z_{n,Q}$ are given in the Online Supplement.

The threshold is applied in a fixed coordinate system and fixed units. The rule is neither affine-equivariant nor invariant to changes of origin, although the asymptotic results hold for every fixed such choice. 

\subsection{Exact local minimax risk}
\begin{theorem}[Exact local minimax limit on a KL neighborhood]
\label{thm:div_kl_minimax}
Assume \eqref{eq:fast_general_local_expansion0} and let $P$ satisfy
\begin{equation*}
    \mathbb E_P\|\bm x\|_2^2<\infty,\qquad
    \bm\mu=\mathbb E_P\bm x,\qquad
    \Sigma=\operatorname{Cov}_P(\bm x)\succ0.
\end{equation*}
Fix $M>0$.
\begin{enumerate}[label=(\roman*),itemsep=2pt]
\item \emph{Lower bound.}
\begin{equation}
    \liminf_{n\to\infty}
       \frac{\mathcal R_{n,M}^{\mathrm{KL}}(P)}{a(n^{-1/2})}
       \ge V_{M,I_d}(\psi\circ\Sigma^{1/2}).
    \label{eq:div_kl_lower}
\end{equation}
\item \emph{Exact limit.} If, in addition, \eqref{eq:ng_reference_distribution} holds, then the limit exists and
\begin{equation}
    \lim_{n\to\infty}
       \frac{\mathcal R_{n,M}^{\mathrm{KL}}(P)}{a(n^{-1/2})}
       = V_{M,I_d}(\psi\circ\Sigma^{1/2}).
    \label{eq:div_kl_exact_limit}
\end{equation}
\item \emph{Rate-minimaxity of truncated SAA.} Under the conditions of (ii), if $V_{M,I_d}(\psi\circ\Sigma^{1/2})>0$, then the policy $\bm w_*(\widetilde{\bm\mu}_n)$ is rate-minimax:
\begin{equation}
    ca(n^{-1/2})
       \le\mathcal R_{n,M}^{\mathrm{KL}}(P)
       \le\sup_{Q\in\mathcal Q_{n,M}^{\mathrm{KL}}(P)}
              \mathbb E_{Q^{\otimes n}}R_Q(\bm w_*(\widetilde{\bm\mu}_n))
       \le Ca(n^{-1/2}),
    \label{eq:div_kl_plugin_rate}
\end{equation}
where constants $0<c\le C<\infty$ may depend on $P$, $M$ and the fixed geometry.
\end{enumerate}
\end{theorem}

The exact limit in \eqref{eq:div_kl_exact_limit} is identical to the matched Gaussian limit in \eqref{eq:matched_gaussian_limit}. Thus, to first order, the full KL neighborhood and the Gaussian experiment over the covariance ellipsoid have the same information bound. The lower bound (i) is information-theoretic and needs only a finite second moment. The exponential-moment condition is used in (ii) and (iii). The upper bound in (ii) uses deterministic policies based on the truncated mean, obtained by averaging out auxiliary randomization; these policies depend on $P$, $M$, and the local geometry and achieve arbitrarily small excess over the minimax constant in asymptotic normalized maximal regret. The truncated plug-in rule $\bm w_*(\widetilde{\bm\mu}_n)$, which does not use $P$, attains the rate in (iii) uniformly over the full KL neighborhood and, by Theorem~\ref{thm:div_plugin_comparison} below, its leading constant is $C_M^{\mathrm{plug}}(\psi,\Sigma)\ge V_{M,I_d}(\psi\circ\Sigma^{1/2})$.

\subsection{The local price of robustness}
\begin{corollary}[Smooth constants and the geometric consequences]
\label{cor:div_kl_smooth}
Under \eqref{eq:ng_reference_distribution}, if $f_*$ is $C^2$ on a neighborhood of $\bm\mu$, then
\begin{equation}
    \lim_{n\to\infty}n\mathcal R_{n,M}^{\mathrm{KL}}(P)
       =\frac12\mathcal B_{M,I_d}(\Sigma^{1/2}H_{\bm\mu}\Sigma^{1/2})
    \label{eq:exact_limit_for_KL_neighborhood}
\end{equation}
\begin{equation}
    \lim_{M\to\infty}\lim_{n\to\infty}n\mathcal R_{n,M}^{\mathrm{KL}}(P)
       =\frac12\operatorname{tr}(\Sigma H_{\bm\mu}),
    \label{eq:asymp_limit_for_KL_neighborhood}
\end{equation}
\begin{equation}
    \lim_{n\to\infty}\sup_{Q\in\mathcal Q_{n,M}^{\mathrm{KL}}(P)}
       \left|n\mathbb E_{Q^{\otimes n}}R_Q(\bm w_*(\widetilde{\bm\mu}_n))
                    -\frac12\operatorname{tr}(\Sigma H_{\bm\mu})\right|=0.
    \label{eq:div_kl_smooth}
\end{equation}
Thus $\bm w_*(\widetilde{\bm\mu}_n)$ attains the sharp constant after $M\to\infty$. 
The geometric identities for $H_{\bm\mu}$ in Corollaries~\ref{cor:fast_distance_smooth} and \ref{cor:fast_principal_curvatures} remain valid under their respective assumptions. 
\end{corollary}

There are three distinct notions of a robustness cost. First, relative to the matched ellipsoidal Gaussian benchmark, the KL neighborhood is first-order costless: it has the same exact fixed-$M$ minimax limit. Second, comparison with the Euclidean-ball Gaussian benchmark can produce a different fixed-$M$ constant because the two experiments use different local parameter sets. Third, at smooth instances and after $M\to\infty$, all of these formulations have the same constant $\operatorname{tr}(\Sigma H_{\bm\mu})/2$. Local robustness is therefore asymptotically free at the ordinary smooth efficiency scale, even though a tail-stable implementation may be needed to attain the optimal rate or sharp constant uniformly over the KL neighborhood.

\begin{corollary}[Boundary contact rates on KL neighborhoods]
\label{cor:div-boundary-contact-rate}
Under the two-sided local boundary condition of Theorem~\ref{thm:geometry-contact-rate}, and under the moment conditions of Theorem~\ref{thm:div_kl_minimax},
\[
\mathcal R_{n,M}^{\mathrm{KL}}(P)
\asymp
\sup_{Q\in\mathcal Q_{n,M}^{\mathrm{KL}}(P)}
\mathbb E_{Q^{\otimes n}}R_Q(\bm w_*(\widetilde{\bm\mu}_n))
\asymp n^{-(1+\beta)/(2\beta)}.
\]
The lower bound needs only a finite second moment, while the matching truncated-SAA upper bound uses the local exponential-moment condition.
\end{corollary}

\section{Limits of SAA under Distributional Robustness}
\label{sec:saa_limits}

Keep $P,d,\Omega,M$ fixed and define
\[
W_{n,M}^{\mathrm{mean}}(P)
=\sup_{Q\in\mathcal Q_{n,M}^{\mathrm{KL}}(P)}
\mathbb E_{Q^{\otimes n}}R_Q(\bm w_*(\bar{\bm x})),
\]
\[
W_{n,M}^{\mathrm{tr}}(P)
=\sup_{Q\in\mathcal Q_{n,M}^{\mathrm{KL}}(P)}
\mathbb E_{Q^{\otimes n}}R_Q(\bm w_*(\widetilde{\bm\mu}_n)),
\qquad
\bm\Delta_n=\bar{\bm x}-\widetilde{\bm\mu}_n.
\]
The truncated and ordinary means have the same sequential Gaussian limit under local alternatives. However, their expected regrets need not agree, because the faster geometric regimes amplify rare tail events, which may alter the rate of maximal expected regret at faster geometric scales and its leading constant even at smooth quadratic instances; see Example~\ref{ex:div_matching_tail_geometry}.

\begin{table}[htbp]
\centering
\caption{Statistical benchmarks and procedure-specific maximal risks used in the paper.}
\label{tab:benchmark-guide}
\small
\renewcommand{\arraystretch}{1.16}
\begin{tabular}{@{}>{\raggedright\arraybackslash}p{0.20\textwidth}>{\raggedright\arraybackslash}p{0.25\textwidth}>{\raggedright\arraybackslash}p{0.24\textwidth}>{\raggedright\arraybackslash}p{0.22\textwidth}@{}}
\toprule
Quantity & Statistical experiment & Local parameter set & Rule may use \\
\midrule
$R_{n,M}^{\Sigma}(\bm\mu)$ & Gaussian, fixed covariance & $\bm\mu+B_2(M)/\sqrt n$ & benchmark center and fixed $\Sigma$ \\
$R_{n,M}^{\mathrm{ell},\Sigma}(\bm\mu)$ & Matched Gaussian & $\bm\mu+\Sigma^{1/2}B_2(M)/\sqrt n$ & benchmark center and fixed $\Sigma$ \\
$\mathcal R_{n,M}^{\mathrm{KL}}(P)$ & KL neighborhood & $D_{\mathrm{KL}}(Q\|P)\le M^2/(2n)$ & reference law $P$, but not $Q$ \\
$W_{n,M}^{\mathrm{mean}},\ W_{n,M}^{\mathrm{tr}}$ & KL maximal risks & same KL neighborhood & prescribed ordinary or truncated plug-in rule \\
\bottomrule
\end{tabular}
\end{table}

\subsection{Leading plug-in risks}
Let $\psi$ be the limit in \eqref{eq:fast_general_local_expansion0}. As a finite convex function, $\psi$ is locally Lipschitz and differentiable Lebesgue-almost everywhere. For $\bm Z\sim\mathcal N(\bm0,\Sigma)$ with $\Sigma\succ0$, the vector $\bm t+\bm Z$ has a density, so $\nabla\psi(\bm t+\bm Z)$ is almost surely well defined, and we define
\begin{equation}
	C_M^{\mathrm{plug}}(\psi,\Sigma)
	:=\sup_{\|\Sigma^{-1/2}\bm t\|_2\le M}
	\mathbb E\left[\psi(\bm t)-\psi(\bm t+\bm Z)
	+\nabla\psi(\bm t+\bm Z)^\top\bm Z\right].
	\label{eq:div_gaussian_plugin_constant}
\end{equation}

\begin{theorem}[Asymptotic comparison of maximal plug-in regrets]
\label{thm:div_plugin_comparison}
Assume \eqref{eq:ng_reference_distribution} and \eqref{eq:fast_general_local_expansion0}, and fix $M>0$ such that $V_{M,I_d}(\psi\circ\Sigma^{1/2})>0$. Then $0<C_M^{\mathrm{plug}}(\psi,\Sigma)<\infty$, and
\begin{equation}
    \lim_{n\to\infty}\frac{W_{n,M}^{\mathrm{tr}}(P)}{a(n^{-1/2})}
       =C_M^{\mathrm{plug}}(\psi,\Sigma)
       \le\liminf_{n\to\infty}
                    \frac{W_{n,M}^{\mathrm{mean}}(P)}{a(n^{-1/2})}.
    \label{eq:div_asymptotic_plugin_comparison}
\end{equation}
\end{theorem}

In particular, Theorem \ref{thm:div_plugin_comparison} implies that
\begin{equation*}
	\limsup_{n\to\infty}
	\frac{W_{n,M}^{\mathrm{tr}}(P)-W_{n,M}^{\mathrm{mean}}(P)}
	{a(n^{-1/2})}\le0,
	\qquad
	\limsup_{n\to\infty}
	\frac{W_{n,M}^{\mathrm{tr}}(P)}{W_{n,M}^{\mathrm{mean}}(P)}\le1.
\end{equation*}
Strict improvement of the rate is exhibited in Example~\ref{ex:div_matching_tail_geometry}.

The constant $C_M^{\mathrm{plug}}$ is the maximal asymptotic Gaussian plug-in risk instead of the minimax constant $V_{M,I_d}$. The theorem says that truncation never worsens the leading maximal plug-in regret in that ordinary SAA cannot have a smaller first-order maximal risk. The remaining issue is whether the ordinary mean has enough tail control to match the truncated rule.

\subsection{A geometry-weighted tail criterion}
\begin{proposition}[A geometry-weighted tail criterion for the ordinary mean]
\label{prop:div_geometric_tail_criterion}
Under the assumptions of Theorem~\ref{thm:div_plugin_comparison}, there is $C<\infty$ such that, for all sufficiently large $n$,
\begin{equation}
    W_{n,M}^{\mathrm{mean}}(P)
       \le C a(n^{-1/2})
          +\frac12\sup_{Q\in\mathcal Q_{n,M}^{\mathrm{KL}}(P)}\mathbb E_{Q^{\otimes n}}
                                r_{\bm\mu}(4\bm\Delta_n)
       \le C\{a(n^{-1/2})+n^{-1}\},
    \label{eq:div_geometric_tail_upper}
\end{equation}
where the constants may depend on the fixed $P,M,d,\Omega$ and the local expansion.
\end{proposition}

The criterion is directional: it measures the truncation error through $r_{\bm\mu}$ rather than through a generic norm. The conditions below are sufficient. This distinction matters because only the directions to which the optimizer is locally sensitive contribute to regret.

If
\begin{equation}
	\sup_{Q\in\mathcal Q_{n,M}^{\mathrm{KL}}(P)}\mathbb E_{Q^{\otimes n}}r_{\bm\mu}(4\bm\Delta_n)
	=O(a(n^{-1/2})).
	\label{eq:div_geometric_tail_rate_condition}
\end{equation}
then Proposition~\ref{prop:div_geometric_tail_criterion} implies that
\begin{equation*}
	W_{n,M}^{\mathrm{mean}}(P)
	\asymp W_{n,M}^{\mathrm{tr}}(P)
	\asymp\mathcal R_{n,M}^{\mathrm{KL}}(P)
	\asymp a(n^{-1/2}).
\end{equation*}
The stronger condition
\begin{equation}
	\sup_{Q\in\mathcal Q_{n,M}^{\mathrm{KL}}(P)}\mathbb E_{Q^{\otimes n}}r_{\bm\mu}(4\bm\Delta_n)
	=o(a(n^{-1/2}))
	\label{eq:div_geometric_tail_constant_condition}
\end{equation}
on the other hand gives
\begin{equation}
	\lim_{n\to\infty}\frac{W_{n,M}^{\mathrm{mean}}(P)}{a(n^{-1/2})}
	=\lim_{n\to\infty}\frac{W_{n,M}^{\mathrm{tr}}(P)}{a(n^{-1/2})}
	=C_M^{\mathrm{plug}}(\psi,\Sigma).
	\label{eq:div_equal_plugin_constants}
\end{equation}
In particular, $\liminf_{n\to\infty}n a(n^{-1/2})>0$ suffices for rate-optimality, and $n a(n^{-1/2})\to\infty$ suffices for \eqref{eq:div_equal_plugin_constants}. 

\subsection{A phase boundary for power geometry and tails}
\begin{corollary}[Local order and reference tail strength]
\label{cor:div_local_order_tail_strength}
Under the assumptions of Theorem~\ref{thm:div_plugin_comparison}, suppose $a(\varepsilon)=\varepsilon^\alpha$ for some $\alpha\ge1$, and that for some $\zeta\ge1$ and $c>0$,
\begin{equation}
    \mathbb E_P e^{c\|\bm x\|_2^\zeta}<\infty.
    \label{eq:div_reference_tail_power}
\end{equation}
Then
\begin{equation*}
    \sup_{Q\in\mathcal Q_{n,M}^{\mathrm{KL}}(P)}\mathbb E_{Q^{\otimes n}}\|\bm\Delta_n\|_2^\zeta=O(n^{-\zeta}),
\end{equation*}
and
\begin{equation}
    W_{n,M}^{\mathrm{mean}}(P)\le C(n^{-\alpha/2}+n^{-\zeta}).
    \label{eq:div_power_tail_regret_bound}
\end{equation}
\end{corollary}

This corollary indicates: if $\alpha\le2\zeta$, the ordinary-mean policy is rate-minimax:
\begin{equation*}
	W_{n,M}^{\mathrm{mean}}(P)
	\asymp W_{n,M}^{\mathrm{tr}}(P)
	\asymp\mathcal R_{n,M}^{\mathrm{KL}}(P)
	\asymp n^{-\alpha/2}.
\end{equation*}
If $\alpha<2\zeta$, it also has the same asymptotic maximal regret constant as the truncated policy:
\begin{equation*}
	n^{\alpha/2}W_{n,M}^{\mathrm{mean}}(P)
	\to C_M^{\mathrm{plug}}(\psi,\Sigma).
\end{equation*}
For $\alpha>2\zeta$, rate-optimality is not guaranteed. The exponent $\zeta$ describes an exponential moment of the reference distribution $P$.

The boundary $\alpha=2\zeta$ separates geometry that is slow enough for the raw sample mean from geometry that is fast enough to expose rare tail contamination. More precisely, ordinary and truncated SAA have the same leading maximal-regret constant when $\alpha<2\zeta$; ordinary SAA remains rate-minimax but can have strictly larger leading maximal regret when $\alpha=2\zeta$; and ordinary SAA can be rate-suboptimal when $\alpha>2\zeta$.
In the regular quadratic case $\alpha=2$, any subexponential reference law with $\zeta\ge1$ lies on the rate-optimal side. Rate suboptimality becomes possible only when the local optimization problem is easier than the tail behavior of the sample mean.

\begin{table}[t]
\centering
\caption{Geometry-tail regimes for ordinary SAA on shrinking KL neighborhoods.}
\label{tab:geometry-tail-phase}
\small
\renewcommand{\arraystretch}{1.15}
\begin{tabular}{@{}>{\centering\arraybackslash}p{0.20\textwidth}>{\raggedright\arraybackslash}p{0.66\textwidth}@{}}
\toprule
Regime & Consequence for ordinary SAA \\
\midrule
$\alpha<2\zeta$ & Same rate and the same leading maximal-regret constant as truncated SAA. \\
$\alpha=2\zeta$ & Rate-minimax; leading maximal regret can be strictly larger than for truncated SAA (Example~\ref{ex:div_matching_tail_geometry}). \\
$\alpha>2\zeta$ & Rate-minimaxity is not guaranteed and can fail sharply. \\
\bottomrule
\end{tabular}
\end{table}

\begin{corollary}[A smooth upper bound allowing degenerate curvature]
\label{cor:div_raw_smooth_rate}
Assume \eqref{eq:ng_reference_distribution}, and suppose $f_*$ is $C^2$ near $\bm\mu$. For each fixed $M>0$,
\begin{equation}
    \sup_{Q\in\mathcal Q_{n,M}^{\mathrm{KL}}(P)}
          \mathbb E_{Q^{\otimes n}}R_Q(\bm w_*(\bar{\bm x}))
                         =O(n^{-1}).
    \label{eq:div_raw_smooth_rate}
\end{equation}
If $H_{\bm\mu}\ne0$, then the sample-mean policy is rate-minimax on this KL neighborhood.
\end{corollary}

This corollary clarifies the role of degeneracy. Twice differentiability alone gives an $O(n^{-1})$ robust bound for ordinary SAA. If the Hessian is nonzero, this matches the minimax rate. If the Hessian vanishes and a faster rate is attainable, ordinary SAA may remain stuck at $n^{-1}$.

\subsection{A sharp separation example}
The mechanism can be seen from a back-of-the-envelope calculation (heuristic, not used in the proofs). Suppose $P$ has tails of order $e^{-\|\bm x\|^\zeta}$ in some coordinate. A perturbation $Q\in\mathcal Q_{n,M}^{\mathrm{KL}}(P)$ may move mass $q$ to distance of order $s$ at KL cost of about $q\,s^\zeta$, so $q\approx n^{-1}s^{-\zeta}$ is admissible. Taking $s\asymp n$ gives $q\approx n^{-1-\zeta}$. This perturbation changes $\bm\mu_Q$ only by $O(n^{-1-\zeta}\cdot n)=O(n^{-\zeta})$, which is negligible on the $n^{-1/2}$ scale, but with probability about $nq\approx n^{-\zeta}$ one observation lands in the perturbed region and displaces $\bar{\bm x}$ by a constant amount. The resulting regret is of order one, so the expected regret of ordinary SAA carries a term of order $n^{-\zeta}$, which dominates $n^{-\alpha/2}$ when $\alpha>2\zeta$. Truncation at $\tau_n=\log^2(n+1)$ caps the displacement caused by any single observation at $\tau_n/n$ and removes this term.

\begin{example}[Power geometry with exponential-power tails]
\label{ex:div_matching_tail_geometry}
Fix $\zeta\ge1$ and $\alpha>1$ (the restriction $\alpha>1$ ensures that the dual exponent $\alpha/(\alpha-1)$ is finite). Let $d=2$ and let $P$ be the law of $\bm x=(x_1,x_2)$ with independent coordinates, $\mathbb P\{x_1=1/2\}=\mathbb P\{x_1=3/2\}=1/2$ and $x_2$ having density proportional to $e^{-|z|^\zeta}$ (this is \eqref{eq:div_tail_reference_law} in the Online Supplement). Then $P$ has mean $\bm e_1$ and satisfies \eqref{eq:ng_reference_distribution} and \eqref{eq:div_reference_tail_power}. Take the feasible set $\Omega=B_{\alpha/(\alpha-1)}^2$, the unit $\ell_{\alpha/(\alpha-1)}$-ball in $\mathbb R^2$, so that the local support remainder at $\bm e_1$ has order $\varepsilon^\alpha$. Then
\[
\mathcal R_{n,M}^{\mathrm{KL}}(P)
\asymp W_{n,M}^{\mathrm{tr}}(P)
\asymp n^{-\alpha/2},
\qquad
W_{n,M}^{\mathrm{mean}}(P)
\asymp n^{-\alpha/2}+n^{-\zeta}.
\]
Hence the phase boundary in Corollary~\ref{cor:div_local_order_tail_strength} is sharp in rate. In particular, with Laplace tails $(\zeta=1)$ and quartic geometry $(\alpha=4)$,
\[
W_{n,M}^{\mathrm{mean}}(P)\asymp n^{-1},
\qquad
W_{n,M}^{\mathrm{tr}}(P)
\asymp\mathcal R_{n,M}^{\mathrm{KL}}(P)\asymp n^{-2}.
\]
Ordinary SAA is then strictly suboptimal in rate despite sharing the same Gaussian weak limit as the truncated procedure.
At the boundary $(\zeta,\alpha)=(1,2)$, the same construction with $\Omega=B_2^2$ gives, for every fixed $M>0$,
\begin{equation*}
	\liminf_{n\to\infty}nW_{n,M}^{\mathrm{mean}}(P)\ge1+\frac{M^2}{4}(1-2^{-1/2})>1=\lim_{n\to\infty}nW_{n,M}^{\mathrm{tr}}(P).
\end{equation*}
Thus equality of the leading maximal regrets can fail even at a smooth quadratic instance. The construction and proofs are given in Section~\ref{app:tail-geometry-construction} of the Online Supplement.
\end{example}

\begin{remark}[Scope of the separation]
The separation concerns expected regret. By \eqref{eq:div_mean_truncation_lone} and Lemma~\ref{lem:div_kl_truncation} of the Online Supplement, the ordinary and truncated means are asymptotically equivalent at the $n^{-1/2}$ scale uniformly over the KL neighborhood and have the same centered Gaussian limit; the gap arises only because regret is evaluated in expectation uniformly over KL-admissible laws, and rare observations of order-$n$ magnitude contribute $O(1)$ regret with probability of order $n^{-\zeta}$. 
\end{remark}

\section{Numerical Experiments}
\label{sec:numerical_experiments}

We illustrate two consequences of the local minimax theory. First, for the Gaussian experiment, the SAA procedure can attain different regret rates as the local geometry changes, even on a fixed feasible set. Second, under distributional perturbations, this geometry determines whether rare observations change the performance of ordinary SAA. 

Throughout, the regret is $R_{\bm\mu}(\bm w)=f_*(\bm\mu)-\bm w^\top\bm\mu$. We report expected regret along the specified Gaussian local alternatives or under the specified distributions $Q_n$. We compare these pointwise estimates with the minimax rates and constants derived above. Error bars are pointwise $95\%$ Monte Carlo intervals. 

For the Gaussian experiments, we generate the exact sufficient statistic
$\bar{\bm x}=\bm\mu+(\bm t+\bm Z)/\sqrt n$, with
$\bm Z\sim\mathcal N(\bm0,I_2)$, using $10^5$ independent replications per sample size.
For the Laplace experiments, we generate all $n$ observations in each replication,
using $100n$ replications for $M=4$ and $800n$ for $M=1$.
Each reported value is the unweighted average of the realized regrets.
The procedures use common observations at a given $n$;
random streams are independent across sample sizes.
Monte Carlo standard errors use the sample standard deviation divided by the square root of the replication count;
ratio errors account for pairing.
Fitted slopes regress log estimated regret on $\log n$ over the final four sample sizes.
The reported standard errors measure Monte Carlo uncertainty in the slopes fitted over these four sample sizes.

\subsection{One feasible set, three local rates}
\label{subsec:numerical_gaussian}

Consider the feasible set in Example~\ref{ex:fast_center_dependent_rates},
\begin{equation}
    \Omega=B_{4/3}^2=\{\bm w\in\mathbb R^2:\|\bm w\|_{4/3}\le1\},
    \qquad f_*(\bm z)=\|\bm z\|_4,
    \qquad w_{*,j}(\bm z)=\frac{z_j^3}{\|\bm z\|_4^3}\quad(\bm z\ne\bm0),
    \label{eq:4/3_ball_optimal_policy}
\end{equation}
with $\bm w_*(\bm0)=\bm0$. We fix $\Sigma=I_2$ and sample $\bm x_i\overset{\mathrm{i.i.d.}}{\sim}\mathcal N(\bm\mu_{n,\bm t},I_2)$, where $\bm\mu_{n,\bm t}=\bm\mu+\bm t/\sqrt n$. Table~\ref{tab:numerical_gaussian_rates} specifies the centers and directions, all with $\|\bm t\|_2=1$. The quantity estimated is
\begin{equation*}
    \mathbb E_{\bm\mu_{n,\bm t},I_2}R_{\bm\mu_{n,\bm t}}\left(\bm w_*(\bar{\bm x})\right).
\end{equation*}
Example~\ref{ex:fast_center_dependent_rates} gives linear, quadratic, and quartic local support expansions at these centers. Theorem~\ref{thm:fast_general_local_rate} therefore yields the respective minimax orders $n^{-1/2}$, $n^{-1}$, and $n^{-2}$, all attained in rate by SAA. We examine these predictions over $n=2^5,2^6,\ldots,2^{12}$.

\begin{table}[htbp]
    \centering
    \caption{Predicted and fitted slopes for SAA expected regret along Gaussian local alternatives. Fitted slopes use $n=2^9,2^{10},2^{11},2^{12}$; parentheses give their Monte Carlo standard errors. }
    \label{tab:numerical_gaussian_rates}
    \renewcommand{\arraystretch}{1.15}
    \begin{tabular}{@{}llll@{}}
        \toprule
        Center $\bm\mu$ & Direction $\bm t$ & Predicted slope & Fitted slope\\
        \midrule
        $\bm0$ & $\bm e_2$ & $-1/2$ & $-0.5001\;(0.0024)$\\
        $2^{-1/4}(1,1)^\top$ & $2^{-1/2}(1,-1)^\top$ & $-1$ & $-0.9857\;(0.0028)$\\
        $\bm e_1$ & $\bm e_2$ & $-2$ & $-2.0107\;(0.0068)$\\
        \bottomrule
    \end{tabular}
\end{table}

\begin{figure}[htbp]
    \centering
    \includegraphics[width=\textwidth]{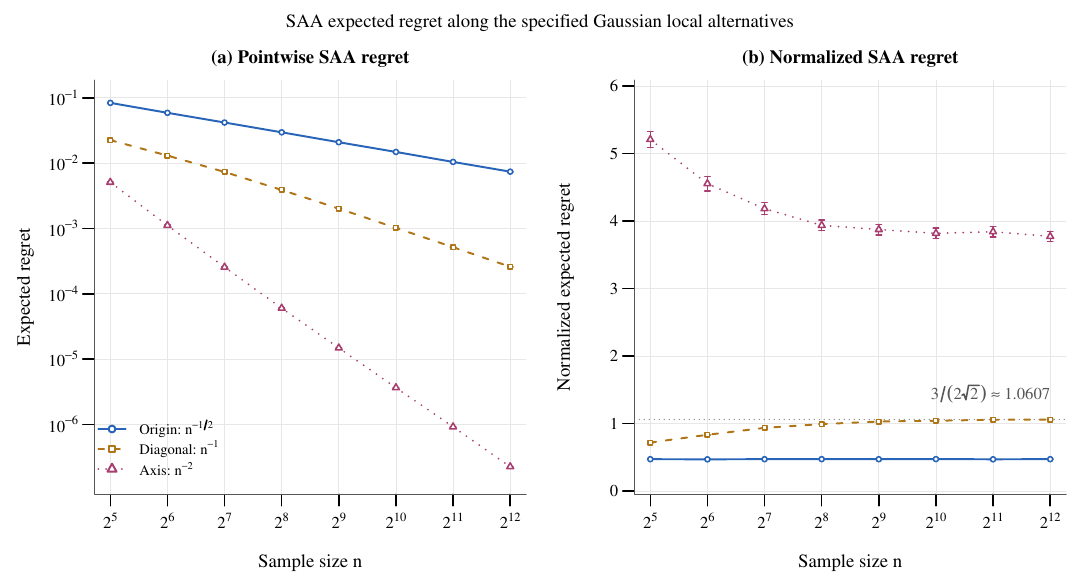}
    \caption{SAA expected regret along the Gaussian local alternatives in Table~\ref{tab:numerical_gaussian_rates}. Panel (a) shows the three rates. Panel (b) normalizes by $\sqrt n$, $n$, and $n^2$ at the origin $\bm\mu=\bm{0}$, diagonal $\bm\mu=2^{-1/4}(1,1)^\top$, and axis $\bm{\mu}=\bm{e}_1$, respectively. The reference level $3/(2\sqrt2)$ is the sharp constant for diagonal center $\bm\mu=2^{-1/4}(1,1)^\top$.}
    \label{fig:numerical_geometry_rates}
\end{figure}

Figure~\ref{fig:numerical_geometry_rates} shows the three distinct rates. The fitted slopes in Table~\ref{tab:numerical_gaussian_rates} are close to $-1/2$, $-1$, and $-2$, while the normalized curves become nearly horizontal at the larger sample sizes. Thus changing the center, with the feasible set and covariance held fixed, changes the statistical difficulty of the decision problem.

For the diagonal center $\bm\mu=2^{-1/4}(1,1)^\top$, the quadratic constant is also sharp in the local asymptotic minimax sense of Theorem~\ref{thm:fast_local_minimax}. The Hessian calculated in Section~\ref{app:l4/3-value-decision} of the Online Supplement gives
\begin{align*}
        H_{\bm\mu}&=\frac{3}{2\sqrt2}\begin{pmatrix}1&-1\\-1&1\end{pmatrix},\qquad
        \lim_{n\to\infty}n\,\mathbb E_{\bm\mu_{n,\bm t},I_2}R_{\bm\mu_{n,\bm t}}\left(\bm w_*(\bar{\bm x})\right)=\frac12\operatorname{tr}(H_{\bm\mu})=\frac3{2\sqrt2}\approx1.0607.
    \end{align*}
This equals the minimax constant. At $n=2^{12}$, the normalized Monte Carlo estimate is approximately $1.0604$, with standard error $4.7083\times10^{-3}$. The other two centers illustrate the rate conclusion of Theorem~\ref{thm:fast_general_local_rate}.

\subsection{Local geometry and the effect of Laplace tails}
\label{subsec:numerical_tail}
\label{subsec:numerical_geometry_threshold}

We next examine how local geometry changes the effect of rare observations on expected regret. Consider the family in Example~\ref{ex:div_matching_tail_geometry} with Laplace tails ($\zeta=1$) and feasible sets $\Omega=B_{\alpha/(\alpha-1)}^2$, where $\alpha\in\{3/2,2,4\}$. 

Under $P$, the coordinates are independent, $P\{x_1=1/2\}=P\{x_1=3/2\}=1/2$, and $x_2\sim\operatorname{Laplace}(0,1)$. Thus $\bm\mu_P=\bm e_1$ and $\Sigma=\operatorname{diag}(1/4,2)$. For a fixed $M>0$, take the perturbation constructed in Section~\ref{app:tail-geometry-construction} of the Online Supplement,
\begin{equation*}
    Q_n=(1-\varepsilon_n)P+\varepsilon_nP(\,\cdot\mid B_n),
\end{equation*}
where $\varepsilon_n=M^2/(4n^2)$, and $B_n=\{n\le|x_2|\le n+1\}$. Symmetry preserves the mean $\bm\mu_{Q_n}=\bm e_1$, and $Q_n\in\mathcal Q_{n,M}^{\mathrm{KL}}(P)$ for large $n$.

We compare the expected regrets of ordinary SAA, $\bm w_*(\bar{\bm x})$, and truncated SAA, $\bm w_*(\widetilde{\bm\mu}_n)$, based on $n$ independent observations from $Q_n$. Here $\widetilde{\bm\mu}_n$ is the coordinatewise truncated mean defined in \eqref{eq:div_truncated_mean}, with $\tau_n=\log^2(n+1)$. Theorem~\ref{thm:div_kl_minimax} gives the local minimax order $n^{-\alpha/2}$, attained by truncated SAA. For ordinary SAA, the additional contribution of rare observations gives maximal expected regret of order $n^{-\alpha/2}+n^{-1}$, as established in Example~\ref{ex:div_matching_tail_geometry}. The comparison therefore changes at $\alpha=2$: rare observations contribute a lower-order term when $\alpha=3/2$, change the leading constant when $\alpha=2$, and determine the regret rate when $\alpha=4$.

At $\alpha=2$, Corollary~\ref{cor:div_kl_smooth} gives the truncated-SAA constant $\operatorname{tr}(\Sigma H_{\bm e_1})/2=1$. Along the present sequence $Q_n$, the construction in Section~\ref{app:tail-geometry-construction} gives
\begin{align*}
        n\mathbb E_{Q_n^{\otimes n}}R_{Q_n}\left(\bm w_*(\widetilde{\bm\mu}_n)\right)\to1,\qquad
        n\mathbb E_{Q_n^{\otimes n}}R_{Q_n}\left(\bm w_*(\bar{\bm x})\right)\to1+\frac{M^2}{4}(1-2^{-1/2}).
\end{align*}
For $M=4$, the ordinary-to-truncated expected-regret ratio tends to $1$ when $\alpha=3/2$, tends to $1+4(1-2^{-1/2})\approx2.1716$ when $\alpha=2$, and grows at order $n$ when $\alpha=4$. These three behaviors isolate equality of leading regret, a constant loss, and a rate loss, respectively.

\begin{figure}[htbp]
    \centering
    \includegraphics[width=\textwidth]{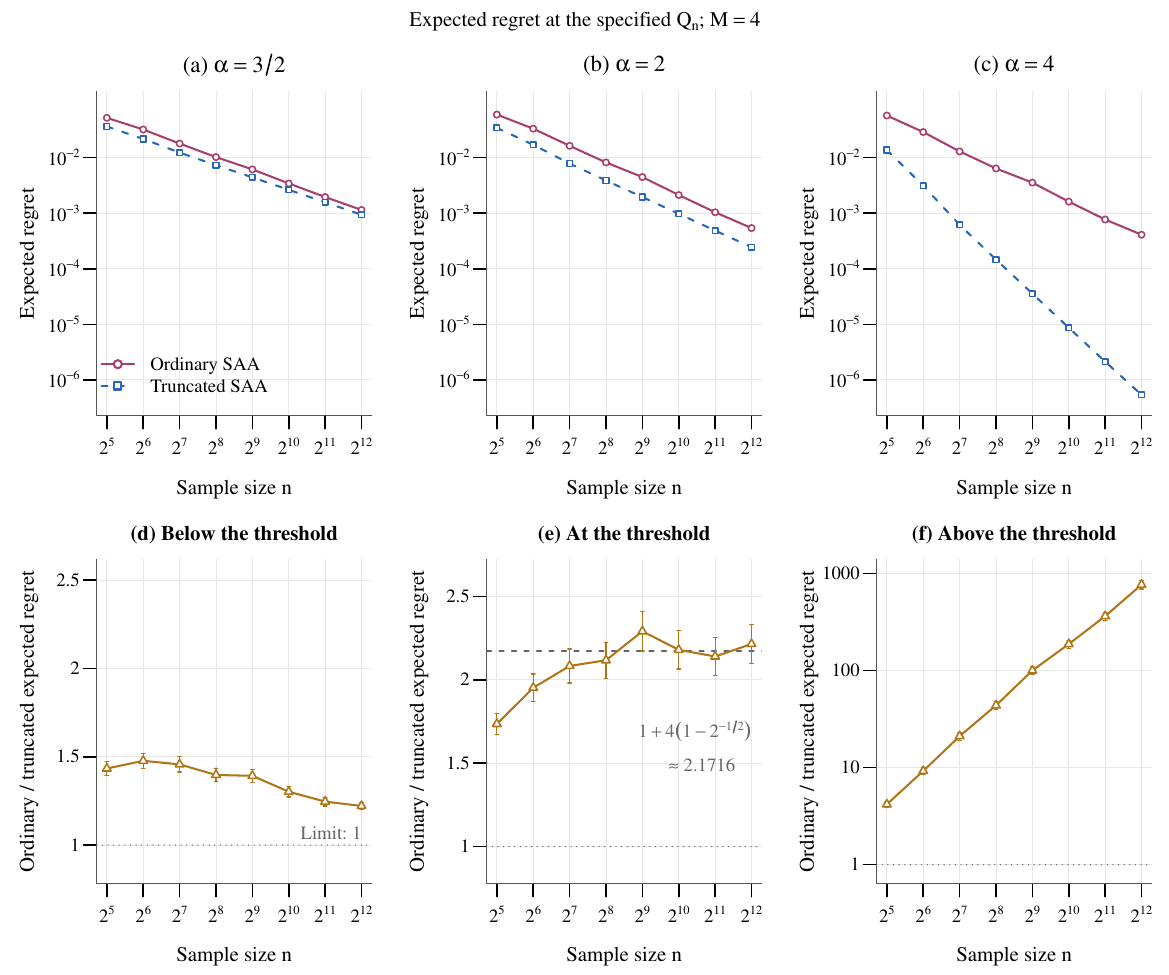}
    \caption{Ordinary and truncated SAA under $Q_n$ with $M=4$. The upper panels show expected regret on a common scale; the lower panels show the ordinary-to-truncated ratio. The reference level $1$ marks equal expected regret and is the limiting ratio for $\alpha=3/2$. For $\alpha=2$, the additional reference is $1+4(1-2^{-1/2})\approx2.1716$. For $\alpha=4$, the ratio grows at order $n$. The ratio scale is linear for $\alpha=3/2,2$ and logarithmic for $\alpha=4$.}
    \label{fig:numerical_geometry_threshold}
\end{figure}

We first set $M=4$ and use the eight sample sizes $n=2^5,2^6,\ldots,2^{12}$. Figure~\ref{fig:numerical_geometry_threshold} displays the expected regrets and their ratios across the three regimes. The fitted slopes in Table~\ref{tab:numerical_geometry_threshold} reflect the transition from a common rate below and at the threshold to different rates above it. The table also reports a quartic comparison at $M=1$ on the same eight-point grid.

\begin{table}[htbp]
    \centering
    \caption{Predicted and fitted slopes for expected regret under $Q_n$. Fits use $n=2^9,2^{10},2^{11},2^{12}$ for both radii; the figures also include $n=2^5,2^6,2^7,2^8$. Parentheses give Monte Carlo standard errors. The pointwise rates match the maximal-regret orders in Example~\ref{ex:div_matching_tail_geometry}.}
    \label{tab:numerical_geometry_threshold}
    \label{tab:numerical_saa_rates}
    \small
    \renewcommand{\arraystretch}{1.15}
    \begin{tabular}{@{}llllll@{}}
        \toprule
        & & \multicolumn{2}{l}{Predicted slope} & \multicolumn{2}{l}{Fitted slope}\\
        \cmidrule(lr){3-4}\cmidrule(l){5-6}
        $\alpha$ & $M$ & Ordinary SAA & Truncated SAA & Ordinary & Truncated\\
        \midrule
        $3/2$ & $4$ & $-3/4$ & $-3/4$ & $-0.8086\;(0.0075)$ & $-0.7457\;(0.0022)$\\
        $2$ & $4$ & $-1$ & $-1$ & $-1.0185\;(0.0174)$ & $-1.0012\;(0.0030)$\\
        $4$ & $4$ & $-1$ & $-2$ & $-1.0421\;(0.0311)$ & $-2.0187\;(0.0071)$\\
        $4$ & $1$ & $-1$ & $-2$ & $-1.1344\;(0.0411)$ & $-2.0034\;(0.0024)$\\
        \bottomrule
    \end{tabular}
\end{table}

At $M=4,n=2^{12}$, the estimated expected-regret ratio for $\alpha=3/2$ is approximately $1.2215$. At larger sample sizes, the estimated ratio declines toward its theoretical limit of $1$. For $\alpha=2$, the ordinary and truncated estimates multiplied by $n$ are approximately $2.2144$ and $0.9999$, with Monte Carlo standard errors $5.9606\times10^{-2}$ and $2.2062\times10^{-3}$. Their separation persists even though both fitted slopes are close to $-1$. For $\alpha=4$, the different slopes display the rate loss of ordinary SAA. For $\alpha=3/2$, the ordinary-SAA slope $-0.8086$ remains below its limit $-3/4$: the rare-event term has order $n^{-1}$ and need not be negligible on this grid. 

Table~\ref{tab:numerical_saa_results} records the quartic endpoint estimates for both radii. At $n=2^{12}$, the estimated expected-regret ratio is approximately $47.3352$ for $M=1$ and $764.6639$ for $M=4$, with Monte Carlo standard errors $3.2736$ and $37.6256$, respectively.

\begin{table}[htbp]
    \centering
    \caption{Quartic expected-regret estimates from direct simulation at $Q_n$. Parentheses give Monte Carlo standard errors; the last column is the ratio of the two estimates. }
    \label{tab:numerical_saa_results}
    \small
    \renewcommand{\arraystretch}{1.2}
    \begin{tabular}{@{}lllll@{}}
        \toprule
        $M$ & $n$ & Ordinary SAA & Truncated SAA & Regret ratio\\
        \midrule
        $1$ & $2^{12}$ & $2.5484\times10^{-5}$ & $5.3838\times10^{-7}$ & $47.3352$\\
            &        & $(1.7621\times10^{-6})$ & $(9.7394\times10^{-10})$ & $(3.2736)$\\
        $4$ & $2^{12}$ & $4.1079\times10^{-4}$ & $5.3721\times10^{-7}$ & $764.6639$\\
            &        & $(2.0132\times10^{-5})$ & $(2.7183\times10^{-9})$ & $(37.6256)$\\
        \bottomrule
    \end{tabular}
\end{table}

The results illustrate how distributional tails and local support geometry jointly affect the performance of ordinary SAA. Under the same Laplace perturbation, ordinary SAA matches the leading asymptotic expected regret of truncated SAA when $\alpha=3/2$, incurs a constant loss when $\alpha=2$, and becomes suboptimal in rate when $\alpha=4$.

\section{Discussion}
\label{sec:discussion}

The results yield a unified but qualified view of SAA. In the nominal Gaussian mean experiment, the support-function remainder completely characterizes the instance-specific statistical difficulty, and SAA automatically attains the resulting rate. At smooth points it is sharp in the standard local asymptotic minimax limit. At nonsmooth or degenerate points it adapts to first-order, higher-order, and non-power geometry without knowing the relevant exponent.

Distributional robustness changes the implementation question more than the geometric benchmark. A KL ball of radius $1/n$ generates the same $n^{-1/2}$ local information scale, transformed by the covariance ellipsoid of the reference distribution. Truncated SAA attains the corresponding minimax rate. Ordinary SAA can fail because weak convergence does not control expected regret uniformly over KL-admissible alternatives. The failure is therefore not a contradiction to Gaussian efficiency: it is a lack of uniform integrability caused by rare tail perturbations.

The distinction suggests three meanings of the price of robustness. Robustness can be asymptotically free in rate; it can change a fixed-radius local constant by changing the shape of the local parameter set; and it can require a tail-stable estimator even when the information bound itself is unchanged. The last phenomenon is most visible at instances with degenerate curvature, where the optimization problem permits a rate faster than $n^{-1}$. These statements compare procedures that do not use $P$ with a benchmark that may depend on the reference law $P$ and the radius $M$; ``robustness is free'' is thus a statement about the attainable rate and smooth constant relative to a $P$-informed benchmark.

The present theory is deliberately local and fixed-dimensional. It assumes a compact convex feasible set, a linear objective whose distributional dependence enters through the mean, pointwise local asymptotics around a fixed reference instance, and one-sided neighborhoods $D_{\mathrm{KL}}(Q\|P)$. The KL lower bound requires a finite second moment and nonsingular covariance, whereas the general matching upper bound uses a local exponential-moment condition. The reference law $P$ is fixed when defining the robust benchmark. These restrictions isolate the geometry-tail mechanism but leave important questions of growing dimension, nonlinear objectives, data-dependent reference models, and other notions of distributional proximity open.

Several extensions are natural. Other $\phi$-divergence neighborhoods should lead to related local ellipsoidal experiments with divergence-specific scaling. Wasserstein neighborhoods would interact differently with the geometry because they constrain transport rather than likelihood ratios. It would also be useful to develop adaptive, coordinate-free truncation, allow $d$ to grow with $n$, and extend the analysis to nonlinear stochastic programs. Finally, combining the present geometry with parametric ETO/IEO models would separate gains from structural dimension reduction from gains due to tail stabilization.

\bibliographystyle{plainnat}
\bibliography{DDOpt_ref}

\newpage
\setcounter{page}{1}

\begin{center}
{\Large Supplement to\\
	\emph{Geometry, Local Minimax Risk, and the Limits of Sample Average Approximation in Linear Stochastic Optimization}}
\end{center}

\medskip
\noindent
This supplement contains all proofs, intermediate technical results, and detailed example calculations for the main paper. Section and equation references without an appendix letter refer to the main paper.

\paragraph{Proof roadmap.}
The Gaussian Bayes identity reduces minimax regret to a support-function Jensen gap (Section~\ref{app:proof-general-rate}).
For the KL limit, bounded-score submodels yield the lower bound under a finite second moment.
Under the exponential-moment condition, entropy bounds confine the local means to the covariance ellipsoid, while clipping gives uniform moment bounds and a Gaussian approximation (Lemma~\ref{lem:div_kl_truncation}).
The upper bound transfers a nearly optimal Gaussian rule after restricting its actions so that its loss is uniformly of order $a(n^{-1/2})$ (Section~\ref{app:proof-kl-minimax}).
Finally, a polynomial envelope and the clipped moment bounds give uniform integrability of truncated-SAA regret; ordinary SAA shares the weak limit but may fail this expectation bound (Section~\ref{app:proof-plugin-comparison}).

\appendix

\section{Proofs for the Global and Local Benchmark Bounds}
\label{app:section2}
\subsection{Proof of Proposition~\ref{prop:global_gaussian_width}}
\label{app:proof-global-width}
\begin{proof}
	Let $\bm w_\mu:=\bm w_*(\bm\mu)$ and
	$\bm w_{\bar x}:=\bm w_*(\bar{\bm x})$. Optimality of
	$\bm w_{\bar x}$ at $\bar{\bm x}$ gives
	$\bar{\bm x}^\top(\bm w_\mu-\bm w_{\bar x})\le0$. Hence
	\begin{align*}
		R_{\bm\mu}(\bm w_{\bar x})
		&=\bm\mu^\top(\bm w_\mu-\bm w_{\bar x})\\
		&\le(\bar{\bm x}-\bm\mu)^\top
		(\bm w_{\bar x}-\bm w_\mu)\\
		&\le f_*(\bar{\bm x}-\bm\mu)
		+f_*(\bm\mu-\bar{\bm x}).
	\end{align*}
	Since $\sqrt n(\bar{\bm x}-\bm\mu)\sim
	\mathcal N(\bm0,\Sigma)$ and this distribution is symmetric,
	\begin{equation*}
		\mathbb E_{\bm\mu,\Sigma}
		R_{\bm\mu}(\bm w_*(\bar{\bm x}))
		\le \frac2{\sqrt n}
		\mathbb E f_*(\Sigma^{1/2}\bm G)
		=2\frac{\mathcal G(\Sigma^{1/2}\Omega)}{\sqrt n}.
	\end{equation*}
	This proves the upper bound in \eqref{eq:global_width_sandwich}.
	
	For the lower bound, place the Gaussian prior
	\begin{equation*}
		\bm\mu\sim\mathcal N\left(\bm0,\frac\Sigma n\right).
	\end{equation*}
	Under this prior, the posterior mean based on the full sample is
	$\bm m=\mathbb E[\bm\mu\mid\bm x_1,\ldots,\bm x_n]
	=\bar{\bm x}/2$. Conditional on the data, a Bayes action maximizes
	$\bm w^\top\bm m$ over $\Omega$. Therefore the minimum Bayes regret is
	\begin{align*}
		\inf_{\widehat{\bm w}_n}
		\mathbb E R_{\bm\mu}(\widehat{\bm w}_n)
		&=\mathbb E f_*(\bm\mu)-\mathbb E f_*(\bm m).
	\end{align*}
	Now
	\[
	\bm\mu\sim\mathcal N\left(\bm0,\frac\Sigma n\right),
	\qquad
	\bm m\sim\mathcal N\left(\bm0,\frac\Sigma{2n}\right),
	\]
	so positive homogeneity of $f_*$ yields
	\begin{equation*}
		\inf_{\widehat{\bm w}_n}
		\mathbb E R_{\bm\mu}(\widehat{\bm w}_n)
		=\left(1-\frac1{\sqrt2}\right)
		\frac{\mathcal G(\Sigma^{1/2}\Omega)}{\sqrt n}.
	\end{equation*}
	The minimax regret dominates every Bayes regret, proving the first
	inequality in \eqref{eq:global_width_sandwich}.
	
	We also record a comparison of Gaussian widths, used in Proposition~\ref{prop:cov_geometric_upper}. For independent standard
	Gaussian vectors $\bm G,\bm G'$, set
	\[
	\bm Y=\Sigma^{1/2}\bm G
	+(\overline\lambda I_d-\Sigma)^{1/2}\bm G'.
	\]
	Then $\bm Y\sim\mathcal N(\bm0,\overline\lambda I_d)$ and
	$\mathbb E[\bm Y\mid\bm G]=\Sigma^{1/2}\bm G$. Convexity of $f_*$ and
	conditional Jensen's inequality give
	\begin{equation}
		\mathcal G(\Sigma^{1/2}\Omega)
		=\mathbb E f_*(\Sigma^{1/2}\bm G)
		\le\mathbb E f_*(\bm Y)
		=\sqrt{\overline\lambda}\,\mathcal G(\Omega).
		\label{eq:width_covariance_monotonicity}
	\end{equation}
	Equality is attained at $\Sigma=\overline\lambda I_d$.
\end{proof}
\subsection{Proof of Proposition~\ref{prop:cov_geometric_upper}}
\label{app:proof-local-fast}
\begin{proof}
	Optimality of $\bm w_*(\bar{\bm x})$ at $\bar{\bm x}$ gives
	\begin{align*}
		R_{\bm\mu}(\bm w_*(\bar{\bm x}))
		&=\bm\mu^\top
		\bigl(\bm w_*(\bm\mu)-\bm w_*(\bar{\bm x})\bigr)\\
		&\le(\bar{\bm x}-\bm\mu)^\top
		\bigl(\bm w_*(\bar{\bm x})-\bm w_*(\bm\mu)\bigr)\\
		&\le\|\bar{\bm x}-\bm\mu\|_K
		\|\bm w_*(\bar{\bm x})-\bm w_*(\bm\mu)\|_\Omega.
	\end{align*}
	On the other hand, \eqref{eq:fast_geometric_condition} at $\bm\mu$ gives
	\begin{equation*}
		R_{\bm\mu}(\bm w_*(\bar{\bm x}))
		\ge c_\Omega \dist(\bm\mu,A)
		\|\bm w_*(\bar{\bm x})-\bm w_*(\bm\mu)\|_\Omega^2.
	\end{equation*}
	Combining the two inequalities yields
	\begin{equation*}
		R_{\bm\mu}(\bm w_*(\bar{\bm x}))
		\le\frac{\|\bar{\bm x}-\bm\mu\|_K^2}
		{c_\Omega \dist(\bm\mu,A)}.
	\end{equation*}
	Since $\bar{\bm x}-\bm\mu\overset d=
	\Sigma^{1/2}\bm G/\sqrt n$, taking expectations proves
	\eqref{eq:cov_geometric_moment_bound}.
	
	Under \eqref{eq:polar_norm_convention},
	$\|\bm y\|_K=f_*(\bm y)$ and hence
	\begin{equation*}
		\mathbb E\|\Sigma^{1/2}\bm G\|_K
		=\mathcal G(\Sigma^{1/2}\Omega).
	\end{equation*}
	Jensen's inequality gives the first inequality in
	\eqref{eq:cov_width_moment_comparison}. Moreover,
	$\bm g\mapsto f_*(\Sigma^{1/2}\bm g)$ is Lipschitz with squared
	Lipschitz constant at most
	$\sup_{\bm w\in\Omega}\bm w^\top\Sigma\bm w$. Gaussian
	Poincare inequality therefore gives
	\begin{equation*}
		\operatorname{Var}\bigl(f_*(\Sigma^{1/2}\bm G)\bigr)
		\le\sup_{\bm w\in\Omega}\bm w^\top\Sigma\bm w.
	\end{equation*}
	Because $\Omega=-\Omega$,
	\[
	f_*(\Sigma^{1/2}\bm G)
	\ge |\bm w^\top\Sigma^{1/2}\bm G|
	\]
	for every $\bm w\in\Omega$. Therefore
	\begin{equation*}
		\mathcal G(\Sigma^{1/2}\Omega)
		\ge\sqrt{\frac2\pi}
		\sup_{\bm w\in\Omega}\sqrt{\bm w^\top\Sigma\bm w}.
	\end{equation*}
	The variance decomposition now gives
	\begin{align*}
		\mathbb E\|\Sigma^{1/2}\bm G\|_K^2
		&=\mathcal G(\Sigma^{1/2}\Omega)^2
		+\operatorname{Var}\bigl(f_*(\Sigma^{1/2}\bm G)\bigr)\\
		&\le\left(1+\frac\pi2\right)
		\mathcal G(\Sigma^{1/2}\Omega)^2,
	\end{align*}
	which proves \eqref{eq:cov_width_moment_comparison}. Finally,
	\eqref{eq:width_covariance_monotonicity} implies the second inequality in
	\eqref{eq:fast_ddopt_bound}.
\end{proof}
\subsection{Strong convexity implies quadratic growth}
Suppose $\Omega=-\Omega$, $\bm0\in\operatorname{int}\Omega$, and $\Omega$ is $c_\Omega$-strongly convex in its gauge in the midpoint sense. For $\bm z\ne\bm0$, choose $\bm v_{\bm z}\in\Omega$ satisfying $\bm z^\top\bm v_{\bm z}=\|\bm z\|_K$, and apply strong convexity with $\bm u=\bm w_*(\bm z)$ and an arbitrary $\bm v=\bm w\in\Omega$. Taking $\bm s=(c_\Omega/2)\|\bm w-\bm w_*(\bm z)\|_\Omega^2\bm v_{\bm z}$ gives $(\bm w_*(\bm z)+\bm w)/2+\bm s\in\Omega$. Optimality of $\bm w_*(\bm z)$ gives
\[
\bm z^\top(\bm w_*(\bm z)-\bm w)
\ge c_\Omega\|\bm z\|_K
\|\bm w-\bm w_*(\bm z)\|_\Omega^2.
\]
In particular, every nonzero $\bm z$ has a unique optimizer, so $A=\{\bm0\}$ and $\dist(\bm z,A)=\|\bm z\|_K$. This proves Remark~\ref{rem:strong-convexity-sufficient}.

\section{Decision-Theoretic Tools and Proofs for the Gaussian Experiment}
\label{app:gaussian-decision}

\subsection{A compact minimax theorem}
\begin{lemma}[Minimax identity on compact local sets]
	\label{lem:cov_compact_minimax}
	Fix $n$, $M>0$, and $\bm\mu\in\mathbb R^d$. Let $\Theta$ be a nonempty compact subset of $B_2(M)\times\mathcal C$, and let $\mathcal W$ be a nonempty compact convex subset of $\Omega$. With samples $(\bm{x}_1,\ldots,\bm{x}_n)$ drawn from the joint distribution
	\begin{equation*}
		P_{\bm t,\Sigma}^{(n)}=\mathcal N\left(\bm\mu+\frac{\bm t}{\sqrt n},\Sigma\right)^{\otimes n},\qquad (\bm t,\Sigma)\in\Theta,
	\end{equation*}
	denote the loss as
	\begin{equation*}
		\ell_n(\bm t,\bm w)=f_*\left(\bm\mu+\frac{\bm t}{\sqrt n}\right)-\bm w^\top\left(\bm\mu+\frac{\bm t}{\sqrt n}\right),\qquad\bm w\in\mathcal W.
	\end{equation*}
	Then
	\begin{equation}
		\inf_{\substack{\hat{\bm w}=\hat{\bm w}(\bm{x}_1,\ldots,\bm{x}_n):(\mathbb R^d)^n\to\mathcal W\\\hat{\bm w}\ \mathrm{measurable}}}\sup_{(\bm t,\Sigma)\in\Theta}\mathbb E_{\bm\mu_{n,\bm t},\Sigma}\ell_n(\bm t,\hat{\bm w})=\sup_{\Pi\in\mathcal P(\Theta)}\inf_{\substack{\hat{\bm w}=\hat{\bm w}(\bm{x}_1,\ldots,\bm{x}_n):(\mathbb R^d)^n\to\mathcal W\\\hat{\bm w}\ \mathrm{measurable}}}\mathbb E_\Pi\ell_n(\bm T,\hat{\bm w}),
		\label{eq:cov_compact_minimax}
	\end{equation}
	where $\mathcal P(\Theta)$ denotes the Borel probability measures on $\Theta$, and $\mathbb E_\Pi$ uses $(\bm T,\Sigma)\sim\Pi$ followed by the sample law $P_{\bm T,\Sigma}^{(n)}$.
	
	Moreover, if for every $\Pi\in\mathcal P(\Theta)$ there exists a Bayes rule for decision space $\Omega$ taking values in $\mathcal W$, then
	\begin{equation*}
		\inf_{\substack{\hat{\bm w}=\hat{\bm w}(\bm{x}_1,\ldots,\bm{x}_n):(\mathbb R^d)^n\to\mathcal W\\\hat{\bm w}\ \mathrm{measurable}}}\sup_{(\bm t,\Sigma)\in\Theta}\mathbb E_{\bm\mu_{n,\bm t},\Sigma}\ell_n(\bm t,\hat{\bm w})
		=\inf_{\substack{\hat{\bm w}=\hat{\bm w}(\bm{x}_1,\ldots,\bm{x}_n):(\mathbb R^d)^n\to\Omega\\\hat{\bm w}\ \mathrm{measurable}}}\sup_{(\bm t,\Sigma)\in\Theta}\mathbb E_{\bm\mu_{n,\bm t},\Sigma}\ell_n(\bm t,\hat{\bm w}).
	\end{equation*}
\end{lemma}
\begin{proof}
	Let $\mathcal D_n^{\mathrm r}(\mathcal W)$ be the class of all Markov kernels from $(\mathbb R^d)^n$ to $\mathcal W$, and define
	\begin{equation*}
		\rho_\delta(\bm t,\Sigma)=\int_{(\mathbb R^d)^n}\int_{\mathcal W}\ell_n(\bm t,\bm w)\,\delta(\dd\bm w\mid\bm x_1,\ldots,\bm x_n)\,P_{\bm t,\Sigma}^{(n)}(\dd\bm x_1\cdots\dd\bm x_n).
		\label{eq:fast_randomized_risk}
	\end{equation*}
	We apply \citet[Theorem~3.78]{LieseMiescke2008} with parameter space $\Theta$, decision space $\mathcal W$, and decision class $\mathcal D_n^{\mathrm r}(\mathcal W)$. Both spaces are compact metric spaces with their Borel $\sigma$-algebras. Since $\Omega$ is compact, $f_*$ is $R_\Omega$-Lipschitz continuous. Thus $((\bm t,\Sigma),\bm w)\mapsto\ell_n(\bm t,\bm w)$ is jointly continuous, and
	\begin{equation*}
		0\le\ell_n(\bm t,\bm w)\le2R_\Omega(\|\bm\mu\|_2+M/\sqrt n),\qquad ((\bm t,\Sigma),\bm w)\in\Theta\times\mathcal W.
	\end{equation*}
	The sample space is Euclidean space $(\mathbb{R}^d)^n$ with its Borel $\sigma$-algebra, and every $P_{\bm t,\Sigma}^{(n)}$ has a strictly positive Lebesgue density. For $(\bm t,\Sigma),(\bm s,\Gamma)\in B_2(M)\times\mathcal C$, the Gaussian KL formula gives
	\begin{equation*}
		D_{\mathrm{KL}}\left(P_{\bm t,\Sigma}^{(n)}\,\middle\|\,P_{\bm s,\Gamma}^{(n)}\right)
		=\frac n2\left[\operatorname{tr}(\Gamma^{-1}\Sigma)-d-\log\det(\Gamma^{-1}\Sigma)\right]
		+\frac12(\bm t-\bm s)^\top\Gamma^{-1}(\bm t-\bm s),
	\end{equation*}
	and therefore, by Pinsker's inequality,
	\begin{equation*}
		\left\|P_{\bm t,\Sigma}^{(n)}-P_{\bm s,\Gamma}^{(n)}\right\|_{\mathrm{TV}}
		\le\sqrt{\frac12D_{\mathrm{KL}}\left(P_{\bm t,\Sigma}^{(n)}\,\middle\|\,P_{\bm s,\Gamma}^{(n)}\right)}\to0
		\quad\text{as }(\bm t,\Sigma)\to(\bm s,\Gamma).
	\end{equation*}
	Here $n$ is fixed, and matrix inversion and the log determinant are continuous on $\mathcal C$, since $\underline\lambda>0$. This verifies condition (A7) of \citet{LieseMiescke2008}, also after restriction to $\Theta$. In particular, $(\bm t,\Sigma)\mapsto P_{\bm t,\Sigma}^{(n)}(B)$ is Borel measurable for every Borel sample event $B$. For fixed $\Sigma$, the same calculation reduces to
	\begin{equation*}
		\left\|P_{\bm t,\Sigma}^{(n)}-P_{\bm s,\Sigma}^{(n)}\right\|_{\mathrm{TV}}\le\frac12\|\Sigma^{-1/2}(\bm t-\bm s)\|_2\le\frac{\|\bm t-\bm s\|_2}{2\sqrt{\underline\lambda}}.
		\label{eq:fast_gaussian_tv_continuity}
	\end{equation*}
	
	The class $\mathcal D_n^{\mathrm r}(\mathcal W)$ is nonempty, since it contains every constant rule in $\mathcal W$, and convex: $a\delta_1+(1-a)\delta_2\in\mathcal D_n^{\mathrm r}(\mathcal W)$ for $a\in[0,1]$ and $\delta_1,\delta_2\in\mathcal D_n^{\mathrm r}(\mathcal W)$. In the weak convergence of \citet[Definition~3.19]{LieseMiescke2008}, $\delta_k\Rightarrow\delta$ means
	\begin{equation*}
		\mathbb E_{\bm\mu_{n,\bm t},\Sigma}\int_{\mathcal W}h(\bm w)\,\delta_k(\dd\bm w\mid\bm x_1,\ldots,\bm x_n)
		\to\mathbb E_{\bm\mu_{n,\bm t},\Sigma}\int_{\mathcal W}h(\bm w)\,\delta(\dd\bm w\mid\bm x_1,\ldots,\bm x_n)
	\end{equation*}
	for every $(\bm t,\Sigma)\in\Theta$ and every bounded continuous $h:\mathcal W\to\mathbb R$, where $\delta$ is a Markov kernel into the decision space $\mathcal W$. The closure condition in \citet[Theorem~3.78]{LieseMiescke2008} is relative to the ambient space of Markov kernels into $\mathcal W$. Our decision class $\mathcal D_n^{\mathrm r}(\mathcal W)$ is that entire ambient class, so it is automatically weakly closed: whenever a convergent net or sequence has a Markov-kernel limit in the topology of \citet[Definition~3.19]{LieseMiescke2008}, that limit again belongs to $\mathcal D_n^{\mathrm r}(\mathcal W)$. All assumptions of the statistical minimax theorem are therefore satisfied, and it gives
	\begin{equation}
		\inf_{\delta\in\mathcal D_n^{\mathrm r}(\mathcal W)}\sup_{(\bm t,\Sigma)\in\Theta}\rho_\delta(\bm t,\Sigma)
		=\sup_{\Pi\in\mathcal P(\Theta)}\inf_{\delta\in\mathcal D_n^{\mathrm r}(\mathcal W)}\int_\Theta\rho_\delta(\bm t,\Sigma)\,\Pi(\dd\bm t,\dd\Sigma).
		\label{eq:minimax_regret_equals_sup_bayes}
	\end{equation}
	
	For $\delta\in\mathcal D_n^{\mathrm r}(\mathcal W)$, define its conditional barycenter
	\begin{equation*}
		\overline{\bm w}_\delta(\bm x_1,\ldots,\bm x_n)=\int_{\mathcal W}\bm w\,\delta(\dd\bm w\mid\bm x_1,\ldots,\bm x_n)\in\mathcal W.
		\label{eq:fast_policy_barycenter}
	\end{equation*}
	This is measurable, and compactness and convexity of $\mathcal W$ give feasibility. Since $\ell_n$ is affine in $\bm w$,
	\begin{equation*}
		\ell_n(\bm t,\overline{\bm w}_\delta)=\int_{\mathcal W}\ell_n(\bm t,\bm w)\,\delta(\dd\bm w\mid\bm x_1,\ldots,\bm x_n).
		\label{eq:fast_derandomization}
	\end{equation*}
	Taking expectation gives, simultaneously for every $(\bm t,\Sigma)\in\Theta$,
	\begin{equation*}
		\rho_\delta(\bm t,\Sigma)=\mathbb E_{\bm\mu_{n,\bm t},\Sigma}\ell_n(\bm t,\overline{\bm w}_\delta).
	\end{equation*}
	Thus every randomized policy has exactly the same risk function as a deterministic policy, which gives 
	\begin{equation}
		\begin{aligned}
			&\inf_{\delta\in\mathcal D_n^{\mathrm r}(\mathcal W)}
			\sup_{(\bm t,\Sigma)\in\Theta}\rho_\delta(\bm t,\Sigma)\\
			&\quad=\inf_{\delta\in\mathcal D_n^{\mathrm r}(\mathcal W)}
			\sup_{(\bm t,\Sigma)\in\Theta}
			\mathbb E_{\bm\mu_{n,\bm t},\Sigma}\ell_n(\bm t,\overline{\bm w}_\delta)\\
			&\quad\ge
			\inf_{\substack{\hat{\bm w}:(\mathbb R^d)^n\to\mathcal W\\
					\hat{\bm w}\ \mathrm{measurable}}}
			\sup_{(\bm t,\Sigma)\in\Theta}
			\mathbb E_{\bm\mu_{n,\bm t},\Sigma}\ell_n(\bm t,\hat{\bm w}).
		\end{aligned}
		\label{eq:randomized_minimax_risk_greater_than_deterministic}
	\end{equation}
	and for any $\Pi\in\mathcal{P}(\Theta)$,
	\begin{equation}
		\begin{aligned}
			\inf_{\delta\in\mathcal D_n^{\mathrm r}(\mathcal W)}\int_\Theta\rho_\delta(\bm t,\Sigma)\,\Pi(\dd\bm t,\dd\Sigma)
			&=\inf_{\delta\in\mathcal D_n^{\mathrm r}(\mathcal W)}\int_\Theta\mathbb E_{\bm\mu_{n,\bm t},\Sigma}\ell_n(\bm t,\overline{\bm w}_\delta)\,\Pi(\dd\bm t,\dd\Sigma)
			\\&\ge\inf_{\substack{\hat{\bm w}:(\mathbb R^d)^n\to\mathcal W\\\hat{\bm w}\ \mathrm{measurable}}}\mathbb E_\Pi\ell_n(\bm T,\hat{\bm w}).
		\end{aligned}
		\label{eq:randomized_bayes_risk_greater_than_deterministic}
	\end{equation}
	Conversely, for every measurable $\hat{\bm w}:(\mathbb R^d)^n\to\mathcal W$, define
	\begin{equation*}
		\delta_{\bm w}(B\mid\bm x_1,\ldots,\bm x_n):=\mathbb I\{\hat{\bm w}(\bm x_1,\ldots,\bm x_n)\in B\},\qquad B\in\mathcal B(\mathcal W).
	\end{equation*}
	For every sample, $\delta_{\bm w}(\cdot\mid\bm x_1,\ldots,\bm x_n)$ is the Dirac probability measure at $\hat{\bm w}(\bm x_1,\ldots,\bm x_n)$; for every Borel $B$, its value is measurable in the sample because $\hat{\bm w}$ is measurable. Hence $\delta_{\bm w}\in\mathcal D_n^{\mathrm r}(\mathcal W)$, and
	\begin{equation*}
		\int_{\mathcal W}\ell_n(\bm t,\bm w)\,\delta_{\bm w}(\dd\bm w\mid\bm x_1,\ldots,\bm x_n)=\ell_n(\bm t,\hat{\bm w}(\bm x_1,\ldots,\bm x_n)),
	\end{equation*}
	which implies $\rho_{\delta_{\bm w}}(\bm t,\Sigma)=\mathbb E_{\bm\mu_{n,\bm t},\Sigma}\ell_n(\bm t,\hat{\bm w})$ for $(\bm t,\Sigma)\in\Theta$. Therefore,
	\begin{equation}
		\begin{aligned}
			&\inf_{\delta\in\mathcal D_n^{\mathrm r}(\mathcal W)}
			\sup_{(\bm t,\Sigma)\in\Theta}\rho_\delta(\bm t,\Sigma)\\
			&\quad\le
			\inf_{\delta_{\bm w}\in\mathcal D_n^{\mathrm r}(\mathcal W)}
			\sup_{(\bm t,\Sigma)\in\Theta}\rho_{\delta_{\bm w}}(\bm t,\Sigma)\\
			&\quad=
			\inf_{\substack{\hat{\bm w}:(\mathbb R^d)^n\to\mathcal W\\
					\hat{\bm w}\ \mathrm{measurable}}}
			\sup_{(\bm t,\Sigma)\in\Theta}
			\mathbb E_{\bm\mu_{n,\bm t},\Sigma}\ell_n(\bm t,\hat{\bm w}).
		\end{aligned}
		\label{eq:randomized_minimax_risk_smaller_than_deterministic}
	\end{equation}
	and for any $\Pi\in\mathcal{P}(\Theta)$, 
	\begin{equation}
		\begin{aligned}
			&\inf_{\delta\in\mathcal D_n^{\mathrm r}(\mathcal W)}
			\int_\Theta\rho_\delta(\bm t,\Sigma)\,
			\Pi(\dd\bm t,\dd\Sigma)\\
			&\quad\le
			\inf_{\delta_{\bm w}\in\mathcal D_n^{\mathrm r}(\mathcal W)}
			\int_\Theta\rho_{\delta_{\bm w}}(\bm t,\Sigma)\,
			\Pi(\dd\bm t,\dd\Sigma)\\
			&\quad=
			\inf_{\substack{\hat{\bm w}:(\mathbb R^d)^n\to\mathcal W\\
					\hat{\bm w}\ \mathrm{measurable}}}
			\mathbb E_\Pi\ell_n(\bm T,\hat{\bm w}).
		\end{aligned}
		\label{eq:randomized_bayes_risk_smaller_than_deterministic}
	\end{equation}
	Combining \eqref{eq:randomized_minimax_risk_greater_than_deterministic} and \eqref{eq:randomized_minimax_risk_smaller_than_deterministic}, 
	\begin{equation*}
		\inf_{\delta\in\mathcal D_n^{\mathrm r}(\mathcal W)}\sup_{(\bm t,\Sigma)\in\Theta}\rho_\delta(\bm t,\Sigma)
		=\inf_{\substack{\hat{\bm w}:(\mathbb R^d)^n\to\mathcal W\\\hat{\bm w}\ \mathrm{measurable}}}\sup_{(\bm t,\Sigma)\in\Theta}\mathbb E_{\bm\mu_{n,\bm t},\Sigma}\ell_n(\bm t,\hat{\bm w}).
	\end{equation*}
	Similarly, for any $\Pi\in\mathcal P(\Theta)$, \eqref{eq:randomized_bayes_risk_greater_than_deterministic} and \eqref{eq:randomized_bayes_risk_smaller_than_deterministic} give
	\begin{equation*}
		\inf_{\delta\in\mathcal D_n^{\mathrm r}(\mathcal W)}\int_\Theta\rho_\delta(\bm t,\Sigma)\,\Pi(\dd\bm t,\dd\Sigma)
		=\inf_{\substack{\hat{\bm w}:(\mathbb R^d)^n\to\mathcal W\\\hat{\bm w}\ \mathrm{measurable}}}\mathbb E_\Pi\ell_n(\bm T,\hat{\bm w}).
	\end{equation*}
	Combining \eqref{eq:minimax_regret_equals_sup_bayes} with the preceding two equations proves \eqref{eq:cov_compact_minimax}. The barycenter step is the affine-loss equality case of \citet[equation~(3.5)]{LieseMiescke2008}.
	
	Finally, under the Bayes-rule assumption in the last assertion, $\mathcal W\subseteq\Omega$ and Bayes optimality imply, for every $\Pi\in\mathcal P(\Theta)$,
	\begin{equation*}
		\inf_{\delta\in\mathcal D_n^{\mathrm r}(\mathcal W)}\int_\Theta\rho_\delta(\bm t,\Sigma)\,\Pi(\dd\bm t,\dd\Sigma)
		=\inf_{\delta\in\mathcal D_n^{\mathrm r}(\Omega)}\int_\Theta\rho_\delta(\bm t,\Sigma)\,\Pi(\dd\bm t,\dd\Sigma).
	\end{equation*}
	Taking the supremum over $\Pi$ and applying \eqref{eq:cov_compact_minimax} to both decision spaces proves the assertion.
\end{proof}

\subsection{Proof of Lemma~\ref{lem:fast_bayes_minimax}}
\label{app:proof-bayes-fixed}
\begin{proof}
	Conditional on any value of $\bm{T}$,  $\bm{Z}\sim\mathcal{N}(\bm{0},\Sigma)$. This conditional distribution does not depend on $\bm{T}$, so $\bm{Z}\perp\bm{T}$ under the joint law. Since $\bar{\bm{x}}$ is a sufficient statistic, $\bm{Y}$ contains all sample information about $\bm{T}$ relevant to the Bayes decision.
	
	Fix $\pi$. Conditional on $\bm{Y}$, the posterior expected regret of a policy $\bm{w}\in\Omega$ is
	\begin{equation*}
		\mathbb{E}_{\pi,\Sigma}\left[f_*(\bm{\mu}_{n,\bm{T}})-\bm{w}^\top\bm{\mu}_{n,\bm{T}}\mid\bm{Y}\right]=\mathbb{E}_{\pi,\Sigma}\left[f_*(\bm{\mu}_{n,\bm{T}})\mid\bm{Y}\right]-\bm{w}^\top\left(\bm{\mu}+\frac{\bm{m}_{\pi,\Sigma}(\bm{Y})}{\sqrt{n}}\right).
		\label{eq:fast_conditional_bayes_regret}
	\end{equation*}
	The first term does not depend on $\bm{w}$. Maximizing the second term over
	$\Omega$ proves \eqref{eq:fast_bayes_policy}. Substituting that policy and
	taking expectation gives the minimum Bayes risk
	\begin{equation}
		\mathbb{E}_{\pi,\Sigma}\left[
		f_*\left(\bm{\mu}+\frac{\bm{T}}{\sqrt{n}}\right)
		-f_*\left(\bm{\mu}+
		\frac{\bm{m}_{\pi,\Sigma}(\bm{Y})}{\sqrt{n}}\right)
		\right].
		\label{eq:fast_bayes_value_gap}
	\end{equation}
	By the tower property, $\mathbb{E}_{\pi,\Sigma}\bm{m}_{\pi,\Sigma}(\bm{Y})=\mathbb{E}_{\pi,\Sigma}\bm{T}$. Hence the constant and affine terms in \eqref{eq:fast_local_remainder} cancel in \eqref{eq:fast_bayes_value_gap}. This proves the second equality in \eqref{eq:fast_exact_jensen_gap}.
	
	It remains to prove the first equality. Apply Lemma~\ref{lem:cov_compact_minimax} with
	\begin{equation*}
		\Theta=B_2(M)\times\{\Sigma\},\qquad\mathcal W=\Omega.
	\end{equation*}
	These are nonempty compact sets, and $\Omega$ is convex. Retain $\mathcal D_n^{\mathrm r}(\mathcal W)$ and $\rho_\delta(\bm t,\Sigma)$ from the proof of Lemma~\ref{lem:cov_compact_minimax}. Identifying a prior on $\Theta$ with its marginal $\pi$ on $B_2(M)$ gives
	\begin{equation}
		\begin{aligned}
			R_{n,M}^{\Sigma}(\bm\mu)
			&=\sup_{\pi\in\mathcal P(B_2(M))}\inf_{\delta\in\mathcal D_n^{\mathrm r}(\Omega)}\int_{B_2(M)}\rho_\delta(\bm t,\Sigma)\,\pi(\dd\bm t)\\
			&=\sup_{\pi\in\mathcal P(B_2(M))}\inf_{\hat{\bm w}}\mathbb E_{\pi,\Sigma}\ell_n(\bm T,\hat{\bm w}).
		\end{aligned}
		\label{eq:fast_statistical_minimax}
	\end{equation}
	The last infimum is over the deterministic sample-only policies in \eqref{eq:fast_local_minimax_definition}; its value is the Bayes risk computed in \eqref{eq:fast_bayes_value_gap}. This proves the first equality in \eqref{eq:fast_exact_jensen_gap}.
\end{proof}
\subsection{A restricted decision set}
\begin{lemma}[Restricted decision set]
	\label{lem:restricted_decision_set}
	Fix $n$, $L>0$, and $\bm\mu\in\mathbb R^d$, and let $b_L(\varepsilon):=\sup_{\bm t\in B_2(L)}r_{\bm\mu}(\varepsilon\bm t)$. Define
	\begin{equation*}
		\Omega_{n,L}
		:=\operatorname{conv}\left(
		\bigcup_{\bm{u}\in B_2(L)}
		\arg\max_{\bm{w}\in\Omega}
		(\bm{\mu}+\bm{u}/\sqrt{n})^\top\bm{w}
		\right).
	\end{equation*}
	Then $\Omega_{n,L}$ is a nonempty compact convex subset of $\Omega$, and
	\begin{equation}
		0\le\ell_n(\bm t,\bm w)\le2b_{3L}(n^{-1/2}),
		\qquad \bm t\in B_2(L),\ \bm w\in\Omega_{n,L},
		\label{eq:local_loss_upper_bound}
	\end{equation}
	where $\ell_n(\bm t,\bm w)=f_*(\bm\mu+\bm t/\sqrt n)-\bm w^\top(\bm\mu+\bm t/\sqrt n)$. Moreover, for every prior $\pi\in\mathcal P(B_2(M))$ with $M\le L$ and every $\Sigma$, the Bayes rule \eqref{eq:fast_bayes_policy} takes values in $\Omega_{n,L}$.
\end{lemma}
\begin{proof}
	The graph of the local maximizers is
	\begin{equation*}
		\left\{(\bm{u},\bm{w})\in B_2(L)\times\Omega:
		(\bm{\mu}+\bm{u}/\sqrt n)^\top\bm{w}
		=f_*(\bm{\mu}+\bm{u}/\sqrt n)\right\}.
	\end{equation*}
	Continuity of $(\bm{\mu}+\bm{u}/\sqrt n)^\top\bm{w}$ and $f_*$ makes this graph closed in the compact space $B_2(L)\times\Omega$, hence compact. Its image under the continuous coordinate projection $(\bm{u},\bm{w})\mapsto\bm{w}$ is
	\begin{equation*}
		\bigcup_{\bm{u}\in B_2(L)}
		\arg\max_{\bm{w}\in\Omega}
		(\bm{\mu}+\bm{u}/\sqrt n)^\top\bm{w},
	\end{equation*}
	which is therefore compact and nonempty. In finite dimensions, its convex hull is compact: by Carath\'eodory's theorem, every point of the convex hull is a convex combination of at most $d+1$ points, and these combinations form the continuous image of a compact product with the closed probability simplex. Since $\Omega$ is convex, $\Omega_{n,L}\subseteq\Omega$. Hence $\Omega_{n,L}$ is a nonempty compact convex set.
	
	For $\bm{t},\bm{u}\in B_2(L)$ and $\bm{w}\in\arg\max_{\Omega}(\bm{\mu}+\bm{u}/\sqrt{n})^\top\bm{w}$, the subgradient inequality for $r_{\bm{\mu}}$ gives
	\begin{equation*}
		\begin{aligned}
			\ell_n(\bm{t},\bm{w})
			&=r_{\bm{\mu}}(\bm{t}/\sqrt{n})
			-r_{\bm{\mu}}(\bm{u}/\sqrt{n})
			+(\bm{w}-\bm{w}_*(\bm{\mu}))^\top
			(\bm{u}-\bm{t})/\sqrt{n}\\
			&\le r_{\bm{\mu}}(\bm{t}/\sqrt{n})
			+r_{\bm{\mu}}((2\bm{u}-\bm{t})/\sqrt{n})
			-2r_{\bm{\mu}}(\bm{u}/\sqrt{n})
			\le2b_{3L}(n^{-1/2}).
		\end{aligned}
	\end{equation*}
	Since $\ell_n(\bm t,\cdot)$ is affine, the same bound holds on the convex hull $\Omega_{n,L}$, and $\ell_n\ge0$ on $\Omega$. Finally, for $\pi\in\mathcal P(B_2(M))$, conditional Jensen's inequality gives $\|\bm m_{\pi,\Sigma}(\bm Y)\|_2\le\mathbb E_{\pi,\Sigma}[\|\bm T\|_2\mid\bm Y]\le M\le L$, so $\bm w_*(\bm\mu+\bm m_{\pi,\Sigma}(\bm Y)/\sqrt n)\in\arg\max_\Omega(\bm\mu+\bm u/\sqrt n)^\top\bm w$ with $\bm u=\bm m_{\pi,\Sigma}(\bm Y)\in B_2(L)$.
\end{proof}

\subsection{Proof of Theorem~\ref{thm:fast_local_minimax}}
\label{app:proof-smooth-gaussian}
\begin{proof}
	Convexity gives $H_{\bm{\mu}}\succeq0$. As argued before, $\bm{w}_*(\bm{\mu})=\nabla f_*(\bm{\mu})$, and Taylor expansion yields
	\begin{equation*}
		\begin{aligned}
			\epsilon_{n,M}&:=\sup_{\bm{t}\in B_2(M)}
			\left|n r_{\bm{\mu}}(\bm{t}/\sqrt{n})
			-\tfrac12\bm{t}^\top H_{\bm{\mu}}\bm{t}\right|\\&=\sup_{\bm{t}\in B_2(M)}\left|nf_*(\bm{\mu}+\bm{t}/\sqrt{n})-nf_*(\bm{\mu})-\sqrt{n}\bm{w}_*(\bm{\mu})^\top\bm{t}-\bm{t}^\top H_{\bm{\mu}}\bm{t}/2\right|\\&\le\frac{M^2}{2}\sup_{\|\bm h\|_2\le M/\sqrt n}\|\nabla^2f_*(\bm\mu+\bm h)-H_{\bm\mu}\|_{\mathrm{op}}\to0\qquad\text{as }n\to\infty.
		\end{aligned}
		\label{eq:fast_uniform_value_expansion}
	\end{equation*}
	Both $\bm{T}$ and $\bm{m}_{\pi,\Sigma}(\bm{Y})$ lie in $B_2(M)$, which yields $\left|n r_{\bm{\mu}}(\bm{T}/\sqrt{n})-\tfrac12\bm{T}^\top H_{\bm{\mu}}\bm{T}\right|\le\epsilon_{n,M}$ and $\left|n r_{\bm{\mu}}(\bm{m}_{\pi,\Sigma}(\bm{Y})/\sqrt{n})-\tfrac12\bm{m}_{\pi,\Sigma}(\bm{Y})^\top H_{\bm{\mu}}\bm{m}_{\pi,\Sigma}(\bm{Y})\right|\le\epsilon_{n,M}$. 
	Therefore, by taking expectation, 
	\begin{equation*}
		\begin{aligned}
			&\left|n\mathbb{E}_{\pi,\Sigma}\left[
			r_{\bm{\mu}}\left(\frac{\bm{T}}{\sqrt{n}}\right)
			-r_{\bm{\mu}}\left(
			\frac{\bm{m}_{\pi,\Sigma}(\bm{Y})}{\sqrt{n}}\right)
			\right]-\frac12\mathbb{E}_{\pi,\Sigma}\left[
			\bm{T}^\top H_{\bm{\mu}}\bm{T}
			-\bm{m}_{\pi,\Sigma}(\bm{Y})^\top H_{\bm{\mu}}\bm{m}_{\pi,\Sigma}(\bm{Y})\right]\right|\\&\le2\epsilon_{n,M}\to0\qquad \text{as }n\to\infty.
		\end{aligned}
	\end{equation*}
	Furthermore, by taking conditional expectation, it is noteworthy that
	\begin{equation*}\mathbb{E}_{\pi,\Sigma}\left[
		\bm{T}^\top H_{\bm{\mu}}\bm{T}
		-\bm{m}_{\pi,\Sigma}(\bm{Y})^\top H_{\bm{\mu}}\bm{m}_{\pi,\Sigma}(\bm{Y})\right]=\mathbb{E}_{\pi,\Sigma}\left[
		(\bm{T}-\bm{m}_{\pi,\Sigma}(\bm{Y}))^\top
		H_{\bm{\mu}}(\bm{T}-\bm{m}_{\pi,\Sigma}(\bm{Y}))\right].
	\end{equation*}
	Combining with equation \eqref{eq:fast_exact_jensen_gap} proves \eqref{eq:fast_fixed_neighborhood_limit}.
	
	We next prove $\lim_{M\to\infty}\mathcal{B}_{M,\Sigma}(H_{\bm{\mu}})=\operatorname{tr}(\Sigma H_{\bm{\mu}})$. Note that $\bm{Y}=\bm{T}+\bm{Z}$, and thus
	\begin{equation}
		\begin{aligned}
			\mathcal{B}_{M,\Sigma}(H_{\bm{\mu}})&=\sup_{\pi\in\mathcal P(B_2(M))}\mathbb{E}_{\pi,\Sigma}\left[
			(\bm{T}-\bm{m}_{\pi,\Sigma}(\bm{Y}))^\top
			H_{\bm{\mu}}(\bm{T}-\bm{m}_{\pi,\Sigma}(\bm{Y}))\right]\\&=\sup_{\pi\in\mathcal P(B_2(M))}\big\{\mathbb{E}_{\pi,\Sigma}\left[
			(\bm{T}-\bm{Y})^\top
			H_{\bm{\mu}}(\bm{T}-\bm{Y})\right]\\&\qquad-\mathbb{E}_{\pi,\Sigma}\left[
			(\bm{m}_{\pi,\Sigma}(\bm{Y})-\bm{Y})^\top
			H_{\bm{\mu}}(\bm{m}_{\pi,\Sigma}(\bm{Y})-\bm{Y})\right]\big\}\\&\le\mathbb{E}_{\pi,\Sigma}\left[
			\bm{Z}^\top
			H_{\bm{\mu}}\bm{Z}\right]=\operatorname{tr}(\Sigma H_{\bm{\mu}}),
		\end{aligned}
		\label{eq:BM_smaller_than_trace}
	\end{equation}
	The first equality follows because the cross term vanishes conditionally on $\bm Y$:
	\begin{equation*}
		\begin{aligned}
			&\mathbb E_{\pi,\Sigma}\!
			\left[(\bm T-\bm m_{\pi,\Sigma}(\bm Y))^\top H_{\bm\mu}
			(\bm m_{\pi,\Sigma}(\bm Y)-\bm Y)\right]\\
			&\quad=\mathbb E_{\pi,\Sigma}\!
			\left[\mathbb E_{\pi,\Sigma}
			[\bm T-\bm m_{\pi,\Sigma}(\bm Y)\mid\bm Y]^\top
			H_{\bm\mu}(\bm m_{\pi,\Sigma}(\bm Y)-\bm Y)\right]=0.
		\end{aligned}
	\end{equation*}
	
	For the reverse inequality, first consider the untruncated prior $\bm{T}\sim\pi_\tau=\mathcal{N}(\bm{0},\tau^2\Sigma)$. The posterior distribution is $\bm{T}|\bm{Y}\sim\mathcal{N}(\tau^2(1+\tau^2)^{-1}\bm{Y},\tau^2(1+\tau^2)^{-1}\Sigma)$, so
	\begin{equation*}
		\begin{aligned}
			B(\pi_\tau)&:=\mathbb{E}_{\pi_\tau,\Sigma}\left[
			(\bm{T}-\bm{m}_{\pi_\tau,\Sigma}(\bm{Y}))^\top
			H_{\bm{\mu}}(\bm{T}-\bm{m}_{\pi_\tau,\Sigma}(\bm{Y}))\right]\\&=\mathbb{E}_{\pi_\tau,\Sigma}\left[\mathbb{E}\left[(\bm{T}-\bm{m}_{\pi_\tau,\Sigma}(\bm{Y}))^\top H_{\bm{\mu}}(\bm{T}-\bm{m}_{\pi_\tau,\Sigma}(\bm{Y}))|\bm{Y}\right]\right]\\&=\mathbb{E}_{\pi_\tau,\Sigma}\left[\operatorname{tr}(H_{\bm{\mu}}\operatorname{Cov}(\bm{T}|\bm{Y}))\right]\\&=\frac{\tau^2}{1+\tau^2}
			\operatorname{tr}(\Sigma H_{\bm{\mu}}).
		\end{aligned}
		\label{eq:fast_gaussian_prior_constant}
	\end{equation*}
	This prior is not supported on $B_2(M)$, so it cannot be used directly for the local lower bound. Let $\pi_{\tau,M}$ be the conditional law of $\pi_\tau=\mathcal{N}(\bm{0},\tau^2\Sigma)$ given $\|\bm{T}\|_2\le M$, $p_{\tau,M}=\mathbb{P}_{\pi_\tau}\{\|\bm{T}\|_2\le M\}$,
	and write $B(\pi_{\tau,M}):=\mathbb{E}_{\pi_{\tau,M},\Sigma}\left[
	(\bm{T}-\bm{m}_{\pi_{\tau,M},\Sigma}(\bm{Y}))^\top
	H_{\bm{\mu}}(\bm{T}-\bm{m}_{\pi_{\tau,M},\Sigma}(\bm{Y}))\right]$.
	By the minimality of $B$, we have
	\begin{equation*}
		\begin{aligned}
			B(\pi_\tau)&\le\mathbb{E}_{\pi_\tau,\Sigma}\left[
			(\bm{T}-\bm{m}_{\pi_{\tau,M},\Sigma}(\bm{Y}))^\top
			H_{\bm{\mu}}(\bm{T}-\bm{m}_{\pi_{\tau,M},\Sigma}(\bm{Y}))\right]\\&=p_{\tau,M}\mathbb{E}_{\pi_\tau,\Sigma}\left[
			(\bm{T}-\bm{m}_{\pi_{\tau,M},\Sigma}(\bm{Y}))^\top
			H_{\bm{\mu}}(\bm{T}-\bm{m}_{\pi_{\tau,M},\Sigma}(\bm{Y}))\big|\|\bm{T}\|_2\le M\right]\\&\quad+\mathbb{E}_{\pi_\tau,\Sigma}\left[
			(\bm{T}-\bm{m}_{\pi_{\tau,M},\Sigma}(\bm{Y}))^\top
			H_{\bm{\mu}}(\bm{T}-\bm{m}_{\pi_{\tau,M},\Sigma}(\bm{Y}))\mathbb{I}\{\|\bm{T}\|_2>M\}\right].
		\end{aligned}
		\label{eq:fast_prior_truncation}
	\end{equation*}
	Here, $\mathbb{E}_{\pi_\tau,\Sigma}\left[(\bm{T}-\bm{m}_{\pi_{\tau,M},\Sigma}(\bm{Y}))^\top H_{\bm{\mu}}(\bm{T}-\bm{m}_{\pi_{\tau,M},\Sigma}(\bm{Y}))\big|\|\bm{T}\|_2\le M\right]=B(\pi_{\tau, M})$. Since the posterior expectation $\|\bm{m}_{\pi_{\tau,M},\Sigma}(\bm{Y})\|_2\le M$, 
	$$\|\bm{T}-\bm{m}_{\pi_{\tau,M},\Sigma}(\bm{Y})\|_2\le\|\bm{T}\|_2+M\le 2\|\bm{T}\|_2$$ on $\{\|\bm{T}\|_2>M\}$, which yields
	$$(\bm{T}-\bm{m}_{\pi_{\tau,M},\Sigma}(\bm{Y}))^\top H_{\bm{\mu}}(\bm{T}-\bm{m}_{\pi_{\tau,M},\Sigma}(\bm{Y}))\le\|H_{\bm{\mu}}\|_{\operatorname{op}}\|\bm{T}-\bm{m}_{\pi_{\tau,M},\Sigma}(\bm{Y})\|_2^2\le4\|H_{\bm{\mu}}\|_{\operatorname{op}}\|\bm{T}\|_2^2.$$
	Therefore,
	\begin{equation*}
		B(\pi_\tau)\le p_{\tau,M}B(\pi_{\tau, M})+4\|H_{\bm{\mu}}\|_{\operatorname{op}}\mathbb{E}_{\pi_\tau,\Sigma}\left[\|\bm{T}\|_2^2
		\mathbb{I}\{\|\bm{T}\|_2>M\}\right].
	\end{equation*}
	With $\tau$ fixed, as $M\to\infty$, we have $p_{\tau,M}\to 1$ and $\mathbb{E}_{\pi_\tau,\Sigma}\left[\|\bm{T}\|_2^2\mathbb{I}\{\|\bm{T}\|_2>M\}\right]\to 0$, which gives
	\begin{equation*}
		\liminf_{M\to\infty}B(\pi_{\tau, M})\ge B(\pi_{\tau})=\frac{\tau^2}{1+\tau^2}\operatorname{tr}(\Sigma H_{\bm{\mu}}).
	\end{equation*}
	Combined with $\mathcal{B}_{M,\Sigma}(H_{\bm{\mu}})\ge B(\pi_{\tau, M})$, we have $\liminf_{M\to\infty}\mathcal{B}_{M,\Sigma}(H_{\bm{\mu}})\ge\frac{\tau^2}{1+\tau^2}\operatorname{tr}(\Sigma H_{\bm{\mu}})$. Letting $\tau\to\infty$ gives 
	\begin{equation}
		\liminf_{M\to\infty}\mathcal{B}_{M,\Sigma}(H_{\bm{\mu}})\ge\operatorname{tr}(\Sigma H_{\bm{\mu}}).
		\label{eq:BM_larger_than_trace}
	\end{equation}
	Combining \eqref{eq:BM_smaller_than_trace} and \eqref{eq:BM_larger_than_trace}, we prove \eqref{eq:fast_sharp_local_minimax}.
	
	Finally, fix $M>0$ and represent the sample mean as $\bar{\bm{x}}=\bm{\mu}_{n,\bm{t}}+\bm{Z}/\sqrt{n}$, where $\bm{Z}\sim\mathcal{N}(\bm{0},\Sigma)$ has the same distribution for every $\bm{t}\in B_2(M)$. Expectations and probabilities without a subscript below are with respect to $\bm{Z}$, and the regret is evaluated at the true mean $\bm{\mu}_{n,\bm{t}}$. Choose $r>0$ such that $f_*$ is twice continuously differentiable on the closed ball $\{\bm{z}:\|\bm{z}-\bm{\mu}\|_2\le r\}$, and choose a finite constant $C$ satisfying
	\begin{equation}
		\sup_{\|\bm{z}-\bm{\mu}\|_2\le r}
		\|\nabla^2f_*(\bm{z})\|_{\mathrm{op}}\le C.
		\label{eq:hessian_operator_norm_bound}
	\end{equation}
	Fix a noise truncation radius $L>0$. For sufficiently large $n$, $M/\sqrt{n}\le r/2$ and $(M+L)/\sqrt{n}\le r$. Thus $\|\bm{Z}\|_2\le L$ implies $\|\bar{\bm{x}}-\bm{\mu}\|_2\le r$, uniformly over $\bm{t}\in B_2(M)$. Since $\mathbb{E}[\bm{Z}^\top H_{\bm{\mu}}\bm{Z}]=\operatorname{tr}(\Sigma H_{\bm{\mu}})$, we first split the expected difference over three disjoint events:
	\begin{equation*}
		\left|n\mathbb{E}_{\bm\mu_{n,\bm{t}},\Sigma}
		R(\bm{w}_*(\bar{\bm{x}}))
		-\tfrac12\operatorname{tr}(\Sigma H_{\bm{\mu}})\right|\le I_{1,n}(\bm{t})+I_{2,n}(\bm{t})+I_{3,n}(\bm{t}),
	\end{equation*}
	where, with $M$ and $L$ fixed,
	\begin{equation*}
		\begin{aligned}
			I_{1,n}(\bm{t})
			&:=\mathbb{E}\left[
			\left|nR(\bm{w}_*(\bar{\bm{x}}))
			-\tfrac12\bm{Z}^\top H_{\bm{\mu}}\bm{Z}\right|
			\mathbb{I}\{\|\bm{Z}\|_2\le L\}\right],\\
			I_{2,n}(\bm{t})
			&:=\mathbb{E}\left[
			\left|nR(\bm{w}_*(\bar{\bm{x}}))
			-\tfrac12\bm{Z}^\top H_{\bm{\mu}}\bm{Z}\right|
			\mathbb{I}\{\|\bm{Z}\|_2>L,\ \|\bar{\bm{x}}-\bm{\mu}\|_2\le r\}
			\right],\\
			I_{3,n}(\bm{t})
			&:=\mathbb{E}\left[
			\left|nR(\bm{w}_*(\bar{\bm{x}}))
			-\tfrac12\bm{Z}^\top H_{\bm{\mu}}\bm{Z}\right|
			\mathbb{I}\{\|\bar{\bm{x}}-\bm{\mu}\|_2>r\}\right].
		\end{aligned}
	\end{equation*}
	
	For the first term, on $\{\|\bm{Z}\|_2\le L\}$, both
	$\bm{\mu}_{n,\bm{t}}$ and $\bar{\bm{x}}$ lie in the smooth ball $\{\bm{z}:\|\bm{z}-\bm{\mu}\|_2\le r\}$.
	Therefore,
	\begin{equation}
		\begin{aligned}
			nR(\bm{w}_*(\bar{\bm{x}}))&=n(f_*(\bm{\mu}_{n,\bm{t}})-\bm{w}_*(\bar{\bm{x}})^\top\bm{\mu}_{n,\bm{t}})
			\\&=n[f_*(\bm{\mu}_{n,\bm{t}})-f_*(\bar{\bm{x}})
			-\nabla f_*(\bar{\bm{x}})^\top
			(\bm{\mu}_{n,\bm{t}}-\bar{\bm{x}})]\\
			&=\int_0^1(1-s)\bm{Z}^\top
			\nabla^2f_*\left(\bm{\mu}
			+\frac{\bm{t}+(1-s)\bm{Z}}{\sqrt{n}}\right)
			\bm{Z}\,\dd s.
		\end{aligned}
		\label{eq:fast_bregman_integral}
	\end{equation}
	Consequently, on $\{\|\bm{Z}\|_2\le L\}$,
	\begin{equation*}
		\begin{aligned}
			\left|nR(\bm{w}_*(\bar{\bm{x}}))-\tfrac12\bm{Z}^\top H_{\bm{\mu}}\bm{Z}\right|&=\left|\int_0^1(1-s)\bm{Z}^\top\left[\nabla^2f_*\left(\bm{\mu}+\frac{\bm{t}+(1-s)\bm{Z}}{\sqrt{n}}\right)-H_{\bm{\mu}}\right]\bm{Z}\,\dd s\right|
			\\&\le\|\bm{Z}\|_2^2\int_0^1(1-s)\left\|\nabla^2f_*\left(\bm{\mu}+\frac{\bm{t}+(1-s)\bm{Z}}{\sqrt{n}}\right)-H_{\bm{\mu}}\right\|_{\operatorname{op}}\,\dd s\\&\le L^2\sup_{\|\bm{z}-\bm{\mu}\|_2\le(M+L)/\sqrt{n}}
			\|\nabla^2f_*(\bm{z})-H_{\bm{\mu}}\|_{\mathrm{op}}\cdot\int_{0}^1(1-s)\dd s\\&=\frac{L^2}{2}
			\sup_{\|\bm{z}-\bm{\mu}\|_2\le(M+L)/\sqrt{n}}
			\|\nabla^2f_*(\bm{z})-H_{\bm{\mu}}\|_{\mathrm{op}},
		\end{aligned}
	\end{equation*}
	which immediately gives
	\begin{equation*}
		\begin{aligned}
			\sup_{\bm{t}\in B_2(M)}I_{1,n}(\bm{t})\le\frac{L^2}{2}\sup_{\|\bm{z}-\bm{\mu}\|_2\le(M+L)/\sqrt{n}}\|\nabla^2f_*(\bm{z})-H_{\bm{\mu}}\|_{\mathrm{op}}
			\to0
		\end{aligned}
	\end{equation*}
	as $n\to\infty$ for each fixed $L$, by continuity of the Hessian.
	
	For the second term, formula \eqref{eq:fast_bregman_integral} remains valid whenever $\|\bar{\bm{x}}-\bm{\mu}\|_2\le r$, since $\|\bm{\mu}_{n,\bm{t}}-\bm{\mu}\|_2\le M/\sqrt{n}\le r/2$. Therefore, \eqref{eq:hessian_operator_norm_bound} implies 
	\begin{equation*}
		|nR(\bm{w}_*(\bar{\bm{x}}))|=\left|\int_0^1(1-s)\bm{Z}^\top\nabla^2f_*\left(\bm{\mu}+\frac{\bm{t}+(1-s)\bm{Z}}{\sqrt{n}}\right)\bm{Z}\,\dd s\right|\le C\|\bm{Z}\|_2^2/2.
	\end{equation*}
	Also, $0\le\bm{Z}^\top H_{\bm{\mu}}\bm{Z}/2\le C\|\bm{Z}\|_2^2/2$. Therefore, for all sufficiently large $n$,
	\begin{equation*}
		\sup_{\bm{t}\in B_2(M)}I_{2,n}(\bm{t})
		\le C\mathbb{E}\left[
		\|\bm{Z}\|_2^2\mathbb{I}\{\|\bm{Z}\|_2>L\}\right].
	\end{equation*}
	This bound tends to zero as $L\to\infty$, since
	$\mathbb{E}\|\bm{Z}\|_2^2=\operatorname{tr}(\Sigma)\le d\overline\lambda<\infty$.
	
	For the third term, $\|\bar{\bm{x}}-\bm{\mu}\|_2>r$ implies
	$\|\bm{Z}\|_2>r\sqrt{n}/2$, uniformly over $\bm{t}\in B_2(M)$.
	By Chernoff's inequality and 
	$\mathbb{E}e^{\|\bm{Z}\|_2^2/(4\overline\lambda)}\le2^{d/2}$,
	\begin{equation*}
		\sup_{\bm{t}\in B_2(M)}
		\mathbb{P}\{\|\bar{\bm{x}}-\bm{\mu}\|_2>r\}
		\le\mathbb{P}\{\|\bm{Z}\|_2>r\sqrt{n}/2\}
		\le e^{-r^2n/(16\overline\lambda)}\mathbb{E}e^{\|\bm{Z}\|_2^2/(4\overline\lambda)}\le2^{d/2}e^{-r^2n/(16\overline\lambda)}.
	\end{equation*}
	Recall $R_\Omega=\sup_{\bm{w}\in\Omega}\|\bm{w}\|_2<\infty$. We have
	\begin{equation*}
		0\le R(\bm{w}_*(\bar{\bm{x}}))=(\bm{w}_*(\bm{\mu}_{n,\bm{t}})-\bm{w}_*(\bar{\bm{x}}))^\top\bm{\mu}_{n,\bm{t}}\le2R_\Omega\left(\|\bm{\mu}\|_2+\frac{M}{\sqrt{n}}\right).
	\end{equation*}
	Noticing that $0\le\bm{Z}^\top H_{\bm{\mu}}\bm{Z}/2\le C\|\bm{Z}\|_2^2/2$, we obtain
	\begin{equation*}
		\begin{aligned}
			\sup_{\bm t\in B_2(M)}I_{3,n}(\bm t)
			&\le2nR_\Omega\left(\|\bm\mu\|_2+\frac{M}{\sqrt n}\right)
			2^{d/2}e^{-r^2n/(16\overline\lambda)}\\
			&\quad+\frac C2\mathbb E\left[
			\|\bm Z\|_2^2
			\mathbb I\{\|\bm Z\|_2>r\sqrt n/2\}\right]
			\longrightarrow0.
		\end{aligned}
	\end{equation*}
	as $n\to\infty$, since $\mathbb{E}\|\bm{Z}\|_2^2=\operatorname{tr}(\Sigma)\le d\overline\lambda<\infty$.
	
	Combining the three bounds and first letting $n\to\infty$ gives
	\begin{equation*}
		\limsup_{n\to\infty}\sup_{\bm{t}\in B_2(M)}
		\left|n\mathbb{E}_{\bm\mu_{n,\bm{t}},\Sigma}
		R(\bm{w}_*(\bar{\bm{x}}))
		-\tfrac12\operatorname{tr}(\Sigma H_{\bm{\mu}})\right|\le C\mathbb{E}\left[
		\|\bm{Z}\|_2^2\mathbb{I}\{\|\bm{Z}\|_2>L\}\right].
	\end{equation*}
	Letting $L\to\infty$ proves \eqref{eq:fast_uniform_plugin_limit} for each fixed $\Sigma$.
	
	The expansion error $\epsilon_{n,M}$ does not depend on $\Sigma$, so
	\begin{equation*}
		\sup_{\Sigma\in\mathcal C}
		\left|nR_{n,M}^{\Sigma}(\bm\mu)
		-\tfrac12\mathcal B_{M,\Sigma}(H_{\bm\mu})\right|
		\le2\epsilon_{n,M}\to0.
	\end{equation*}
	It remains to make the plug-in limit uniform over $\Sigma\in\mathcal C$. Keep the same $r$ and Hessian bound as above; neither depends on $\Sigma$. In $I_{1,n}(\bm t)$, $I_{2,n}(\bm t)$, and $I_{3,n}(\bm t)$, the expectation is now evaluated under $\bm Z=\Sigma^{1/2}\bm G$, where $\bm G\sim\mathcal N(\bm0,I_d)$; the dependence on $\Sigma$ remains implicit in this notation. For each fixed $L$ and all sufficiently large $n$, the first bound gives
	\begin{equation*}
		\sup_{\Sigma\in\mathcal C}\sup_{\bm t\in B_2(M)}I_{1,n}(\bm t)
		\le\frac{L^2}{2}\sup_{\|\bm z-\bm\mu\|_2\le(M+L)/\sqrt n}\|\nabla^2f_*(\bm z)-H_{\bm\mu}\|_{\mathrm{op}}\to0.
	\end{equation*}
	For the second bound, $\Sigma\preceq\overline\lambda I_d$ implies, pointwise in $\bm G$, $\|\bm Z\|_2^2=\bm G^\top\Sigma\bm G\le\overline\lambda\|\bm G\|_2^2$, implying $\{\|\bm Z\|_2>L\}\subseteq\{\|\bm G\|_2>L/\sqrt{\overline\lambda}\}$. Consequently,
	\begin{equation*}
		\begin{aligned}
			\sup_{\Sigma\in\mathcal C}\sup_{\bm t\in B_2(M)}I_{2,n}(\bm t)
			&\le C\sup_{\Sigma\in\mathcal C}
			\mathbb E\left[\|\bm Z\|_2^2\mathbb I\{\|\bm Z\|_2>L\}\right]
			\\&\le C\overline\lambda\mathbb E\left[
			\|\bm G\|_2^2\mathbb I\{\|\bm G\|_2>L/\sqrt{\overline\lambda}\}
			\right].
		\end{aligned}
	\end{equation*}
	The right-hand side tends to zero as $L\to\infty$ by DCT, since $\mathbb E\|\bm G\|_2^2=d<\infty$.
	
	For the third bound, if $M/\sqrt n\le r/2$, then
	\begin{equation*}
		\{\|\bar{\bm x}-\bm\mu\|_2>r\}
		\subseteq\{\|\bm Z\|_2>r\sqrt n/2\}
		\subseteq\{\|\bm G\|_2>r\sqrt n/(2\sqrt{\overline\lambda})\}.
	\end{equation*}
	The preceding exponential probability bound and the same pointwise comparison therefore yield
	\begin{equation*}
		\begin{aligned}
			\sup_{\Sigma\in\mathcal C}\sup_{\bm t\in B_2(M)}I_{3,n}(\bm t)
			&\le2nR_\Omega\left(\|\bm\mu\|_2+\frac M{\sqrt n}\right)2^{d/2}e^{-r^2n/(16\overline\lambda)}\\
			&\quad+\frac{C\overline\lambda}{2}\mathbb E\left[\|\bm G\|_2^2\mathbb I\{\|\bm G\|_2>r\sqrt n/(2\sqrt{\overline\lambda})\}\right]\to0.
		\end{aligned}
	\end{equation*}
	The first term tends to zero because exponential decay dominates its polynomial factor, and the second tends to zero by DCT. All three bounds hold for the same sufficiently large $n$, independently of $\bm t$ and $\Sigma$. Taking the supremum in the three-term decomposition and then letting $n\to\infty$, with $M$ and $L$ fixed, gives
	\begin{equation*}
		\limsup_{n\to\infty}\sup_{\Sigma\in\mathcal C}\sup_{\bm t\in B_2(M)}
		\left|n\mathbb E_{\bm\mu_{n,\bm t},\Sigma}R(\bm w_*(\bar{\bm x}))-\tfrac12\operatorname{tr}(\Sigma H_{\bm\mu})\right|\le C\overline\lambda\mathbb E\left[\|\bm G\|_2^2\mathbb I\{\|\bm G\|_2>L/\sqrt{\overline\lambda}\}\right].
	\end{equation*}
	Letting $L\to\infty$ proves \eqref{eq:fast_uniform_plugin_limit}, uniformly over $\Sigma\in\mathcal C$.
\end{proof}

\section{Proofs for General Local Scales and the Optimization Modulus}
\label{app:general-scales}

\subsection{From pointwise to uniform convergence}
\begin{proposition}
	\label{pr:pointtouniform}
	The pointwise convergence \eqref{eq:fast_general_local_expansion0} implies convergence uniformly over $\bm{t}\in B_2(M)$:
	\begin{equation}
		\lim_{\varepsilon\downarrow0}
		\sup_{\bm{t}\in B_2(M)}
		\left|\frac{r_{\bm{\mu}}(\varepsilon\bm{t})}{a(\varepsilon)}
		-\psi(\bm{t})\right|=0.
		\label{eq:fast_general_local_expansion}
	\end{equation}
\end{proposition}

\begin{proof}
	Indeed, $r_{\bm{\mu}}$ is nonnegative and convex, with
	$r_{\bm{\mu}}(\bm{0})=0$, and $a(\varepsilon)>0$. Passing to pointwise
	limits gives $\psi(\bm{t})\ge0$, $\psi(\bm{0})=0$, and
	\begin{equation*}
		\psi(\lambda\bm{s}+(1-\lambda)\bm{t})
		\le\lambda\psi(\bm{s})+(1-\lambda)\psi(\bm{t}),
		\qquad \bm{s},\bm{t}\in\mathbb{R}^d,\quad \lambda\in[0,1].
	\end{equation*}
	We next prove uniform convergence on $B_2(M)$ for a fixed $M>0$.
	
	Let $\mathcal{V}_M=\{-(M+1),M+1\}^d$ be the vertices of
	$[-(M+1),M+1]^d$, and set
	$C_M=1+\sup_{\bm{z}\in\mathcal{V}_M}\psi(\bm{z})<\infty$.
	Since $\mathcal{V}_M$ is finite, pointwise convergence gives
	$0<\varepsilon_M<\varepsilon_0$ such that
	\begin{equation*}
		\sup_{\bm{z}\in\mathcal{V}_M}
		\frac{r_{\bm{\mu}}(\varepsilon\bm{z})}{a(\varepsilon)}
		\le C_M,
		\qquad 0<\varepsilon<\varepsilon_M.
	\end{equation*}
	For every $\bm{t}\in B_2(M+1)\subseteq\operatorname{conv}(\mathcal{V}_M)$,
	convexity therefore yields the uniform bound
	\begin{equation*}
		0\le\frac{r_{\bm{\mu}}(\varepsilon\bm{t})}{a(\varepsilon)}
		\le\sup_{\bm{z}\in\mathcal{V}_M}
		\frac{r_{\bm{\mu}}(\varepsilon\bm{z})}{a(\varepsilon)}
		\le C_M,
		\qquad 0<\varepsilon<\varepsilon_M.
	\end{equation*}
	For distinct $\bm{s},\bm{t}\in B_2(M)$, write
	\begin{equation*}
		\begin{aligned}
			h&=\|\bm{t}-\bm{s}\|_2,\qquad
			\bm{z}=\bm{t}+\frac{\bm{t}-\bm{s}}{h}\in B_2(M+1),\qquad
			\bm{t}=\frac{1}{1+h}\bm{s}+\frac{h}{1+h}\bm{z}.
		\end{aligned}
	\end{equation*}
	Applying convexity to this combination gives
	\begin{equation*}
		\frac{r_{\bm{\mu}}(\varepsilon\bm{t})
			-r_{\bm{\mu}}(\varepsilon\bm{s})}{a(\varepsilon)}
		\le\frac{h}{1+h}
		\frac{r_{\bm{\mu}}(\varepsilon\bm{z})
			-r_{\bm{\mu}}(\varepsilon\bm{s})}{a(\varepsilon)}\le\frac{h}{1+h}C_M\le C_M\|\bm{t}-\bm{s}\|_2.
	\end{equation*}
	Interchanging $\bm{s}$ and $\bm{t}$ proves, for all
	$\bm{s},\bm{t}\in B_2(M)$ and $0<\varepsilon<\varepsilon_M$,
	\begin{equation*}
		\left|\frac{r_{\bm{\mu}}(\varepsilon\bm{t})}{a(\varepsilon)}
		-\frac{r_{\bm{\mu}}(\varepsilon\bm{s})}{a(\varepsilon)}\right|
		\le C_M\|\bm{t}-\bm{s}\|_2.
	\end{equation*}
	Taking pointwise limits also gives
	\begin{equation*}
		|\psi(\bm{t})-\psi(\bm{s})|
		\le C_M\|\bm{t}-\bm{s}\|_2,
		\qquad \bm{s},\bm{t}\in B_2(M).
	\end{equation*}
	Thus $\psi$ is locally Lipschitz, hence continuous.
	
	Finally, fix $\delta>0$. By compactness, there are finitely many points
	$\bm{t}_1,\ldots,\bm{t}_N\in B_2(M)$ satisfying
	\begin{equation*}
		\sup_{\bm{t}\in B_2(M)}\min_{1\le j\le N}
		\|\bm{t}-\bm{t}_j\|_2\le\delta.
	\end{equation*}
	For any $\bm{t}\in B_2(M)$, choose $j$ with
	$\|\bm{t}-\bm{t}_j\|_2\le\delta$. The two Lipschitz bounds imply
	\begin{equation*}
		\begin{aligned}
			\left|\frac{r_{\bm{\mu}}(\varepsilon\bm{t})}{a(\varepsilon)}
			-\psi(\bm{t})\right|
			&\le\left|\frac{r_{\bm{\mu}}(\varepsilon\bm{t})
				-r_{\bm{\mu}}(\varepsilon\bm{t}_j)}{a(\varepsilon)}\right|+\left|\frac{r_{\bm{\mu}}(\varepsilon\bm{t}_j)}
			{a(\varepsilon)}-\psi(\bm{t}_j)\right|
			+|\psi(\bm{t}_j)-\psi(\bm{t})|\\
			&\le2C_M\delta+\sup_{1\le j\le N}
			\left|\frac{r_{\bm{\mu}}(\varepsilon\bm{t}_j)}{a(\varepsilon)}
			-\psi(\bm{t}_j)\right|.
		\end{aligned}
	\end{equation*}
	For this fixed finite set of points, the maximum tends to zero as
	$\varepsilon\downarrow0$. Hence
	\begin{equation*}
		\limsup_{\varepsilon\downarrow0}
		\sup_{\bm{t}\in B_2(M)}
		\left|\frac{r_{\bm{\mu}}(\varepsilon\bm{t})}{a(\varepsilon)}
		-\psi(\bm{t})\right|
		\le2C_M\delta,\qquad \delta>0.
	\end{equation*}
	Letting $\delta\downarrow0$ proves
	\eqref{eq:fast_general_local_expansion}. Thus finiteness of all the pointwise limits suffices.
\end{proof}

\subsection{Proof of Theorem~\ref{thm:fast_general_local_rate}}
\label{app:proof-general-rate}
\label{app:proof-general-scale}
\begin{proof}
	For every $\pi\in\mathcal{P}(B_2(M))$, both $\bm{T}$ and
	$\bm{m}_{\pi,\Sigma}(\bm{Y})$ belong to $B_2(M)$ almost surely; for the latter,
	$\|\mathbb{E}_{\pi,\Sigma}[\bm{T}\mid\bm{Y}]\|_2
	\le\mathbb{E}_{\pi,\Sigma}[\|\bm{T}\|_2\mid\bm{Y}]\le M$.
	Consequently,
	\begin{equation*}
		\begin{aligned}
			&\quad\left|\mathbb{E}_{\pi,\Sigma}\left[
			\frac{r_{\bm{\mu}}(\bm{T}/\sqrt{n})
				-r_{\bm{\mu}}(\bm{m}_{\pi,\Sigma}(\bm{Y})/\sqrt{n})}
			{a(n^{-1/2})}
			\right]
			-\mathbb{E}_{\pi,\Sigma}\left[
			\psi(\bm{T})-\psi(\bm{m}_{\pi,\Sigma}(\bm{Y}))\right]\right|\\
			&\le 2\sup_{\bm{t}\in B_2(M)}
			\left|\frac{r_{\bm{\mu}}(\bm{t}/\sqrt{n})}{a(n^{-1/2})}
			-\psi(\bm{t})\right|.
		\end{aligned}
	\end{equation*}
	The bound is independent of $\pi$. Taking suprema over $\pi$ and using
	\eqref{eq:fast_exact_jensen_gap} proves
	\eqref{eq:fast_general_local_error}; assumption
	\eqref{eq:fast_general_local_expansion} then proves
	\eqref{eq:fast_general_local_minimax}.
	
	It remains to determine when $V_{M,\Sigma}(\psi)$ is positive. For every prior
	$\pi\in\mathcal{P}(B_2(M))$, the tower property and conditional Jensen's
	inequality give
	\begin{equation*}
		\mathbb{E}_{\pi,\Sigma}\left[
		\psi(\bm{T})-\psi(\bm{m}_{\pi,\Sigma}(\bm{Y}))\right]=\mathbb{E}_{\pi,\Sigma}\left[
		\mathbb{E}_{\pi,\Sigma}[\psi(\bm{T})\mid\bm{Y}]
		-\psi\left(\mathbb{E}_{\pi,\Sigma}[\bm{T}\mid\bm{Y}]\right)\right]\ge0.
	\end{equation*}
	Since $\psi$ is continuous and both arguments belong to the compact
	ball $B_2(M)$, we also have
	\begin{equation*}
		0\le V_{M,\Sigma}(\psi)
		\le\sup_{\bm{t}\in B_2(M)}\psi(\bm{t})
		-\min_{\bm{t}\in B_2(M)}\psi(\bm{t})<\infty.
	\end{equation*}
	If $\psi$ is affine on this ball, write
	$\psi(\bm{t})=b+\bm{c}^\top\bm{t}$ there. For every such prior,
	\begin{equation*}
		\mathbb{E}_{\pi,\Sigma}\left[
		\psi(\bm{T})-\psi(\bm{m}_{\pi,\Sigma}(\bm{Y}))\right]
		=\bm{c}^\top\left(
		\mathbb{E}_{\pi,\Sigma}[\bm{T}]
		-\mathbb{E}_{\pi,\Sigma}[\mathbb{E}_{\pi,\Sigma}[\bm{T}\mid\bm{Y}]]\right)
		=0,
	\end{equation*}
	so $V_{M,\Sigma}(\psi)=0$.
	
	Conversely, suppose that $\psi$ is not affine on $B_2(M)$. There exist
	$\bm{t}_1,\bm{t}_2\in B_2(M)$ such that
	\begin{equation*}
		\Delta:=\frac{\psi(\bm{t}_1)+\psi(\bm{t}_2)}{2}
		-\psi\left(\frac{\bm{t}_1+\bm{t}_2}{2}\right)>0.
	\end{equation*}
	Indeed, if midpoint equality held for every pair, $\psi$ is affine on every segment in
	the ball, hence affine on the ball. Choose the two-point prior
	$\pi$ with
	$\mathbb{P}_\pi\{\bm{T}=\bm{t}_1\}
	=\mathbb{P}_\pi\{\bm{T}=\bm{t}_2\}=1/2$.
	In the experiment $\bm{Y}=\bm{T}+\bm{Z}$, Bayes' formula gives the
	posterior probability
	\begin{equation*}
		q(\bm{y}):=\mathbb{P}_{\pi,\Sigma}\{\bm{T}=\bm{t}_1\mid\bm{Y}=\bm{y}\}
		=\frac{\exp(-\|\Sigma^{-1/2}(\bm{y}-\bm{t}_1)\|_2^2/2)}
		{\exp(-\|\Sigma^{-1/2}(\bm{y}-\bm{t}_1)\|_2^2/2)
			+\exp(-\|\Sigma^{-1/2}(\bm{y}-\bm{t}_2)\|_2^2/2)}.
	\end{equation*}
	Both likelihoods are strictly positive, so $0<q(\bm{y})<1$ for every
	$\bm{y}\in\mathbb{R}^d$, and
	\begin{equation*}
		\begin{aligned}
			\bm{m}_{\pi,\Sigma}(\bm{y})
			&=q(\bm{y})\bm{t}_1+(1-q(\bm{y}))\bm{t}_2,\\
			\mathbb{E}_{\pi,\Sigma}[\psi(\bm{T})\mid\bm{Y}=\bm{y}]
			&=q(\bm{y})\psi(\bm{t}_1)
			+(1-q(\bm{y}))\psi(\bm{t}_2).
		\end{aligned}
	\end{equation*}
	To see explicitly why the conditional Jensen gap is positive, define
	\begin{equation*}
		g(q):=q\psi(\bm{t}_1)+(1-q)\psi(\bm{t}_2)
		-\psi(q\bm{t}_1+(1-q)\bm{t}_2),\qquad 0\le q\le1.
	\end{equation*}
	Convexity of $\psi$ makes $g$ concave, with
	$g(0)=g(1)=0$ and $g(1/2)=\Delta$. Applying concavity on $[0,1/2]$
	and $[1/2,1]$, respectively, yields
	\begin{equation*}
		g(q)\ge
		\begin{cases}
			2q\Delta,&0\le q\le1/2,\\
			2(1-q)\Delta,&1/2\le q\le1.
		\end{cases}
	\end{equation*}
	Consequently, this particular prior satisfies
	\begin{equation*}
		\begin{aligned}
			V_{M,\Sigma}(\psi)
			&\ge\mathbb E_{\pi,\Sigma}
			[\psi(\bm T)-\psi(\bm m_{\pi,\Sigma}(\bm Y))]\\
			&=\mathbb E_{\pi,\Sigma}[g(q(\bm Y))]\\
			&\ge2\Delta\,\mathbb E_{\pi,\Sigma}
			[\min\{q(\bm Y),1-q(\bm Y)\}]>0.
		\end{aligned}
	\end{equation*}
For the equally weighted two-point prior, the posterior testing error is
\[
    \mathbb E_{\pi,\Sigma}\min\{q(\bm Y),1-q(\bm Y)\}
    =\Phi\!\left(-\frac12\|\Sigma^{-1/2}(\bm t_1-\bm t_2)\|_2\right)
    \ge\Phi(-M/\sqrt{\underline\lambda})>0.
\]
Indeed, the log likelihood ratio is
$(\bm t_1-\bm t_2)^\top\Sigma^{-1}(\bm Y-(\bm t_1+\bm t_2)/2)$;
under the two hypotheses it has means $\pm D^2/2$ and variance $D^2$, where
$D=\|\Sigma^{-1/2}(\bm t_1-\bm t_2)\|_2\le2M/\sqrt{\underline\lambda}$.
	Thus the positive lower bound for $V_{M,\Sigma}(\psi)$ is uniform over
	$\mathcal C$, where $\Phi$ denotes the standard normal distribution
	function. Finally, \eqref{eq:fast_general_local_minimax} gives $R_{n,M}^{\Sigma}(\bm{\mu})\asymp a(n^{-1/2})$ when $V_{M,\Sigma}(\psi)>0$, and $R_{n,M}^{\Sigma}(\bm{\mu})=o(a(n^{-1/2}))$ when $V_{M,\Sigma}(\psi)=0$.
	
	To prove the plug-in bound, fix $M>0$ with $V_{M,\Sigma}(\psi)>0$. We first derive a growth bound for the remainder that will control unbounded Gaussian perturbations. Recall $b_M(\varepsilon):=\sup_{\bm{t}\in B_2(M)}r_{\bm{\mu}}(\varepsilon\bm{t})$ for $\varepsilon\ge0$. Since $r_{\bm{\mu}}$ is nonnegative and convex and vanishes at zero, $b_M$ is nonnegative, convex, and nondecreasing, with $b_M(0)=0$. $f_*$ is $R_\Omega$-Lipschitz, where $R_\Omega=\sup_{\bm{w}\in\Omega}\|\bm{w}\|_2$, so
	\begin{equation}
		0\le b_M(\varepsilon)\le2R_\Omega M\varepsilon,\qquad \varepsilon\ge0.
		\label{eq:fast_radial_global_bound}
	\end{equation}
	\eqref{eq:fast_general_local_expansion} implies $\psi\ge0$ and $\psi(\bm{0})=0$. Since $V_{M,\Sigma}(\psi)>0$, $\sup_{\bm{t}\in B_2(M)}\psi(\bm{t})>0$. \eqref{eq:fast_general_local_expansion} gives
	\begin{equation}
		\left|\frac{b_M(\varepsilon)}{a(\varepsilon)}
		-\sup_{\bm{t}\in B_2(M)}\psi(\bm{t})\right|\le\sup_{\bm{t}\in B_2(M)}\left|\frac{r_{\bm{\mu}}(\varepsilon\bm{t})}{a(\varepsilon)}-\psi(\bm{t})\right|\to0.
		\label{eq:fast_radial_scale_limits}
	\end{equation}
	Applying the same argument to $b_M(2\varepsilon)$ gives
	\begin{equation}
		\frac{b_M(2\varepsilon)}{b_M(\varepsilon)}
		\to\frac{\sup_{\bm{t}\in B_2(2M)}\psi(\bm{t})}{\sup_{\bm{t}\in B_2(M)}\psi(\bm{t})}<\infty.
	\end{equation}
	In particular, $b_M(\varepsilon)>0$ for small positive $\varepsilon$, hence for every $\varepsilon>0$ by monotonicity. Choose $\varepsilon_1>0$ sufficiently small that the ratio $b_M(2\varepsilon)/b_M(\varepsilon)$ is bounded for $0<\varepsilon\le\varepsilon_1$. For $\varepsilon\ge\varepsilon_1$, convexity and \eqref{eq:fast_radial_global_bound} imply
	\begin{equation*}
		b_M(\varepsilon)
		\ge\frac{\varepsilon}{\varepsilon_1}b_M(\varepsilon_1),
		\qquad
		\frac{b_M(2\varepsilon)}{b_M(\varepsilon)}
		\le\frac{4R_\Omega M\varepsilon_1}{b_M(\varepsilon_1)}.
	\end{equation*}
	Consequently, a finite constant $C_b\ge2$ can be chosen so that $b_M(2\varepsilon)\le C_b b_M(\varepsilon)$ for every $\varepsilon>0$. Set $c_b=\log_2C_b$. For $\lambda\ge1$, take $k=\lceil\log_2\lambda\rceil$ and iterate this inequality to obtain
	\begin{equation}
		b_M(\lambda\varepsilon)
		\le b_M(2^k\varepsilon)
		\le C_b^k b_M(\varepsilon)
		\le C_b\lambda^{c_b} b_M(\varepsilon),
		\qquad \varepsilon>0.
		\label{eq:fast_radial_polynomial_growth}
	\end{equation}
	
	Under the true mean
	$\bm{\mu}_{n,\bm{t}}=\bm{\mu}+\bm{t}/\sqrt{n}$, write
	$\bar{\bm{x}}=\bm{\mu}+(\bm{t}+\bm{Z})/\sqrt{n}$ with
	$\bm{Z}\sim\mathcal{N}(\bm{0},\Sigma)$. Expanding the remainder and using
	$f_*(\bar{\bm{x}})=\bm{w}_*(\bar{\bm{x}})^\top\bar{\bm{x}}$ gives
	\begin{equation*}
		R(\bm{w}_*(\bar{\bm{x}}))
		=r_{\bm{\mu}}(\bm{t}/\sqrt{n})
		-r_{\bm{\mu}}((\bm{t}+\bm{Z})/\sqrt{n})+
		(\bm{w}_*(\bar{\bm{x}})-\bm{w}_*(\bm{\mu}))^\top\bm{Z}/\sqrt{n}.
	\end{equation*}
	The maximizer $\bm{w}_*(\bar{\bm{x}})$ is feasible at
	$\bar{\bm{x}}+\bm{Z}/\sqrt{n}$, so
	\begin{equation*}
		r_{\bm{\mu}}((\bm{t}+2\bm{Z})/\sqrt{n})\ge r_{\bm{\mu}}((\bm{t}+\bm{Z})/\sqrt{n})+(\bm{w}_*(\bar{\bm{x}})-\bm{w}_*(\bm{\mu}))^\top\bm{Z}/\sqrt{n}.
	\end{equation*}
	Combining these two formulas and using $r_{\bm{\mu}}\ge0$ yields
	\begin{equation}
		\begin{aligned}
			0\le R(\bm{w}_*(\bar{\bm{x}}))
			&\le r_{\bm{\mu}}(\bm{t}/\sqrt{n})
			+r_{\bm{\mu}}((\bm{t}+2\bm{Z})/\sqrt{n})-2r_{\bm{\mu}}((\bm{t}+\bm{Z})/\sqrt{n})\\
			&\le r_{\bm{\mu}}(\bm{t}/\sqrt{n})
			+r_{\bm{\mu}}((\bm{t}+2\bm{Z})/\sqrt{n}).
		\end{aligned}
		\label{eq:fast_general_plugin_remainder_bound}
	\end{equation}
	For $\bm{t}\in B_2(M)$, we have
	$\|\bm{t}+2\bm{Z}\|_2\le M+2\|\bm{Z}\|_2$. The definition of $b_M$
	and \eqref{eq:fast_radial_polynomial_growth} give
	\begin{equation*}
		r_{\bm{\mu}}((\bm{t}+2\bm{Z})/\sqrt{n})\le b_M\left(\left(1+\frac{2\|\bm{Z}\|_2}{M}\right)\Big/\sqrt{n}\right)\le C_b\left(1+\frac{2\|\bm{Z}\|_2}{M}\right)^{c_b} b_M(n^{-1/2}).
	\end{equation*}
	Gaussian random vectors have finite moments of every fixed order. Taking expectations in \eqref{eq:fast_general_plugin_remainder_bound}, we obtain
	\begin{equation*}
		\sup_{\bm{t}\in B_2(M)}
		\mathbb{E}_{\bm\mu_{n,\bm{t}},\Sigma}R(\bm{w}_*(\bar{\bm{x}}))\le b_M(n^{-1/2})
		\left[1+C_b\mathbb{E}
		\left(1+\frac{2\|\bm{Z}\|_2}{M}\right)^{c_b}\right].
	\end{equation*}
	By \eqref{eq:fast_radial_scale_limits}, for sufficiently large $n$, $b_M(n^{-1/2})\le2a(n^{-1/2})\sup_{\bm{t}\in B_2(M)}\psi(\bm{t})$. Therefore, taking $C_M=2\sup_{\bm t\in B_2(M)}\psi(\bm t)\left[1+C_b\mathbb E\left(1+2\|\bm Z\|_2/M\right)^{c_b}\right]<\infty$, we have
	\begin{equation*}
		\begin{aligned}
			\sup_{\bm{t}\in B_2(M)}\mathbb{E}_{\bm\mu_{n,\bm{t}},\Sigma}R(\bm{w}_*(\bar{\bm{x}}))\le C_M a(n^{-1/2}).
		\end{aligned}
	\end{equation*}
	Finally, \eqref{eq:fast_general_local_minimax} gives the lower bound with $c_M=V_{M,\Sigma}(\psi)/2$ for large $n$, and the definition of $R_{n,M}^{\Sigma}$ gives the middle inequality in \eqref{eq:fast_general_plugin_rate}.
	
	The right side of \eqref{eq:fast_general_local_error} is independent of
	$\Sigma$, proving the asserted uniform convergence. Writing
	$\bm Z=\Sigma^{1/2}\bm G$ with $\bm G\sim\mathcal N(\bm0,I_d)$ gives
	\begin{equation*}
		\sup_{\Sigma\in\mathcal C}
		\mathbb E(1+2\|\bm Z\|_2/M)^{c_b}
		\le\mathbb E(1+2\sqrt{\overline\lambda}\|\bm G\|_2/M)^{c_b}<\infty.
	\end{equation*}
	Together with the uniform positive lower bound just proved, this permits
	$c_M,C_M$ in \eqref{eq:fast_general_plugin_rate} to be chosen uniformly over $\Sigma\in\mathcal C$.
\end{proof}

\subsection{Proof of Corollary~\ref{cor:fast_first_order_limit}}
\label{app:proof-first-order}
\begin{proof}
	For each fixed $\bm{t}$ and $\varepsilon>0$, $\forall \bm{w}\in \mathcal W_*(\bm{\mu})$, optimality gives
	\begin{equation*}
		f_*(\bm{\mu}+\varepsilon\bm{t})-f_*(\bm{\mu})\ge(\bm{\mu}+\varepsilon\bm{t})^\top\bm{w}-\bm{\mu}^\top\bm{w}=\varepsilon\bm{t}^\top\bm{w},
	\end{equation*}
	and thus
	\begin{equation}
		\sup_{\bm{w}\in \mathcal W_*(\bm{\mu})}\bm{w}^\top\bm{t}
		\le\frac{f_*(\bm{\mu}+\varepsilon\bm{t})-f_*(\bm{\mu})}
		{\varepsilon}
		\le\bm{w}_*(\bm{\mu}+\varepsilon\bm{t})^\top\bm{t}.
		\label{eq:two_bounds}
	\end{equation}
	Take any
	sequence $\varepsilon_k\downarrow0$. By compactness of $\Omega$, there
	is a subsequence, still indexed by $k$, such that
	\begin{equation*}
		\bm{w}_*(\bm{\mu}+\varepsilon_k\bm{t})
		\to\overline{\bm{w}}\in\Omega.
	\end{equation*}
	Optimality at $\bm{\mu}+\varepsilon_k\bm{t}$ gives
	\begin{equation*}
		\bm{\mu}^\top\bm{w}_*(\bm{\mu}+\varepsilon_k\bm{t})
		=f_*(\bm{\mu}+\varepsilon_k\bm{t})
		-\varepsilon_k\bm{t}^\top
		\bm{w}_*(\bm{\mu}+\varepsilon_k\bm{t})\to f_*(\bm{\mu}),
	\end{equation*}
	where the first term converges by continuity of $f_*$ and the second term tends to zero because $\Omega$ is bounded. Thus
	$\bm{\mu}^\top\overline{\bm{w}}=f_*(\bm{\mu})$, so
	$\overline{\bm{w}}\in \mathcal W_*(\bm{\mu})$. This argument applies to every
	convergent subsequence of the perturbed maximizers. In particular,
	\begin{equation*}
		\limsup_{\varepsilon\downarrow0}
		\bm{w}_*(\bm{\mu}+\varepsilon\bm{t})^\top\bm{t}
		\le\sup_{\bm{w}\in \mathcal W_*(\bm{\mu})}\bm{w}^\top\bm{t}.
	\end{equation*}
	Indeed, choose a sequence realizing the limsup and extract a convergent
	subsequence of its maximizers; its limit belongs to $\mathcal W_*(\bm{\mu})$.
	Combining this inequality with \eqref{eq:two_bounds}, we have
	\begin{equation*}
		\begin{aligned}
			\sup_{\bm{w}\in \mathcal W_*(\bm{\mu})}\bm{w}^\top\bm{t}&\le \liminf_{\varepsilon\downarrow0}\frac{f_*(\bm{\mu}+\varepsilon\bm{t})-f_*(\bm{\mu})}{\varepsilon}\le\limsup_{\varepsilon\downarrow0}\frac{f_*(\bm{\mu}+\varepsilon\bm{t})-f_*(\bm{\mu})}{\varepsilon}\\&\le\limsup_{\varepsilon\downarrow0}
			\bm{w}_*(\bm{\mu}+\varepsilon\bm{t})^\top\bm{t}
			\le\sup_{\bm{w}\in \mathcal W_*(\bm{\mu})}\bm{w}^\top\bm{t},
		\end{aligned}
	\end{equation*}
	and thus 
	\begin{equation*}
		\lim_{\varepsilon\downarrow0}\frac{f_*(\bm{\mu}+\varepsilon\bm{t})-f_*(\bm{\mu})}{\varepsilon}=\sup_{\bm{w}\in \mathcal W_*(\bm{\mu})}\bm{w}^\top\bm{t}.
	\end{equation*}
	Subtracting
	$\bm{w}_*(\bm{\mu})^\top\bm{t}$ proves
	\begin{equation*}
		\lim_{\varepsilon\downarrow0}
		\frac{r_{\bm{\mu}}(\varepsilon\bm{t})}{\varepsilon}
		=\lim_{\varepsilon\downarrow0}\frac{f_*(\bm{\mu}+\varepsilon\bm{t})-f_*(\bm{\mu})}{\varepsilon}-\bm{w}_*(\bm{\mu})^\top\bm{t}=\sup_{\bm{w}\in \mathcal W_*(\bm{\mu})}\bm{w}^\top\bm{t}-\bm{w}_*(\bm{\mu})^\top\bm{t}=\psi_1(\bm{t}).
	\end{equation*}
	Proposition~\ref{pr:pointtouniform} makes this convergence uniform on each $B_2(M)$, and Theorem \ref{thm:fast_general_local_rate} gives \eqref{eq:fast_first_order_minimax}. If $\mathcal W_*(\bm{\mu})$ is a singleton, then $\psi_1\equiv0$. Otherwise, choosing $\bm{w}\in \mathcal W_*(\bm{\mu})\setminus\{\bm{w}_*(\bm{\mu})\}$ gives $\psi_1(\bm{w}-\bm{w}_*(\bm{\mu}))>0$.
	Since $\psi_1$ is nonnegative, vanishes at zero, and is positively homogeneous, it cannot be affine on any $B_2(M)$ with $M>0$. The positivity assertion follows from the theorem. 

\end{proof}
\subsection{Proof of Proposition~\ref{prop:fast_symmetric_difference}}
\label{app:proof-finite-modulus}
\begin{proof}
	By the convexity of $f_*$,
	\begin{equation*}
		\frac{f_*(\bm{\mu}+\varepsilon\bm{t})+f_*(\bm{\mu}-\varepsilon\bm{t})}{2}\ge f_*(\bm{\mu}),
	\end{equation*}
	which yields $D_{\bm{\mu}}(\varepsilon,M)\ge0$. By continuity of $f_*$, since $B_2(M)$ is compact, a maximizer exists. Since $r_{\bm{\mu}}\ge0$,
	\begin{equation*}
		0\le r_{\bm{\mu}}(\varepsilon\bm{t})
		\le r_{\bm{\mu}}(\varepsilon\bm{t})
		+r_{\bm{\mu}}(-\varepsilon\bm{t})
		\le2D_{\bm{\mu}}(\varepsilon,M),\qquad
		\bm{t}\in B_2(M).
	\end{equation*}
	Lemma \ref{lem:fast_bayes_minimax} gives
	\begin{equation*}
		R_{n,M}^{\Sigma}(\bm{\mu})\le\sup_{\pi\in\mathcal{P}(B_2(M))}\mathbb{E}_{\pi,\Sigma} r_{\bm{\mu}}(\bm{T}/\sqrt n)\le\sup_{\bm{t}\in B_2(M)}r_{\bm{\mu}}(\bm{t}/\sqrt n)\le2D_{\bm{\mu}}(n^{-1/2},M).
	\end{equation*}
	This proves the upper bound. 
	
	For the lower bound, suppose $D_{\bm{\mu}}(n^{-1/2},M)>0$ and choose
	a maximizer $\bm{t}_0\in B_2(M)$ in
	\eqref{eq:fast_symmetric_difference_definition}. Then $\bm{t}_0\ne0$.
	Use the prior
	\begin{equation*}
		\mathbb{P}_\pi\{\bm{T}=\bm{t}_0\}
		=\mathbb{P}_\pi\{\bm{T}=-\bm{t}_0\}=\frac12,\qquad
		\bm{Y}=\bm{T}+\bm{Z},\qquad
		q(\bm{Y}):=\mathbb{P}_{\pi,\Sigma}\{\bm{T}=\bm{t}_0\mid\bm{Y}\}.
	\end{equation*}
	For $q\in[0,1]$, set
	\begin{equation*}
		\begin{aligned}
			g(q)&:=q r_{\bm{\mu}}(\bm{t}_0/\sqrt n)
			+(1-q)r_{\bm{\mu}}(-\bm{t}_0/\sqrt n)
			-r_{\bm{\mu}}((2q-1)\bm{t}_0/\sqrt n),\\
			g(0)&=g(1)=0,\qquad
			g(1/2)=D_{\bm{\mu}}(n^{-1/2},M).
		\end{aligned}
	\end{equation*}
	Convexity of $r_{\bm{\mu}}$ makes $g$ concave. Applying concavity on $[0,1/2]$ and $[1/2,1]$ gives
	\begin{equation*}
		g(q)\ge2\min\{q,1-q\}D_{\bm{\mu}}(n^{-1/2},M).
	\end{equation*}
	The likelihood ratio between $\mathcal{N}(\bm{t}_0,\Sigma)$ and $\mathcal{N}(-\bm{t}_0,\Sigma)$ is $\exp(2\bm{t}_0^\top\Sigma^{-1}\bm{Y})$. Thus the Bayes-optimal test chooses $\bm{t}_0$ when $\bm{t}_0^\top\Sigma^{-1}\bm{Y}>0$, and
	\begin{equation*}
		\begin{aligned}
			\mathbb{E}_{\pi,\Sigma}\min\{q(\bm{Y}),1-q(\bm{Y})\}
			&=\frac12\mathbb{P}_{\pi,\Sigma}\{\bm{t}_0^\top\Sigma^{-1}\bm{Y}\le0\mid\bm T=\bm t_0\}
			+\frac12\mathbb{P}_{\pi,\Sigma}\{\bm{t}_0^\top\Sigma^{-1}\bm{Y}>0\mid\bm T=-\bm t_0\}\\&=\Phi(-\|\Sigma^{-1/2}\bm{t}_0\|_2)\ge\Phi(-M/\sqrt{\lambda_{\min}(\Sigma)}).
		\end{aligned}
	\end{equation*}
	Since
	\[
	\bm m_{\pi,\Sigma}(\bm Y)=(2q(\bm Y)-1)\bm t_0
	\]
	and
	\[
	\begin{aligned}
		\mathbb E_{\pi,\Sigma}
		[r_{\bm\mu}(\bm T/\sqrt n)\mid\bm Y]
		&=q(\bm Y)r_{\bm\mu}(\bm t_0/\sqrt n)\\
		&\quad+(1-q(\bm Y))r_{\bm\mu}(-\bm t_0/\sqrt n),
	\end{aligned}
	\]
	Lemma~\ref{lem:fast_bayes_minimax} yields
	\begin{equation*}
		\begin{aligned}
			R_{n,M}^{\Sigma}(\bm{\mu})&\ge\mathbb E_{\pi,\Sigma}\left[r_{\bm\mu}(\bm T/\sqrt n)-r_{\bm\mu}(\bm m_{\pi,\Sigma}(\bm Y)/\sqrt n)\right]\\&=\mathbb{E}_{\pi,\Sigma} g(q(\bm{Y}))
			\ge2\Phi(-M/\sqrt{\lambda_{\min}(\Sigma)})D_{\bm{\mu}}(n^{-1/2},M).
		\end{aligned}
	\end{equation*}
	If $D_{\bm{\mu}}(n^{-1/2},M)=0$, the result follows directly from nonnegativity of regret.
	
	Since $\lambda_{\min}(\Sigma)\ge\underline\lambda$, the two bounds give, for every $\Sigma\in\mathcal C$,
	\begin{equation*}
		2\Phi(-M/\sqrt{\underline\lambda})D_{\bm\mu}(n^{-1/2},M)
		\le R_{n,M}^{\Sigma}(\bm\mu)\le2D_{\bm\mu}(n^{-1/2},M),
	\end{equation*}
	so that the comparison factor between the two sides is finite, independent of $n$ and $\Sigma$, for fixed $M$ and $\underline\lambda>0$. The inequalities also apply when both risks vanish.
\end{proof}
\subsection{Distance normalization for general scales}
\begin{corollary}[Distance normalization under a general local scale]
	\label{cor:fast_distance_general}
	Suppose that $\Omega$ contains at least two points, $\|\cdot\|_K$ is a fixed norm, and $\dist(\bm{\mu},A)>0$. Use $\bm{\nu}_{\bm{\mu}}=\bm{\mu}/\dist(\bm{\mu},A)$ as in \eqref{eq:fast_distance_normalization}. For every $\varepsilon,M>0$,
	\begin{equation}
		\begin{aligned}
			r_{\bm{\mu}}(\varepsilon\bm{t})
			&=\dist(\bm{\mu},A)
			r_{\bm{\nu}_{\bm{\mu}}}
			\left(\frac{\varepsilon}{\dist(\bm{\mu},A)}\bm{t}\right),\\
			D_{\bm{\mu}}(\varepsilon,M)
			&=\dist(\bm{\mu},A)
			D_{\bm{\nu}_{\bm{\mu}}}
			\left(\frac{\varepsilon}{\dist(\bm{\mu},A)},M\right).
		\end{aligned}
		\label{eq:fast_distance_remainder_scaling}
	\end{equation}
	In particular, the finite-sample bounds become
	\begin{equation}
		\begin{aligned}
			2\Phi(-M/\sqrt{\lambda_{\min}(\Sigma)})\dist(\bm{\mu},A)
			D_{\bm{\nu}_{\bm{\mu}}}
			\left(\frac{1}{\sqrt n \dist(\bm{\mu},A)},M\right)
			&\le R_{n,M}^{\Sigma}(\bm{\mu})\\
			&\le2\dist(\bm{\mu},A)
			D_{\bm{\nu}_{\bm{\mu}}}
			\left(\frac{1}{\sqrt n \dist(\bm{\mu},A)},M\right).
		\end{aligned}
		\label{eq:fast_distance_modulus_bounds}
	\end{equation}
	Suppose, in addition, that there exist $a(\varepsilon)>0$ with
	$a(\varepsilon)\to0$ and a finite continuous convex
	$\psi_{\bm{\nu}_{\bm{\mu}}}:\mathbb{R}^d\to\mathbb{R}$ such that,
	for every fixed $M>0$,
	\begin{equation}
		\sup_{\bm{t}\in B_2(M)}
		\left|\frac{r_{\bm{\nu}_{\bm{\mu}}}(\varepsilon\bm{t})}
		{a(\varepsilon)}
		-\psi_{\bm{\nu}_{\bm{\mu}}}(\bm{t})\right|
		\to0\qquad\text{as }\varepsilon\downarrow0.
		\label{eq:fast_distance_normalized_expansion}
	\end{equation}
	For fixed $d,\Omega,\bm{\mu},M$,
	\begin{equation}
		\lim_{n\to\infty}
		\frac{R_{n,M}^{\Sigma}(\bm{\mu})}
		{\dist(\bm{\mu},A)a((\sqrt n \dist(\bm{\mu},A))^{-1})}
		=V_{M,\Sigma}(\psi_{\bm{\nu}_{\bm{\mu}}}).
		\label{eq:fast_distance_general_limit}
	\end{equation}
	If $\psi_{\bm{\nu}_{\bm{\mu}}}$ is nonaffine on $B_2(M)$, then
	$V_{M,\Sigma}(\psi_{\bm{\nu}_{\bm{\mu}}})>0$, and the plug-in policy attains
	the corresponding rate. In particular, for $a(\varepsilon)=\varepsilon^\alpha$,
	\begin{equation}
		\lim_{n\to\infty}
		n^{\alpha/2}R_{n,M}^{\Sigma}(\bm{\mu})
		=\dist(\bm{\mu},A)^{1-\alpha}
		V_{M,\Sigma}(\psi_{\bm{\nu}_{\bm{\mu}}}).
		\label{eq:fast_distance_power_limit}
	\end{equation}
\end{corollary}

\begin{proof}
	Positive homogeneity gives $\dist(\bm{\nu}_{\bm{\mu}},A)=1$ and
	$\bm{w}_*(\bm{\nu}_{\bm{\mu}})=\bm{w}_*(\bm{\mu})$. Hence
	\begin{equation*}
		\begin{aligned}
			r_{\bm{\mu}}(\varepsilon\bm{t})
			&=\dist(\bm{\mu},A)
			\left[f_*\left(\bm{\nu}_{\bm{\mu}}
			+\frac{\varepsilon\bm{t}}{\dist(\bm{\mu},A)}\right)
			-f_*(\bm{\nu}_{\bm{\mu}})
			-\frac{\varepsilon\bm{w}_*(\bm{\nu}_{\bm{\mu}})^\top\bm{t}}
			{\dist(\bm{\mu},A)}\right]\\
			&=\dist(\bm{\mu},A)r_{\bm{\nu}_{\bm{\mu}}}
			(\varepsilon\bm{t}/\dist(\bm{\mu},A)).
		\end{aligned}
	\end{equation*}
	Taking symmetric differences proves \eqref{eq:fast_distance_remainder_scaling}, and \eqref{eq:fast_symmetric_difference_bounds} yields \eqref{eq:fast_distance_modulus_bounds}.
	
	By \eqref{eq:fast_distance_remainder_scaling} and
	\eqref{eq:fast_distance_normalized_expansion}, for every fixed $M>0$,
	\begin{equation*}
		\sup_{\bm{t}\in B_2(M)}
		\left|\frac{r_{\bm{\mu}}(\varepsilon\bm{t})}
		{\dist(\bm{\mu},A)a(\varepsilon/\dist(\bm{\mu},A))}
		-\psi_{\bm{\nu}_{\bm{\mu}}}(\bm{t})\right|
		\to0.
	\end{equation*}
	Apply Theorem \ref{thm:fast_general_local_rate} with scale
	$\varepsilon\mapsto \dist(\bm{\mu},A)a(\varepsilon/\dist(\bm{\mu},A))$
	to obtain \eqref{eq:fast_distance_general_limit}, positivity, and
	attainment of the rate by the plug-in policy.
	For $a(\varepsilon)=\varepsilon^\alpha$,
	$\dist(\bm{\mu},A)a((\sqrt n \dist(\bm{\mu},A))^{-1})
	=\dist(\bm{\mu},A)^{1-\alpha}n^{-\alpha/2}$,
	which proves \eqref{eq:fast_distance_power_limit}.
\end{proof}

\begin{remark}[What the distance does and does not determine]
	\label{rem:fast_distance_scope}
	Let $C_K:=\sup_{\|\bm{t}\|_2\le1}\|\bm{t}\|_K$.
	The distance is $1$-Lipschitz in its defining norm, so
	\begin{equation*}
		\sup_{\bm{t}\in B_2(M)}
		|\dist(\bm{\mu}_{n,\bm{t}},A)-\dist(\bm{\mu},A)|
		\le\frac{C_K M}{\sqrt n}.
	\end{equation*}
	For fixed $\dist(\bm{\mu},A)>0$, the relative variation tends to zero, but the normalized curvature in \eqref{eq:fast_distance_sharp_constant} and the limiting function in \eqref{eq:fast_distance_general_limit} still depend on the direction of $\bm{\mu}$ and on $\Omega$.
	
	Without a local expansion, \eqref{eq:fast_distance_modulus_bounds}
	provides the following criterion for each fixed $M>0$:
	\begin{equation*}
		\liminf_{\varepsilon\downarrow0}
		\frac{D_{\bm{\nu}_{\bm{\mu}}}(\varepsilon,M)}{\varepsilon^2}>0
		\quad\Longrightarrow\quad
		\liminf_{n\to\infty}n\,\dist(\bm{\mu},A)R_{n,M}^{\Sigma}(\bm{\mu})>0.
	\end{equation*}
	If this ratio tends to zero, then
	$R_{n,M}^{\Sigma}(\bm{\mu})=o(1/[n\,\dist(\bm{\mu},A)])$.
	If its liminf is positive and its limsup is finite, then
	$R_{n,M}^{\Sigma}(\bm{\mu})\asymp1/[n\,\dist(\bm{\mu},A)]$ for fixed
	$d,\Omega,\bm{\mu},M$. 
	
	When $\dist(\bm{\mu},A)=0$, reciprocal-distance formulas do not apply.
	Corollary \ref{cor:fast_first_order_limit} instead gives
	\begin{equation*}
		\lim_{n\to\infty}\sqrt n R_{n,M}^{\Sigma}(\bm{\mu})=V_{M,\Sigma}(\psi_1),
		\qquad V_{M,\Sigma}(\psi_1)>0\ \Longleftrightarrow\ \bm{\mu}\in A,
		\quad M>0.
	\end{equation*}
	In particular, $\bm{\mu}\in\overline A\setminus A$ implies
	$R_{n,M}^{\Sigma}(\bm{\mu})=o(n^{-1/2})$.
	
	The lower spectral bound is substantive. With singular covariance, $\operatorname{tr}(\Sigma H_{\bm\mu})$ can vanish even when $H_{\bm\mu}\ne0$, and the full-support Gaussian testing argument for positive $V_{M,\Sigma}$ can fail. 
\end{remark}

\subsection{Finite-sample plug-in bounds}
\begin{lemma}[A finite-sample plug-in bound through symmetric differences]
	\label{lem:fast_plugin_modulus_upper}
	Let $\Omega$ be compact and convex, let $\bm\mu\in\mathbb R^d$, and let $M>0$. For $\bm G\sim\mathcal N(\bm0,I_d)$, every $n\ge2$, and every $\Sigma\in\mathcal C$,
	\begin{equation}
		\begin{aligned}
			\sup_{\bm t\in B_2(M)}
			\mathbb E_{\bm\mu_{n,\bm t},\Sigma}R(\bm w_*(\bar{\bm x}))
			&\le2\mathbb E D_{\bm\mu}\!\left(n^{-1/2},M+2\|\Sigma^{1/2}\bm G\|_2\right)\\
			&\le2\mathbb E D_{\bm\mu}\!\left(n^{-1/2},M+2\sqrt{\overline\lambda}\|\bm G\|_2\right).
		\end{aligned}
		\label{eq:fast_plugin_integrated_modulus}
	\end{equation}
	The last expression also bounds the supremum over $\Sigma\in\mathcal C$ on the left. 
\end{lemma}

\begin{proof}
	Denote $\bm Z=\Sigma^{1/2}\bm G$, so $\bar{\bm x}=\bm\mu+(\bm t+\bm Z)/\sqrt{n}$. Optimality at $\bar{\bm x}$ gives
	\begin{equation*}
		R(\bm w_*(\bar{\bm x}))
		=r_{\bm\mu}(\bm t/\sqrt{n})-r_{\bm\mu}((\bm t+\bm Z)/\sqrt{n})
		+(\bm w_*(\bar{\bm x})-\bm w_*(\bm\mu))^{\top}\bm Z/\sqrt{n},
	\end{equation*}
	and the subgradient inequality for $r_{\bm\mu}$ gives
	\begin{equation*}
		r_{\bm\mu}((\bm t+2\bm Z)/\sqrt{n})\ge r_{\bm\mu}((\bm t+\bm Z)/\sqrt{n})+(\bm w_*(\bar{\bm x})-\bm w_*(\bm\mu))^{\top}\bm Z/\sqrt{n}.
	\end{equation*}
	Consequently,
	\begin{equation}
		R(\bm w_*(\bar{\bm x}))
		\le r_{\bm\mu}(\bm t/\sqrt{n})
		+r_{\bm\mu}((\bm t+2\bm Z)/\sqrt{n})
		-2r_{\bm\mu}((\bm t+\bm Z)/\sqrt{n}).
		\label{eq:plug_in_upper_bound_via_perterb}
	\end{equation}
	By Jensen's inequality, convexity of $r_{\bm{\mu}}$ and $\mathbb E\bm Z=\bm0$ yield
	\begin{equation*}
		\mathbb E r_{\bm\mu}((\bm t+\bm Z)/\sqrt{n})
		\ge r_{\bm\mu}(\bm t/\sqrt{n}),
		\qquad r_{\bm\mu}(\bm t/\sqrt{n})\ge0.
	\end{equation*}
	Taking expectations in \eqref{eq:plug_in_upper_bound_via_perterb} therefore gives
	\begin{equation*}
		\begin{aligned}
			\mathbb E_{\bm\mu_{n,\bm t},\Sigma}R(\bm w_*(\bar{\bm x}))
			&\le\mathbb E r_{\bm\mu}((\bm t+2\bm Z)/\sqrt{n})-r_{\bm\mu}(\bm t/\sqrt{n})\le\mathbb E r_{\bm\mu}((\bm t+2\bm Z)/\sqrt{n})\\&\le2\mathbb E D_{\bm\mu}(n^{-1/2},M+2\|\bm Z\|_2)
		\end{aligned}
	\end{equation*}
	for $\bm t\in B_2(M)$. The last inequality uses $0\le r_{\bm\mu}(\varepsilon\bm u)\le2D_{\bm\mu}(\varepsilon,L)$ for $\bm u\in B_2(L)$, as proved in Proposition~\ref{prop:fast_symmetric_difference}. All expectations are finite because $r_{\bm\mu}(\bm h)\le2R_\Omega\|\bm h\|_2$ and, by Lipschitz continuity of $f_*$,
	\begin{equation*}
		0\le D_{\bm\mu}(\varepsilon,L)\le R_\Omega\varepsilon L.
	\end{equation*}
	Finally, $D_{\bm\mu}(\varepsilon,L)$ is nondecreasing in $L$, and $\|\Sigma^{1/2}\bm G\|_2\le\sqrt{\overline\lambda}\|\bm G\|_2$. Taking suprema over $\bm t$ proves \eqref{eq:fast_plugin_integrated_modulus} and its uniform-covariance assertion.
\end{proof}

\begin{lemma}[Rate minimaxity under a local doubling condition]
	\label{cor:fast_plugin_modulus_rate}
	Retain the assumptions of Lemma~\ref{lem:fast_plugin_modulus_upper}. Suppose that, for some $\varepsilon_0>0$ and $\Gamma\ge2$,
	\begin{equation}
		0<D_{\bm\mu}(\varepsilon,M),
		\qquad D_{\bm\mu}(2\varepsilon,M)\le\Gamma D_{\bm\mu}(\varepsilon,M),
		\qquad 0<\varepsilon\le\varepsilon_0.
		\label{eq:fast_modulus_local_doubling}
	\end{equation}
	Then
	\begin{equation}
		\limsup_{n\to\infty}
		\frac{\displaystyle\sup_{\substack{\bm t\in B_2(M)\\\Sigma\in\mathcal C}}
			\mathbb E_{\bm\mu_{n,\bm t},\Sigma}R(\bm w_*(\bar{\bm x}))}
		{D_{\bm\mu}(n^{-1/2},M)}
		\le2\Gamma\mathbb E\left(1+\frac{2\sqrt{\overline\lambda}\|\bm G\|_2}{M}\right)^{\log_2\Gamma}<\infty.
		\label{eq:fast_plugin_modulus_ratio}
	\end{equation}
	In particular, as $n\to\infty$ with $M$ fixed,
	\begin{equation}
		\begin{aligned}
			R_{n,M}^{\Sigma}(\bm\mu)
			\asymp D_{\bm\mu}(n^{-1/2},M)
			\asymp\sup_{\bm t\in B_2(M)}
			\mathbb E_{\bm\mu_{n,\bm t},\Sigma}R(\bm w_*(\bar{\bm x})).
		\end{aligned}
		\label{eq:fast_plugin_modulus_rate}
	\end{equation}
	The comparison constants are uniform over $\Sigma\in\mathcal C$. 
\end{lemma}

\begin{proof}
	Changing variables in \eqref{eq:fast_symmetric_difference_definition} gives
	\begin{equation*}
		D_{\bm\mu}(\varepsilon,L)
		=D_{\bm\mu}(\varepsilon L/M,M),\qquad \varepsilon,L>0.
	\end{equation*}
	The function $\varepsilon\mapsto D_{\bm\mu}(\varepsilon,M)$ is convex and nondecreasing, and vanishes at zero: for every fixed $\bm t$, the symmetric difference is an even convex function of $\varepsilon$, and taking a supremum preserves these properties.
	
	For $k\ge1$ and $k\varepsilon\le\varepsilon_0$, applying \eqref{eq:fast_modulus_local_doubling} at $\varepsilon,2\varepsilon,\ldots,2^{\lceil\log_2k\rceil-1}\varepsilon$ gives
	\begin{equation*}
		D_{\bm\mu}(k\varepsilon,M)
		\le\Gamma^{\lceil\log_2k\rceil}D_{\bm\mu}(\varepsilon,M)
		\le\Gamma k^{\log_2\Gamma}D_{\bm\mu}(\varepsilon,M).
	\end{equation*}
	Let
	\[
	H(\bm G):=1+\frac{2\sqrt{\overline\lambda}\|\bm G\|_2}{M}.
	\]
	Applying this estimate on $\{\varepsilon H(\bm G)\le\varepsilon_0\}$ and the global bound $D_{\bm\mu}(\varepsilon,L)\le R_\Omega\varepsilon L$ on its complement gives
	\begin{equation*}
		\begin{aligned}
			\mathbb E D_{\bm\mu}(\varepsilon,M+2\sqrt{\overline\lambda}\|\bm G\|_2)
			&\le\Gamma D_{\bm\mu}(\varepsilon,M)
			\mathbb E[H(\bm G)^{\log_2\Gamma}]\\
			&\quad+R_\Omega M\varepsilon\,
			\mathbb E[H(\bm G)\mathbb I\{\varepsilon H(\bm G)>\varepsilon_0\}].
		\end{aligned}
	\end{equation*}
	To compare the last term with $D_{\bm\mu}(\varepsilon,M)$, take $c_{\varepsilon}=\lfloor\log_2(\varepsilon_0/\varepsilon)\rfloor$ for $0<\varepsilon\le\varepsilon_0/2$. Then $\varepsilon_0/2<2^{c_{\varepsilon}}\varepsilon\le\varepsilon_0$, and the same doubling condition gives
	\begin{equation*}
		D_{\bm\mu}(\varepsilon,M)
		\ge\Gamma^{-c_{\varepsilon}}D_{\bm\mu}(2^{c_{\varepsilon}}\varepsilon,M)
		\ge D_{\bm\mu}(\varepsilon_0/2,M)
		(\varepsilon/\varepsilon_0)^{\log_2\Gamma}>0.
	\end{equation*}
	For these $\varepsilon$, the exceptional event implies
	$\|\bm G\|_2>M\varepsilon_0/(4\sqrt{\overline\lambda}\varepsilon)$. Cauchy--Schwarz and the Gaussian exponential-moment bound imply
	\begin{equation*}
		\begin{aligned}
			&\quad\mathbb E\left[\left(1+\frac{2\sqrt{\overline\lambda}\|\bm G\|_2}{M}\right)
			\mathbb I\left\{\varepsilon\left(1+\frac{2\sqrt{\overline\lambda}\|\bm G\|_2}{M}\right)>\varepsilon_0\right\}\right]\\
			&\le
			\left[\mathbb E\left(1+\frac{2\sqrt{\overline\lambda}\|\bm G\|_2}{M}\right)^2\right]^{1/2}
			2^{d/4}\exp\left(-\frac{M^2\varepsilon_0^2}{128\overline\lambda\varepsilon^2}\right).
		\end{aligned}
	\end{equation*}
	After division by $D_{\bm\mu}(\varepsilon,M)$, the exceptional contribution is bounded by
	\begin{equation*}
		C\varepsilon^{1-\log_2\Gamma}
		\exp\left(-\frac{M^2\varepsilon_0^2}{128\overline\lambda\varepsilon^2}\right)
		\to0.
	\end{equation*}
	Lemma~\ref{lem:fast_plugin_modulus_upper} now proves \eqref{eq:fast_plugin_modulus_ratio}. Finally, Proposition~\ref{prop:fast_symmetric_difference} gives the uniform positive lower factor $2\Phi(-M/\sqrt{\underline\lambda})$, while the minimax definition compares $R_{n,M}^{\Sigma}(\bm\mu)$ with the risk of the same sample-only policy. These inequalities prove \eqref{eq:fast_plugin_modulus_rate}.
\end{proof}

\section{Convex-Geometric Proofs}
\label{app:boundary-geometry}
\subsection{Proof of Corollary~\ref{cor:fast_distance_smooth}}
\label{app:proof-distance-smooth}
\begin{proof}
	For every $b>0$ and $\bm{x}\in\mathbb{R}^d$,
	\begin{equation*}
		f_*(b\bm{x})=b f_*(\bm{x}),\qquad
		\arg\max_{\bm{w}\in\Omega}\bm{w}^\top(b\bm{x})
		=\arg\max_{\bm{w}\in\Omega}\bm{w}^\top\bm{x}.
	\end{equation*}
	It follows that $bA=A$. Also $\bm{0}\in A$, since $\Omega$ is
	not a singleton. Therefore
	\begin{equation}
		\dist(b\bm{x},A)
		=\inf_{\bm{y}\in A}\|b\bm{x}-\bm{y}\|_K
		=b\inf_{\bm{y}\in A}\|\bm{x}-\bm{y}\|_K
		=b \dist(\bm{x},A).
		\label{eq:fast_distance_homogeneity}
	\end{equation}
	We first verify $\partial f_*(\bm{x})=\arg\max_{\bm{w}\in\Omega}\bm{w}^\top\bm{x}$. If $\bm{w}\in\Omega$ maximizes $\bm{w}^\top\bm{x}$, then
	\begin{equation*}
		f_*(\bm{y})\ge\bm{w}^\top\bm{y}
		=f_*(\bm{x})+\bm{w}^\top(\bm{y}-\bm{x})
		\qquad\forall\bm{y}\in\mathbb{R}^d,
	\end{equation*}
	so $\bm{w}\in\partial f_*(\bm{x})$.
	Conversely, if $\bm{w}\in\partial f_*(\bm{x})$, substituting
	$\bm{y}=\bm{0}$ and $\bm{y}=2\bm{x}$ into the subgradient
	inequality yields
	\begin{equation*}
		0\ge f_*(\bm{x})-\bm{w}^\top\bm{x},\qquad
		2f_*(\bm{x})\ge f_*(\bm{x})+\bm{w}^\top\bm{x}.
	\end{equation*}
	Hence $\bm{w}^\top\bm{x}=f_*(\bm{x})$, and the same inequality gives
	\begin{equation*}
		\bm{w}^\top\bm{y}\le f_*(\bm{y})
		=\sup_{\bm{v}\in\Omega}\bm{v}^\top\bm{y}
		\qquad\forall\bm{y}\in\mathbb{R}^d.
	\end{equation*}
	Since $\Omega$ is closed and convex, by the separating hyperplane theorem, if $\bm{w}\notin\Omega$, there exists $\bm{y}\in\mathbb{R}^d$, such that $\bm{w}^\top\bm{y}> \sup_{\bm{v}\in\Omega}\bm{v}^\top\bm{y}$. Hence $\bm{w}\in\Omega$, and $\partial f_*(\bm{x})=\arg\max_{\bm{w}\in\Omega}\bm{w}^\top\bm{x}$.
	
	A finite convex function on $\mathbb{R}^d$ is differentiable at
	$\bm{x}$ if and only if its subdifferential is a singleton. Therefore
	\begin{equation*}
		f_*\text{ is differentiable at }\bm{x}
		\quad\Longleftrightarrow\quad
		\partial f_*(\bm{x})=\{\nabla f_*(\bm{x})\}
		\quad\Longleftrightarrow\quad \bm{x}\notin A.
	\end{equation*}
	Differentiability throughout a neighborhood of $\bm{\mu}$ gives
	$r>0$ such that
	\begin{equation*}
		\{\bm{x}:\|\bm{x}-\bm{\mu}\|_2<r\}\cap A=\varnothing.
	\end{equation*}
	By equivalence of norms, there is $c>0$ such that
	$\|\bm{v}\|_K\ge c\|\bm{v}\|_2$ for all $\bm{v}$, whence
	\begin{equation*}
		\dist(\bm{\mu},A)
		=\inf_{\bm{a}\in A}\|\bm{\mu}-\bm{a}\|_K
		\ge c\inf_{\bm{a}\in A}\|\bm{\mu}-\bm{a}\|_2
		\ge cr>0.
	\end{equation*}
	Also $\dist(\bm{\mu},A)\le\|\bm{\mu}\|_K<\infty$ since $\bm{0}\in A$.
	Thus $\bm{\nu}_{\bm{\mu}}$ is well-defined, and
	\eqref{eq:fast_distance_homogeneity} gives
	\begin{equation*}
		\dist(\bm{\nu}_{\bm{\mu}},A)
		=\frac{\dist(\bm{\mu},A)}{\dist(\bm{\mu},A)}=1.
	\end{equation*}
	
	Let $U$ be an open neighborhood of $\bm{\mu}$ on which $f_*\in C^2$.
	For each fixed $b>0$, the identity
	$f_*(\bm{y})=b f_*(\bm{y}/b)$ for $\bm{y}\in bU$
	implies $f_*\in C^2(bU)$ as well.
	Differentiating $f_*(b\bm{x})=b f_*(\bm{x})$ with respect to
	$\bm{x}\in U$, with $b$ fixed, the chain rule gives
	\begin{equation*}
		\begin{aligned}
			b\nabla f_*(b\bm{x})
			&=\nabla_{\bm{x}}[f_*(b\bm{x})]
			=\nabla_{\bm{x}}[b f_*(\bm{x})]
			=b\nabla f_*(\bm{x}),\\
			b^2\nabla^2 f_*(b\bm{x})
			&=\nabla_{\bm{x}}^2[f_*(b\bm{x})]
			=\nabla_{\bm{x}}^2[b f_*(\bm{x})]
			=b\nabla^2 f_*(\bm{x}).
		\end{aligned}
	\end{equation*}
	Dividing by $b$ and $b^2$, respectively, yields
	\begin{equation*}
		\nabla f_*(b\bm{x})=\nabla f_*(\bm{x}),\qquad
		\nabla^2f_*(b\bm{x})=b^{-1}\nabla^2f_*(\bm{x}).
	\end{equation*}
	Taking $\bm{x}=\bm{\mu}$ and $b=1/\dist(\bm{\mu},A)$ gives
	$H_{\bm{\nu}_{\bm{\mu}}}=\dist(\bm{\mu},A)H_{\bm{\mu}}$.
	Since $\mathcal{B}_{M,\Sigma}(bH)=b\mathcal{B}_{M,\Sigma}(H)$ for $b\ge0$,
	Theorem \ref{thm:fast_local_minimax} proves \eqref{eq:fast_distance_fixed_M} and \eqref{eq:fast_distance_sharp_constant}.
\end{proof}

\subsection{Proof of Corollary~\ref{cor:fast_principal_curvatures}}
\label{app:proof-principal-curvature}
\begin{proof}
	Since $\bm{\mu}\ne\bm{0}$, if $\bm{w}_*(\bm{\mu})\in\operatorname{int}\Omega$, then $\bm{w}_*(\bm{\mu})+t\bm{n}_{\bm{\mu}}\in\Omega$ for some $t>0$, and $\bm{\mu}^\top(\bm{w}_*(\bm{\mu})+t\bm{n}_{\bm{\mu}})=f_*(\bm{\mu})+t\|\bm{\mu}\|_2>f_*(\bm{\mu})$, which contradicts the definition of $\bm{w}_*(\bm{\mu})$. Hence $\bm{w}_*(\bm{\mu})\in\partial\Omega$. 
	
	By optimality of $\bm{w}_*(\bm{\mu})$, $\bm{\mu}^\top(\bm{w}-\bm{w}_*(\bm{\mu}))\le0$ for all $\bm{w}\in\Omega$. Thus the support hyperplane of $\Omega$ at $\bm{w}_*(\bm{\mu})$ is $\{\bm{w}:\bm{\mu}^\top\bm{w}=f_*(\bm{\mu})\}$, with orthogonal unit normal vector $\bm{N}(\bm{w}_*(\bm{\mu}))=\bm{n}_{\bm{\mu}}$.
	
	To justify the local graph representation, the $C^2$ hypersurface assumption provides an open neighborhood $W$ of $\bm{w}_*(\bm{\mu})$ and a defining function $F\in C^2(W)$ such that
	\begin{equation*}
		\partial\Omega\cap W=\{\bm{w}\in W:F(\bm{w})=0\},
		\qquad \nabla F(\bm{w}_*(\bm{\mu}))\ne\bm{0}.
	\end{equation*}
	Since this gradient is normal to the boundary, multiplying $F$ by a nonzero constant allows us to assume
	$\nabla F(\bm{w}_*(\bm{\mu}))=\bm{n}_{\bm{\mu}}$.
	In tangent and normal coordinates, define
	\begin{equation*}
		G(\bm{z},a)
		:=F\bigl(\bm{w}_*(\bm{\mu})+U_{\bm{\mu}}\bm{z}
		+\bm{n}_{\bm{\mu}}a\bigr),
		\qquad (\bm{z},a)\in\mathbb{R}^{d-1}\times\mathbb{R}
		\text{ near }(\bm{0},0).
	\end{equation*}
	The chain rule gives
	\begin{equation*}
		\begin{aligned}
			G(\bm{0},0)&=0,\\
			\partial_aG(\bm{0},0)
			&=\nabla F(\bm{w}_*(\bm{\mu}))^\top\bm{n}_{\bm{\mu}}=1,\\
			D_{\bm{z}}G(\bm{0},0)
			&=\nabla F(\bm{w}_*(\bm{\mu}))^\top U_{\bm{\mu}}
			=\bm{0}^\top.
		\end{aligned}
	\end{equation*}
	By the implicit function theorem, there are neighborhoods
	$V\subset\mathbb{R}^{d-1}$ of $\bm{0}$ and $I\subset\mathbb{R}$
	of $0$, and $\varrho_{\bm\mu}\in C^2(V)$ with $\varrho_{\bm\mu}(\bm{0})=0$, such that
	\begin{equation*}
		G(\bm{z},a)=0\quad\Longleftrightarrow\quad a=-\varrho_{\bm\mu}(\bm{z}),
		\qquad (\bm{z},a)\in V\times I.
	\end{equation*}
	Differentiating $G(\bm{z},-\varrho_{\bm\mu}(\bm{z}))=0$ at $\bm{z}=\bm{0}$ gives
	\begin{equation*}
		D_{\bm{z}}G(\bm{0},0)
		-\partial_aG(\bm{0},0)D\varrho_{\bm\mu}(\bm{0})=\bm{0}^\top,
		\qquad D\varrho_{\bm\mu}(\bm{0})=\bm{0}^\top.
	\end{equation*}
	Consequently, in this neighborhood the boundary is parametrized by
	\begin{equation}
		\begin{aligned}
			\bm{w}(\bm{z})=\bm{w}_*(\bm{\mu})+U_{\bm{\mu}}\bm{z}-\bm{n}_{\bm{\mu}}\varrho_{\bm\mu}(\bm{z})\in\partial\Omega,
		\end{aligned}
		\label{eq:fast_boundary_local_graph}
	\end{equation}
	where $\varrho_{\bm\mu}\in C^2$ near $\bm{0}$, $\varrho_{\bm\mu}(\bm{0})=0$, and $\nabla \varrho_{\bm\mu}(\bm{0})=\bm{0}$. For every $\bm{a}\in\mathbb{R}^{d-1}$, by \eqref{eq:fast_boundary_local_graph}, the tangent vector takes the form
	\begin{equation*}
		D\bm{w}(\bm{z})[\bm{a}]=U_{\bm{\mu}}\bm{a}-\bm{n}_{\bm{\mu}}\nabla \varrho_{\bm\mu}(\bm{z})^\top\bm{a}.
	\end{equation*}
	Note that $(\bm{n}_{\bm{\mu}}+U_{\bm{\mu}}\nabla \varrho_{\bm\mu}(\bm{z}))^\top D\bm{w}(\bm{z})[\bm{a}]=0$, and therefore the unit normal vector is
	\begin{equation*}
		\bm{N}(\bm{w}(\bm{z}))=\frac{\bm{n}_{\bm{\mu}}+U_{\bm{\mu}}\nabla \varrho_{\bm\mu}(\bm{z})}{\sqrt{1+\|\nabla \varrho_{\bm\mu}(\bm{z})\|_2^2}}.
	\end{equation*}
	Again, by \eqref{eq:fast_boundary_local_graph}, $\bm{z}=U_{\bm{\mu}}^\top(\bm{w}(\bm{z})-\bm{w}_*(\bm{\mu}))$, which gives $\bm{N}\in C^1$ on the neighborhood, since $\varrho_{\bm\mu}\in C^2$, and
	\begin{equation*}
		D\bm{N}(\bm{w}_*(\bm{\mu}))[U_{\bm{\mu}}\bm{a}]
		=U_{\bm{\mu}}\nabla^2\varrho_{\bm\mu}(\bm{0})\bm{a}
		=U_{\bm{\mu}}\mathscr{S}_{\bm{\mu}}\bm{a},
		\qquad\bm{a}\in\mathbb{R}^{d-1}
	\end{equation*}
	on $T_{\bm{n}_{\bm{\mu}}}S^{d-1}=\bm{n}_{\bm{\mu}}^\perp=\mathcal{T}_{\bm{\mu}}$, the tangent space of $S^{d-1}$ at $\bm{n}_{\bm{\mu}}$. Thus
	\begin{equation*}
		D\bm{N}(\bm{w}_*(\bm{\mu})):
		T_{\bm{w}_*(\bm{\mu})}\partial\Omega
		\to T_{\bm{n}_{\bm{\mu}}}S^{d-1}
	\end{equation*}
	is a linear isomorphism, since its matrix in the orthonormal basis $U_{\bm{\mu}}$ is $\mathscr{S}_{\bm{\mu}}\succ0$.
	Therefore, by the $C^1$ inverse function theorem, $\bm N$ has a $C^1$ local inverse $\bm N^{-1}$, defined on a neighborhood of $\bm n_{\bm\mu}$ in $S^{d-1}$ and taking values in a neighborhood of $\bm w_*(\bm\mu)$ in $\partial\Omega$.
	
	To identify this inverse with the optimizer, let $\bm{x}_k\to\bm{\mu}$ and $\bm{w}_k\in\arg\max_{\bm{w}\in\Omega}\bm{x}_k^\top\bm{w}$. Compactness and continuity imply, for every convergent subsequence,
	\begin{equation*}
		\begin{aligned}
			\bm{w}_{k_j}\to\bm{w}\in\Omega
			&\ \Longrightarrow\
			\bm{\mu}^\top\bm{w}
			=\lim_{j\to\infty}\bm{x}_{k_j}^\top\bm{w}_{k_j}
			=\lim_{j\to\infty}f_*(\bm{x}_{k_j})=f_*(\bm{\mu})\\
			&\ \Longrightarrow\ \bm{w}=\bm{w}_*(\bm{\mu}).
		\end{aligned}
	\end{equation*}
	Consequently,
	\begin{equation*}
		\lim_{\bm{x}\to\bm{\mu}}\sup_{\bm{w}\in\arg\max_{\bm{v}\in\Omega}\bm{x}^\top\bm{v}}
		\|\bm{w}-\bm{w}_*(\bm{\mu})\|_2=0.
	\end{equation*}
	All maximizers thus lie in the boundary neighborhood on which
	$\bm{N}$ is injective for $\bm{x}$ sufficiently close to $\bm{\mu}$.
	At each such maximizer,
	$\bm{N}(\bm{w})=\bm{x}/\|\bm{x}\|_2$; local injectivity yields
	\begin{equation}
		\arg\max_{\bm{w}\in\Omega}\bm{x}^\top\bm{w}
		=\left\{\bm{N}^{-1}\left(\frac{\bm{x}}{\|\bm{x}\|_2}\right)\right\},
		\qquad
		\bm{N}(\bm{w}_*(\bm{x}))=\frac{\bm{x}}{\|\bm{x}\|_2}.
		\label{eq:fast_inverse_normal_map}
	\end{equation}
	Here the objective neighborhood is also chosen so that
	$\bm{x}/\|\bm{x}\|_2$ belongs to the domain of $\bm{N}^{-1}$.
	The support-function subdifferential identity gives
	\begin{equation*}
		\begin{aligned}
			\partial f_*(\bm{x})&=\{\bm{w}_*(\bm{x})\},\qquad
			\nabla f_*(\bm{x})=\bm{w}_*(\bm{x})
			=\bm{N}^{-1}\left(\frac{\bm{x}}{\|\bm{x}\|_2}\right),
		\end{aligned}
	\end{equation*}
	and therefore $\bm{w}_*\in C^1\ \Longrightarrow\ \nabla f_*\in C^1\ \Longrightarrow\ f_*\in C^2.$
	
	Since $\bm{w}_*=\nabla f_*$ locally, for every $\bm{v}\in\mathbb{R}^d$,
	\begin{equation*}
		D\bm{w}_*(\bm{\mu})[\bm{v}]
		=D(\nabla f_*)(\bm{\mu})[\bm{v}]
		=H_{\bm{\mu}}\bm{v}
		=\left.\frac{\mathrm{d}}{\mathrm{d}t}
		\bm{w}_*(\bm{\mu}+t\bm{v})\right|_{t=0}.
	\end{equation*}
	The curve $t\mapsto\bm{w}_*(\bm{\mu}+t\bm{v})$ lies in
	$\partial\Omega$ for sufficiently small $|t|$, so $H_{\bm{\mu}}\bm{v}\in\mathcal{T}_{\bm{\mu}}$.
	Also, the columns of $(U_{\bm{\mu}},\bm{n}_{\bm{\mu}})$ form an
	orthonormal basis of $\mathbb{R}^d$, i.e., $U_{\bm{\mu}}U_{\bm{\mu}}^\top+\bm{n}_{\bm{\mu}}\bm{n}_{\bm{\mu}}^\top=I_d$. Hence $U_{\bm{\mu}}U_{\bm{\mu}}^\top H_{\bm{\mu}}\bm{v}=H_{\bm{\mu}}\bm{v}$.
	
	To differentiate the right-hand side of
	\eqref{eq:fast_inverse_normal_map}, first compute
	\begin{equation*}
		\left.\frac{\mathrm{d}}{\mathrm{d}t}
		\frac{\bm{\mu}+t\bm{v}}{\|\bm{\mu}+t\bm{v}\|_2}
		\right|_{t=0}
		=\frac{\bm{v}}{\|\bm{\mu}\|_2}
		-\frac{\bm{\mu}(\bm{\mu}^\top\bm{v})}{\|\bm{\mu}\|_2^3}=\frac{1}{\|\bm{\mu}\|_2}
		(I_d-\bm{n}_{\bm{\mu}}\bm{n}_{\bm{\mu}}^\top)\bm{v}
		=\frac{1}{\|\bm{\mu}\|_2}U_{\bm{\mu}}U_{\bm{\mu}}^\top\bm{v}.
	\end{equation*}
	For the left-hand side, the chain rule gives
	\begin{equation*}
		\left.\frac{\mathrm{d}}{\mathrm{d}t}
		\bm{N}(\bm{w}_*(\bm{\mu}+t\bm{v}))\right|_{t=0}
		=D\bm{N}(\bm{w}_*(\bm{\mu}))
		[D\bm{w}_*(\bm{\mu})[\bm{v}]]=D\bm{N}(\bm{w}_*(\bm{\mu}))
		[H_{\bm{\mu}}\bm{v}].
	\end{equation*}
	Equating these derivatives yields
	\begin{equation*}
		D\bm{N}(\bm{w}_*(\bm{\mu}))[H_{\bm{\mu}}\bm{v}]
		=\frac{1}{\|\bm{\mu}\|_2}
		U_{\bm{\mu}}U_{\bm{\mu}}^\top\bm{v}.
	\end{equation*}
	Substitute $H_{\bm{\mu}}\bm{v}
	=U_{\bm{\mu}}(U_{\bm{\mu}}^\top H_{\bm{\mu}}\bm{v})$
	and multiply by $U_{\bm{\mu}}^\top$. By
	\eqref{eq:fast_boundary_curvature_matrix},
	\begin{equation*}
		\mathscr{S}_{\bm{\mu}}U_{\bm{\mu}}^\top H_{\bm{\mu}}\bm{v}
		=U_{\bm{\mu}}^\top D\bm{N}(\bm{w}_*(\bm{\mu}))U_{\bm{\mu}}
		(U_{\bm{\mu}}^\top H_{\bm{\mu}}\bm{v})
		=\frac{1}{\|\bm{\mu}\|_2}
		U_{\bm{\mu}}^\top U_{\bm{\mu}}U_{\bm{\mu}}^\top\bm{v}
		=\frac{1}{\|\bm{\mu}\|_2}U_{\bm{\mu}}^\top\bm{v}.
	\end{equation*}
	Since $\mathscr{S}_{\bm{\mu}}\succ0$, multiplication by its inverse gives $U_{\bm{\mu}}^\top H_{\bm{\mu}}\bm{v}=\frac{1}{\|\bm{\mu}\|_2}\mathscr{S}_{\bm{\mu}}^{-1}U_{\bm{\mu}}^\top\bm{v}$, and thus
	\begin{equation*}
		H_{\bm{\mu}}\bm{v}
		=U_{\bm{\mu}}U_{\bm{\mu}}^\top H_{\bm{\mu}}\bm{v}
		=\frac{1}{\|\bm{\mu}\|_2}U_{\bm{\mu}}\mathscr{S}_{\bm{\mu}}^{-1}U_{\bm{\mu}}^\top\bm{v}.
	\end{equation*}
	As this holds for every $\bm{v}\in\mathbb{R}^d$,
	\begin{equation*}
		H_{\bm{\mu}}=\frac{1}{\|\bm{\mu}\|_2}
		U_{\bm{\mu}}\mathscr{S}_{\bm{\mu}}^{-1}U_{\bm{\mu}}^\top.
	\end{equation*}
	Moreover,
	\begin{equation*}
		\bm{v}^\top H_{\bm{\mu}}\bm{v}
		=\frac{(U_{\bm{\mu}}^\top\bm{v})^\top
			\mathscr{S}_{\bm{\mu}}^{-1}(U_{\bm{\mu}}^\top\bm{v})}
		{\|\bm{\mu}\|_2}\ge0.
	\end{equation*}
	Positive definiteness of $\mathscr{S}_{\bm{\mu}}^{-1}$ gives
	\begin{equation*}
		\begin{aligned}
			\bm{v}^\top H_{\bm{\mu}}\bm{v}=0\ \Longleftrightarrow\ U_{\bm{\mu}}^\top\bm{v}=0\ \Longleftrightarrow\ \bm{v}\in\mathcal{T}_{\bm{\mu}}^\perp
			=\operatorname{span}\{\bm{\mu}\}.
		\end{aligned}
	\end{equation*}
	Thus
	\begin{equation*}
		\ker(H_{\bm{\mu}})=\operatorname{span}\{\bm{\mu}\},\qquad
		\operatorname{rank}(H_{\bm{\mu}})
		=d-\dim\ker(H_{\bm{\mu}})=d-1,
	\end{equation*}
	since $\bm{\mu}\ne\bm{0}$. This proves
	\eqref{eq:fast_hessian_principal_curvatures} and the rank and kernel assertions.
	
	In the principal basis, we have
	\begin{equation*}
		\operatorname{tr}(\Sigma H_{\bm\mu})
		=\frac{\operatorname{tr}(\mathscr{S}_{\bm\mu}^{-1}
			U_{\bm\mu}^\top\Sigma U_{\bm\mu})}{\|\bm\mu\|_2}
		=\frac1{\|\bm\mu\|_2}
		\sum_{i=1}^{d-1}\frac{\bm u_i^\top\Sigma\bm u_i}{\kappa_i(\bm\mu)}.
	\end{equation*}
	The established local $C^2$ smoothness permits the application of
	Theorem \ref{thm:fast_local_minimax}, yielding
	\begin{equation*}
		\lim_{M\to\infty}\lim_{n\to\infty}nR_{n,M}^{\Sigma}(\bm{\mu})
		=\frac12\operatorname{tr}(\Sigma H_{\bm{\mu}})
		=\frac{1}{2\|\bm{\mu}\|_2}
		\sum_{i=1}^{d-1}\frac{\bm u_i^\top\Sigma\bm u_i}{\kappa_i(\bm{\mu})}.
	\end{equation*}
	Corollary \ref{cor:fast_distance_smooth} also gives $\dist(\bm{\mu},A)>0$.
	Since $\bm{\mu}=\|\bm{\mu}\|_2\bm{n}_{\bm{\mu}}$,
	\eqref{eq:fast_distance_homogeneity} implies
	\begin{equation*}
		\dist(\bm{\mu},A)=\dist(\|\bm{\mu}\|_2\bm{n}_{\bm{\mu}},A)
		=\|\bm{\mu}\|_2\dist(\bm{n}_{\bm{\mu}},A).
	\end{equation*}
	Substituting into the minimax constant gives
	\begin{equation*}
		\lim_{M\to\infty}\lim_{n\to\infty}nR_{n,M}^{\Sigma}(\bm{\mu})
		=\frac{\dist(\bm{n}_{\bm{\mu}},A)}{2\dist(\bm{\mu},A)}
		\sum_{i=1}^{d-1}\frac{\bm u_i^\top\Sigma\bm u_i}{\kappa_i(\bm{\mu})},
	\end{equation*}
	which proves \eqref{eq:fast_principal_curvature_minimax}.
\end{proof}

\subsection{Quadratic growth implies positive curvature}
\label{app:growth-implies-curvature}
\begin{lemma}
	\label{lem:growth_implies_curvature}
	Let $d\ge2$, let $\Omega$ be compact and convex with nonempty interior, and fix $\bm{\mu}\ne\bm{0}$ with $\dist(\bm{\mu},A)>0$. Assume that $\partial\Omega$ is a $C^2$ hypersurface near $\bm{w}_*(\bm{\mu})$ and that \eqref{eq:fast_geometric_condition} holds at $\bm z=\bm\mu$. Then $\mathscr{S}_{\bm\mu}\succ0$, so Corollary~\ref{cor:fast_principal_curvatures} applies without a separate positive-curvature assumption.
\end{lemma}
\begin{proof}
	Since $\bm{\mu}\notin A$, the maximizer $\bm{w}_*(\bm{\mu})$ is unique. The local parametrization \eqref{eq:fast_boundary_local_graph} and the identity $D\bm N(\bm w_*(\bm\mu))[U_{\bm\mu}\bm a]=U_{\bm\mu}\nabla^2\varrho_{\bm\mu}(\bm0)\bm a$ were derived in the proof of Corollary~\ref{cor:fast_principal_curvatures} using only the $C^2$ property of $\partial\Omega$. In particular $\mathscr{S}_{\bm\mu}=\nabla^2\varrho_{\bm\mu}(\bm0)$ is symmetric.
	
	For each $\bm{a}\in\mathbb{R}^{d-1}$ and sufficiently small $t>0$, \eqref{eq:fast_geometric_condition} with $\bm{z}=\bm{\mu}$ and $\bm{w}=\bm{w}(t\bm{a})$ gives
	\begin{equation*}
		\|\bm{\mu}\|_2\varrho_{\bm\mu}(t\bm{a})
		=\bm{\mu}^\top(\bm{w}_*(\bm{\mu})-\bm{w}(t\bm{a}))\ge c_\Omega \dist(\bm{\mu},A)
		\|tU_{\bm{\mu}}\bm{a}-\bm{n}_{\bm{\mu}}\varrho_{\bm\mu}(t\bm{a})\|_\Omega^2.
	\end{equation*}
	Since $\varrho_{\bm\mu}\in C^2$, $\varrho_{\bm\mu}(\bm{0})=0$, and $\nabla \varrho_{\bm\mu}(\bm{0})=\bm{0}$,
	\begin{equation*}
		\varrho_{\bm\mu}(t\bm{a})=\frac{t^2}{2}\bm{a}^\top \mathscr{S}_{\bm{\mu}}\bm{a}+o(t^2),
		\qquad \frac{\varrho_{\bm\mu}(t\bm{a})}{t}\to0.
	\end{equation*}
	Dividing the inequality by $t^2$ and letting $t\downarrow0$ yields
	\begin{equation*}
		\frac{\|\bm{\mu}\|_2}{2}\bm{a}^\top \mathscr{S}_{\bm{\mu}}\bm{a}
		\ge c_\Omega \dist(\bm{\mu},A)
		\lim_{t\downarrow0}
		\left\|U_{\bm{\mu}}\bm{a}
		-\bm{n}_{\bm{\mu}}\frac{\varrho_{\bm\mu}(t\bm{a})}{t}\right\|_\Omega^2=c_\Omega \dist(\bm{\mu},A)\|U_{\bm{\mu}}\bm{a}\|_\Omega^2.
	\end{equation*}
	Consequently,
	\begin{equation*}
		\min_{1\le i\le d-1}\kappa_i(\bm{\mu})
		=\min_{\|\bm{a}\|_2=1}\bm{a}^\top \mathscr{S}_{\bm{\mu}}\bm{a}\ge\frac{2c_\Omega \dist(\bm{\mu},A)}{\|\bm{\mu}\|_2}
		\min_{\|\bm{a}\|_2=1}\|U_{\bm{\mu}}\bm{a}\|_\Omega^2>0.
	\end{equation*}
	The final minimum is positive because $U_{\bm{\mu}}$ is injective and $\{\bm{a}:\|\bm{a}\|_2=1\}$ is compact.
\end{proof}

\subsection{Technical boundary-growth theorem}
The following result supplies the explicit lower bound, matching upper rate, and homogeneous-profile conclusion summarized in Theorem~\ref{thm:geometry-contact-rate}.
\begin{theorem}[Boundary growth and distance-based lower bounds without $C^2$]
	\label{thm:fast_boundary_growth_lower}
	\label{thm:fast_two_sided_growth}
	Let $d\ge2$, let $\Omega$ be compact and convex, and fix $\bm{\mu}\ne\bm{0}$ with $\dist(\bm{\mu},A)>0$, where the distance is computed in the fixed norm $\|\cdot\|_K$. Retain $\bm{n}_{\bm{\mu}}=\bm{\mu}/\|\bm{\mu}\|_2$ and the orthonormal tangent basis $U_{\bm{\mu}}$ used in Corollary~\ref{cor:fast_principal_curvatures}. Suppose there exist $\rho>0$, $\bar{\kappa}>0$, $\beta>0$, and a boundary profile $\varrho_{\bm\mu}$ on $\{\bm{z}\in\mathbb{R}^{d-1}:\|\bm{z}\|_2\le\rho\}$ such that
	\begin{equation}
		\bm{w}_*(\bm{\mu})+U_{\bm{\mu}}\bm{z}-\bm{n}_{\bm{\mu}}\varrho_{\bm\mu}(\bm{z})\in\Omega,
		\label{eq:fast_boundary_growth_assumption}
	\end{equation}
	where $\varrho_{\bm\mu}(\bm{0})=0$, and $0\le \varrho_{\bm\mu}(\bm{z})\le\frac{\bar{\kappa}}{1+\beta}\|\bm{z}\|_2^{1+\beta}$ for $\|\bm{z}\|_2\le\rho$. This uses the same parametrization as \eqref{eq:fast_boundary_local_graph}, but does not require $\varrho_{\bm\mu}$ nor $f_*$ to be twice differentiable. For any $M>0$, set $h=\min\{1,M\}$. Whenever
	\begin{equation}
		\frac{h}{\sqrt n}
		\le\bar{\kappa}\|\bm{\mu}\|_2\rho^{\beta},
		\label{eq:fast_boundary_growth_sample_condition}
	\end{equation}
	for every unit $\bm v\in\mathbb R^{d-1}$, the local minimax risk satisfies
	\begin{equation}
		\begin{aligned}
			R_{n,M}^{\Sigma}(\bm{\mu})
			&\ge 2\Phi\left(-h\sqrt{\bm v^\top U_{\bm\mu}^\top\Sigma^{-1}U_{\bm\mu}\bm v}\right)\frac{\beta}{1+\beta}\left(\frac{h^{1+\beta}}{\bar{\kappa}\|\bm{\mu}\|_2n^{(1+\beta)/2}}\right)^{1/\beta}
			\\&=2\Phi\left(-h\sqrt{\bm v^\top U_{\bm\mu}^\top\Sigma^{-1}U_{\bm\mu}\bm v}\right)\frac{\beta}{1+\beta}
			h^{(1+\beta)/\beta}n^{-(1+\beta)/(2\beta)}
			\left(
			\frac{\dist(\bm{n}_{\bm{\mu}},A)}
			{\bar{\kappa}\,\dist(\bm{\mu},A)}
			\right)^{1/\beta}.
		\end{aligned}
		\label{eq:fast_boundary_growth_distance_lower}
	\end{equation}
	For the matching upper bound, consider $\bm\mu\ne\bm0$ with a unique maximizer $\bm w_*(\bm\mu)$. Suppose that, for some $\rho>0$, the continuous convex boundary profile $\varrho_{\bm\mu}$ is defined on a neighborhood of $\{\bm z\in\mathbb R^{d-1}:\|\bm z\|_2\le\rho\}$, satisfies $\varrho_{\bm\mu}(\bm z)<\rho$ there, and describes the local feasible set by
	\begin{equation}
		\bm w_*(\bm\mu)+U_{\bm\mu}\bm z-\bm n_{\bm\mu}b\in\Omega
		\quad\Longleftrightarrow\quad \varrho_{\bm\mu}(\bm z)\le b,
		\qquad \|\bm z\|_2\le\rho,\quad |b|\le\rho.
		\label{eq:fast_two_sided_local_graph}
	\end{equation}
	Assume that, for $\beta>0$ and $0<\underline\kappa\le\bar\kappa<\infty$,
	\begin{equation}
		\frac{\underline\kappa}{1+\beta}\|\bm z\|_2^{1+\beta}
		\le \varrho_{\bm\mu}(\bm z)\le
		\frac{\bar\kappa}{1+\beta}\|\bm z\|_2^{1+\beta},
		\qquad \|\bm z\|_2\le\rho.
		\label{eq:fast_two_sided_boundary_growth}
	\end{equation}
	For every fixed $M>0$, there are $0<c_M\le C_M<\infty$ such that, for all sufficiently large $n$ and all $\Sigma\in\mathcal C$,
	\begin{equation}
		\begin{aligned}
			c_M\|\bm\mu\|_2^{-1/\beta}n^{-(1+\beta)/(2\beta)}
			&\le R_{n,M}^{\Sigma}(\bm\mu)
			\le\sup_{\bm t\in B_2(M)}
			\mathbb E_{\bm\mu_{n,\bm t},\Sigma}R(\bm w_*(\bar{\bm x}))
			\le C_M\|\bm\mu\|_2^{-1/\beta}n^{-(1+\beta)/(2\beta)}.
		\end{aligned}
		\label{eq:fast_two_sided_regret_rate}
	\end{equation}
	The constants can depend on $d,M,\beta,\underline\kappa,\bar\kappa,\underline\lambda,\overline\lambda$.
	
	If, in addition, a function $\varrho_0:\mathbb R^{d-1}\to\mathbb R$ satisfies
	\begin{equation}
		\lim_{\delta\downarrow0}\sup_{\|\bm z\|_2\le L}
		\left|\frac{\varrho_{\bm\mu}(\delta\bm z)}{\delta^{1+\beta}}-\varrho_0(\bm z)\right|=0,
		\qquad L>0,
		\label{eq:fast_boundary_profile_limit}
	\end{equation}
	then $\varrho_0$ is finite, convex, and positively homogeneous of degree $1+\beta$. Define its convex conjugate by
	\begin{equation*}
		\varrho_0^*(\bm y):=\sup_{\bm z\in\mathbb R^{d-1}}
		\{\bm y^{\top}\bm z-\varrho_0(\bm z)\}.
	\end{equation*}
	The expansion \eqref{eq:fast_general_local_expansion} holds with
	\begin{equation}
		a(\varepsilon)=\varepsilon^{(1+\beta)/\beta},
		\qquad
		\psi(\bm t)=\|\bm\mu\|_2^{-1/\beta}
		\varrho_0^*(U_{\bm\mu}^{\top}\bm t).
		\label{eq:fast_boundary_conjugate_profile}
	\end{equation}
\end{theorem}

\begin{proof}
	Fix any $\bm{v}\in\mathbb{R}^{d-1}$ with $\|\bm{v}\|_2=1$ and take $\bm{t}=hU_{\bm{\mu}}\bm{v}\in B_2(M)$. For $0\le b\le\rho$, $\bm{w}_*(\bm{\mu})+bU_{\bm{\mu}}\bm{v}-\bm{n}_{\bm{\mu}}\varrho_{\bm\mu}(b\bm{v})\in\Omega$. Since $\bm{\mu}^\top U_{\bm{\mu}}\bm{v}=\bm{n}_{\bm{\mu}}^\top U_{\bm{\mu}}\bm{v}=0$, and
	$U_{\bm{\mu}}^\top U_{\bm{\mu}}=I_{d-1}$, the support-function
	definition gives
	\begin{equation*}
		\begin{aligned}
			r_{\bm{\mu}}(\bm{t}/\sqrt n)
			&=\sup_{\bm{w}\in\Omega}
			(\bm{\mu}+\bm{t}/\sqrt n)^\top
			(\bm{w}-\bm{w}_*(\bm{\mu}))\\
			&\ge(\bm{\mu}+hU_{\bm{\mu}}\bm{v}/\sqrt n)^\top
			(bU_{\bm{\mu}}\bm{v}-\bm{n}_{\bm{\mu}}\varrho_{\bm\mu}(b\bm{v}))\\
			&=\frac{hb}{\sqrt n}-\|\bm{\mu}\|_2\varrho_{\bm\mu}(b\bm{v})\\
			&\ge\frac{hb}{\sqrt n}
			-\frac{\bar{\kappa}\|\bm{\mu}\|_2}{1+\beta}b^{1+\beta}.
		\end{aligned}
	\end{equation*}
	The same calculation with $-b\bm{v}$ gives the identical bound for
	$r_{\bm{\mu}}(-\bm{t}/\sqrt n)$.
	The concave function $b\mapsto\frac{hb}{\sqrt n}-\frac{\bar{\kappa}\|\bm{\mu}\|_2}{1+\beta}b^{1+\beta}$ on the last line is maximized over
	$b\ge0$ at
	\begin{equation*}
		b=\left(\frac{h}{\bar{\kappa}\|\bm{\mu}\|_2\sqrt n}
		\right)^{1/\beta}.
	\end{equation*}
	By \eqref{eq:fast_boundary_growth_sample_condition}, this value belongs to $[0,\rho]$. Substitution yields
	\begin{equation}
		\min\{r_{\bm{\mu}}(\bm{t}/\sqrt n),
		r_{\bm{\mu}}(-\bm{t}/\sqrt n)\}
		\ge\frac{\beta}{1+\beta}
		\left(\frac{h}{\sqrt n}\right)^{(1+\beta)/\beta}
		(\bar{\kappa}\|\bm{\mu}\|_2)^{-1/\beta}.
		\label{eq:fast_boundary_growth_remainder_lower}
	\end{equation}
	
	Use the equal-probability prior on $\bm{T}=\bm{t}$ and
	$\bm{T}=-\bm{t}$. The same binary-prior argument in the proof of Proposition~\ref{prop:fast_symmetric_difference} gives, for this particular
	$\bm{t}$,
	\begin{equation*}
		\begin{aligned}
			R_{n,M}^{\Sigma}(\bm{\mu})
			&\ge 2\Phi(-\|\Sigma^{-1/2}\bm{t}\|_2)
			\left\{\frac{f_*(\bm{\mu}+\bm{t}/\sqrt n)
				+f_*(\bm{\mu}-\bm{t}/\sqrt n)}2
			-f_*(\bm{\mu})\right\}\\
			&=\Phi\left(-h\sqrt{\bm v^\top U_{\bm\mu}^\top\Sigma^{-1}U_{\bm\mu}\bm v}\right)
			\{r_{\bm{\mu}}(\bm{t}/\sqrt n)
			+r_{\bm{\mu}}(-\bm{t}/\sqrt n)\}.
		\end{aligned}
	\end{equation*}
	Combining with \eqref{eq:fast_boundary_growth_remainder_lower} and $\|\bm{\mu}\|_2^{-1}=\dist(\bm{n}_{\bm{\mu}},A)/\dist(\bm{\mu},A)$
	proves \eqref{eq:fast_boundary_growth_distance_lower}.
	
	We next prove the two-sided bounds and the conjugate limit. Uniqueness at $\bm\mu$ and compactness localize every perturbed maximizer. Indeed, if $\bm h_k\to\bm0$ and $\bm w_k\in\arg\max_{\bm w\in\Omega}(\bm\mu+\bm h_k)^{\top}\bm w$, every convergent subsequence satisfies
	\begin{equation*}
		\bm\mu^{\top}\lim_k\bm w_k
		=\lim_k\{f_*(\bm\mu+\bm h_k)-\bm h_k^{\top}\bm w_k\}
		=f_*(\bm\mu),
		\qquad \lim_k\bm w_k=\bm w_*(\bm\mu).
	\end{equation*}
	Hence all maximizers lie in the interior of the coordinate cylinder in \eqref{eq:fast_two_sided_local_graph} for sufficiently small $\bm h$. Since $f_*(\bm\mu)=\bm\mu^{\top}\bm w_*(\bm\mu)$, the definition of the remainder gives
	\begin{equation*}
		\begin{aligned}
			r_{\bm\mu}(\bm h)
			&=f_*(\bm\mu+\bm h)-f_*(\bm\mu)-\bm w_*(\bm\mu)^{\top}\bm h\\
			&=\max_{\bm w\in\Omega}(\bm\mu+\bm h)^{\top}
			(\bm w-\bm w_*(\bm\mu))\\
			&=\max_{\substack{\|\bm z\|_2\le\rho,\ |b|\le\rho\\\varrho_{\bm\mu}(\bm z)\le b}}
			(\bm\mu+\bm h)^{\top}(U_{\bm\mu}\bm z-\bm n_{\bm\mu}b).
		\end{aligned}
	\end{equation*}
	The last equality follows from the localization of the maximizers and \eqref{eq:fast_two_sided_local_graph}. Using $U_{\bm\mu}^{\top}\bm\mu=\bm0$ and $\bm\mu^{\top}\bm n_{\bm\mu}=\|\bm\mu\|_2$, we obtain
	\begin{equation*}
		(\bm\mu+\bm h)^{\top}(U_{\bm\mu}\bm z-\bm n_{\bm\mu}b)
		=(U_{\bm\mu}^{\top}\bm h)^{\top}\bm z
		-(\|\bm\mu\|_2+\bm n_{\bm\mu}^{\top}\bm h)b.
	\end{equation*}
	For $\|\bm h\|_2<\|\bm\mu\|_2$,
	\begin{equation*}
		\|\bm\mu\|_2+\bm n_{\bm\mu}^{\top}\bm h
		\ge\|\bm\mu\|_2-\|\bm h\|_2>0.
	\end{equation*}
	Moreover, \eqref{eq:fast_two_sided_boundary_growth} and the assumption $\varrho_{\bm\mu}(\bm z)<\rho$ give $0\le \varrho_{\bm\mu}(\bm z)<\rho$ for $\|\bm z\|_2\le\rho$. Thus, for each fixed $\bm z$, the feasible values of $b$ form $[\varrho_{\bm\mu}(\bm z),\rho]$, and the objective is strictly decreasing in $b$. Consequently,
	\begin{equation*}
		\begin{aligned}
			&\max_{b\in[\varrho_{\bm\mu}(\bm z),\rho]}
			\left\{(U_{\bm\mu}^{\top}\bm h)^{\top}\bm z
			-(\|\bm\mu\|_2+\bm n_{\bm\mu}^{\top}\bm h)b\right\}\\
			&\qquad=(U_{\bm\mu}^{\top}\bm h)^{\top}\bm z
			-(\|\bm\mu\|_2+\bm n_{\bm\mu}^{\top}\bm h)\varrho_{\bm\mu}(\bm z).
		\end{aligned}
	\end{equation*}
	Maximizing over $\bm z$ therefore gives
	\begin{equation}
		r_{\bm\mu}(\bm h)
		=\max_{\|\bm z\|_2\le\rho}
		\left\{(U_{\bm\mu}^{\top}\bm h)^{\top}\bm z
		-(\|\bm\mu\|_2+\bm n_{\bm\mu}^{\top}\bm h)\varrho_{\bm\mu}(\bm z)\right\}.
		\label{eq:fast_local_conjugate_representation}
	\end{equation}
	For $c>0$, direct maximization gives
	\begin{equation*}
		\sup_{\bm z\in\mathbb R^{d-1}}
		\left\{\bm y^{\top}\bm z-\frac{c}{1+\beta}\|\bm z\|_2^{1+\beta}\right\}
		=\frac{\beta}{1+\beta}c^{-1/\beta}\|\bm y\|_2^{(1+\beta)/\beta}.
	\end{equation*}
	For $\bm y\ne\bm0$, the maximizing vector has direction $\bm y$ and norm $(\|\bm y\|_2/c)^{1/\beta}$; for $\bm y=\bm0$, it is zero. Since $\|\bm\mu\|_2+\bm n_{\bm\mu}^{\top}\bm h>0$, the lower bound on $\varrho_{\bm\mu}$ in \eqref{eq:fast_two_sided_boundary_growth} gives
	\begin{equation*}
		\begin{aligned}
			r_{\bm\mu}(\bm h)
			&\le\max_{\|\bm z\|_2\le\rho}
			\left\{(U_{\bm\mu}^{\top}\bm h)^{\top}\bm z
			-\frac{\underline\kappa(\|\bm\mu\|_2+\bm n_{\bm\mu}^{\top}\bm h)}{1+\beta}
			\|\bm z\|_2^{1+\beta}\right\}\\
			&\le\sup_{\bm z\in\mathbb R^{d-1}}
			\left\{(U_{\bm\mu}^{\top}\bm h)^{\top}\bm z
			-\frac{\underline\kappa(\|\bm\mu\|_2+\bm n_{\bm\mu}^{\top}\bm h)}{1+\beta}
			\|\bm z\|_2^{1+\beta}\right\}\\
			&=\frac{\beta}{1+\beta}
			\frac{\|U_{\bm\mu}^{\top}\bm h\|_2^{(1+\beta)/\beta}}
			{[\underline\kappa(\|\bm\mu\|_2+\bm n_{\bm\mu}^{\top}\bm h)]^{1/\beta}}.
		\end{aligned}
	\end{equation*}
	Here the domain is enlarged only after replacing $\varrho_{\bm\mu}$ by its lower bound on $B_2^{d-1}(\rho):=\{\bm z\in\mathbb R^{d-1}:\|\bm z\|_2\le\rho\}$.
	
	Conversely, the upper bound on $\varrho_{\bm\mu}$ gives
	\begin{equation*}
		r_{\bm\mu}(\bm h)
		\ge\max_{\|\bm z\|_2\le\rho}
		\left\{(U_{\bm\mu}^{\top}\bm h)^{\top}\bm z
		-\frac{\bar\kappa(\|\bm\mu\|_2+\bm n_{\bm\mu}^{\top}\bm h)}{1+\beta}
		\|\bm z\|_2^{1+\beta}\right\}.
	\end{equation*}
	To evaluate this maximum, the maximizer over $\mathbb R^{d-1}$ must belong to $B_2^{d-1}(\rho)$. If $U_{\bm\mu}^{\top}\bm h=\bm0$, it is $\bm0$. Otherwise, its norm is
	\begin{equation*}
		\left(\frac{\|U_{\bm\mu}^{\top}\bm h\|_2}
		{\bar\kappa(\|\bm\mu\|_2+\bm n_{\bm\mu}^{\top}\bm h)}\right)^{1/\beta}
		\le\left(\frac{2\|\bm h\|_2}{\bar\kappa\|\bm\mu\|_2}\right)^{1/\beta}
		\longrightarrow0,
		\qquad \bm h\to\bm0,
	\end{equation*}
	where the inequality holds for $\|\bm h\|_2\le\|\bm\mu\|_2/2$. Thus, for sufficiently small $\bm h$, the maximum over $B_2^{d-1}(\rho)$ equals the supremum over $\mathbb R^{d-1}$, and
	\begin{equation*}
		r_{\bm\mu}(\bm h)
		\ge\frac{\beta}{1+\beta}
		\frac{\|U_{\bm\mu}^{\top}\bm h\|_2^{(1+\beta)/\beta}}
		{[\bar\kappa(\|\bm\mu\|_2+\bm n_{\bm\mu}^{\top}\bm h)]^{1/\beta}}.
	\end{equation*}
	Combining the two inequalities gives, for all sufficiently small $\|\bm h\|_2$,
	\begin{equation}
		\frac{\beta}{1+\beta}
		\frac{\|U_{\bm\mu}^{\top}\bm h\|_2^{(1+\beta)/\beta}}
		{[\bar\kappa(\|\bm\mu\|_2+\bm n_{\bm\mu}^{\top}\bm h)]^{1/\beta}}
		\le r_{\bm\mu}(\bm h)\le
		\frac{\beta}{1+\beta}
		\frac{\|U_{\bm\mu}^{\top}\bm h\|_2^{(1+\beta)/\beta}}
		{[\underline\kappa(\|\bm\mu\|_2+\bm n_{\bm\mu}^{\top}\bm h)]^{1/\beta}}.
		\label{eq:fast_two_sided_remainder}
	\end{equation}
	For fixed $M>0$, apply \eqref{eq:fast_two_sided_remainder} to $\bm h=\varepsilon\bm t$. For the lower bound choose $\bm t=MU_{\bm\mu}\bm v$ with $\|\bm v\|_2=1$, so $\bm n_{\bm\mu}^{\top}\bm t=0$. For the upper bound use $\|U_{\bm\mu}^{\top}\bm t\|_2\le M$ and $\|\bm\mu\|_2\pm\varepsilon\bm n_{\bm\mu}^{\top}\bm t\ge\|\bm\mu\|_2/2$. The definition of $D_{\bm\mu}$ then gives, for all sufficiently small $\varepsilon>0$,
	\begin{equation*}
		\frac{\beta}{1+\beta}(\bar\kappa\|\bm\mu\|_2)^{-1/\beta}
		(\varepsilon M)^{(1+\beta)/\beta}
		\le D_{\bm\mu}(\varepsilon,M)
		\le\frac{\beta}{1+\beta}
		(\underline\kappa\|\bm\mu\|_2/2)^{-1/\beta}
		(\varepsilon M)^{(1+\beta)/\beta}.
	\end{equation*}
	Consequently,
	\begin{equation*}
		D_{\bm\mu}(\varepsilon,M)>0,
		\qquad
		\frac{D_{\bm\mu}(2\varepsilon,M)}{D_{\bm\mu}(\varepsilon,M)}
		\le2^{(1+\beta)/\beta}
		\left(\frac{2\bar\kappa}{\underline\kappa}\right)^{1/\beta}.
	\end{equation*}
	Lemma~\ref{cor:fast_plugin_modulus_rate} therefore proves \eqref{eq:fast_two_sided_regret_rate}. The constant in \eqref{eq:fast_plugin_modulus_ratio} depends only on this doubling factor, $d,M,\overline\lambda$; the lower bound in Proposition~\ref{prop:fast_symmetric_difference} is uniform over $\mathcal C$. 
	
	Under \eqref{eq:fast_boundary_profile_limit}, convexity and \eqref{eq:fast_two_sided_boundary_growth} pass to $\varrho_0$. For $b>0$,
	\begin{equation*}
		\varrho_0(b\bm z)
		=\lim_{\delta\downarrow0}\frac{\varrho_{\bm\mu}(\delta b\bm z)}{\delta^{1+\beta}}
		=b^{1+\beta}\varrho_0(\bm z),
		\qquad
		\frac{\underline\kappa}{1+\beta}\|\bm z\|_2^{1+\beta}
		\le \varrho_0(\bm z)\le\frac{\bar\kappa}{1+\beta}\|\bm z\|_2^{1+\beta}.
	\end{equation*}
	In \eqref{eq:fast_local_conjugate_representation}, substitute $\bm h=\varepsilon\bm t$ and rescale the tangent coordinate by $\varepsilon^{1/\beta}$. Then
	\begin{equation*}
		\frac{r_{\bm\mu}(\varepsilon\bm t)}{\varepsilon^{(1+\beta)/\beta}}
		=\max_{\|\bm z\|_2\le\rho\varepsilon^{-1/\beta}}
		\left\{(U_{\bm\mu}^{\top}\bm t)^{\top}\bm z
		-(\|\bm\mu\|_2+\varepsilon\bm n_{\bm\mu}^{\top}\bm t)
		\frac{\varrho_{\bm\mu}(\varepsilon^{1/\beta}\bm z)}{\varepsilon^{(1+\beta)/\beta}}\right\}.
	\end{equation*}
	The value at $\bm z=\bm0$ is zero. Every nonzero maximizer, uniformly for $\|\bm t\|_2\le M$ and small $\varepsilon$, therefore satisfies
	\begin{equation*}
		0\le M\|\bm z\|_2
		-\frac{\underline\kappa\|\bm\mu\|_2}{2(1+\beta)}\|\bm z\|_2^{1+\beta},
		\qquad
		\|\bm z\|_2\le
		\left(\frac{2(1+\beta)M}{\underline\kappa\|\bm\mu\|_2}\right)^{1/\beta}.
	\end{equation*}
	The same bound contains the maximizers of the limiting objective. Uniform convergence in \eqref{eq:fast_boundary_profile_limit} on this fixed ball consequently gives, uniformly for $\bm t\in B_2(M)$,
	\begin{equation*}
		\frac{r_{\bm\mu}(\varepsilon\bm t)}{\varepsilon^{(1+\beta)/\beta}}
		\to
		\sup_{\bm z}\{(U_{\bm\mu}^{\top}\bm t)^{\top}\bm z-\|\bm\mu\|_2\varrho_0(\bm z)\}
		=\|\bm\mu\|_2\varrho_0^*(U_{\bm\mu}^{\top}\bm t/\|\bm\mu\|_2).
	\end{equation*}
	Rescaling $\bm z$ in the definition of $\varrho_0^*$ gives $\varrho_0^*(b\bm y)=b^{(1+\beta)/\beta}\varrho_0^*(\bm y)$, proving \eqref{eq:fast_boundary_conjugate_profile}. The two-sided bounds imply that $\psi$ is nonaffine on every $B_2(M)$, $M>0$. Theorem~\ref{thm:fast_general_local_rate} and the distance homogeneity identity \eqref{eq:fast_distance_homogeneity} complete the proof.
\end{proof}

\subsection{Proof of Proposition~\ref{prop:boundary-local-parametrization}}
\label{app:proof-boundary-local-parametrization}
\begin{proof}
First, we claim that $\bm w_*(\bm\mu)-b_0\bm n_{\bm\mu}\in\operatorname{int}\Omega$ for some $b_0>0$. Indeed, otherwise there would be $b_k\downarrow0$, such that $\bm w_*(\bm\mu)-b_k\bm n_{\bm\mu}\notin \operatorname{int}\Omega$ for all $k$. By hyperplane separation theorem, for each $k$, there exists a unit vector $\bm{v}_k$ satisfying
\begin{equation}
\bm v_k^\top(\bm w-\bm w_*(\bm\mu)+b_k\bm n_{\bm\mu})\le0\quad\forall\bm w\in\Omega,
\label{eq:v_k_separating_hyperplane}
\end{equation}
and in particular, letting $\bm{w}=\bm w_*(\bm\mu)$ gives $\bm v_k^\top\bm n_{\bm\mu}\le0$. Suppose $\bm{v}$ is a limit point of $\bm v_k$. Then $\bm{v}$ is also a unit vector, and letting $k\to\infty$ in \eqref{eq:v_k_separating_hyperplane} gives $\forall\bm w\in\Omega$, $\bm v^\top(\bm w-\bm w_*(\bm\mu))\le0$ and $\bm v^\top\bm n_{\bm\mu}\le0$, which implies that $\bm v$ is a unit supporting normal at $\bm w_*(\bm\mu)$ distinct from $\bm n_{\bm\mu}$, a contradiction to the uniqueness of supporting hyperplane. 

For each $\bm z\in\bigl\{U_{\bm\mu}^\top(\bm w-\bm w_*(\bm\mu)):\bm w\in\Omega\bigr\}
\subseteq\mathbb R^{d-1}$, define
\begin{equation}
\varrho_{\bm\mu}(\bm z)=\min\{b:\bm w_*(\bm\mu)+U_{\bm\mu}\bm z-\bm n_{\bm\mu}b\in\Omega\}.
\label{eq:defn_of_varrho}
\end{equation}
To verify that the minimum exists, fix such a $\bm z$ and choose $\bm w\in\Omega$ with $U_{\bm\mu}^\top(\bm w-\bm w_*(\bm\mu))=\bm z$. Since $U_{\bm\mu}U_{\bm\mu}^\top+\bm n_{\bm\mu}\bm n_{\bm\mu}^\top=I_d$,
\begin{equation*}
\bm w=\bm w_*(\bm\mu)+U_{\bm\mu}\bm z-\bm n_{\bm\mu}\bigl[\bm n_{\bm\mu}^\top(\bm w_*(\bm\mu)-\bm w)\bigr],
\end{equation*}
so $b=\bm n_{\bm\mu}^\top(\bm w_*(\bm\mu)-\bm w)$ is feasible. Moreover, every feasible $b$ satisfies
\begin{equation*}
0\ge\bm\mu^\top(\bm w-\bm w_*(\bm\mu))=\bm\mu^\top(U_{\bm\mu}\bm z-\bm n_{\bm\mu}b)=-\|\bm\mu\|_2b,
\end{equation*}
and hence $b\ge0$. The feasible set is therefore nonempty and bounded below, with a finite infimum
\begin{equation*}
b_*:=\inf\{b:\bm w_*(\bm\mu)+U_{\bm\mu}\bm z-\bm n_{\bm\mu}b\in\Omega\}\ge0.
\end{equation*}
Choose feasible $b_j\to b_*$. Since $\Omega$ is closed,
\begin{equation*}
\bm w_*(\bm\mu)+U_{\bm\mu}\bm z-\bm n_{\bm\mu}b_j
\longrightarrow\bm w_*(\bm\mu)+U_{\bm\mu}\bm z-\bm n_{\bm\mu}b_*\in\Omega.
\end{equation*}
Thus $b_*$ is feasible, so the minimum is attained and $\varrho_{\bm\mu}(\bm z)=b_*$.

To prove convexity, take $\bm z_1,\bm z_2$ in the domain of $\varrho_{\bm\mu}$ and $t\in[0,1]$. Attainment of the minima gives
\begin{equation*}
\bm w_i=\bm w_*(\bm\mu)+U_{\bm\mu}\bm z_i-\bm n_{\bm\mu}\varrho_{\bm\mu}(\bm z_i)\in\Omega,
\qquad i=1,2.
\end{equation*}
By convexity of $\Omega$,
\begin{equation*}
(1-t)\bm w_1+t\bm w_2
=\bm w_*(\bm\mu)+U_{\bm\mu}\bigl((1-t)\bm z_1+t\bm z_2\bigr)-\bm n_{\bm\mu}\bigl((1-t)\varrho_{\bm\mu}(\bm z_1)+t\varrho_{\bm\mu}(\bm z_2)\bigr)\in\Omega.
\end{equation*}
By definition of $\varrho_{\bm{\mu}}$ (Equation \eqref{eq:defn_of_varrho}),
\begin{equation*}
\varrho_{\bm\mu}\bigl((1-t)\bm z_1+t\bm z_2\bigr)
\le(1-t)\varrho_{\bm\mu}(\bm z_1)+t\varrho_{\bm\mu}(\bm z_2),
\end{equation*}
which proves convexity.

The interior point $\bm w_*(\bm\mu)-b_0\bm n_{\bm\mu}$ gives some $\rho>0$ such that
\begin{equation*}
\bm w_*(\bm\mu)+U_{\bm\mu}\bm z-b_0\bm n_{\bm\mu}\in\Omega,
\qquad 0\le\varrho_{\bm\mu}(\bm z)\le b_0,
\qquad \|\bm z\|_2\le\rho.
\end{equation*}
Hence $\bm0$ is in the interior of the domain of $\varrho_{\bm\mu}$. A finite convex function is continuous on the interior of its domain, so $\varrho_{\bm\mu}$ is continuous near $\bm0$. Since $\bm w_*(\bm\mu)\in\Omega$, $b=0$ is feasible at $\bm z=\bm0$; together with nonnegativity, this gives $\varrho_{\bm\mu}(\bm0)=0$.

For every $\bm a\in\partial\varrho_{\bm\mu}(\bm0)$, the subgradient inequality gives $\varrho_{\bm\mu}(\bm z)\ge\varrho_{\bm\mu}(\bm0)+\bm a^\top\bm z=\bm a^\top\bm z$ throughout the domain. For any $\bm w\in\Omega$, write $\bm z=U_{\bm\mu}^\top(\bm w-\bm w_*(\bm\mu))$, $b=\bm n_{\bm\mu}^\top(\bm w_*(\bm\mu)-\bm w)$, and $\bm w-\bm w_*(\bm\mu)=U_{\bm\mu}\bm z-\bm n_{\bm\mu}b$. Here $b$ is feasible, so $b\ge\varrho_{\bm\mu}(\bm z)\ge\bm a^\top\bm z$. Consequently,
\begin{equation*}
(U_{\bm\mu}\bm a+\bm n_{\bm\mu})^\top(\bm w-\bm w_*(\bm\mu))
=\bm a^\top\bm z-b\le0,
\qquad\forall\bm w\in\Omega.
\end{equation*}
Thus $U_{\bm\mu}\bm a+\bm n_{\bm\mu}$ is a supporting normal at $\bm w_*(\bm\mu)$. Uniqueness of the supporting hyperplane implies $U_{\bm\mu}\bm a+\bm n_{\bm\mu}=\lambda\bm n_{\bm\mu}$ for some $\lambda>0$. Taking inner products with $\bm n_{\bm\mu}$ gives $\lambda=1$, and then $U_{\bm\mu}\bm a=\bm0$ implies $\bm a=\bm0$. Conversely, $\varrho_{\bm\mu}(\bm z)\ge0=\varrho_{\bm\mu}(\bm0)$ shows that $\bm0\in\partial\varrho_{\bm\mu}(\bm0)$. Therefore $\partial\varrho_{\bm\mu}(\bm0)=\{\bm0\}$. A convex function is differentiable at an interior point of its domain if its subdifferential there is a singleton, so
\begin{equation*}
\nabla\varrho_{\bm\mu}(\bm0)=\bm0,\qquad
\varrho_{\bm\mu}(\bm z)=\varrho_{\bm\mu}(\bm0)+\nabla\varrho_{\bm\mu}(\bm0)^\top\bm z+o(\|\bm z\|_2)=o(\|\bm z\|_2).
\end{equation*}
Since $\bm w_*(\bm\mu)-b_0\bm n_{\bm\mu}\in\operatorname{int}\Omega$, for sufficiently small $\rho\in(0,b_0)$ we have
\begin{equation*}
\bm w_*(\bm\mu)+U_{\bm\mu}\bm z-b_0\bm n_{\bm\mu}\in\Omega,
\qquad 0\le\varrho_{\bm\mu}(\bm z)<\rho,
\qquad \|\bm z\|_2\le\rho.
\end{equation*}
For each fixed $\bm z$ with $\|\bm z\|_2\le\rho$, convexity of $\Omega$ implies that the feasible $b$ form an interval. Its minimum is $\varrho_{\bm\mu}(\bm z)$, and it contains $b_0>\rho$. Thus every $b\in[\varrho_{\bm\mu}(\bm z),\rho]$ is feasible, while every $b<\varrho_{\bm\mu}(\bm z)$ is infeasible. Since $\varrho_{\bm\mu}(\bm z)\ge0$, the feasible heights within $|b|\le\rho$ are exactly $[\varrho_{\bm\mu}(\bm z),\rho]$, proving \eqref{eq:boundary-local-parametrization}.

For (i), let $\beta>0$.

	For every boundary point $\bm w=\bm w_*(\bm\mu)+U_{\bm\mu}\bm z-\bm n_{\bm\mu}\varrho_{\bm\mu}(\bm z)$ in this neighborhood,
	\begin{equation*}
		\bm n_{\bm\mu}^\top(\bm w_*(\bm\mu)-\bm w)=\varrho_{\bm\mu}(\bm z),\qquad
		U_{\bm\mu}^\top(\bm w-\bm w_*(\bm\mu))=\bm z.
	\end{equation*}
	Consequently, the local support-distance growth condition
	\begin{equation}
		\bm n_{\bm\mu}^\top(\bm w_*(\bm\mu)-\bm w)
		\ge\frac{\underline\kappa}{1+\beta}
		\left\|U_{\bm\mu}^\top(\bm w-\bm w_*(\bm\mu))\right\|_2^{1+\beta},
		\qquad \bm w\in\partial\Omega\text{ near }\bm w_*(\bm\mu),
		\label{eq:fast_boundary_support_distance_growth}
	\end{equation}
	is equivalent to $\varrho_{\bm\mu}(\bm z)\ge\underline\kappa\|\bm z\|_2^{1+\beta}/(1+\beta)$ on a sufficiently small coordinate ball. The stronger condition
	\begin{equation*}
		\bm n_{\bm\mu}^\top(\bm w_*(\bm\mu)-\bm w)
		\ge\frac{\underline\kappa}{1+\beta}
		\|\bm w-\bm w_*(\bm\mu)\|_2^{1+\beta}
	\end{equation*}
	also suffices, since $\|U_{\bm\mu}^\top(\bm w-\bm w_*(\bm\mu))\|_2\le\|\bm w-\bm w_*(\bm\mu)\|_2$. 

For (ii), assume $0<\beta\le1$ and that $\partial\Omega$ is locally $C^{1,\beta}$.

	Following the same implicit-function argument as in the proof of
	Corollary \ref{cor:fast_principal_curvatures}, with a $C^{1,\beta}$
	defining function $F$, we obtain the parametrization similar to \eqref{eq:fast_boundary_local_graph}:
	\begin{equation*}
		\bm{w}(\bm{z})=\bm{w}_*(\bm{\mu})+U_{\bm{\mu}}\bm{z}-\bm{n}_{\bm{\mu}}\varrho_{\bm\mu}(\bm{z})\in\partial\Omega,
	\end{equation*}
	with $\varrho_{\bm\mu}\in C^1$, $\varrho_{\bm\mu}(\bm{0})=0$, and $\nabla\varrho_{\bm\mu}(\bm{0})=\bm{0}$.
	Using $G$ as defined there, differentiation gives
	\begin{equation*}
		D\varrho_{\bm\mu}(\bm{z})=\frac{D_{\bm{z}}G(\bm{z},-\varrho_{\bm\mu}(\bm{z}))}{\partial_aG(\bm{z},-\varrho_{\bm\mu}(\bm{z}))}.
	\end{equation*}
	On a sufficiently small neighborhood, the denominator is bounded away from zero. Since $DG$ is locally $\beta$-H\"older continuous and $\bm{z}\mapsto(\bm{z},-\varrho_{\bm\mu}(\bm{z}))$ is locally Lipschitz, this formula gives $\varrho_{\bm\mu}\in C^{1,\beta}$. Choose $\rho>0$ small enough that $\varrho_{\bm\mu}$ is defined and belongs to $C^{1,\beta}$ on a neighborhood of $\{\bm z\in\mathbb R^{d-1}:\|\bm z\|_2\le\rho\}$, and choose
	$\bar{\kappa}>0$ such that
	\begin{equation*}
		\|\nabla\varrho_{\bm\mu}(\bm{z})-\nabla\varrho_{\bm\mu}(\bm{z}')\|_2
		\le\bar{\kappa}\|\bm{z}-\bm{z}'\|_2^\beta,
		\qquad \|\bm{z}\|_2,\|\bm{z}'\|_2\le\rho.
	\end{equation*}
	Optimality at $\bm{\mu}$ gives
	\begin{equation*}
		-\|\bm{\mu}\|_2\varrho_{\bm\mu}(\bm{z})
		=\bm{\mu}^\top(\bm{w}(\bm{z})-\bm{w}_*(\bm{\mu}))
		\le0,
	\end{equation*}
	so $\varrho_{\bm\mu}(\bm{z})\ge0$. For $\|\bm{z}\|_2\le\rho$, integration
	along the segment from $\bm{0}$ to $\bm{z}$ yields
	\begin{equation*}
		\begin{aligned}
			\varrho_{\bm\mu}(\bm{z})
			&=\int_0^1
			(\nabla\varrho_{\bm\mu}(s\bm{z})-\nabla\varrho_{\bm\mu}(\bm{0}))^\top\bm{z}\,\dd s\\
			&\le\int_0^1
			\bar{\kappa}s^\beta\|\bm{z}\|_2^{1+\beta}\,\dd s
			=\frac{\bar{\kappa}}{1+\beta}\|\bm{z}\|_2^{1+\beta}.
		\end{aligned}
	\end{equation*}
	Thus \eqref{eq:fast_boundary_growth_assumption} and its growth bound hold with exponent $1+\beta$. 
	
\end{proof}

\subsection{Sufficient conditions for quadratic boundary growth}
\label{app:quadratic-boundary-growth}
\begin{lemma}[Quadratic boundary growth from convexity and curvature]
\label{lem:quadratic_boundary_growth}
Let $d\ge2$, let $\Omega$ be compact and convex with nonempty interior, and fix $\bm\mu\ne\bm0$ with a unique optimizer $\bm w_*(\bm\mu)$. Retain $\bm n_{\bm\mu}$ and $U_{\bm\mu}$ from Proposition~\ref{prop:boundary-local-parametrization}. Each of the following conditions implies the local representation \eqref{eq:boundary-local-parametrization} and, for some $\rho>0$ and $0<\underline\kappa\le\bar\kappa<\infty$,
\begin{equation*}
\frac{\underline\kappa}{2}\|\bm z\|_2^2
\le\varrho_{\bm\mu}(\bm z)\le\frac{\bar\kappa}{2}\|\bm z\|_2^2,
\qquad \|\bm z\|_2\le\rho.
\end{equation*}
\begin{enumerate}[label=(\roman*),itemsep=2pt]
\item $\partial\Omega$ is $C^{1,1}$ near $\bm w_*(\bm\mu)$, and $\Omega$ is strongly convex, i.e., for some $c>0$,
	\begin{equation*}
		(1-s)\bm w+s\bm w'
		+c\,s(1-s)\|\bm w-\bm w'\|_2^2B_2(1)\subseteq\Omega,
		\qquad \bm w,\bm w'\in\Omega,\quad s\in[0,1].
	\end{equation*}
In this case, one may take $\underline\kappa=2c$.
\item $\partial\Omega$ is $C^2$ near $\bm w_*(\bm\mu)$ and $\mathscr S_{\bm\mu}\succ0$.
\end{enumerate}
In case (ii), $\varrho_{\bm\mu}(\bm z)=\frac12\bm z^\top\mathscr S_{\bm\mu}\bm z+o(\|\bm z\|_2^2)$. \eqref{eq:fast_boundary_profile_limit} holds with $\beta=1$ and $\varrho_0(\bm z)=\bm z^\top\mathscr S_{\bm\mu}\bm z/2$. Its convex conjugate and the resulting support-remainder function in Theorem~\ref{thm:geometry-contact-rate} are
\begin{equation*}
\varrho_0^*(\bm y)=\frac12\bm y^\top\mathscr S_{\bm\mu}^{-1}\bm y,\qquad
\psi(\bm t)=\frac{1}{2\|\bm\mu\|_2}\bm t^\top U_{\bm\mu}\mathscr S_{\bm\mu}^{-1}U_{\bm\mu}^\top\bm t
=\frac12\bm t^\top H_{\bm\mu}\bm t.
\end{equation*}
\end{lemma}

\begin{proof}
In either case, $\partial\Omega$ is locally $C^1$, so the supporting hyperplane at $\bm w_*(\bm\mu)$ is unique and Proposition~\ref{prop:boundary-local-parametrization} supplies the local boundary profile.

	For (i), fix $\bm w'\in\Omega$ and $s\in(0,1)$. Apply the strong-convexity condition with $\bm w=\bm w_*(\bm\mu)$. Since $\|\bm n_{\bm\mu}\|_2=1$, the point on the resulting ball in the direction $\bm n_{\bm\mu}$ satisfies
	\begin{equation*}
		\bm v:=(1-s)\bm w_*(\bm\mu)+s\bm w'
		+c\,s(1-s)\|\bm w'-\bm w_*(\bm\mu)\|_2^2\bm n_{\bm\mu}\in\Omega.
	\end{equation*}
	Optimality of $\bm w_*(\bm\mu)$ gives $\bm\mu^\top(\bm v-\bm w_*(\bm\mu))\le0$. Dividing by $\|\bm\mu\|_2$ and using $\bm n_{\bm\mu}^\top\bm n_{\bm\mu}=1$, we obtain
	\begin{equation*}
		\begin{aligned}
			0&\ge\bm n_{\bm\mu}^\top(\bm v-\bm w_*(\bm\mu))\\
			&=s\bm n_{\bm\mu}^\top(\bm w'-\bm w_*(\bm\mu))
			+c\,s(1-s)\|\bm w'-\bm w_*(\bm\mu)\|_2^2.
		\end{aligned}
	\end{equation*}
	Divide by $s$ and let $s\downarrow0$. Then
	\begin{equation*}
		\bm n_{\bm\mu}^\top(\bm w_*(\bm\mu)-\bm w')
		\ge c\|\bm w'-\bm w_*(\bm\mu)\|_2^2,
		\qquad \varrho_{\bm\mu}(\bm z)\ge c(\|\bm z\|_2^2+\varrho_{\bm\mu}(\bm z)^2)\ge c\|\bm z\|_2^2.
	\end{equation*}
	This proves the lower bound with $\underline\kappa=2c$. Proposition~\ref{prop:boundary-local-parametrization}(ii), applied with $\beta=1$, supplies the upper bound from the local $C^{1,1}$ regularity.
	
	For (ii), the calculation in the proof of Corollary~\ref{cor:fast_principal_curvatures} gives $\varrho_{\bm\mu}\in C^2$ locally and $\nabla^2\varrho_{\bm\mu}(\bm0)=\mathscr S_{\bm\mu}$. Since $\varrho_{\bm\mu}(\bm0)=0$ and $\nabla\varrho_{\bm\mu}(\bm0)=\bm0$, Taylor expansion gives
	\begin{equation*}
		\varrho_{\bm\mu}(\bm z)=\frac12\bm z^\top \mathscr{S}_{\bm\mu}\bm z+o(\|\bm z\|_2^2).
	\end{equation*}
	The eigenvalues of $\mathscr{S}_{\bm\mu}$ are positive and finite, so \eqref{eq:fast_two_sided_boundary_growth} holds with $\beta=1$ and suitable $\underline\kappa,\bar\kappa>0$. The Taylor remainder also implies, for every $L>0$,
	\begin{equation*}
		\sup_{\|\bm z\|_2\le L}\left|\frac{\varrho_{\bm\mu}(\delta\bm z)}{\delta^2}-\frac12\bm z^\top\mathscr S_{\bm\mu}\bm z\right|\longrightarrow0,
		\qquad \delta\downarrow0,
	\end{equation*}
	which proves \eqref{eq:fast_boundary_profile_limit}. Completing the square gives
\begin{equation*}
\begin{aligned}
\bm y^\top\bm z-\frac12\bm z^\top\mathscr S_{\bm\mu}\bm z
&=\frac12\bm y^\top\mathscr S_{\bm\mu}^{-1}\bm y\\
&\quad-\frac12(\bm z-\mathscr S_{\bm\mu}^{-1}\bm y)^\top\mathscr S_{\bm\mu}(\bm z-\mathscr S_{\bm\mu}^{-1}\bm y).
\end{aligned}
\end{equation*}
The maximum over $\bm z$ is attained at $\mathscr S_{\bm\mu}^{-1}\bm y$, yielding the stated formula for $\varrho_0^*$. Substitution into \eqref{eq:fast_boundary_conjugate_profile} and the identity \eqref{eq:fast_hessian_principal_curvatures} give
\begin{equation*}
\psi(\bm t)=\frac{1}{2\|\bm\mu\|_2}\bm t^\top U_{\bm\mu}\mathscr S_{\bm\mu}^{-1}U_{\bm\mu}^\top\bm t
=\frac12\bm t^\top H_{\bm\mu}\bm t.
\end{equation*}
Thus the conjugate limit recovers the principal-curvature formula in the smooth case.
\end{proof}

\section{Detailed Calculations for the Gaussian Examples}
\label{app:examples}

\subsection{Rectangular constraints: details for Example~\ref{ex:fast_rectangle}}
\label{app:rectangle-details}
For the rectangle $\Omega=[0,1]^d$, we have
$f_*(\bm{\mu})=\sum_{j=1}^d(\mu_j)_+$. Write
$\mathcal{S}_0=\{j\in[d]:\mu_j=0\}$ and
$s=\|\bm{\mu}\|_0=d-|\mathcal{S}_0|$, so that $s$ is the number of
nonzero coordinates. For every fixed $M$, when $\varepsilon$ is sufficiently small, 
$\varepsilon M<\min_{j\notin\mathcal{S}_0}|\mu_j|$, with the minimum of the
empty set interpreted as $+\infty$, every $\bm{t}\in B_2(M)$ satisfies
\begin{equation*}
	|\varepsilon t_j|
	\le\varepsilon\|\bm{t}\|_2
	\le\varepsilon M
	<\min_{\ell\notin\mathcal{S}_0}|\mu_\ell|
	\le|\mu_j|,
	\qquad j\notin\mathcal{S}_0.
\end{equation*}
Thus $\mu_j+\varepsilon t_j$ has the same sign as $\mu_j$ for every nonzero coordinate, uniformly over $\bm{t}\in B_2(M)$. Coordinatewise optimality gives $w_{*,j}(\bm{\mu})=1$ when $\mu_j>0$ and $w_{*,j}(\bm{\mu})=0$ when $\mu_j<0$; when $\mu_j=0$, the selected maximizer has $w_{*,j}(\bm{\mu})\in[0,1]$. Expanding the definition of $r_{\bm{\mu}}$ yields
\begin{equation*}
	\begin{aligned}
		r_{\bm{\mu}}(\varepsilon\bm{t})
		&=f_*(\bm{\mu}+\varepsilon\bm{t})-f_*(\bm{\mu})
		-\varepsilon\bm{w}_*(\bm{\mu})^\top\bm{t}\\
		&=\sum_{j=1}^d\left[
		(\mu_j+\varepsilon t_j)_+-(\mu_j)_+
		-\varepsilon w_{*,j}(\bm{\mu})t_j\right].
	\end{aligned}
\end{equation*}
When $\varepsilon$ is sufficiently small, the coordinate contributions are
\begin{equation*}
		(\mu_j+\varepsilon t_j)_+-(\mu_j)_+
		-\varepsilon w_{*,j}(\bm{\mu})t_j=
		\begin{cases}
			(\mu_j+\varepsilon t_j)-\mu_j-\varepsilon t_j=0,
			&\mu_j>0,\\
			0-0-0=0,&\mu_j<0,\\
			\varepsilon(t_j)_+-\varepsilon w_{*,j}(\bm{\mu})t_j,
			&\mu_j=0,
		\end{cases}
\end{equation*}
where $(\varepsilon t_j)_+=\varepsilon(t_j)_+$ since $\varepsilon>0$.
All nonzero coordinates therefore contribute zero, leaving exactly
\begin{equation}
	r_{\bm{\mu}}(\varepsilon\bm{t})
	=\varepsilon\sum_{j\in\mathcal{S}_0}
	\left[(t_j)_+-w_{*,j}(\bm{\mu})t_j\right],
	\qquad \bm{t}\in B_2(M).
	\label{eq:fast_rectangle_remainder}
\end{equation}
This is an exact identity on the whole local ball. If $\bm{\mu}=\bm{0}$,
then $\mathcal{S}_0=[d]$ and the minimum above is $+\infty$, so the
identity holds for every $\varepsilon>0$.
As in Example~\ref{ex:fast_rectangle}, write $\psi(\bm{t})=\sum_{j\in\mathcal{S}_0}[(t_j)_+-w_{*,j}(\bm{\mu})t_j]$, so that $r_{\bm{\mu}}(\varepsilon\bm{t})=\varepsilon\psi(\bm{t})$ on the local ball.
Suppose that $M/\sqrt n<\min_{j\notin\mathcal{S}_0}|\mu_j|$.
By \eqref{eq:fast_exact_jensen_gap}, the local minimax risk is
\begin{equation*}
	R_{n,M}^{\Sigma}(\bm{\mu})
	=\sup_{\pi\in\mathcal{P}(B_2(M))}
	\mathbb{E}_{\pi,\Sigma}\left[
	r_{\bm{\mu}}\left(\frac{\bm{T}}{\sqrt n}\right)
	-r_{\bm{\mu}}\left(
	\frac{\bm{m}_{\pi,\Sigma}(\bm{Y})}{\sqrt n}\right)\right].
\end{equation*}
Here $\bm{T}\sim\pi$, $\bm{Y}=\bm{T}+\bm{Z}$ with independent
$\bm{Z}\sim\mathcal{N}(\bm{0},\Sigma)$, and
$\bm{m}_{\pi,\Sigma}(\bm{Y})=\mathbb{E}_{\pi,\Sigma}[\bm{T}\mid\bm{Y}]$.
For every such prior, $\bm{T}\in B_2(M)$ almost surely, and
\begin{equation*}
	\|\bm{m}_{\pi,\Sigma}(\bm{Y})\|_2
	\le\mathbb{E}_{\pi,\Sigma}[\|\bm{T}\|_2\mid\bm{Y}]
	\le M.
\end{equation*}
Applying \eqref{eq:fast_rectangle_remainder} at scale $n^{-1/2}$ gives
\begin{equation*}
	r_{\bm{\mu}}\left(\frac{\bm{t}}{\sqrt n}\right)
	=\frac1{\sqrt n}\psi(\bm{t}),
	\qquad \bm{t}\in B_2(M).
\end{equation*}
This formula therefore applies to both $\bm{T}$ and
$\bm{m}_{\pi,\Sigma}(\bm{Y})$. Substituting it into their expected remainder
difference yields
\begin{equation*}
	\begin{aligned}
		\mathbb{E}_{\pi,\Sigma}\left[
		r_{\bm{\mu}}\left(\frac{\bm{T}}{\sqrt n}\right)
		-r_{\bm{\mu}}\left(\frac{\bm{m}_{\pi,\Sigma}(\bm{Y})}{\sqrt n}\right)
		\right]&=\frac1{\sqrt n}\mathbb{E}_{\pi,\Sigma}\left[
		\psi(\bm{T})-\psi(\bm{m}_{\pi,\Sigma}(\bm{Y}))\right]\\
		&=\frac1{\sqrt n}\mathbb{E}_{\pi,\Sigma}\left[
		\sum_{j\in\mathcal S_0}\bigl((T_j)_+-(m_{\pi,\Sigma,j}(\bm Y))_+\bigr)\right]\\
		&\quad-\frac1{\sqrt n}\mathbb{E}_{\pi,\Sigma}\left[
		\sum_{j\in\mathcal{S}_0}w_{*,j}(\bm{\mu})
		(T_j-m_{\pi,\Sigma,j}(\bm{Y}))\right].
	\end{aligned}
\end{equation*}
The last expectation vanishes by the tower property: for every $j$,
\begin{equation*}
	\mathbb{E}_{\pi,\Sigma}[m_{\pi,\Sigma,j}(\bm{Y})]
	=\mathbb{E}_{\pi,\Sigma}[\mathbb{E}_{\pi,\Sigma}[T_j\mid\bm{Y}]]
	=\mathbb{E}_{\pi,\Sigma}[T_j],
\end{equation*}
and the coefficients $w_{*,j}(\bm{\mu})$ are fixed. Therefore, by definition of $V_{M,\Sigma}$,
\begin{equation}
		R_{n,M}^{\Sigma}(\bm{\mu})
		=\frac1{\sqrt n}
		\sup_{\pi\in\mathcal{P}(B_2(M))}
		\mathbb{E}_{\pi,\Sigma}\left[
		\psi(\bm{T})-\psi(\bm{m}_{\pi,\Sigma}(\bm{Y}))\right]=\frac1{\sqrt n}V_{M,\Sigma}(\psi).
	\label{eq:fast_rectangle_exact_minimax}
\end{equation}
This equality is exact whenever the sign condition above holds. If $s<d$, then for every fixed $M>0$ and all
sufficiently large $n$,
\begin{equation}
	c\,\frac{\min\{M\sqrt{d-s},d-s\}}{\sqrt{n}}
	\le R_{n,M}^{\Sigma}(\bm{\mu})
	\le C\,\frac{\min\{M\sqrt{d-s},d-s\}}{\sqrt{n}}.
	\label{eq:fast_rectangle_dimension_bound}
\end{equation}
where the constants can be chosen as $c=\Phi(-1)\min\{1,\sqrt{\underline\lambda}\}$ and $C=\max\{1,\sqrt{\overline\lambda}\}$. Here $\Phi$ and $\phi$ denote the standard normal distribution function and density. We prove the lower bound by evaluating one prior in $V_{M,\Sigma}(\psi)$, and the upper bound using the sample-only plug-in policy.

For the lower bound, first note that $V_{M,\Sigma}(\psi)\ge V_{M,\underline\lambda I_d}(\psi)$ because $\Sigma\succeq\underline\lambda I_d$. Indeed, write $\bm Y_\Sigma=\bm Y_1+\bm Z'$ with $\bm Y_1=\bm T+\bm Z_1$, $\bm Z_1\sim\mathcal N(\bm0,\underline\lambda I_d)$, and $\bm Z'\sim\mathcal N(\bm0,\Sigma-\underline\lambda I_d)$ independent of $(\bm T,\bm Z_1)$. Since $\bm T\to\bm Y_1\to\bm Y_\Sigma$ is a Markov chain, $\bm m_{\pi,\Sigma}(\bm Y_\Sigma)=\mathbb E[\bm m_{\pi,\underline\lambda I_d}(\bm Y_1)\mid\bm Y_\Sigma]$, and Jensen's inequality for the convex function $\psi$ gives $\mathbb E\psi(\bm m_{\pi,\Sigma}(\bm Y_\Sigma))\le\mathbb E\psi(\bm m_{\pi,\underline\lambda I_d}(\bm Y_1))$. Hence the Bayes gap under $\Sigma$ is at least that under $\underline\lambda I_d$ for every prior $\pi$, and taking suprema over $\pi$ proves the claim.
Evaluate the latter constant in the experiment
$\bm Y=\bm T+\bm Z$ with $\bm Z\sim\mathcal N(\bm0,\underline\lambda I_d)$.
Set
\begin{equation*}
	h=\min\left\{\sqrt{\underline\lambda},\frac{M}{\sqrt{d-s}}\right\},
	\qquad
	\mathbb{P}_{\pi,\underline\lambda I_d}\{T_j=h\}=\mathbb{P}_{\pi,\underline\lambda I_d}\{T_j=-h\}=\frac12,
	\quad j\in\mathcal{S}_0,
\end{equation*}
with these coordinates independent and $T_j=0$ for
$j\notin\mathcal{S}_0$. Every point in the support satisfies
\begin{equation*}
	\|\bm{T}\|_2^2=(d-s)h^2\le M^2,
\end{equation*}
so $\pi\in\mathcal{P}(B_2(M))$. When $\bm{Y}=\bm{T}+\bm{Z}$, independence across coordinates implies that,
for $j\in\mathcal{S}_0$, the posterior of $T_j$ depends only on $Y_j$.
Bayes' formula gives
\begin{equation*}
	\begin{aligned}
		\frac{\mathbb{P}_{\pi,\underline\lambda I_d}\{T_j=h\mid\bm{Y}\}}
		{\mathbb{P}_{\pi,\underline\lambda I_d}\{T_j=-h\mid\bm{Y}\}}
		&=\frac{\phi((Y_j-h)/\sqrt{\underline\lambda})}{\phi((Y_j+h)/\sqrt{\underline\lambda})}=e^{2hY_j/\underline\lambda},\\
		m_{\pi,\underline\lambda I_d,j}(\bm{Y})
		&=h\frac{e^{hY_j/\underline\lambda}-e^{-hY_j/\underline\lambda}}{e^{hY_j/\underline\lambda}+e^{-hY_j/\underline\lambda}}.
	\end{aligned}
\end{equation*}
In particular, $m_{\pi,\underline\lambda I_d,j}(\bm{Y})$ has the same sign as $Y_j$, and the
tower property yields
\begin{equation*}
	\begin{aligned}
		\mathbb{E}_{\pi,\underline\lambda I_d}[(m_{\pi,\underline\lambda I_d,j}(\bm{Y}))_+]
		=\mathbb{E}_{\pi,\underline\lambda I_d}[m_{\pi,\underline\lambda I_d,j}(\bm{Y})\mathbb{I}\{Y_j>0\}]=\mathbb{E}_{\pi,\underline\lambda I_d}[T_j\mathbb{I}\{Y_j>0\}].
	\end{aligned}
\end{equation*}
Since $Y_j\mid (T_j=\pm h)\sim\mathcal{N}(\pm h,\underline\lambda)$, the expected gap in coordinate $j$ is exactly
\begin{equation*}
	\begin{aligned}
		\mathbb{E}_{\pi,\underline\lambda I_d}[(T_j)_+-(m_{\pi,\underline\lambda I_d,j}(\bm{Y}))_+]&=\mathbb{E}_{\pi,\underline\lambda I_d}[(T_j)_+-T_j\mathbb{I}\{Y_j>0\}]\\
		&=h\mathbb{P}_{\pi,\underline\lambda I_d}\{T_j=h,Y_j\le0\}
		+h\mathbb{P}_{\pi,\underline\lambda I_d}\{T_j=-h,Y_j>0\}\\
		&=\frac h2\Phi(-h/\sqrt{\underline\lambda})+\frac h2\Phi(-h/\sqrt{\underline\lambda})=h\Phi(-h/\sqrt{\underline\lambda}).
	\end{aligned}
\end{equation*}
Summing over $\mathcal{S}_0$, using the cancellation of the linear terms established above and $h\le\sqrt{\underline\lambda}$, we obtain
\begin{equation*}
	\begin{aligned}
		V_{M,\Sigma}(\psi)
		&\ge\mathbb{E}_{\pi,\underline\lambda I_d}[\psi(\bm{T})-\psi(\bm{m}_{\pi,\underline\lambda I_d}(\bm{Y}))]\\
		&=(d-s)h\Phi(-h/\sqrt{\underline\lambda})\\
		&\ge\Phi(-1)\min\{M\sqrt{d-s},\sqrt{\underline\lambda}(d-s)\}.
	\end{aligned}
	\label{eq:fast_rectangle_jensen_lower}
\end{equation*}
Since $\min\{M\sqrt{d-s},\sqrt{\underline\lambda}(d-s)\}
\ge\min\{1,\sqrt{\underline\lambda}\}\min\{M\sqrt{d-s},d-s\}$,
together with \eqref{eq:fast_rectangle_exact_minimax}, this proves the
lower bound in \eqref{eq:fast_rectangle_dimension_bound}.

For the upper bound, use the plug-in policy
$\hat{\bm{w}}=\bm{w}_*(\bar{\bm{x}})$, whose coordinates are $\hat w_j=\mathbb{I}\{\bar x_j>0\}$. At a true mean $\bm{\mu}_{n,\bm{t}}$, choosing the wrong sign in coordinate $j$ costs $|\mu_{n,\bm{t},j}|$ and occurs with probability
$\Phi(-\sqrt n|\mu_{n,\bm{t},j}|/\sqrt{\Sigma_{jj}})$. Thus
\begin{equation*}
	\begin{aligned}
		\mathbb{E}_{\bm\mu_{n,\bm{t}},\Sigma}R(\hat{\bm{w}})
		&=\sum_{j=1}^d|\mu_{n,\bm{t},j}|
		\Phi(-\sqrt n|\mu_{n,\bm{t},j}|/\sqrt{\Sigma_{jj}})\\
		&=\frac1{\sqrt n}\sum_{j\in\mathcal{S}_0}|t_j|\Phi(-|t_j|/\sqrt{\Sigma_{jj}})+\sum_{j\notin\mathcal{S}_0}
		\left|\mu_j+\frac{t_j}{\sqrt n}\right|
		\Phi\left(-\frac{\sqrt n}{\sqrt{\Sigma_{jj}}}\left|\mu_j+\frac{t_j}{\sqrt n}\right|\right).
	\end{aligned}
	\label{eq:fast_rectangle_plugin_risk}
\end{equation*}
For the first sum, Cauchy--Schwarz gives
\begin{equation*}
	\sum_{j\in\mathcal{S}_0}|t_j|\Phi(-|t_j|/\sqrt{\Sigma_{jj}})
	\le\frac12\sum_{j\in\mathcal{S}_0}|t_j|
	\le\frac12\sqrt{d-s}\|\bm{t}\|_2
	\le\frac12 M\sqrt{d-s}.
\end{equation*}
Also, for every $x\ge0$, by the Mills' ratio,
\begin{equation*}
	x\Phi(-x)\le\phi(x)\le\frac1{\sqrt{2\pi}}<\frac12.
\end{equation*}
In particular,
$|t_j|\Phi(-|t_j|/\sqrt{\Sigma_{jj}})
\le\sqrt{\Sigma_{jj}}/\sqrt{2\pi}\le\sqrt{\overline\lambda}/2$.
There are $d-s$ terms, so these two bounds imply
\begin{equation*}
	\sup_{\bm{t}\in B_2(M)}
	\sum_{j\in\mathcal{S}_0}|t_j|\Phi(-|t_j|/\sqrt{\Sigma_{jj}})
	\le\frac12\min\{M\sqrt{d-s},\sqrt{\overline\lambda}(d-s)\}.
\end{equation*}
If $s>0$, write $\delta=\min_{j\notin\mathcal{S}_0}|\mu_j|>0$.
Whenever $M/\sqrt n\le\delta/2$, every $\bm{t}\in B_2(M)$ satisfies
\begin{equation*}
	\left|\mu_j+\frac{t_j}{\sqrt n}\right|
	\ge|\mu_j|-\frac{M}{\sqrt n}\ge\frac\delta2,
	\qquad j\notin\mathcal{S}_0.
\end{equation*}
Using $x\Phi(-x)\le\phi(x)$ again bounds each term in the second sum by
\begin{equation*}
	\left|\mu_j+\frac{t_j}{\sqrt n}\right|
	\Phi\left(-\frac{\sqrt n}{\sqrt{\Sigma_{jj}}}\left|\mu_j+\frac{t_j}{\sqrt n}\right|\right)
	\le\frac{\sqrt{\overline\lambda}}{\sqrt{2\pi n}}e^{-n\delta^2/(8\overline\lambda)}.
\end{equation*}
Consequently,
\begin{equation*}
	R_{n,M}^{\Sigma}(\bm{\mu})\le\sup_{\bm{t}\in B_2(M)}
	\mathbb{E}_{\bm\mu_{n,\bm{t}},\Sigma}R(\bm{w}_*(\bar{\bm{x}}))\le\frac{\min\{M\sqrt{d-s},\sqrt{\overline\lambda}(d-s)\}}{2\sqrt n}
	+\frac{s\sqrt{\overline\lambda}}{\sqrt{2\pi n}}e^{-n\delta^2/(8\overline\lambda)}.
	\label{eq:fast_rectangle_plugin_upper}
\end{equation*}
For $0<s<d$ and fixed $d,\bm\mu,M$, the exponential term is eventually at most $\min\{M\sqrt{d-s},\sqrt{\overline\lambda}(d-s)\}/(2\sqrt n)$. If $s=0$, that term is absent. The upper bound is uniform over $\Sigma\in\mathcal C$ and uses the same plug-in rule. Combining it with the preceding lower bound gives, for $s<d$ and all sufficiently large $n$,
\begin{equation}
	\frac{\Phi(-1)\min\{M\sqrt{d-s},\sqrt{\underline\lambda}(d-s)\}}{\sqrt n}
	\le R_{n,M}^{\Sigma}(\bm\mu)
	\le\frac{\min\{M\sqrt{d-s},\sqrt{\overline\lambda}(d-s)\}}{\sqrt n}.
	\label{eq:fast_rectangle_rates_explicit}
\end{equation}
Thus \eqref{eq:fast_rectangle_dimension_bound} holds with the stated constants, and the same plug-in policy attains the local minimax rate.

If $s<d$, then for any $j\in\mathcal S_0$ and $M>0$, $\psi(M\bm e_j)+\psi(-M\bm e_j)=M>2\psi(\bm0)=0$. Thus $\psi$ is nonaffine on $B_2(M)$, and $V_{M,\Sigma}(\psi)>0$ by Theorem~\ref{thm:fast_general_local_rate}. The exact identity \eqref{eq:fast_rectangle_exact_minimax} in particular gives
\begin{equation*}
	\lim_{n\to\infty}\sqrt n R_{n,M}^{\Sigma}(\bm\mu)
	=V_{M,\Sigma}(\psi).
\end{equation*}

Now we consider the case $s=d$. In this case $\mathcal S_0=\varnothing$ and $\delta=\min_{1\le j\le d}|\mu_j|>0$. If $M/\sqrt n<\delta$, then $\mu_j+t_j/\sqrt n$ shares the same sign as $\mu_j$ for every $j$ and every $\bm t\in B_2(M)$. Consequently,
\begin{equation*}
	w_{*,j}(\bm\mu_{n,\bm t})=w_{*,j}(\bm\mu)=\mathbb I\{\mu_j>0\},\qquad j=1,\ldots,d,\quad\bm t\in B_2(M).
\end{equation*}
The constant policy $\bm w_*(\bm\mu)$ is admissible in the local benchmark and is optimal throughout the local ball. Therefore,
\begin{equation*}
	0\le R_{n,M}^{\Sigma}(\bm\mu)\le\sup_{\bm t\in B_2(M)}\mathbb E_{\bm\mu_{n,\bm t},\Sigma}R(\bm w_*(\bm\mu))=0.
\end{equation*}
Equivalently, $\psi\equiv0$ and \eqref{eq:fast_rectangle_exact_minimax} gives the same conclusion.

For SAA, coordinatewise optimality gives $w_{*,j}(\bar{\bm x})=\mathbb I\{\bar x_j>0\}$. If $\mu_{n,\bm t,j}>0$, choosing $0$ instead of $1$ loses $\mu_{n,\bm t,j}$; if $\mu_{n,\bm t,j}<0$, choosing $1$ instead of $0$ loses $-\mu_{n,\bm t,j}$. Thus
\begin{equation*}
	R(\bm w_*(\bar{\bm x}))
	=\sum_{j=1}^d|\mu_{n,\bm t,j}|\mathbb I\{\bar x_j\mu_{n,\bm t,j}<0\}.
\end{equation*}
Under the true mean $\bm\mu_{n,\bm t}$, $\bar x_j\sim\mathcal N\left(\mu_{n,\bm t,j},\Sigma_{jj}/n\right)$. Hence
\begin{equation*}
	\begin{aligned}
		\mathbb P_{\bm\mu_{n,\bm t},\Sigma}\{\bar x_j\mu_{n,\bm t,j}<0\}
		&=\begin{cases}
			\Phi(-\sqrt n\,\mu_{n,\bm t,j}/\sqrt{\Sigma_{jj}}),&\mu_{n,\bm t,j}>0,\\
			1-\Phi(-\sqrt n\,\mu_{n,\bm t,j}/\sqrt{\Sigma_{jj}}),&\mu_{n,\bm t,j}<0
		\end{cases}\\
		&=\Phi\left(-\frac{\sqrt n|\mu_{n,\bm t,j}|}{\sqrt{\Sigma_{jj}}}\right).
	\end{aligned}
\end{equation*}
Moreover, if $M/\sqrt n\le\delta/2$, then for every $\bm t\in B_2(M)$ and every $j$,
\begin{equation*}
	|\mu_{n,\bm t,j}|=\left|\mu_j+\frac{t_j}{\sqrt n}\right|
	\ge|\mu_j|-\frac{|t_j|}{\sqrt n}
	\ge\delta-\frac M{\sqrt n}\ge\frac\delta2.
\end{equation*}
For $x\ge0$, the Mills' ratio gives $x\Phi(-x)\le\phi(x)$. Applying this inequality with $x=\sqrt n|\mu_{n,\bm t,j}|/\sqrt{\Sigma_{jj}}$ gives
\begin{equation*}
|\mu_{n,\bm t,j}|\Phi\left(-\frac{\sqrt n|\mu_{n,\bm t,j}|}{\sqrt{\Sigma_{jj}}}\right)
\le\frac{\sqrt{\Sigma_{jj}}}{\sqrt{2\pi n}}\exp\left(-\frac{n|\mu_{n,\bm t,j}|^2}{2\Sigma_{jj}}\right)\le\frac{\sqrt{\overline\lambda}}{\sqrt{2\pi n}}\exp\left(-\frac{n\delta^2}{8\overline\lambda}\right).
\end{equation*}
Here the second inequality uses $\Sigma_{jj}\le\overline\lambda$ and
\begin{equation*}
	\frac{n|\mu_{n,\bm t,j}|^2}{2\Sigma_{jj}}
	\ge\frac{n(\delta/2)^2}{2\overline\lambda}
	=\frac{n\delta^2}{8\overline\lambda}.
\end{equation*}
By linearity of expectation, summing the coordinate losses yields the following bound, regardless of dependence between coordinates:
\begin{equation*}
	\sup_{\bm t\in B_2(M)}\mathbb E_{\bm\mu_{n,\bm t},\Sigma}R(\bm w_*(\bar{\bm x}))
	\le\frac{d\sqrt{\overline\lambda}}{\sqrt{2\pi n}}\exp\left(-\frac{n\delta^2}{8\overline\lambda}\right),\qquad\Sigma\in\mathcal C.
\end{equation*}
For fixed $d$, $\bm\mu$, and $M$, this bound holds for every $n\ge4M^2/\delta^2$, uniformly over $\Sigma\in\mathcal C$. Since $\delta^2/(8\overline\lambda)>0$, SAA has exponentially decaying maximal expected regret on each fixed local ball.

\subsection{A strongly convex set with center-dependent rates: details for Example~\ref{ex:fast_center_dependent_rates}}
\label{app:l4/3-value-decision}
Consider the strongly convex set
\begin{equation*}
	\Omega=B_{4/3}^{\,2}
	=\{\bm{w}:|w_1|^{4/3}+|w_2|^{4/3}\le1\},
\end{equation*}
which gives $f_*(\bm{\mu})=(\mu_1^4+\mu_2^4)^{1/4}$.
The local minimax rate depends on the reference center $\bm{\mu}$ even though the feasible set remains fixed. We first take
$\bm{\mu}=\bm{e}_1$, where $\bm{w}_*(\bm{e}_1)=\bm{e}_1$ and
$H_{\bm{e}_1}=\nabla^2f_*(\bm{e}_1)=0$. Uniformly for bounded $\bm{t}$,
\begin{equation*}
	f_*(\bm{e}_1+\varepsilon\bm{t})
	=1+\varepsilon t_1+\frac{\varepsilon^4t_2^4}{4}
	+o(\varepsilon^4).
\end{equation*}
Theorem \ref{thm:fast_general_local_rate} therefore applies with
$a(\varepsilon)=\varepsilon^4$ and $\psi(\bm{t})=t_2^4/4$, yielding
\begin{equation*}
	R_{n,M}^{\Sigma}(\bm{e}_1)\asymp_M n^{-2}
	\qquad\text{for every fixed }M>0.
	\label{eq:fast_quartic_example}
\end{equation*}
More generally, at any fixed nonzero point on either coordinate axis,
the first nonzero local remainder is quartic in the transverse
coordinate, giving the same $n^{-2}$ rate.

Next take $\bm{\mu}=2^{-1/4}(1,1)^\top$, so that
$f_*(\bm{\mu})=1$ as before. Direct differentiation gives
\begin{equation*}
	H_{\bm{\mu}}
	=\frac{3}{2\sqrt{2}}
	\begin{pmatrix}1&-1\\-1&1\end{pmatrix},
	\qquad
	r_{\bm{\mu}}(\varepsilon\bm{t})
	=\varepsilon^2\frac{3(t_1-t_2)^2}{4\sqrt{2}}
	+o(\varepsilon^2),
\end{equation*}
uniformly for bounded $\bm{t}$. The limit
$\psi(\bm{t})=3(t_1-t_2)^2/(4\sqrt{2})$ is nonaffine on every
$B_2(M)$ with $M>0$. Therefore,
\begin{equation*}
	R_{n,M}^{\Sigma}\bigl(2^{-1/4}(1,1)^\top\bigr)\asymp_M n^{-1}.
	\label{eq:fast_quadratic_center_example}
\end{equation*}
Theorem \ref{thm:fast_local_minimax} also gives the sharp constant when $M\to\infty$:
\begin{equation*}
	\lim_{M\to\infty}\lim_{n\to\infty}
	nR_{n,M}^{\Sigma}\bigl(2^{-1/4}(1,1)^\top\bigr)
	=\frac12\operatorname{tr}(\Sigma H_{\bm{\mu}})
	=\frac{3}{4\sqrt{2}}(\Sigma_{11}+\Sigma_{22}-2\Sigma_{12}).
\end{equation*}
At any fixed center with $\mu_1\mu_2\ne0$, the Hessian is nonzero and
the same $n^{-1}$ rate holds. Here the ambiguity set is $A=\{\bm{0}\}$
and, under the polar-norm convention, $\dist(\bm{\mu},A)=f_*(\bm{\mu})$.
In particular, the two displayed nonzero centers have the same distance
$\dist(\bm{\mu},A)=1$, yet their local minimax rates differ.

Finally, take $\bm{\mu}=\bm{0}$. With the convention
$\bm{w}_*(\bm{0})=\bm{0}$, positive homogeneity gives the exact identity
\begin{equation*}
	r_{\bm{0}}(\varepsilon\bm{t})
	=\varepsilon\|\bm{t}\|_4,
	\qquad
	R_{n,M}^{\Sigma}(\bm{0})
	=\frac{1}{\sqrt{n}}V_{M,\Sigma}(\bm{t}\mapsto\|\bm{t}\|_4)
	\asymp_M n^{-1/2}.
	\label{eq:fast_origin_center_example}
\end{equation*}
The positivity follows because $\|\bm{t}\|_4$ is nonaffine on every
$B_2(M)$ with $M>0$. Thus this single strongly convex set exhibits
$n^{-2}$, $n^{-1}$, and $n^{-1/2}$ local rates at different centers.

All of these rates concern neighborhoods of radius $M/\sqrt{n}$ with
$M$ fixed. In contrast, every fixed-radius neighborhood of $\bm{e}_1$
contains interior points with nonzero Hessian and consequently has an
$n^{-1}$ minimax lower bound. If its closure excludes the origin,
boundedness of the Hessian on a larger neighborhood and Gaussian tail
bounds also give an $O(n^{-1})$ upper bound for the plug-in policy.

More generally, for $d\ge2$, consider
\begin{equation*}
	\Omega=B_{4/3}^d
	=\left\{\bm{w}:\sum_{j=1}^d|w_j|^{4/3}\le1\right\},
\end{equation*}
which gives $f_*(\bm{\mu})=\|\bm{\mu}\|_4$.
For $\bm{\mu}\ne\bm{0}$, denote
\begin{equation*}
	\mathcal{S}_0=\{j\in[d]:\mu_j=0\},\qquad
	s=\|\bm{\mu}\|_0=d-|\mathcal{S}_0|.
\end{equation*}
The maximizer is unique, with $w_{*,j}(\bm{\mu})=\mu_j^3/\|\bm{\mu}\|_4^3$.
Thus $A=\{\bm{0}\}$ and the polar-norm convention gives
$\dist(\bm{\mu},A)=\|\bm{\mu}\|_4$.

The strong-convexity constant can be chosen independently of $d$. To verify this, set $F(\bm{w})=\|\bm{w}\|_{4/3}^2/2$. At points with nonzero coordinates,
differentiation and H\"older's inequality give
\begin{equation*}
	\bm{z}^\top\nabla^2F(\bm{w})\bm{z}
	\ge\frac13\|\bm{w}\|_{4/3}^{2/3}
	\sum_{j=1}^d|w_j|^{-2/3}z_j^2
	\ge\frac13\|\bm{z}\|_{4/3}^2.
\end{equation*}
Integrating along line segments and extending to zero coordinates by
continuity proves that $F$ is $1/3$-strongly convex in the
$\ell_{4/3}$ norm. Since
$F(\bm{w}_*(\bm{\mu}))=1/2$ and
$\nabla F(\bm{w}_*(\bm{\mu}))=\bm{\mu}/\|\bm{\mu}\|_4$, it follows that
\begin{equation*}
	F(\bm{w})\ge\frac12
	+\frac{\bm{\mu}^\top(\bm{w}-\bm{w}_*(\bm{\mu}))}{\|\bm{\mu}\|_4}
	+\frac16\|\bm{w}-\bm{w}_*(\bm{\mu})\|_{4/3}^2.
\end{equation*}
For $\bm{w}\in\Omega$, we have $F(\bm{w})\le1/2$, hence
\begin{equation*}
	\bm{\mu}^\top(\bm{w}_*(\bm{\mu})-\bm{w})
	\ge\frac {\|\bm{\mu}\|_4} 6\|\bm{w}-\bm{w}_*(\bm{\mu})\|_{4/3}^2.
	\label{eq:fast_highdim_strong_convexity}
\end{equation*}
Thus \eqref{eq:fast_geometric_condition} holds with $c_\Omega=1/6$.
Also, for $\bm{G}\sim\mathcal{N}(\bm{0},I_d)$, norm comparison and Jensen's inequality yield
\begin{equation*}
	\sqrt{\frac2\pi}\,d^{1/4}
	=d^{-3/4}\mathbb{E}\|\bm{G}\|_1
	\le\mathcal{G}(\Omega)=\mathbb{E}\|\bm{G}\|_4\le\left(\mathbb{E}\sum_{j=1}^d G_j^4\right)^{1/4}
	=(3d)^{1/4}.
\end{equation*}
Conditional Jensen's inequality, applied after adding independent
Gaussian noise, also gives
\begin{equation*}
	\sqrt{\underline\lambda}\,\mathcal G(\Omega)
	\le\mathcal G(\Sigma^{1/2}\Omega)
	\le\sqrt{\overline\lambda}\,\mathcal G(\Omega).
\end{equation*}
Consequently, the right-hand side of \eqref{eq:fast_ddopt_bound} has
order $\sqrt d/(n\|\bm{\mu}\|_4)$, with constants depending only on the fixed spectral bounds.

For every nonzero $\bm{\mu}$, the support function is smooth, and direct differentiation gives
\begin{equation*}
	\operatorname{tr}(H_{\bm{\mu}})
	=\frac3{\|\bm{\mu}\|_4^3}\sum_{j=1}^d\mu_j^2
	-\frac3{\|\bm{\mu}\|_4^7}\sum_{j=1}^d\mu_j^6.
	\label{eq:fast_highdim_trace}
\end{equation*}
For general $\Sigma$, the corresponding weighted trace is
\begin{equation*}
	\operatorname{tr}(\Sigma H_{\bm\mu})
	=\frac3{\|\bm\mu\|_4^3}\sum_j\Sigma_{jj}\mu_j^2
	-\frac3{\|\bm\mu\|_4^7}\sum_{i,j}\Sigma_{ij}\mu_i^3\mu_j^3,
\end{equation*}
because $(H_{\bm\mu})_{ij}
=3\mu_i^2\mathbb I\{i=j\}/\|\bm\mu\|_4^3
-3\mu_i^3\mu_j^3/\|\bm\mu\|_4^7$.
Moreover,
\begin{equation*}
	\underline\lambda\operatorname{tr}(H_{\bm\mu})
	\le\operatorname{tr}(\Sigma H_{\bm\mu})
	\le\overline\lambda\operatorname{tr}(H_{\bm\mu}).
\end{equation*}
Suppose now that the $s$ nonzero coordinates have equal absolute value,
so $|\mu_j|=\|\bm{\mu}\|_4 s^{-1/4}$ for $j\notin\mathcal{S}_0$. Then
\begin{equation*}
	\begin{aligned}
		&\sum_{j=1}^d\mu_j^2=\|\bm{\mu}\|_4^2\sqrt s,\qquad
		\sum_{j=1}^d\mu_j^6=\frac{\|\bm{\mu}\|_4^6}{\sqrt s},\\
		&\operatorname{tr}(H_{\bm{\mu}})
		=\frac3{\|\bm{\mu}\|_4}\left(\sqrt s-\frac1{\sqrt s}\right).
	\end{aligned}
	\label{eq:fast_highdim_sparse_trace}
\end{equation*}
For $s\ge2$, Theorem \ref{thm:fast_local_minimax} therefore gives
\begin{equation*}
	\lim_{M\to\infty}\lim_{n\to\infty}nR_{n,M}^{\Sigma}(\bm{\mu})
	=\frac12\operatorname{tr}(\Sigma H_{\bm\mu})
	\asymp\frac{\sqrt s}{\|\bm{\mu}\|_4}.
	\label{eq:fast_highdim_sparse_constant}
\end{equation*}
The comparison constants in the last relation depend only on the fixed spectral bounds and are independent of
$d$ and $s\ge2$. The plug-in policy has the same asymptotic risk
constant uniformly on every fixed local ball by
\eqref{eq:fast_uniform_plugin_limit}. 
Thus the fast rate matches the sharp local asymptotic scale when $s\asymp d$. For $2\le s\ll d$, its $n^{-1}$ power is optimal, but its geometric coefficient exceeds the sharp one by a factor of order $\sqrt{d/s}$.

When $s=1$, the Hessian vanishes. Coordinate permutations and sign changes preserve $\Omega$ and $B_2(M)$; the same transformation must be applied to the covariance matrix, which remains in $\mathcal C$. Since the following calculation holds for every $\Sigma\in\mathcal C$, it suffices to take $\bm\mu=\mu_1\bm e_1$ with $\mu_1>0$, where
$\bm{w}_*(\bm{\mu})=\bm{e}_1$. Uniformly for bounded $\bm{t}$,
\begin{equation*}
	\begin{aligned}
		f_*(\mu_1\bm{e}_1+\varepsilon\bm{t})
		&=\left((\mu_1+\varepsilon t_1)^4+\varepsilon^4\sum_{j=2}^d t_j^4\right)^{1/4}\\
		&=\mu_1+\varepsilon t_1+\frac{\varepsilon^4}{4\mu_1^3}\sum_{j=2}^d t_j^4+o(\varepsilon^4).
	\end{aligned}
\end{equation*}
Hence Theorem \ref{thm:fast_general_local_rate} applies with
$a(\varepsilon)=\varepsilon^4$ and
$\psi(\bm{t})=\sum_{j=2}^d t_j^4/(4\mu_1^3)$, which is nonaffine on
every $B_2(M)$ with $M>0$. For fixed $d,\mu_1,M$,
\begin{equation*}
	\begin{aligned}
		n^2R_{n,M}^{\Sigma}(\mu_1\bm e_1)
		&\to V_{M,\Sigma}
		\left(\bm t\mapsto\frac1{4\mu_1^3}\sum_{j=2}^d t_j^4\right).
	\end{aligned}
\end{equation*}
The limit is positive, and the plug-in conclusion of the same theorem gives
\begin{equation*}
	R_{n,M}^{\Sigma}(\mu_1\bm{e}_1)
	\asymp\sup_{\bm{t}\in B_2(M)}
	\mathbb{E}_{\bm\mu_{n,\bm{t}},\Sigma}
	R(\bm{w}_*(\bar{\bm{x}}))
	\asymp n^{-2}.
	\label{eq:fast_highdim_axis_rate}
\end{equation*}

At zero, Corollary~\ref{cor:fast_first_order_limit} and Theorem~\ref{thm:fast_general_local_rate}
give the $n^{-1/2}$ rate. The plug-in rule attains all three rates,
uniformly over $\Sigma\in\mathcal C$, without using $\Sigma$. 

\subsection{Optimal value versus optimal decision on \texorpdfstring{$\ell_{p^\star}$}{lp-star} balls}
\label{app:lq-value-decision}
Let $d\ge2$, $1<p<\infty$, $p^\star=p/(p-1)$, and $\Omega=B_{p^\star}^d$. By H\"older duality,
\begin{equation*}
f_*(\bm z)=\|\bm z\|_p,
\end{equation*}
and, for $\bm z\ne\bm0$, the unique optimizer has coordinates
\begin{equation*}
w_{*,j}(\bm z)
=\frac{\operatorname{sgn}(z_j)|z_j|^{p-1}}
{\|\bm z\|_p^{p-1}}.
\end{equation*}
Fix $\bm\mu=\mu_1\bm e_1$, $\mu_1>0$. For $\bm t$ in a fixed bounded set and sufficiently small $\varepsilon$, $\mu_1+\varepsilon t_1>0$, and
\begin{equation*}
\begin{aligned}
	f_*(\bm\mu+\varepsilon\bm t)
	&=(\mu_1+\varepsilon t_1)
	\left\{1+
	\frac{\varepsilon^p\sum_{j=2}^d|t_j|^p}
	{(\mu_1+\varepsilon t_1)^p}\right\}^{1/p}\\
	&=\mu_1+\varepsilon t_1
	+\frac{\varepsilon^p}{p(\mu_1+\varepsilon t_1)^{p-1}}
	\sum_{j=2}^d|t_j|^p+o(\varepsilon^p)\\
	&=\mu_1+\varepsilon t_1
	+\frac{\varepsilon^p}{p\mu_1^{p-1}}
	\sum_{j=2}^d|t_j|^p+o(\varepsilon^p),
\end{aligned}
\end{equation*}
uniformly on that bounded set. Since $\bm w_*(\bm\mu)=\bm e_1$,
\begin{equation*}
r_{\bm\mu}(\varepsilon\bm t)
=\frac{\varepsilon^p}{p\mu_1^{p-1}}
\sum_{j=2}^d|t_j|^p+o(\varepsilon^p).
\end{equation*}
Thus $a(\varepsilon)=\varepsilon^p$ and $\psi(\bm t)=(p\mu_1^{p-1})^{-1}\sum_{j=2}^d|t_j|^p$. Since $d\ge2$, the profile is nonaffine on every $B_2(M)$ with $M>0$. Theorem~\ref{thm:fast_general_local_rate} therefore yields
\begin{equation*}
R_{n,M}^{\Sigma}(\bm\mu)
\asymp\sup_{\bm t\in B_2(M)}\mathbb E_{\bm\mu_{n,\bm t},\Sigma}R(\bm w_*(\bar{\bm x}))
\asymp n^{-p/2}.
\end{equation*}

Under the reference mean $\bm\mu$, write $\bar{\bm x}=\bm\mu+\bm Z/\sqrt n$, where $\bm Z\sim\mathcal N(\bm0,\Sigma)$. Since $\nabla f_*(\bm\mu)=\bm e_1$, by the delta method,
\begin{equation*}
	\sqrt n\left(f_*(\bm\mu+\bm Z/\sqrt n)-f_*(\bm\mu)\right)\Longrightarrow Z_1,
\end{equation*}
where $Z_1\sim\mathcal N(0,\Sigma_{11})$. Since $f_*(\cdot)=\|\cdot\|_p$, the reverse triangle inequality also gives
\begin{equation*}
	\sqrt n\,|f_*(\bm\mu+\bm Z/\sqrt n)-f_*(\bm\mu)|\le\|\bm Z\|_p.
\end{equation*}
Both $\mathbb E\|\bm Z\|_p$ and $\mathbb E\|\bm Z\|_p^2$ are finite. DCT therefore yields
\begin{equation*}
\begin{aligned}
	\sqrt n\,\mathbb E_{\bm\mu,\Sigma}|f_*(\bar{\bm x})-f_*(\bm\mu)|
	&\longrightarrow\mathbb E|Z_1|=\sqrt{2\Sigma_{11}/\pi},\\
	n\,\mathbb E_{\bm\mu,\Sigma}(f_*(\bar{\bm x})-f_*(\bm\mu))^2
	&\longrightarrow\mathbb E Z_1^2=\Sigma_{11}.
\end{aligned}
\end{equation*}

Now we discuss the local minimax value-estimation rates. Fix $M>0$. Following the local minimax framework of \citet[Chapter~8]{vanderVaart1998}, let
\begin{equation*}
	\widehat f_n=\widehat f_n(\bm x_1,\ldots,\bm x_n),\qquad
	\widehat f_n:(\mathbb R^d)^n\to\mathbb R,
\end{equation*}
range over all Borel measurable estimators of the scalar value $f_*(\bm\mu_{n,\bm t})$. An estimator may depend on the fixed benchmark center $\bm\mu$, covariance $\Sigma$, and radius $M$, but must be the same rule for every unknown $\bm t\in B_2(M)$. The minimax risk takes the supremum over $\bm t$ for each such rule and then the infimum over the rules. The plug-in estimator is $f_*(\bar{\bm x})$.

For the upper bound, the reverse triangle inequality gives
\begin{equation*}
\begin{aligned}
	\sup_{\bm t\in B_2(M)}\mathbb E_{\bm\mu_{n,\bm t},\Sigma}|f_*(\bar{\bm x})-f_*(\bm\mu_{n,\bm t})|
	&\le\frac{\mathbb E\|\bm Z\|_p}{\sqrt n},\\
	\sup_{\bm t\in B_2(M)}\mathbb E_{\bm\mu_{n,\bm t},\Sigma}(f_*(\bar{\bm x})-f_*(\bm\mu_{n,\bm t}))^2
	&\le\frac{\mathbb E\|\bm Z\|_p^2}{n}.
\end{aligned}
\end{equation*}
Both right-hand sides are finite. Moreover, $\|\bm Z\|_p\le d^{\max\{1/p-1/2,0\}}\|\bm Z\|_2$ and $\mathbb E\|\bm Z\|_2^2=\operatorname{tr}(\Sigma)\le d\overline\lambda$, so the upper bounds are uniform over $\Sigma\in\mathcal C$.

For the lower bound, fix $h=\min\{M,1\}>0$ and use the prior $\pi$ assigning probability $1/2$ to each of $\bm T=\pm h\bm e_1\in B_2(M)$. Since $\bm\mu=\mu_1\bm e_1$ and $\mu_1>0$, for $h/\sqrt n<\mu_1$,
\begin{equation*}
		f_*(\bm\mu+h\bm e_1/\sqrt n)-f_*(\bm\mu-h\bm e_1/\sqrt n)=\frac{2h}{\sqrt n}.
\end{equation*}
For every estimator $\widehat f_n$, its maximal risk dominates its average risk at these two alternatives:
\begin{equation}
	\begin{aligned}
		\sup_{\bm t\in B_2(M)}\mathbb E_{\bm\mu_{n,\bm t},\Sigma}
		|\widehat f_n-f_*(\bm\mu_{n,\bm t})|&\ge\frac12\mathbb E_{\bm\mu+h\bm e_1/\sqrt n,\Sigma}
		|\widehat f_n-\mu_1-h/\sqrt n|\\
		&\quad+\frac12\mathbb E_{\bm\mu-h\bm e_1/\sqrt n,\Sigma}
		|\widehat f_n-\mu_1+h/\sqrt n|\\
		&=\mathbb E_{\pi,\Sigma}|\widehat f_n-f_*(\bm\mu_{n,\bm T})|.
	\end{aligned}
\label{eq:maximal_risk_dominates_Bayes_risk}
\end{equation}
We compute the minimum of this Bayes risk. The sample mean is sufficient, so with $\bm Y=\sqrt n(\bar{\bm x}-\bm\mu)=\bm T+\bm Z$, $q(\bm Y):=\mathbb P_{\pi,\Sigma}\{\bm T=h\bm e_1\mid\bm Y\}=\mathbb P_{\pi,\Sigma}\{\bm T=h\bm e_1\mid\bm x_1,\ldots,\bm x_n\}$.
Conditional on the full sample, the expected absolute loss of a scalar output $a\in\mathbb R$ is
\begin{equation*}
\begin{aligned}
	&\quad q(\bm Y)|a-\mu_1-h/\sqrt n|+(1-q(\bm Y))|a-\mu_1+h/\sqrt n|\\
		&\ge\min\{q(\bm Y),1-q(\bm Y)\}
		\bigl(|a-\mu_1-h/\sqrt n|+|a-\mu_1+h/\sqrt n|\bigr)\\
		&\ge\frac{2h}{\sqrt n}\min\{q(\bm Y),1-q(\bm Y)\}.
	\end{aligned}
\end{equation*}
Equality is attained by $a=\mu_1+h/\sqrt n$ when $q(\bm Y)\ge1/2$ and by $a=\mu_1-h/\sqrt n$ when $q(\bm Y)<1/2$. These choices are measurable functions of $\bm Y$. Thus minimizing the conditional loss and taking expectations gives
\begin{equation}
	\inf_{\widehat f_n}\mathbb E_{\pi,\Sigma}|\widehat f_n-f_*(\bm\mu_{n,\bm T})|
	=\frac{2h}{\sqrt n}\mathbb E_{\pi,\Sigma}\min\{q(\bm Y),1-q(\bm Y)\}.
    \label{eq:exact_expression_of_Bayes_risk}
\end{equation}
To evaluate the last expectation, the Gaussian densities give
\begin{equation*}
	\begin{aligned}
		\log\frac{q(\bm Y)}{1-q(\bm Y)}
		&=-\frac12(\bm Y-h\bm e_1)^\top\Sigma^{-1}(\bm Y-h\bm e_1)+\frac12(\bm Y+h\bm e_1)^\top\Sigma^{-1}(\bm Y+h\bm e_1)\\
		&=2h\bm e_1^\top\Sigma^{-1}\bm Y.
	\end{aligned}
\end{equation*}
Note that $\bm e_1^\top\Sigma^{-1}\bm Y\mid\bm T=\pm h\bm e_1\sim\mathcal N\!\left(\pm h\bm e_1^\top\Sigma^{-1}\bm e_1,\ \bm e_1^\top\Sigma^{-1}\bm e_1\right)$. Hence
\begin{equation}
	\begin{aligned}
		\mathbb E_{\pi,\Sigma}\min\{q(\bm Y),1-q(\bm Y)\}
		&=\frac12\mathbb P_{\pi,\Sigma}
		\left\{\bm e_1^\top\Sigma^{-1}\bm Y<0\mid\bm T=h\bm e_1\right\}\\
		&\quad+\frac12\mathbb P_{\pi,\Sigma}
		\left\{\bm e_1^\top\Sigma^{-1}\bm Y\ge0\mid\bm T=-h\bm e_1\right\}\\
		&=\Phi\!\left(-h\sqrt{\bm e_1^\top\Sigma^{-1}\bm e_1}\right)\ge\Phi(-h/\sqrt{\underline\lambda}),
	\end{aligned}
\label{eq:expectation_of_posterior_prob}
\end{equation}
where we use $\bm e_1^\top\Sigma^{-1}\bm e_1\le1/\underline\lambda$. Combining Equations~\eqref{eq:maximal_risk_dominates_Bayes_risk}, \eqref{eq:exact_expression_of_Bayes_risk} and \eqref{eq:expectation_of_posterior_prob} gives
\begin{equation*}
		\inf_{\widehat f_n}\sup_{\bm t\in B_2(M)}
		\mathbb E_{\bm\mu_{n,\bm t},\Sigma}|\widehat f_n-f_*(\bm\mu_{n,\bm t})|
		\ge\inf_{\widehat f_n}
		\mathbb E_{\pi,\Sigma}|\widehat f_n-f_*(\bm\mu_{n,\bm T})|\ge\frac{2h}{\sqrt n}\Phi(-h/\sqrt{\underline\lambda}).
\end{equation*}
The minimax risk is at most the maximal risk of $f_*(\bar{\bm x})$, so the upper and lower bounds prove
\begin{equation*}
    \inf_{\widehat f_n}\sup_{\bm t\in B_2(M)}\mathbb E_{\bm\mu_{n,\bm t},\Sigma}|\widehat f_n-f_*(\bm\mu_{n,\bm t})|\asymp n^{-1/2}.
\end{equation*}

For squared loss, by Jensen's inequality, for every $\widehat f_n$,
\begin{equation*}
	\begin{aligned}
		\sup_{\bm t\in B_2(M)}
		\mathbb E_{\bm\mu_{n,\bm t},\Sigma}(\widehat f_n-f_*(\bm\mu_{n,\bm t}))^2&\ge\sup_{\bm t\in B_2(M)}
		\left(\mathbb E_{\bm\mu_{n,\bm t},\Sigma}|\widehat f_n-f_*(\bm\mu_{n,\bm t})|\right)^2\\
		&=\left(\sup_{\bm t\in B_2(M)}
		\mathbb E_{\bm\mu_{n,\bm t},\Sigma}|\widehat f_n-f_*(\bm\mu_{n,\bm t})|\right)^2.
	\end{aligned}
\end{equation*}
Taking the infimum over $\widehat f_n$ and using monotonicity of squaring on $[0,\infty)$ gives
\begin{equation*}
	\begin{aligned}
		\inf_{\widehat f_n}\sup_{\bm t\in B_2(M)}
		\mathbb E_{\bm\mu_{n,\bm t},\Sigma}(\widehat f_n-f_*(\bm\mu_{n,\bm t}))^2&\ge\left(\inf_{\widehat f_n}\sup_{\bm t\in B_2(M)}
		\mathbb E_{\bm\mu_{n,\bm t},\Sigma}|\widehat f_n-f_*(\bm\mu_{n,\bm t})|\right)^2\\
		&\ge\frac{4h^2}{n}\Phi(-h/\sqrt{\underline\lambda})^2.
	\end{aligned}
\end{equation*}
Together with the squared-loss upper bound for $f_*(\bar{\bm x})$, this proves 
\begin{equation*}
    \inf_{\widehat f_n}\sup_{\bm t\in B_2(M)}
		\mathbb E_{\bm\mu_{n,\bm t},\Sigma}(\widehat f_n-f_*(\bm\mu_{n,\bm t}))^2\asymp n^{-1}.
\end{equation*}
Both comparisons hold uniformly over $\Sigma\in\mathcal C$ for fixed $M>0$.

\subsection{A non-power local rate: details for Example~\ref{ex:fast_nonpower_rate}}
\label{app:non_power_local_rate}
A non-power scale can be necessary even when the optimal decision is unique. Choose $0<\delta<e^{-2}$ and define, for $|s|\le\delta$,
\begin{equation*}
	g(0)=0,\qquad
	g(s)=\frac{|s|}{\log(1/|s|)}\quad(s\ne0).
\end{equation*}
The choice of $g$ separates differentiability from higher-order power growth:
\begin{equation*}
	\begin{aligned}
		\frac{g(s)}{|s|}&=\frac1{\log(1/|s|)}\to0,\\
		\frac{g(s)}{|s|^\alpha}
		&=\frac{|s|^{1-\alpha}}{\log(1/|s|)}\to+\infty
		\quad(\alpha>1),\qquad s\to0.
	\end{aligned}
\end{equation*}
The first limit gives $g'(0)=0$, while the second rules out an $O(|s|^\alpha)$ remainder for any $\alpha>1$. The function $g$ is even and continuously differentiable, with
\begin{equation*}
	g'(s)=\operatorname{sgn}(s)
	\left\{\frac1{\log(1/|s|)}+\frac1{[\log(1/|s|)]^2}\right\}
	\to0\qquad(s\to0,\ s\ne0).
\end{equation*}
$g$ is convex because its derivative is continuous at zero and
\begin{equation*}
	g''(s)=\frac{\log(1/s)+2}{s[\log(1/s)]^3}>0,
	\qquad 0<s<\delta.
\end{equation*}
We realize $1+g$ as a section of a support function by using the intercept and slope of each tangent line as the coordinates of a feasible point:
\begin{equation*}
	\Omega=\operatorname{conv}\left(
	\{\bm{0}\}\cup
	\left\{\left(1+g(s)-sg'(s),\,g'(s)\right):|s|\le\delta\right\}
	\right),\qquad \bm{\mu}=\bm{e}_1.
\end{equation*}
The generating curve $\left\{\left(1+g(s)-sg'(s),\,g'(s)\right):|s|\le\delta\right\}$ is continuous on $[-\delta,\delta]$, so its union with $\{\bm0\}$ is compact. Its convex hull is compact in $\mathbb R^2$, and $\Omega$ is therefore compact, convex, and contains $\bm0$. For $s\ne0$,
\begin{equation*}
	1+g(s)-sg'(s)=1-\frac{|s|}{[\log(1/|s|)]^2}<1,
\end{equation*}
whereas $s=0$ gives the point $(1,0)^\top=\bm e_1$. Every point of $\Omega$ is a finite convex combination of generating points, i.e., the origin or points of the form $\left(1+g(s)-sg'(s),\,g'(s)\right)$ for some $|s|\le\delta$; its first coordinate can equal $1$ only if all points receiving positive weight are $\bm e_1$. Hence
\begin{equation*}
	f_*(\bm e_1)=1,\qquad
	\mathcal W_*(\bm e_1)=\{\bm e_1\},\qquad
	\nabla f_*(\bm e_1)=\bm e_1.
\end{equation*}

For $|u|\le\delta$, convexity gives $g(s)+g'(s)(u-s)\le g(u)$, with equality at $s=u$. Taking the supremum of the corresponding linear functionals therefore yields
\begin{equation*}
	\begin{aligned}
		f_*((1,u)^\top)
		&=\max\left\{0,\sup_{|s|\le\delta}
		[1+g(s)-sg'(s)+ug'(s)]\right\}\\
		&=1+g(u).
	\end{aligned}
\end{equation*}
Here $1+g(u)\ge1$, so the feasible point $\bm0$ does not affect the maximum. Fix $M>0$ and take $\varepsilon M\le\min\{1/2,\delta/2\}$. For every $\bm t\in B_2(M)$,
\begin{equation*}
	1+\varepsilon t_1\ge\frac12,\qquad
	\left|\frac{\varepsilon t_2}{1+\varepsilon t_1}\right|
	\le2\varepsilon M\le\delta.
\end{equation*}
Positive homogeneity of $f_*$ and the section identity give
\begin{equation*}
	\begin{aligned}
		f_*(\bm e_1+\varepsilon\bm t)
		&=(1+\varepsilon t_1)
		f_*\left(\left(1,\frac{\varepsilon t_2}{1+\varepsilon t_1}\right)^\top\right)\\
		&=(1+\varepsilon t_1)
		\left[1+g\left(\frac{\varepsilon t_2}{1+\varepsilon t_1}\right)\right].
	\end{aligned}
\end{equation*}
Hence
\begin{equation*}
\begin{aligned}
    r_{\bm e_1}(\varepsilon\bm t)
		&=f_*(\bm e_1+\varepsilon\bm t)-1-\varepsilon t_1\\
		&=(1+\varepsilon t_1)
		g\left(\frac{\varepsilon t_2}{1+\varepsilon t_1}\right).
\end{aligned}
\end{equation*}

To identify the local scale, set
$a(\varepsilon)=\varepsilon/\log(1/\varepsilon)$ for
$0<\varepsilon<1$. If $t_2=0$, the preceding identity and $g(0)=0$
give $r_{\bm{e}_1}(\varepsilon\bm{t})=0$ for all sufficiently small
$\varepsilon$.
If $t_2\ne0$, substituting the definition of $g$ gives
\begin{equation*}
	\begin{aligned}
		r_{\bm{e}_1}(\varepsilon\bm{t})
		&=(1+\varepsilon t_1)
		\frac{\varepsilon|t_2|/(1+\varepsilon t_1)}
		{\log\bigl((1+\varepsilon t_1)/(\varepsilon|t_2|)\bigr)}\\
		&=\frac{\varepsilon|t_2|}
		{\log(1/\varepsilon)+\log(1+\varepsilon t_1)-\log|t_2|}.
	\end{aligned}
\end{equation*}
Consequently,
\begin{equation*}
	\frac{r_{\bm{e}_1}(\varepsilon\bm{t})}{a(\varepsilon)}
	=\frac{|t_2|}
	{1+\dfrac{\log(1+\varepsilon t_1)-\log|t_2|}
		{\log(1/\varepsilon)}}
	\to|t_2|.
\end{equation*}
Indeed, for fixed $\bm{t}$ with $t_2\ne0$,
$\log(1+\varepsilon t_1)\to0$, $\log|t_2|$ is finite, and
$\log(1/\varepsilon)\to+\infty$.
Combining the two cases proves $\psi(\bm{t})=|t_2|$. Proposition~\ref{pr:pointtouniform} upgrades this pointwise limit to
\begin{equation*}
	\lim_{\varepsilon\downarrow0}
	\sup_{\bm{t}\in B_2(M)}
	\left|\frac{r_{\bm{e}_1}(\varepsilon\bm{t})}{a(\varepsilon)}
	-|t_2|\right|=0,
	\qquad M>0.
\end{equation*}
Since $\psi(\bm t)=|t_2|$ is not affine on any $B_2(M)$ with $M>0$ and
$a(n^{-1/2})=2/(\sqrt{n}\log n)$for $n>1$,
Theorem \ref{thm:fast_general_local_rate} gives
\begin{equation*}
	\frac{\sqrt n\log n}{2}R_{n,M}^{\Sigma}(\bm e_1)
	\to V_{M,\Sigma}(\psi)>0.
\end{equation*}
The same theorem gives the matching SAA rate:
\begin{equation*}
	R_{n,M}^{\Sigma}(\bm e_1)
	\asymp\sup_{\bm t\in B_2(M)}
	\mathbb E_{\bm\mu_{n,\bm t},\Sigma}R(\bm w_*(\bar{\bm x}))
	\asymp\frac1{\sqrt n\log n},\qquad M>0.
	\label{eq:fast_nonpower_rate}
\end{equation*}
The comparison constants can be chosen uniformly over $\Sigma\in\mathcal C$, and SAA uses neither the center nor $\Sigma$. Along the transverse direction, the remainder satisfies, for every $\alpha>1$,
\begin{equation*}
	\begin{aligned}
		\varepsilon^{-1}r_{\bm e_1}(\varepsilon\bm e_2)
		&=\frac1{\log(1/\varepsilon)}\to0,\\
		\varepsilon^{-\alpha}r_{\bm e_1}(\varepsilon\bm e_2)
		&=\frac{\varepsilon^{1-\alpha}}{\log(1/\varepsilon)}
		\to+\infty.
	\end{aligned}
\end{equation*}
Thus the first-order remainder vanishes, but no pure-power scale $\varepsilon^\alpha$ with $\alpha>1$ yields a finite local profile. The logarithmic factor captures the intermediate scale and produces the corresponding non-power regret rate.

\section{Technical Results for KL Neighborhoods}
\label{app:kl-technical}

\subsection{Proof of Corollary~\ref{cor:matched_gaussian_benchmark}}
\label{app:proof-matched-gaussian}
\begin{proof}
	Write the local alternative as
	$\bm\mu+\Sigma^{1/2}\bm u/\sqrt n$ with $\bm u\in B_2(M)$ and set
	\[
	\bm Y=\Sigma^{-1/2}\sqrt n(\bar{\bm x}-\bm\mu)
	=\bm u+\bm G,
	\qquad \bm G\sim\mathcal N(\bm0,I_d).
	\]
	The transformed experiment is therefore the standard Gaussian shift experiment on $B_2(M)$. Under \eqref{eq:fast_general_local_expansion}, its normalized support-function remainder is
	\[
	\frac{r_{\bm\mu}(\Sigma^{1/2}\bm u/\sqrt n)}{a(n^{-1/2})}
	\longrightarrow (\psi\circ\Sigma^{1/2})(\bm u)
	\]
	uniformly on bounded sets. Applying Theorem~\ref{thm:fast_general_local_rate} to this transformed experiment proves \eqref{eq:matched_gaussian_limit} and the rate-optimality assertion for SAA. In the smooth case,
	$\psi(\bm t)=\bm t^\top H_{\bm\mu}\bm t/2$, so
	$\psi(\Sigma^{1/2}\bm u)=\bm u^\top
	\Sigma^{1/2}H_{\bm\mu}\Sigma^{1/2}\bm u/2$.
	The smooth formulas follow from Theorem~\ref{thm:fast_local_minimax}.
\end{proof}

\subsection{Uniform properties of the truncated mean}
\begin{lemma}[Means and truncated observations in a KL neighborhood]
	\label{lem:div_kl_truncation}
    Assume \eqref{eq:ng_reference_distribution}.
	Write $\bm x_i=(x_{i,1},\ldots,x_{i,d})^\top$, $i=1,\ldots,n$, for the i.i.d. observations with law $Q$, and set
	\begin{equation}
		\begin{aligned}
			\tau_n&:=\log^2(n+1),&
			x_{i,j}^{(n)}&:=\operatorname{sgn}(x_{i,j})\min\{|x_{i,j}|,\tau_n\},
			\quad i=1,\ldots,n,\ j=1,\ldots,d,\\
			\widetilde{\bm\mu}_n&:=\frac1n\sum_{i=1}^n\bm x_i^{(n)},&
			\bm Z_{n,Q}&:=\sqrt n(\widetilde{\bm\mu}_n-\bm\mu_Q).
		\end{aligned}
	\end{equation}
	Here $\bm x_i^{(n)}=(x_{i,1}^{(n)},\ldots,x_{i,d}^{(n)})^\top$. Expectations $\mathbb E_Q$ involving $\bm x_1$ are with respect to its one-observation law $Q$; $\mathbb E_P$ denotes integration of the same one-observation function under the reference law $P$. Expectations of sample statistics are denoted by $\mathbb E_{Q^{\otimes n}}$.
	For each fixed $M>0$, all $Q\in\mathcal Q_{n,M}^{\mathrm{KL}}(P)$ have finite first moments, and
	\begin{equation}
		\sup_{Q\in\mathcal Q_{n,M}^{\mathrm{KL}}(P)}
		\sqrt n\|\bm\mu_Q-\bm\mu\|_2=O(1).
		\label{eq:div_mean_shift_bound}
	\end{equation}
	\begin{equation}
		\sup_{Q\in\mathcal Q_{n,M}^{\mathrm{KL}}(P)}
		n(\bm\mu_Q-\bm\mu)^\top\Sigma^{-1}(\bm\mu_Q-\bm\mu)
		\le M^2+o(1).
		\label{eq:div_mean_ellipsoid_bound}
	\end{equation}
	\begin{equation}
		\sup_{Q\in\mathcal Q_{n,M}^{\mathrm{KL}}(P)}
		\|\mathbb E_Q\bm x_1^{(n)}-\bm\mu_Q\|_2=O(n^{-1}).
		\label{eq:div_truncation_bias_bound}
	\end{equation}
	\begin{equation}
		\sup_{Q\in\mathcal Q_{n,M}^{\mathrm{KL}}(P)}
		\|\operatorname{Cov}_Q(\bm x_1^{(n)})-\Sigma\|_{\mathrm{op}}\to0.
		\label{eq:div_truncation_covariance_limit}
	\end{equation}
	\begin{equation}
		\sup_{n\ \mathrm{large}}\sup_{Q\in\mathcal Q_{n,M}^{\mathrm{KL}}(P)}
		\mathbb E_{Q^{\otimes n}}\|\bm Z_{n,Q}\|_2^q<\infty,
		\qquad q>0.
		\label{eq:div_truncation_moment_bound}
	\end{equation}
	For every sequence $Q_n\in\mathcal Q_{n,M}^{\mathrm{KL}}(P)$,
	$\bm Z_{n,Q_n}\xrightarrow{d}\mathcal N(\bm0,\Sigma)$. Moreover, for each fixed $\epsilon>0$ and an independent $\bm G\sim\mathcal N(\bm0,I_d)$,
	\begin{equation}
		\sup_{Q\in\mathcal Q_{n,M}^{\mathrm{KL}}(P)}
		\left\|\mathcal L_{Q^{\otimes n}}(\bm Z_{n,Q}+\epsilon\bm G)
		-\mathcal N(\bm0,\Sigma+\epsilon^2I_d)\right\|_{\mathrm{TV}}
		\to0.
		\label{eq:div_smoothed_clt}
	\end{equation}
	The law in \eqref{eq:div_smoothed_clt} includes the independent $\bm G$.
\end{lemma}

\begin{proof}
	\eqref{eq:ng_reference_distribution} implies, for some $\gamma>0$,
	\begin{equation*}
		\mathbb E_P e^{\gamma\|\bm x_1\|_2}<\infty.
	\end{equation*}
	Indeed, for sufficiently small $\gamma>0$, $e^{\gamma\|\bm x_1\|_2}\le e^{\gamma\sum_j|x_{1,j}|}\le\sum_{\bm v\in\{-1,1\}^d}e^{\gamma\bm v^\top\bm x_1}$. For bounded measurable $g$, nonnegativity of
$D_{\mathrm{KL}}(Q\|e^gP/\mathbb E_Pe^g)$ gives the entropy inequality
\begin{equation}
    \mathbb E_Qg\le D_{\mathrm{KL}}(Q\|P)+\log\mathbb E_Pe^g.
    \label{eq:div_entropy_inequality}
\end{equation}
For nonnegative $g$ with $\mathbb E_Pe^g<\infty$, apply this inequality to $g\wedge K$ and let $K\to\infty$ by monotone convergence.
Taking $g(\bm x_1)=\gamma\|\bm x_1\|_2$ proves 
	\begin{equation*}
		\mathbb{E}_Q\|\bm x_1\|_2\le\frac{D_{\mathrm{KL}}(Q\|P)+\log\mathbb E_Pe^{g}}{\gamma}<\infty.
	\end{equation*}
	For $g(\bm x_1)=\bm v^\top(\bm x_1-\bm\mu)$ with $\|\bm v\|_2<\eta$, we have
	\begin{equation*}
		\mathbb E_Q|g|\le\|\bm v\|_2(\mathbb E_Q\|\bm x_1\|_2+\|\bm\mu\|_2)<\infty,
		\qquad \log\mathbb E_Pe^g=\Lambda_P(\bm v)<\infty.
	\end{equation*}
	Thus the same log-density identity is integrable under $Q$ and gives
	\begin{equation*}
		\bm v^\top(\bm\mu_Q-\bm\mu)
		\le D_{\mathrm{KL}}(Q\|P)+\Lambda_P(\bm v).
	\end{equation*}
	There exist $C_P,\rho_P>0$ such that $\Lambda_P(\bm v)\le C_P\|\bm v\|_2^2/2$ for $\|\bm v\|_2\le\rho_P$. For sufficiently large $n$, when $\bm\mu_Q\ne\bm\mu$, $\bm v=n^{-1/2}(\bm\mu_Q-\bm\mu)/\|\bm\mu_Q-\bm\mu\|_2$ satisfies $\|\bm v\|_2\le\rho_P$, and therefore
	\begin{equation*}
		\sqrt n\|\bm\mu_Q-\bm\mu\|_2
		\le\frac{M^2+C_P}{2},
	\end{equation*}
	which proves \eqref{eq:div_mean_shift_bound}. Thus, for sufficiently large $n$, $\bm h=\bm\mu_Q-\bm\mu$ belongs to the domain of $\bm\theta(\bm h)$ in Lemma~\ref{lem:ng_mean_parametrization}, uniformly for large $n$. Since $\mathbb E_Q(\bm x_1-\bm\mu)=\bm h$, and $\log \dd P_{\bm h}/\dd P$ is $Q$-integrable,
	\begin{equation}
		D_{\mathrm{KL}}(Q\|P)
		=D_{\mathrm{KL}}(Q\|P_{\bm h})
		+\bm\theta(\bm h)^\top\bm h-\Lambda_P(\bm\theta(\bm h)),
		\label{eq:div_entropy_projection}
	\end{equation}
	and \eqref{eq:ng_local_kl_expansion} gives
	\begin{equation*}
		\bm\theta(\bm h)^\top\bm h-\Lambda_P(\bm\theta(\bm h))=\frac12\bm h^\top\Sigma^{-1}\bm h+O(\|\bm h\|_2^3).
	\end{equation*}
	Thus,
	\begin{equation*}
		\frac12\bm h^\top\Sigma^{-1}\bm h\le D_{\mathrm{KL}}(Q\|P)-D_{\mathrm{KL}}(Q\|P_{\bm h})+C\|\bm h\|_2^3\le\frac{M^2}{2n}+C\|\bm{h}\|_2^3
	\end{equation*}
	for some constant $C>0$ and small $\bm{h}$. By \eqref{eq:div_mean_shift_bound}, $n\|\bm{h}\|_2^3=o(1)$, which proves \eqref{eq:div_mean_ellipsoid_bound}.
	
	For \eqref{eq:div_truncation_bias_bound}, choose $0<\gamma_1<\gamma$ and let $g(\bm x_1)=\gamma_1\|\bm x_1\|_2\mathbb I\{\|\bm x_1\|_2>\tau\}$ for some constant $\tau>0$. Then
	\begin{equation*}
		\mathbb E_Pe^g\le1+\mathbb E_P[e^{\gamma_1\|\bm x_1\|_2}\mathbb I\{\|\bm x_1\|_2>\tau\}]\le1+e^{-(\gamma-\gamma_1)\tau}\mathbb E_Pe^{\gamma\|\bm x_1\|_2},
	\end{equation*}
	and therefore
	\begin{equation*}
		\begin{aligned}
			\mathbb E_Qg=\gamma_1\mathbb E_Q[\|\bm x_1\|_2\mathbb I\{\|\bm x_1\|_2>\tau\}]&\le D_{\mathrm{KL}}(Q\|P)+\log\left(1+e^{-(\gamma-\gamma_1)\tau}\mathbb E_Pe^{\gamma\|\bm x_1\|_2}\right)\\&\le\frac{M^2}{2n}+e^{-(\gamma-\gamma_1)\tau}\mathbb E_Pe^{\gamma\|\bm x_1\|_2}.
		\end{aligned}
	\end{equation*}
	Since the observations are identically distributed,
	\begin{equation*}
		\mathbb E_{Q^{\otimes n}}\widetilde{\bm\mu}_n
		=\frac1n\sum_{i=1}^n\mathbb E_Q\bm x_i^{(n)}
		=\mathbb E_Q\bm x_1^{(n)},
		\qquad \bm\mu_Q=\mathbb E_Q\bm x_1.
	\end{equation*}
	Moreover, $\|\bm x_1-\bm x_1^{(n)}\|_2\le\|\bm x_1\|_2\mathbb I\{\|\bm x_1\|_2>\tau_n\}$. Taking $\tau=\tau_n$ gives
	\begin{equation}
		\begin{aligned}
			\sup_{Q\in\mathcal Q_{n,M}^{\mathrm{KL}}(P)}
			\mathbb E_Q\|\bm x_1^{(n)}-\bm x_1\|_2
			&\le\sup_{Q\in\mathcal Q_{n,M}^{\mathrm{KL}}(P)}
			\mathbb E_Q[\|\bm x_1\|_2\mathbb I\{\|\bm x_1\|_2>\tau_n\}]\\
			&\le\gamma_1^{-1}\left[\frac{M^2}{2n}
			+e^{-(\gamma-\gamma_1)\tau_n}\mathbb E_Pe^{\gamma\|\bm x_1\|_2}\right]
			=O(n^{-1}).
		\end{aligned}
		\label{eq:div_truncation_lone_bound}
	\end{equation}
	Here $\tau_n=\log^2(n+1)$. Since
	\begin{equation*}
		\|\mathbb E_{Q^{\otimes n}}\widetilde{\bm\mu}_n-\bm\mu_Q\|_2
		=\|\mathbb E_Q(\bm x_1^{(n)}-\bm x_1)\|_2
		\le\mathbb E_Q\|\bm x_1^{(n)}-\bm x_1\|_2,
	\end{equation*}
	this also proves \eqref{eq:div_truncation_bias_bound}.
	
	Pinsker's inequality gives
	\begin{equation*}
		\sup_{Q\in\mathcal Q_{n,M}^{\mathrm{KL}}(P)}
		\|Q-P\|_{\mathrm{TV}}
		\le\sup_{Q\in\mathcal Q_{n,M}^{\mathrm{KL}}(P)}
		\sqrt{\frac12D_{\mathrm{KL}}(Q\|P)}\le\frac{M}{2\sqrt n}.
	\end{equation*}
	For every fixed $k>0$, $\|\bm x_1^{(n)}\|_2\le\sqrt d\,\tau_n$, and hence
	\begin{equation}
		\begin{aligned}
			\sup_{Q\in\mathcal Q_{n,M}^{\mathrm{KL}}(P)}
			\left|\mathbb E_Q\|\bm x_1^{(n)}\|_2^k-\mathbb E_P\|\bm x_1^{(n)}\|_2^k\right|
			&\le2(\sqrt d\,\tau_n)^k
			\sup_{Q\in\mathcal Q_{n,M}^{\mathrm{KL}}(P)}\|Q-P\|_{\mathrm{TV}}\\
			&\le\frac{M(\sqrt d\,\tau_n)^k}{\sqrt n}=o(1).
		\end{aligned}
		\label{eq:div_clipped_moments}
	\end{equation}
	Since $\|\bm x_1^{(n)}\|_2^k\uparrow\|\bm x_1\|_2^k$ pointwise, monotone convergence and the exponential-moment bound under $P$ give
	\begin{equation*}
		\mathbb E_P\|\bm x_1^{(n)}\|_2^k\to\mathbb E_P\|\bm x_1\|_2^k<\infty.
	\end{equation*}
	Combining this limit with \eqref{eq:div_clipped_moments} gives, for every fixed $k>0$,
	\begin{equation}
		\sup_{n\ \mathrm{large}}\sup_{Q\in\mathcal Q_{n,M}^{\mathrm{KL}}(P)}
		\mathbb E_Q\|\bm x_1^{(n)}\|_2^k<\infty.
		\label{eq:div_clipped_uniform_moments}
	\end{equation}
	
For the first and second clipped coordinate moments, Pinsker's inequality and the bounds $|x_{1,j}^{(n)}|\le\tau_n$ and $|x_{1,j}^{(n)}x_{1,k}^{(n)}|\le\tau_n^2$ give
\begin{align*}
    \sup_Q|\mathbb E_Qx_{1,j}^{(n)}-\mathbb E_Px_{1,j}^{(n)}|
       &\le M\tau_n/\sqrt n=o(1),\\
    \sup_Q|\mathbb E_Q[x_{1,j}^{(n)}x_{1,k}^{(n)}]
          -\mathbb E_P[x_{1,j}^{(n)}x_{1,k}^{(n)}]|
       &\le M\tau_n^2/\sqrt n=o(1).
\end{align*}
Dominated convergence under $P$, using $|x_{1,j}^{(n)}|\le|x_{1,j}|$ and
$|x_{1,j}^{(n)}x_{1,k}^{(n)}|\le(x_{1,j}^2+x_{1,k}^2)/2$, identifies the limits as $\mu_j$ and $\mathbb E_Px_{1,j}x_{1,k}$.
Expanding each covariance entry and using fixed $d$ proves \eqref{eq:div_truncation_covariance_limit} uniformly over $Q$.

	For an integer $k\ge1$, expand the $2k$th moment of
	$n^{-1/2}\sum_{i=1}^n(x_{i,j}^{(n)}-\mathbb E_Qx_{1,j}^{(n)})$. By independence and centering, any term with an index occurring only once is zero. The remaining terms contain at most $k$ distinct indices, and their number is at most $C_kn^k$. H\"older's inequality and \eqref{eq:div_clipped_uniform_moments} with exponent $2k$ therefore give
	\begin{equation*}
		\mathbb E_{Q^{\otimes n}}\left|\frac1{\sqrt n}\sum_{i=1}^n(x_{i,j}^{(n)}-\mathbb E_Qx_{1,j}^{(n)})\right|^{2k}\le C_k\mathbb E_Q|x_{1,j}^{(n)}-\mathbb E_Qx_{1,j}^{(n)}|^{2k}\le C_{k,P,M},
	\end{equation*}
	Combined with \eqref{eq:div_truncation_bias_bound}, summing over coordinates and using $\|\bm z\|_2^{2k}\le d^{k-1}\sum_j|z_j|^{2k}$ proves \eqref{eq:div_truncation_moment_bound}; every fixed $q>0$ is bounded by a sufficiently large even moment.
	
Fix $\bm v\in\mathbb R^d$ and put $\bm U_{n,i}=\bm x_i^{(n)}-\mathbb E_Q\bm x_1^{(n)}$.
The third moments of $\bm U_{n,1}$ are uniformly bounded by \eqref{eq:div_clipped_uniform_moments}.
Taylor's inequality $|e^{iu}-1-iu+u^2/2|\le|u|^3/6$ and independence give
\begin{equation}
    \mathbb E_{Q^{\otimes n}}e^{i\bm v^\top n^{-1/2}\sum_i\bm U_{n,i}}
    =\left[1-\frac{\bm v^\top\operatorname{Cov}_Q(\bm x_1^{(n)})\bm v}{2n}
             +O_{\bm v,P,M}(n^{-3/2})\right]^n
    \longrightarrow e^{-\bm v^\top\Sigma\bm v/2}
    \label{eq:cf_of_z_nQ}
\end{equation}
uniformly over $Q$, by \eqref{eq:div_truncation_covariance_limit} and the expansion $\log(1+z)=z+O(|z|^2)$ near $1$.
The centering correction $\sqrt n(\mathbb E_Q\bm x_1^{(n)}-\bm\mu_Q)$ is uniformly $O(n^{-1/2})$ by \eqref{eq:div_truncation_bias_bound}; hence
\[
    \sup_Q\left|\mathbb E_{Q^{\otimes n}}e^{i\bm v^\top\bm Z_{n,Q}}
                    -e^{-\bm v^\top\Sigma\bm v/2}\right|\to0.
\]
L\'evy's continuity theorem proves the stated weak limit along every sequence $Q_n$.

	For completeness, fix $\epsilon>0$ and let $\bm G\sim\mathcal N(\bm0,I_d)$ be independent of the sample. All suprema over $Q$ in this proof are over $\mathcal Q_{n,M}^{\mathrm{KL}}(P)$. Independence gives
	\begin{equation*}
		\mathbb E e^{i\bm v^\top(\bm Z_{n,Q}+\epsilon\bm G)}
		=e^{-\epsilon^2\|\bm v\|_2^2/2}
		\mathbb E_{Q^{\otimes n}}e^{i\bm v^\top\bm Z_{n,Q}}.
	\end{equation*}
The Gaussian factor is integrable. Fourier inversion gives continuous densities $p_{n,Q,\epsilon}$ and $p_\epsilon$ for the smoothed statistic and $\mathcal N(\bm0,\Sigma+\epsilon^2I_d)$, respectively, with
\begin{equation}
    \sup_Q\|p_{n,Q,\epsilon}-p_\epsilon\|_\infty
    \le(2\pi)^{-d}\int e^{-\epsilon^2\|\bm v\|_2^2/2}
       \sup_Q\left|\mathbb E e^{i\bm v^\top\bm Z_{n,Q}}
                        -e^{-\bm v^\top\Sigma\bm v/2}\right|\,\dd\bm v
    \to0.
    \label{eq:uniform_density_convergence}
\end{equation}
	For every fixed $\bm v$, the supremum inside the integral tends to zero and is at most $2$. DCT with the integrable bound $2e^{-\epsilon^2\|\bm v\|_2^2/2}$ proves the displayed limit. The supremum is measurable: the uniform first-moment bound from \eqref{eq:div_truncation_moment_bound} makes these characteristic functions uniformly Lipschitz in $\bm v$, and hence their supremum of absolute differences is continuous for sufficiently large $n$.
	
For $L>0$, split the density integral over $B_2(L)$ and its complement. The uniform second-moment bound in \eqref{eq:div_truncation_moment_bound} and Markov's inequality give
\[
    \sup_Q\int|p_{n,Q,\epsilon}-p_\epsilon|
    \le\operatorname{vol}(B_2(L))\sup_Q\|p_{n,Q,\epsilon}-p_\epsilon\|_\infty
       +C_{P,M,\epsilon}L^{-2}.
\]
First let $n\to\infty$, then $L\to\infty$, and use the density formula for total variation.
	This proves \eqref{eq:div_smoothed_clt} for every fixed $\epsilon>0$.
\end{proof}

\subsection{Proof of Theorem~\ref{thm:div_kl_minimax}}
\label{app:proof-kl-minimax}
\begin{proof}
	We first prove \eqref{eq:div_kl_lower}. Fix $0<\delta<M$ and use $\bm g_L,\Gamma_L,J_L$ and $Q_{n,\bm u}^{(L)}$ from Lemma~\ref{lem:ng_bounded_score_tilt}. Since $J_L\to\Sigma^{1/2}$ as $L\to\infty$, fix $L$ sufficiently large that $\Gamma_L\succ0$ and $J_L$ is invertible. Throughout the sampling limit below, $L$ and $\delta$ are fixed. Equation~\eqref{eq:ng_bounded_score_divergence}, with $\phi(s)=s\log s-s+1$, gives
	\begin{equation*}
		\sup_{\bm u\in B_2(M-\delta)}nD_{\mathrm{KL}}(Q_{n,\bm u}^{(L)}\|P)
		\le\frac{(M-\delta)^2}{2}+o(1)<\frac{M^2}{2}
	\end{equation*}
	for all sufficiently large $n$. The density ratios $\dd Q_{n,\bm{u}}^{(L)}/\dd P$ in \eqref{eq:div_bounded_score_submodel} belong to $[1/2,3/2]$ for large $n$, so
	\begin{equation*}
		\sup_{\bm u\in B_2(M-\delta)}
		\mathbb E_{Q_{n,\bm u}^{(L)}}\|\bm x\|_2
		\le\frac32\mathbb E_P\|\bm x\|_2<\infty.
	\end{equation*}
	Consequently,
	\begin{equation*}
		\{Q_{n,\bm u}^{(L)}:\bm u\in B_2(M-\delta)\}
		\subseteq\mathcal Q_{n,M}^{\mathrm{KL}}(P),\qquad
		\mathbb E_{Q_{n,\bm u}^{(L)}}\bm x
		=\bm\mu+\frac{J_L\bm u}{\sqrt n}.
	\end{equation*}
	The second identity is from \eqref{eq:div_bounded_score_submodel}. Under $P^{\otimes n}$, the sample likelihoods $L_n^{(L)}$ have strictly positive continuous paths, and Lemma~\ref{lem:ng_bounded_score_tilt} gives
	\begin{equation*}
		\begin{aligned}
			&L_n^{(L)}\xrightarrow[n\to\infty]{d}
			\left[\bm u\mapsto e^{\bm u^\top\bm Y-\|\bm u\|_2^2/2}\right]
			\quad\text{in }C(B_2(M-\delta)),\qquad
			\bm Y\sim\mathcal N(\bm0,I_d),\\
			&\sup_{n\ \mathrm{large}}\mathbb E_{P^{\otimes n}}
			\sup_{\bm u\in B_2(M-\delta)}|L_n^{(L)}(\bm u)|^2<\infty.
		\end{aligned}
	\end{equation*}
	Thus \eqref{eq:ng_likelihood_process_limit} and its required moment bound hold on $B_2(M-\delta)$.
	
	For every fixed policy $\hat{\bm w}_n$, the submodel inclusion gives
	\begin{equation*}
		\sup_{Q\in\mathcal Q_{n,M}^{\mathrm{KL}}(P)}
		\mathbb E_{Q^{\otimes n}}R_Q(\hat{\bm w}_n)
		\ge\sup_{\bm u\in B_2(M-\delta)}
		\mathbb E_{(Q_{n,\bm u}^{(L)})^{\otimes n}}
		\ell_n(J_L\bm u,\hat{\bm w}_n).
	\end{equation*}
	Taking the infimum over the same class of sample-based policies on both sides gives
	\begin{equation*}
		\mathcal R_{n,M}^{\mathrm{KL}}(P)
		\ge\inf_{\hat{\bm w}_n}\sup_{\bm u\in B_2(M-\delta)}
		\mathbb E_{(Q_{n,\bm u}^{(L)})^{\otimes n}}
		\ell_n(J_L\bm u,\hat{\bm w}_n).
	\end{equation*}
	All assumptions of Lemma~\ref{lem:ng_likelihood_lower_bound} have now been verified with radius $M-\delta$ and $\Gamma=J_L$, so
	\begin{equation*}
		\liminf_{n\to\infty}
		\frac{\mathcal R_{n,M}^{\mathrm{KL}}(P)}{a(n^{-1/2})}
		\ge V_{M-\delta,I_d}(\psi\circ J_L).
	\end{equation*}
	This inequality holds for every sufficiently large fixed $L$. By $J_L\to\Sigma^{1/2}$ and continuity in Lemma~\ref{lem:div_gaussian_comparison}, letting first $L\to\infty$ and then $\delta\downarrow0$ gives
	\begin{equation*}
		\liminf_{n\to\infty}\frac{\mathcal R_{n,M}^{\mathrm{KL}}(P)}{a(n^{-1/2})}
		\ge V_{M-\delta,I_d}(\psi\circ\Sigma^{1/2})
		\xrightarrow{\delta\downarrow0}V_{M,I_d}(\psi\circ\Sigma^{1/2}).
	\end{equation*}
	This proves \eqref{eq:div_kl_lower} under the finite-second-moment conditions.
	
	For the upper bound, now additionally assume \eqref{eq:ng_reference_distribution}, and fix $\delta,\epsilon,\upsilon>0$, with $\delta<1$. Lemma~\ref{lem:div_kl_truncation} implies
	\begin{equation*}
		\bm t_Q:=\sqrt n(\bm\mu_Q-\bm\mu)
		\in\Sigma^{1/2}B_2(M+\delta)
		\subset B_2(L),\qquad
		L:=\|\Sigma^{1/2}\|_{\mathrm{op}}(M+1),
	\end{equation*}
	uniformly over $Q\in\mathcal Q_{n,M}^{\mathrm{KL}}(P)$, for all sufficiently large $n$. All suprema over $Q$ below are over $\mathcal Q_{n,M}^{\mathrm{KL}}(P)$. The invertible transformation $\bm y\mapsto\Sigma^{-1/2}\bm y$ identifies the experiment $\mathcal N(\bm t,\Sigma)$, $\bm t=\Sigma^{1/2}\bm u$, with $\mathcal N(\bm u,I_d)$, and preserves the loss $\ell_n(\Sigma^{1/2}\bm u,\bm w)$. Thus Lemma~\ref{lem:div_gaussian_comparison} applies. In the experiment $\bm Y\sim\mathcal N(\bm t,\Sigma)$, $\bm t\in\Sigma^{1/2}B_2(M+\delta)$, there exists a kernel $\delta_n^\upsilon(\dd\bm w\mid\bm y)$ with $\delta_n^\upsilon(\Omega_{n,L}\mid\bm y)=1$ for every $\bm y\in\mathbb R^d$, satisfying
	\begin{equation*}
		\begin{aligned}
			&\quad\sup_{\bm t\in\Sigma^{1/2}B_2(M+\delta)}
			\mathbb E_{\mathcal N(\bm t,\Sigma)}
			\int\ell_n(\bm t,\bm w)\,\delta_n^\upsilon(\dd\bm w\mid\bm Y)\\
			&\le\inf_{\delta'}\sup_{\bm u\in B_2(M+\delta)}
			\mathbb E_{\mathcal N(\bm u,I_d)}
			\int_\Omega\ell_n(\Sigma^{1/2}\bm u,\bm w)\,
			\delta'(\dd\bm w\mid\bm Y)
			+\upsilon a(n^{-1/2}).
		\end{aligned}
	\end{equation*}
	Let $\bm G\sim\mathcal N(\bm0,I_d)$ be independent of the sample and apply $\delta_n^\upsilon$ to
	\begin{equation*}
		\bm Y_n'=\sqrt n(\widetilde{\bm\mu}_n-\bm\mu)+\epsilon\bm G
		=\bm t_Q+\bm Z_{n,Q}+\epsilon\bm G.
	\end{equation*}
For fixed $0<\epsilon\le1$, the Gaussian KL formula and Pinsker's inequality yield
\[
    \|\mathcal N(\bm0,\Sigma+\epsilon^2I_d)-\mathcal N(\bm0,\Sigma)\|_{\mathrm{TV}}
    \le C_\Sigma\epsilon^2,\qquad
    C_\Sigma=\sqrt{\operatorname{tr}(\Sigma^{-2})/8}.
\]
Thus \eqref{eq:div_smoothed_clt} and translation invariance imply
\begin{equation}
    \sup_Q\|\mathcal L(\bm Y_n')-\mathcal N(\bm t_Q,\Sigma)\|_{\mathrm{TV}}
       \le o_n(1)+C_\Sigma\epsilon^2.
    \label{eq:bound_on_TV}
\end{equation}
The laws include the independent $\bm G$. Here $o_n(1)$ is uniform in $Q$ for each fixed $\epsilon$; the limit $n\to\infty$ precedes $\epsilon\downarrow0$.
	By \eqref{eq:local_loss_upper_bound},
	\begin{equation*}
		0\le\int\ell_n(\bm t_Q,\bm w)\,\delta_n^\upsilon(\dd\bm w\mid\bm y)
		\le2b_{3L}(n^{-1/2}),\qquad\bm y\in\mathbb R^d.
	\end{equation*}
	For $0\le g\le C$, 
	\begin{equation}
		|\mathbb E_{P_1}g-\mathbb E_{P_2}g|
		=\left|\int_0^C(P_1\{g>s\}-P_2\{g>s\})\,\dd s\right|
		\le C\|P_1-P_2\|_{\mathrm{TV}}.
		\label{eq:bounded_function_derive_bounded_TV}
	\end{equation}
Apply \eqref{eq:bounded_function_derive_bounded_TV} with
$g_Q(\bm y)=\int\ell_n(\bm t_Q,\bm w)\delta_n^\upsilon(\dd\bm w\mid\bm y)$.
The loss bound is uniform in $Q$, $\bm y$, and $n$ after normalization:
\begin{equation}
\begin{aligned}
    \sup_Q\mathbb E g_Q(\bm Y_n')
       &\le\sup_Q\mathbb E_{\mathcal N(\bm t_Q,\Sigma)}g_Q(\bm Y)\\
       &\quad+2b_{3L}(n^{-1/2})\{o_n(1)+C_\Sigma\epsilon^2\}.
\end{aligned}
\label{eq:first_ineq}
\end{equation}
In particular, $b_{3L}(n^{-1/2})/a(n^{-1/2})=O(1)$ by the local expansion.
This restriction of the action space is essential: the global loss bound on $\Omega$ alone need not be of order $a(n^{-1/2})$.
	For all sufficiently large $n$, the mean bound above gives
	\begin{equation*}
		\{\bm t_Q:Q\in\mathcal Q_{n,M}^{\mathrm{KL}}(P)\}
		\subseteq\Sigma^{1/2}B_2(M+\delta).
	\end{equation*}
	Consequently,
	\begin{equation}
		\sup_Q\mathbb E_{\mathcal N(\bm t_Q,\Sigma)}
		\int\ell_n(\bm t_Q,\bm w)\,\delta_n^\upsilon(\dd\bm w\mid\bm Y)\le\sup_{\bm t\in\Sigma^{1/2}B_2(M+\delta)}
		\mathbb E_{\mathcal N(\bm t,\Sigma)}
		\int\ell_n(\bm t,\bm w)\,\delta_n^\upsilon(\dd\bm w\mid\bm Y).
		\label{eq:second_ineq}
	\end{equation}
Combining \eqref{eq:first_ineq}--\eqref{eq:second_ineq} bounds the randomized rule by the nearly optimal Gaussian risk plus the displayed transfer error.
	The expectation on the left includes the independent $\bm G$. For each fixed sample, average the selected action over both $\bm G$ and the kernel $\delta_n^\upsilon$. This defines the deterministic policy
	\begin{equation*}
		\hat{\bm w}_n(\bm x_1,\ldots,\bm x_n)
		:=\int_{\mathbb R^d}\int_{\Omega_{n,L}}\bm w\,
		\delta_n^\upsilon\left(\dd\bm w\,\middle|\,
		\sqrt n(\widetilde{\bm\mu}_n-\bm\mu)+\epsilon\bm g\right)
		\mathcal N(\bm0,I_d)(\dd\bm g).
	\end{equation*}
The bounded Markov-kernel integral is measurable, and compactness and convexity give $\hat{\bm w}_n\in\Omega_{n,L}$ for every sample.
Since regret is affine in the action, averaging preserves risk exactly:
\[
    \mathbb E_{Q^{\otimes n}}R_Q(\hat{\bm w}_n)
       =\mathbb E\int\ell_n(\bm t_Q,\bm w)
                          \delta_n^\upsilon(\dd\bm w\mid\bm Y_n').
\]
Consequently, the maximal risk of this deterministic policy bounds $\mathcal R_{n,M}^{\mathrm{KL}}(P)$ from above.

	Also, the uniform local expansion implies
	\begin{equation*}
		\left|\frac{b_{3L}(n^{-1/2})}{a(n^{-1/2})}
		-\sup_{\bm t\in B_2(3L)}\psi(\bm t)\right|
		\le\sup_{\bm t\in B_2(3L)}
		\left|\frac{r_{\bm\mu}(\bm t/\sqrt n)}{a(n^{-1/2})}-\psi(\bm t)\right|\to0.
	\end{equation*}
	Dividing the risk bound by $a(n^{-1/2})$, taking $n\to\infty$ with $\delta,\epsilon,\upsilon$ fixed, and applying \eqref{eq:div_gaussian_ellipsoid} gives
	\begin{equation*}
		\limsup_{n\to\infty}
		\frac{\mathcal R_{n,M}^{\mathrm{KL}}(P)}{a(n^{-1/2})}
		\le V_{M+\delta,I_d}(\psi\circ\Sigma^{1/2})+\upsilon+2C_\Sigma\epsilon^2
		\sup_{\bm t\in B_2(3L)}\psi(\bm t).
	\end{equation*}
	Letting $\upsilon\downarrow0$, $\epsilon\downarrow0$, and $\delta\downarrow0$ and using continuity in Lemma~\ref{lem:div_gaussian_comparison} proves the upper bound and hence \eqref{eq:div_kl_exact_limit}. 
	
	If $\psi$ is affine on $\Sigma^{1/2}B_2(M)$, then $\psi\circ\Sigma^{1/2}$ is affine on $B_2(M)$, and the conditional-mean identity makes every Bayes Jensen gap zero. Thus $V_{M,I_d}(\psi\circ\Sigma^{1/2})=0$. If $\psi$ is nonaffine on $\Sigma^{1/2}B_2(M)$, Theorem~\ref{thm:fast_general_local_rate} gives $V_{M,I_d}(\psi\circ\Sigma^{1/2})>0$. Choose $L$ large enough that $\bm t_Q\in B_2(L)$ for large $n$ and $\sup_{B_2(L)}\psi>0$. Note that the deterministic argument leading to \eqref{eq:fast_general_plugin_remainder_bound} applies with $\bar{\bm x}$ replaced by $\widetilde{\bm\mu}_n$. With a similar bound, \eqref{eq:fast_radial_polynomial_growth} gives
	\begin{equation*}
		\begin{aligned}
			R_Q(\bm w_*(\widetilde{\bm\mu}_n))
			&\le r_{\bm\mu}(\bm t_Q/\sqrt n)
			+r_{\bm\mu}((\bm t_Q+2\bm Z_{n,Q})/\sqrt n)\\
			&\le b_L(n^{-1/2})
			\left[1+C_b\left(1+\frac{2\|\bm Z_{n,Q}\|_2}{L}\right)^{c_b}\right].
		\end{aligned}
	\end{equation*}
	Choose an integer $q\ge\max\{1,c_b\}$. By \eqref{eq:div_truncation_moment_bound},
	\begin{equation*}
		\sup_Q\mathbb E_{Q^{\otimes n}}
		\left(1+\frac{2\|\bm Z_{n,Q}\|_2}{L}\right)^{c_b}
		\le C_L\left(1+\sup_Q\mathbb E_{Q^{\otimes n}}\|\bm Z_{n,Q}\|_2^q\right)
		\le C_{P,M,L}.
	\end{equation*}
	Moreover, \eqref{eq:fast_radial_scale_limits} gives
	\begin{equation*}
		\frac{b_L(n^{-1/2})}{a(n^{-1/2})}
		\to\sup_{\bm t\in B_2(L)}\psi(\bm t)\in(0,\infty).
	\end{equation*}
	Taking expectations in the deterministic bound therefore yields
	\begin{equation}
		\sup_Q\mathbb E_{Q^{\otimes n}}R_Q(\bm w_*(\widetilde{\bm\mu}_n))
		\le C_{P,M}b_L(n^{-1/2})\le C_{P,M}'a(n^{-1/2}).
		\label{eq:div_truncated_uniform_risk_upper}
	\end{equation}
	Since the limiting constant is positive, \eqref{eq:div_kl_exact_limit} gives, for all sufficiently large $n$,
	\begin{equation*}
		\frac12 V_{M,I_d}(\psi\circ\Sigma^{1/2})a(n^{-1/2})
		\le\mathcal R_{n,M}^{\mathrm{KL}}(P)
		\le\sup_Q\mathbb E_{Q^{\otimes n}}
		R_Q(\bm w_*(\widetilde{\bm\mu}_n)).
	\end{equation*}
	The second inequality follows from the definition of the minimax risk and feasibility of the truncated plug-in policy. Together with the upper bound just proved, this gives \eqref{eq:div_kl_plugin_rate}. If the limiting constant is zero, \eqref{eq:div_kl_exact_limit} instead gives $\mathcal R_{n,M}^{\mathrm{KL}}(P)=o(a(n^{-1/2}))$.
\end{proof}
\subsection{Proof of Corollary~\ref{cor:div_kl_smooth}}
\label{app:proof-kl-smooth}
\begin{proof}
The expansion has $a(\varepsilon)=\varepsilon^2$ and
$\psi(\bm t)=\bm t^\top H_{\bm\mu}\bm t/2$.
Theorem~\ref{thm:div_kl_minimax} and the quadratic Gaussian limit give
\eqref{eq:exact_limit_for_KL_neighborhood}--\eqref{eq:asymp_limit_for_KL_neighborhood}.

For \eqref{eq:div_kl_smooth}, fix $M$ and take all suprema over
$\mathcal Q_{n,M}^{\mathrm{KL}}(P)$.
Write $\bm t_Q=\sqrt n(\bm\mu_Q-\bm\mu)$ and $\bm Z=\bm Z_{n,Q}$.
By \eqref{eq:div_mean_shift_bound}, $\sup_Q\|\bm t_Q\|_2\le L_0$ for large $n$.
Choose $r>0$ so that $f_*$ is $C^2$ near $\bm\mu+B_2(r)$, and let $C_r$ bound its Hessian there.
For large $n$, $\bm\mu_Q\in\bm\mu+B_2(r/2)$.
On $A_{n,Q}=\{\|\widetilde{\bm\mu}_n-\bm\mu\|_2\le r\}$,
the support identity and Taylor's integral formula yield
\[
    nR_Q(\bm w_*(\widetilde{\bm\mu}_n))
      =\int_0^1(1-s)\bm Z^\top
         \nabla^2f_*\!\left(\bm\mu+
                   \frac{\bm t_Q+(1-s)\bm Z}{\sqrt n}\right)\bm Z\,\dd s.
\]
Put $D_{n,Q}=nR_Q(\bm w_*(\widetilde{\bm\mu}_n))
                    -\bm Z^\top H_{\bm\mu}\bm Z/2$.
For each fixed $L$ and sufficiently large $n$, $\{\|\bm Z\|_2\le L\}\subseteq A_{n,Q}$ uniformly over $Q$. Continuity of the Hessian gives
\[
    \sup_Q\mathbb E[|D_{n,Q}|\mathbb I\{\|\bm Z\|_2\le L\}]
       \le\frac{L^2}{2}
         \sup_{\|\bm h\|_2\le(L_0+L)/\sqrt n}
           \|\nabla^2f_*(\bm\mu+\bm h)-H_{\bm\mu}\|_{\mathrm{op}}
       \to0.
\]
On $A_{n,Q}$, $|D_{n,Q}|\le(C_r+\|H_{\bm\mu}\|_{\mathrm{op}})\|\bm Z\|_2^2/2$.
The uniform fourth-moment bound in \eqref{eq:div_truncation_moment_bound} therefore implies
\[
    \sup_Q\mathbb E[|D_{n,Q}|\mathbb I_{A_{n,Q}}\mathbb I\{\|\bm Z\|_2>L\}]
        \le CL^{-2}.
\]
On $A_{n,Q}^c$, $\|\bm Z\|_2>r\sqrt n-L_0\ge r\sqrt n/2$.
Since every feasible action has regret at most $2R_\Omega\|\bm\mu_Q\|_2\le C$,
the same fourth-moment bound gives
\begin{align*}
    \sup_Q\mathbb E[|D_{n,Q}|\mathbb I_{A_{n,Q}^c}]
       &\le Cn\sup_Q\mathbb P\{\|\bm Z\|_2>r\sqrt n/2\}\\
       &\quad+C\sup_Q\mathbb E[\|\bm Z\|_2^2
                                  \mathbb I\{\|\bm Z\|_2>r\sqrt n/2\}]\\
       &=O(n^{-1}).
\end{align*}
First taking $n\to\infty$ and then $L\to\infty$ proves
$\sup_Q\mathbb E|D_{n,Q}|\to0$.

Finally, independence gives
\[
    \mathbb E\bm Z\bm Z^\top
      =\operatorname{Cov}_Q(\bm x_1^{(n)})+n\bm b_{n,Q}\bm b_{n,Q}^\top,
      \qquad \bm b_{n,Q}=\mathbb E_Q\bm x_1^{(n)}-\bm\mu_Q.
\]
Equations~\eqref{eq:div_truncation_covariance_limit} and
\eqref{eq:div_truncation_bias_bound} imply uniform convergence of this matrix to $\Sigma$.
Taking traces against $H_{\bm\mu}$ and combining with the $L^1$ approximation proves \eqref{eq:div_kl_smooth}.
\end{proof}

\subsection{Boundary growth on KL neighborhoods}
\begin{theorem}[Boundary growth and matching rates on a KL neighborhood]
	\label{thm:div_boundary_growth}
	Let $d\ge2$, let $\Omega$ be compact and convex with nonempty interior, and suppose $\bm\mu\ne\bm0$ has a unique maximizer $\bm w_*(\bm\mu)$. Retain $\bm n_{\bm\mu}=\bm\mu/\|\bm\mu\|_2$ and the orthonormal basis matrix $U_{\bm\mu}$ of $\bm n_{\bm\mu}^{\perp}$. For the lower bound, assume only $\mathbb E_P\|\bm x\|_2^2<\infty$ and $\Sigma=\operatorname{Cov}_P(\bm x)\succ0$. Suppose that, for some $\rho,\bar\kappa,\beta>0$, a boundary profile $\varrho_{\bm\mu}$ satisfies
	\begin{equation}
		\begin{gathered}
			\bm w_*(\bm\mu)+U_{\bm\mu}\bm z-\bm n_{\bm\mu}\varrho_{\bm\mu}(\bm z)\in\Omega,\qquad
			\varrho_{\bm\mu}(\bm0)=0,\\
			0\le \varrho_{\bm\mu}(\bm z)\le\frac{\bar\kappa}{1+\beta}\|\bm z\|_2^{1+\beta},
			\qquad \|\bm z\|_2\le\rho.
		\end{gathered}
		\label{eq:div_boundary_feasible_growth}
	\end{equation}
	For every $M>0$, every unit $\bm v\in\mathbb R^{d-1}$, and every $h>0$ with $h\|\Sigma^{-1/2}U_{\bm\mu}\bm v\|_2<M$,
	\begin{equation}
		\begin{aligned}
			\liminf_{n\to\infty}n^{(1+\beta)/(2\beta)}
			\mathcal R_{n,M}^{\mathrm{KL}}(P)
			&\ge 2\Phi\!\left(-h\|\Sigma^{-1/2}U_{\bm\mu}\bm v\|_2\right)
			\frac{\beta}{1+\beta}
			(\bar\kappa\|\bm\mu\|_2)^{-1/\beta}h^{(1+\beta)/\beta}>0.
		\end{aligned}
		\label{eq:div_boundary_rate_lower}
	\end{equation}
	
	For the matching upper bound, assume additionally \eqref{eq:ng_reference_distribution}. Suppose that $\varrho_{\bm\mu}$ describes the entire local feasible set and has a corresponding lower growth bound: for some $0<\underline\kappa\le\bar\kappa$, after decreasing $\rho$ if necessary,
	\begin{equation}
		\begin{gathered}
			\bm w_*(\bm\mu)+U_{\bm\mu}\bm z-\bm n_{\bm\mu}b\in\Omega
			\quad\Longleftrightarrow\quad \varrho_{\bm\mu}(\bm z)\le b,
			\qquad \|\bm z\|_2\le\rho,\quad |b|\le\rho,\\
			\frac{\underline\kappa}{1+\beta}\|\bm z\|_2^{1+\beta}
			\le \varrho_{\bm\mu}(\bm z)\le\frac{\bar\kappa}{1+\beta}\|\bm z\|_2^{1+\beta},
			\qquad \|\bm z\|_2\le\rho.
		\end{gathered}
		\label{eq:div_boundary_two_sided_growth}
	\end{equation}
	Then, for every fixed $M>0$ and all sufficiently large $n$,
	\begin{equation}
		c_{P,M}n^{-(1+\beta)/(2\beta)}
		\le\mathcal R_{n,M}^{\mathrm{KL}}(P)\le\sup_{Q\in\mathcal Q_{n,M}^{\mathrm{KL}}(P)}
		\mathbb E_{Q^{\otimes n}}R_Q(\bm w_*(\widetilde{\bm\mu}_n))
		\le C_{P,M}n^{-(1+\beta)/(2\beta)},
		\label{eq:div_boundary_matching_rates}
	\end{equation}
	where $0<c_{P,M}\le C_{P,M}<\infty$ may also depend on the fixed geometry. 
\end{theorem}

\begin{proof}
	We first prove the lower bound. Fix any $\bm{v}\in\mathbb{R}^{d-1}$ with $\|\bm{v}\|_2=1$, fix $h>0$ with $h\|\Sigma^{-1/2}U_{\bm\mu}\bm v\|_2<M$, and take $\bm{t}=hU_{\bm{\mu}}\bm{v}$. For $0\le b\le\rho$, \eqref{eq:div_boundary_feasible_growth} gives $\bm{w}_*(\bm{\mu})+bU_{\bm{\mu}}\bm{v}-\bm{n}_{\bm{\mu}}\varrho_{\bm\mu}(b\bm{v})\in\Omega$. Since $\bm{\mu}^\top U_{\bm{\mu}}\bm{v}=\bm{n}_{\bm{\mu}}^\top U_{\bm{\mu}}\bm{v}=0$, and
	$U_{\bm{\mu}}^\top U_{\bm{\mu}}=I_{d-1}$, the support-function
	definition gives
	\begin{equation*}
		\begin{aligned}
			r_{\bm{\mu}}(\bm{t}/\sqrt n)
			&=\sup_{\bm{w}\in\Omega}
			(\bm{\mu}+\bm{t}/\sqrt n)^\top
			(\bm{w}-\bm{w}_*(\bm{\mu}))\\
			&\ge(\bm{\mu}+hU_{\bm{\mu}}\bm{v}/\sqrt n)^\top
			(bU_{\bm{\mu}}\bm{v}-\bm{n}_{\bm{\mu}}\varrho_{\bm\mu}(b\bm{v}))\\
			&=\frac{hb}{\sqrt n}-\|\bm{\mu}\|_2\varrho_{\bm\mu}(b\bm{v})\\
			&\ge\frac{hb}{\sqrt n}
			-\frac{\bar{\kappa}\|\bm{\mu}\|_2}{1+\beta}b^{1+\beta}.
		\end{aligned}
	\end{equation*}
	The same calculation with $-b\bm{v}$ gives the identical bound for
	$r_{\bm{\mu}}(-\bm{t}/\sqrt n)$.
	The concave function $b\mapsto\frac{hb}{\sqrt n}-\frac{\bar{\kappa}\|\bm{\mu}\|_2}{1+\beta}b^{1+\beta}$ on the last line is maximized over
	$b\ge0$ at
	\begin{equation*}
		b=\left(\frac{h}{\bar{\kappa}\|\bm{\mu}\|_2\sqrt n}
		\right)^{1/\beta}.
	\end{equation*}
	For sufficiently large $n$, $h/\sqrt n\le\bar\kappa\|\bm\mu\|_2\rho^\beta$, so this value belongs to $[0,\rho]$. Substitution yields
	\begin{equation}
		\min\{r_{\bm{\mu}}(\bm{t}/\sqrt n),
		r_{\bm{\mu}}(-\bm{t}/\sqrt n)\}
		\ge\frac{\beta}{1+\beta}
		\left(\frac{h}{\sqrt n}\right)^{(1+\beta)/\beta}
		(\bar{\kappa}\|\bm{\mu}\|_2)^{-1/\beta}.
		\label{eq:div_boundary_directional_remainder}
	\end{equation}
	
	Use the tilted construction of distributions in Lemma~\ref{lem:ng_bounded_score_tilt}. Since $J_L\to\Sigma^{1/2}$, for every sufficiently large fixed $L$ the vector $\bm u=J_L^{-1}\bm t$ satisfies $\|\bm u\|_2<M$. Equations~\eqref{eq:div_bounded_score_submodel} and~\eqref{eq:ng_bounded_score_divergence}, with $\phi(s)=s\log s-s+1$, give
	\begin{equation*}
		\bm\mu_{Q_{n,\pm\bm u}^{(L)}}=\bm\mu\pm\frac{\bm t}{\sqrt n},\qquad
		nD_{\mathrm{KL}}(Q_{n,\pm\bm u}^{(L)}\|P)
		\to\frac{\|\bm u\|_2^2}{2}<\frac{M^2}{2}.
	\end{equation*}
	Both distributions belong to $\mathcal Q_{n,M}^{\mathrm{KL}}(P)$ for large $n$.
	
	Use the equal-probability prior on $Q_{n,\bm u}^{(L)}$ and $Q_{n,-\bm u}^{(L)}$, whose mean perturbations are $\bm t$ and $-\bm t$. The same binary-prior argument in the proof of Proposition~\ref{prop:fast_symmetric_difference} applies. Since both distributions belong to $\mathcal Q_{n,M}^{\mathrm{KL}}(P)$, every measurable sample-based policy $\hat{\bm w}_n$ satisfies
	\begin{equation*}
		\sup_{Q\in\mathcal Q_{n,M}^{\mathrm{KL}}(P)}
		\mathbb E_{Q^{\otimes n}}R_Q(\hat{\bm w}_n)\ge\frac12\mathbb E_{(Q_{n,\bm u}^{(L)})^{\otimes n}}
		\ell_n(\bm t,\hat{\bm w}_n)
		+\frac12\mathbb E_{(Q_{n,-\bm u}^{(L)})^{\otimes n}}
		\ell_n(-\bm t,\hat{\bm w}_n).
	\end{equation*}
	Taking the infimum over the same class of policies gives
	\begin{equation*}
		\mathcal R_{n,M}^{\mathrm{KL}}(P)
		\ge\inf_{\hat{\bm w}_n}\left\{
		\frac12\mathbb E_{(Q_{n,\bm u}^{(L)})^{\otimes n}}
		\ell_n(\bm t,\hat{\bm w}_n)
		+\frac12\mathbb E_{(Q_{n,-\bm u}^{(L)})^{\otimes n}}
		\ell_n(-\bm t,\hat{\bm w}_n)\right\}.
	\end{equation*}
	The right side is the minimum Bayes risk for this two-point prior.
	
	For $0\le\lambda\le1$, by convexity of $r_{\bm{\mu}}$, note that
	\begin{equation}
		\begin{aligned}
			&\quad\lambda r_{\bm\mu}(\bm t/\sqrt n)
			+(1-\lambda)r_{\bm\mu}(-\bm t/\sqrt n)
			-r_{\bm\mu}((2\lambda-1)\bm t/\sqrt n)\\
			&\ge\min\{\lambda,1-\lambda\}
			\{r_{\bm\mu}(\bm t/\sqrt n)+r_{\bm\mu}(-\bm t/\sqrt n)\}.
		\end{aligned}
		\label{eq:jensen_gap_lb}
	\end{equation}
	Now let $\lambda$ be the posterior probability of $Q_{n,\bm u}^{(L)}$ given the sample. Bayes' formula gives
	\begin{equation*}
		\lambda
		=\frac{\tfrac12 L_n^{(L)}(\bm u)}
		{\tfrac12 L_n^{(L)}(\bm u)+\tfrac12 L_n^{(L)}(-\bm u)}
		=\frac{L_n^{(L)}(\bm u)}
		{L_n^{(L)}(\bm u)+L_n^{(L)}(-\bm u)}.
	\end{equation*}
	The denominator is positive. Conditional on the sample, the expected regret of an action $\bm w\in\Omega$ is
	\begin{equation*}
		\begin{aligned}
			&\quad\lambda\ell_n(\bm t,\bm w)+(1-\lambda)\ell_n(-\bm t,\bm w)\\
			&=\lambda f_*(\bm\mu+\bm t/\sqrt n)
			+(1-\lambda)f_*(\bm\mu-\bm t/\sqrt n)
			-\bm w^\top\left(\bm\mu+\frac{(2\lambda-1)\bm t}{\sqrt n}\right)\\
			&\ge\lambda f_*(\bm\mu+\bm t/\sqrt n)
			+(1-\lambda)f_*(\bm\mu-\bm t/\sqrt n)
			-f_*(\bm\mu+(2\lambda-1)\bm t/\sqrt n)\\
			&=\lambda r_{\bm\mu}(\bm t/\sqrt n)
			+(1-\lambda)r_{\bm\mu}(-\bm t/\sqrt n)
			-r_{\bm\mu}((2\lambda-1)\bm t/\sqrt n).
		\end{aligned}
	\end{equation*}
	where the equality is attained by the Bayes policy $\bm{w}_*(\bm\mu+(2\lambda-1)\bm t/\sqrt n)$. The marginal sample law under the prior is $\tfrac12(Q_{n,\bm u}^{(L)})^{\otimes n}+\tfrac12(Q_{n,-\bm u}^{(L)})^{\otimes n}$, whose density relative to $P^{\otimes n}$ is $(L_n^{(L)}(\bm u)+L_n^{(L)}(-\bm u))/2$. Since the Bayes policy $\bm{w}_*(\bm\mu+(2\lambda-1)\bm t/\sqrt n)$ attains the conditional minimum for each sample, its integrated risk equals
	\begin{equation}
		\begin{aligned}
			&\quad\inf_{\hat{\bm w}_n}\left\{
			\frac12\mathbb E_{(Q_{n,\bm u}^{(L)})^{\otimes n}}
			\ell_n(\bm t,\hat{\bm w}_n)
			+\frac12\mathbb E_{(Q_{n,-\bm u}^{(L)})^{\otimes n}}
			\ell_n(-\bm t,\hat{\bm w}_n)\right\}\\
			&=\mathbb E_{P^{\otimes n}}\Bigg[
			\frac{L_n^{(L)}(\bm u)+L_n^{(L)}(-\bm u)}2
			\Big\{\lambda r_{\bm\mu}(\bm t/\sqrt n)
			+(1-\lambda)r_{\bm\mu}(-\bm t/\sqrt n)-r_{\bm\mu}((2\lambda-1)\bm t/\sqrt n)\Big\}\Bigg],
		\end{aligned}
		\label{eq:inf_of_expect}
	\end{equation}
	where
	\begin{equation*}
		\begin{aligned}
			\frac{L_n^{(L)}(\bm u)+L_n^{(L)}(-\bm u)}2\min\{\lambda,1-\lambda\}
			&=\frac{L_n^{(L)}(\bm u)+L_n^{(L)}(-\bm u)}2
			\frac{\min\{L_n^{(L)}(\bm u),L_n^{(L)}(-\bm u)\}}
			{L_n^{(L)}(\bm u)+L_n^{(L)}(-\bm u)}\\
			&=\frac12\min\{L_n^{(L)}(\bm u),L_n^{(L)}(-\bm u)\}.
		\end{aligned}
	\end{equation*}
	Applying \eqref{eq:jensen_gap_lb} inside \eqref{eq:inf_of_expect}  gives
	\begin{equation}
		\mathcal R_{n,M}^{\mathrm{KL}}(P)
		\ge \frac12\left(r_{\bm\mu}(\bm t/\sqrt n)+r_{\bm\mu}(-\bm t/\sqrt n)\right)\mathbb E_{P^{\otimes n}}
		\min\{L_n^{(L)}(\bm u),L_n^{(L)}(-\bm u)\}.
		\label{eq:div_boundary_binary_bayes_bound}
	\end{equation}
	By \eqref{eq:ng_bounded_score_likelihood} and the continuous mapping theorem, for $\bm{Y}\sim\mathcal{N}(0,I_d)$,
	\begin{equation*}
		\left(L_n^{(L)}(\bm u),L_n^{(L)}(-\bm u)\right)\xrightarrow{d}\left(e^{\bm u^\top\bm Y-\|\bm u\|_2^2/2},e^{-\bm u^\top\bm Y-\|\bm u\|_2^2/2}\right),
	\end{equation*}
	which, again by the continuous mapping theorem, yields
	\begin{equation*}
		\min\left\{L_n^{(L)}(\bm u),L_n^{(L)}(-\bm u)\right\}\xrightarrow{d}\min\left\{e^{\bm u^\top\bm Y-\|\bm u\|_2^2/2},e^{-\bm u^\top\bm Y-\|\bm u\|_2^2/2}\right\}=e^{-|\bm u^\top\bm Y|-\|\bm u\|_2^2/2}.
	\end{equation*}
	Note that $\min\{L_n^{(L)}(\bm u),L_n^{(L)}(-\bm u)\}$ is uniformly integrable, since
	\begin{equation*}
		\sup_{n\ \mathrm{large}}\mathbb E_{P^{\otimes n}}
		\min\{L_n^{(L)}(\bm u),L_n^{(L)}(-\bm u)\}^2
		\le\sup_{n\ \mathrm{large}}\mathbb E_{P^{\otimes n}}
		\|L_n^{(L)}\|_\infty^2<\infty.
	\end{equation*}
	Hence
	\begin{equation}
		\begin{aligned}
			\frac12\mathbb E_{P^{\otimes n}}
			\min\{L_n^{(L)}(\bm u),L_n^{(L)}(-\bm u)\}
			&\to\frac12\mathbb E_{\bm Y\sim\mathcal N(\bm0,I_d)} e^{-|\bm u^\top\bm Y|-\|\bm u\|_2^2/2}\\
			&=\frac1{\sqrt{2\pi}}\int_0^\infty
			e^{-(s+\|\bm u\|_2)^2/2}\,\dd s
			=\Phi(-\|\bm u\|_2).
		\end{aligned}
		\label{eq:lb_of_minimum_likelihoods}
	\end{equation}
	By \eqref{eq:div_boundary_directional_remainder},
	\begin{equation*}
		r_{\bm\mu}(\bm t/\sqrt n)+r_{\bm\mu}(-\bm t/\sqrt n)
		\ge 2\frac{\beta}{1+\beta}
		(\bar\kappa\|\bm\mu\|_2)^{-1/\beta}
		h^{(1+\beta)/\beta}n^{-(1+\beta)/(2\beta)}.
	\end{equation*}
	Substituting this into \eqref{eq:div_boundary_binary_bayes_bound} and multiplying by $n^{(1+\beta)/(2\beta)}$ gives, for each sufficiently large fixed $L$ and all sufficiently large $n$,
	\begin{equation*}
		\begin{aligned}
			n^{(1+\beta)/(2\beta)}\mathcal R_{n,M}^{\mathrm{KL}}(P)
			&\ge2\frac{\beta}{1+\beta}
			(\bar\kappa\|\bm\mu\|_2)^{-1/\beta}h^{(1+\beta)/\beta}\left[\frac12\mathbb E_{P^{\otimes n}}
			\min\{L_n^{(L)}(\bm u),L_n^{(L)}(-\bm u)\}\right],
		\end{aligned}
	\end{equation*}
	where $\bm u=J_L^{-1}\bm t$. For this fixed $L$, \eqref{eq:lb_of_minimum_likelihoods} therefore yields
	\begin{equation}
		\begin{aligned}
			\liminf_{n\to\infty}n^{(1+\beta)/(2\beta)}
			\mathcal R_{n,M}^{\mathrm{KL}}(P)
			&\ge2\Phi(-\|J_L^{-1}\bm t\|_2)
			\frac{\beta}{1+\beta}
			(\bar\kappa\|\bm\mu\|_2)^{-1/\beta}h^{(1+\beta)/\beta}.
		\end{aligned}
		\label{eq:fix_l_limit_for_regret}
	\end{equation}
	This inequality holds for every sufficiently large fixed $L$, while its left side is independent of $L$. Since $J_L\to\Sigma^{1/2}$ and $\Sigma\succ0$, continuity of matrix inversion gives $J_L^{-1}\to\Sigma^{-1/2}$ in operator norm as $L\to\infty$. Consequently,
	\begin{equation*}
		\left|\|J_L^{-1}\bm t\|_2-\|\Sigma^{-1/2}\bm t\|_2\right|
		\le\|(J_L^{-1}-\Sigma^{-1/2})\bm t\|_2\le\|J_L^{-1}-\Sigma^{-1/2}\|_{\mathrm{op}}\|\bm t\|_2
		\xrightarrow{L\to\infty}0,
	\end{equation*}
	where $\|\Sigma^{-1/2}\bm t\|_2=h\|\Sigma^{-1/2}U_{\bm\mu}\bm v\|_2$. Letting $L\to\infty$ on the right side of \eqref{eq:fix_l_limit_for_regret}, and using continuity of $\Phi$ gives
	\begin{equation*}
		\begin{aligned}
			\liminf_{n\to\infty}n^{(1+\beta)/(2\beta)}
			\mathcal R_{n,M}^{\mathrm{KL}}(P)
			&\ge2\Phi\!\left(-h\|\Sigma^{-1/2}U_{\bm\mu}\bm v\|_2\right)
			\frac{\beta}{1+\beta}
			(\bar\kappa\|\bm\mu\|_2)^{-1/\beta}h^{(1+\beta)/\beta}.
		\end{aligned}
	\end{equation*}
	This proves \eqref{eq:div_boundary_rate_lower}, with $n\to\infty$ taken first for each fixed $L$, followed by $L\to\infty$.
	
	For the upper bound, uniqueness and compactness imply that every maximizer at $\bm\mu+\bm h$ converges to $\bm w_*(\bm\mu)$ as $\bm h\to\bm0$. Otherwise, a sequence of maximizers outside a fixed neighborhood would have a cluster point $\bm w\ne\bm w_*(\bm\mu)$ satisfying $\bm\mu^\top\bm w=f_*(\bm\mu)$.
	Decrease $\rho$ so that $\bar\kappa\rho^\beta/(1+\beta)\le1$. The upper growth bound in \eqref{eq:div_boundary_two_sided_growth} then gives
	\begin{equation*}
		0\le \varrho_{\bm\mu}(\bm z)\le\frac{\bar\kappa}{1+\beta}\|\bm z\|_2^{1+\beta}
		\le\frac{\bar\kappa}{1+\beta}\rho^{1+\beta}\le\rho,
		\qquad \|\bm z\|_2\le\rho.
	\end{equation*}
	For this smaller fixed $\rho$, all maximizers still belong to the corresponding coordinate neighborhood for sufficiently small $\bm h$. The local feasible points have the form
	\begin{equation*}
		\bm w=\bm w_*(\bm\mu)+U_{\bm\mu}\bm z-\bm n_{\bm\mu}b,
		\qquad \|\bm z\|_2\le\rho,\quad \varrho_{\bm\mu}(\bm z)\le b\le\rho.
	\end{equation*}
	Here $b$ is the inward normal displacement. Since $\bm\mu^\top U_{\bm\mu}=\bm0^\top$ and $\bm\mu^\top\bm n_{\bm\mu}=\|\bm\mu\|_2$,
	\begin{equation*}
		\begin{aligned}
			(\bm\mu+\bm h)^\top(\bm w-\bm w_*(\bm\mu))
			&=(\bm\mu+\bm h)^\top(U_{\bm\mu}\bm z-\bm n_{\bm\mu}b)\\
			&=(U_{\bm\mu}^\top\bm h)^\top\bm z
			-(\|\bm\mu\|_2+\bm n_{\bm\mu}^\top\bm h)b.
		\end{aligned}
	\end{equation*}
	For $\|\bm h\|_2<\|\bm\mu\|_2/2$,
	\begin{equation*}
		\|\bm\mu\|_2+\bm n_{\bm\mu}^\top\bm h
		\ge\|\bm\mu\|_2-\|\bm h\|_2>\frac12\|\bm\mu\|_2>0.
	\end{equation*}
	For each fixed $\bm z$, the objective $(\bm\mu+\bm h)^\top(\bm w-\bm w_*(\bm\mu))$ is therefore decreasing in $b$ and is maximized at $b=\varrho_{\bm\mu}(\bm z)$. By definition of $r_{\bm{\mu}}$, for sufficiently small $\bm h$,
	\begin{equation*}
		\begin{aligned}
			r_{\bm\mu}(\bm h)
			&=\sup_{\bm w\in\Omega}(\bm\mu+\bm h)^\top
			(\bm w-\bm w_*(\bm\mu))\\
			&=\sup_{\substack{\|\bm z\|_2\le\rho\\\varrho_{\bm\mu}(\bm z)\le b\le\rho}}
			\{(U_{\bm\mu}^\top\bm h)^\top\bm z
			-(\|\bm\mu\|_2+\bm n_{\bm\mu}^\top\bm h)b\}\\
			&=\sup_{\|\bm z\|_2\le\rho}
			\{(U_{\bm\mu}^\top\bm h)^\top\bm z
			-(\|\bm\mu\|_2+\bm n_{\bm\mu}^\top\bm h)\varrho_{\bm\mu}(\bm z)\}\\
			&\le\sup_{\|\bm z\|_2\le\rho}
			\left\{(U_{\bm\mu}^\top\bm h)^\top\bm z
			-\frac{\underline\kappa(\|\bm\mu\|_2+\bm n_{\bm\mu}^\top\bm h)}
			{1+\beta}\|\bm z\|_2^{1+\beta}\right\}\\
			&\le\sup_{\bm z\in\mathbb R^{d-1}}
			\left\{(U_{\bm\mu}^\top\bm h)^\top\bm z
			-\frac{\underline\kappa(\|\bm\mu\|_2+\bm n_{\bm\mu}^\top\bm h)}
			{1+\beta}\|\bm z\|_2^{1+\beta}\right\}\\
			&=\sup_{s\ge0}\left\{
			\|U_{\bm\mu}^\top\bm h\|_2s
			-\frac{\underline\kappa(\|\bm\mu\|_2+\bm n_{\bm\mu}^\top\bm h)}
			{1+\beta}s^{1+\beta}\right\}
		\end{aligned}
	\end{equation*}
	where for nonzero $U_{\bm\mu}^\top\bm h$, the last term is maximized at
	\begin{equation*}
		s=\left(
		\frac{\|U_{\bm\mu}^\top\bm h\|_2}
		{\underline\kappa(\|\bm\mu\|_2+\bm n_{\bm\mu}^\top\bm h)}
		\right)^{1/\beta}.
	\end{equation*}
	At this value,
	\begin{equation*}
		\begin{aligned}
			r_{\bm\mu}(\bm h)
			&\le\frac{\beta}{1+\beta}\|U_{\bm\mu}^\top\bm h\|_2
			\left(\frac{\|U_{\bm\mu}^\top\bm h\|_2}
			{\underline\kappa(\|\bm\mu\|_2+\bm n_{\bm\mu}^\top\bm h)}
			\right)^{1/\beta}\\
			&=\frac{\beta}{1+\beta}
			[\underline\kappa(\|\bm\mu\|_2+\bm n_{\bm\mu}^\top\bm h)]^{-1/\beta}
			\|U_{\bm\mu}^\top\bm h\|_2^{(1+\beta)/\beta}.
		\end{aligned}
	\end{equation*}
	If $U_{\bm\mu}^\top\bm h=\bm0$, the supremum is attained at $s=0$ with value zero, and the same bound holds. Consequently, for some $r,C>0$,
	$r_{\bm\mu}(\bm h)\le C\|\bm h\|_2^{(1+\beta)/\beta}$ on $B_2(r)$. Outside this ball, the global Lipschitz bound gives
	\begin{equation*}
		r_{\bm\mu}(\bm h)\le2R_\Omega\|\bm h\|_2
		\le2R_\Omega r^{-1/\beta}\|\bm h\|_2^{(1+\beta)/\beta}.
	\end{equation*}
	After enlarging $C$ we therefore have
	\begin{equation}
		0\le r_{\bm\mu}(\bm h)
		\le C\|\bm h\|_2^{(1+\beta)/\beta},
		\qquad \bm h\in\mathbb R^d.
		\label{eq:div_boundary_global_remainder}
	\end{equation}
	Write $\bm t_Q=\sqrt n(\bm\mu_Q-\bm\mu)$ and retain $\bm Z_{n,Q}$ from \eqref{eq:div_truncated_mean}. The deterministic subgradient argument in \eqref{eq:fast_general_plugin_remainder_bound} applies with the selected maximizer at $\widetilde{\bm\mu}_n=\bm\mu+(\bm t_Q+\bm Z_{n,Q})/\sqrt n$. Together with \eqref{eq:div_boundary_global_remainder}, it gives
	\begin{equation*}
		\begin{aligned}
			R_Q(\bm w_*(\widetilde{\bm\mu}_n))
			&\le r_{\bm\mu}(\bm t_Q/\sqrt n)
			+r_{\bm\mu}((\bm t_Q+2\bm Z_{n,Q})/\sqrt n)\\
			&\le Cn^{-(1+\beta)/(2\beta)}
			\{\|\bm t_Q\|_2^{(1+\beta)/\beta}
			+\|\bm t_Q+2\bm Z_{n,Q}\|_2^{(1+\beta)/\beta}\}.
		\end{aligned}
	\end{equation*}
	The uniform mean bound \eqref{eq:div_mean_shift_bound} and moment bound \eqref{eq:div_truncation_moment_bound} imply
	\begin{equation*}
		\sup_{Q\in\mathcal Q_{n,M}^{\mathrm{KL}}(P)}
		\mathbb E_{Q^{\otimes n}}R_Q(\bm w_*(\widetilde{\bm\mu}_n))
		\le C_{P,M}n^{-(1+\beta)/(2\beta)}.
	\end{equation*}
	Together with \eqref{eq:div_boundary_rate_lower} and the definition of $\mathcal R_{n,M}^{\mathrm{KL}}(P)$, this proves \eqref{eq:div_boundary_matching_rates}.
\end{proof}

\begin{remark}[Boundary regularity and distance]
	The local $C^{1,\beta}$ argument in the proof of Proposition~\ref{prop:boundary-local-parametrization}, for $0<\beta\le1$, verifies \eqref{eq:div_boundary_feasible_growth} and hence the KL lower bound. The matching upper bound additionally uses the lower growth and local feasible-set description in \eqref{eq:div_boundary_two_sided_growth}. For the boundary graph in \eqref{eq:div_boundary_two_sided_growth}, the identities
	\begin{equation*}
		\bm n_{\bm\mu}^\top(\bm w_*(\bm\mu)-\bm w)=\varrho_{\bm\mu}(\bm z),\qquad
		U_{\bm\mu}^\top(\bm w-\bm w_*(\bm\mu))=\bm z,
		\quad \bm w=\bm w_*(\bm\mu)+U_{\bm\mu}\bm z-\bm n_{\bm\mu}\varrho_{\bm\mu}(\bm z),
	\end{equation*}
	show that the local support-distance condition \eqref{eq:fast_boundary_support_distance_growth} is equivalent to
	\begin{equation*}
		\varrho_{\bm\mu}(\bm z)\ge\frac{\underline\kappa}{1+\beta}\|\bm z\|_2^{1+\beta}
	\end{equation*}
	on a sufficiently small coordinate ball. Thus, for $0<\beta\le1$, a local $C^{1,\beta}$ boundary together with \eqref{eq:fast_boundary_support_distance_growth} implies \eqref{eq:div_boundary_two_sided_growth} and the matching KL rates in \eqref{eq:div_boundary_matching_rates}, under the moment assumptions of Theorem~\ref{thm:div_boundary_growth}. In particular, the strong-convexity and positive-principal-curvature conditions verified in Lemma~\ref{lem:quadratic_boundary_growth} give the corresponding quadratic growth conclusions.
	
	If $\dist(\bm\mu,A)>0$, homogeneity \eqref{eq:fast_distance_homogeneity} rewrites the geometric factor in \eqref{eq:div_boundary_rate_lower} as
	\begin{equation*}
		(\bar\kappa\|\bm\mu\|_2)^{-1/\beta}
		=\left(\frac{\dist(\bm n_{\bm\mu},A)}
		{\bar\kappa\,\dist(\bm\mu,A)}\right)^{1/\beta}.
	\end{equation*}
\end{remark}

\begin{remark}[First-order rates and the role of exponential perturbations]
	
	For every compact convex $\Omega$, Corollary~\ref{cor:fast_first_order_limit} and Theorem~\ref{thm:div_kl_minimax} give
	\begin{equation*}
		\lim_{n\to\infty}\sqrt n\,\mathcal R_{n,M}^{\mathrm{KL}}(P)
		=V_{M,I_d}(\psi_1\circ\Sigma^{1/2}).
	\end{equation*}
	Because $\Sigma\succ0$, this constant is positive for every $M>0$ if $\mathcal W_*(\bm\mu)$ contains more than one point; it is zero if $\mathcal W_*(\bm\mu)$ is a singleton. The geometric expansions in the preceding examples therefore give their respective nonzero rates on KL neighborhoods as well, with Gaussian constant $V_{M,I_d}(\psi\circ\Sigma^{1/2})$ and the truncated plug-in policy attaining the rate.
	
	The neighborhood \eqref{eq:div_kl_benchmark} allows all dominated $Q$ satisfying its entropy bound, rather than only the exponential family \eqref{eq:ng_linear_tilt}. Its mean perturbations satisfy
	\begin{equation*}
		\sqrt n(\bm\mu_Q-\bm\mu)\in
		\Sigma^{1/2}B_2(M+o(1)),
	\end{equation*}
	so its fixed-$M$ constant is generally different from $V_{M,\Sigma}(\psi)$ for the Gaussian experiment on a Euclidean mean ball in Section~\ref{sect:fast_rate}. The bounded-score tilts in Lemma~\ref{lem:ng_bounded_score_tilt} supply the lower-bound submodel in the proof of \eqref{eq:div_kl_exact_limit}; the upper bound additionally controls distributions outside that submodel.
\end{remark}

\section{Proofs for Ordinary and Truncated SAA}
\label{app:plugin-comparison}
\subsection{Proof of Theorem~\ref{thm:div_plugin_comparison}}
\label{app:proof-plugin-comparison}
\begin{proof}
	For this proof, write
	\begin{equation*}
		\ell_\psi(\bm t,\bm z)
		:=\psi(\bm t)-\psi(\bm t+\bm z)
		+\nabla\psi(\bm t+\bm z)^\top\bm z
	\end{equation*}
	where $\psi$ is differentiable at $\bm t+\bm z$, and assign value zero on the remaining set. By Rademacher's theorem, the finite convex function $\psi$ on $\mathbb R^d$ is locally Lipschitz and differentiable Lebesgue almost everywhere. Since $\Sigma\succ0$, the assigned values do not affect $\mathbb E\ell_\psi(\bm t,\bm Z)$ for any fixed $\bm t$.
	
	For $\bm u\in\mathbb R^d$, the selected maximizer gives the subgradient
	\begin{equation*}
		\frac{\bm w_*(\bm\mu+\bm u/\sqrt n)-\bm w_*(\bm\mu)}
		{\sqrt n\,a(n^{-1/2})}
		\in\partial\left(\frac{r_{\bm\mu}(\,\cdot\,/\sqrt n)}
		{a(n^{-1/2})}\right)(\bm u).
	\end{equation*}
	Suppose $\bm u_n\to\bm u$ and $\psi$ is differentiable at $\bm u$. For every unit vector $\bm v$ and every $h>0$, the subgradient inequalities give
	\begin{equation*}
		\begin{aligned}
			\frac{r_{\bm\mu}(\bm u_n/\sqrt n)
				-r_{\bm\mu}((\bm u_n-h\bm v)/\sqrt n)}
			{h\,a(n^{-1/2})}
			&\le\bm v^\top
			\frac{\bm w_*(\bm\mu+\bm u_n/\sqrt n)-\bm w_*(\bm\mu)}
			{\sqrt n\,a(n^{-1/2})}\\
			&\le\frac{r_{\bm\mu}((\bm u_n+h\bm v)/\sqrt n)
				-r_{\bm\mu}(\bm u_n/\sqrt n)}
			{h\,a(n^{-1/2})}.
		\end{aligned}
	\end{equation*}
	For fixed $h>0$, $\bm{u}_n$, $\bm u_n-h\bm v$ and $\bm u_n+h\bm v$ belong to the compact ball $B_2(\|\bm{u}\|_2+h+1)$ for sufficiently large $n$. \eqref{eq:fast_general_local_expansion} and continuity of $\psi$ therefore imply
	\begin{equation*}
		\begin{aligned}
			\frac{\psi(\bm u)-\psi(\bm u-h\bm v)}h
			&\le\liminf_{n\to\infty}\bm v^\top
			\frac{\bm w_*(\bm\mu+\bm u_n/\sqrt n)-\bm w_*(\bm\mu)}
			{\sqrt n\,a(n^{-1/2})}\\
			&\le\limsup_{n\to\infty}\bm v^\top
			\frac{\bm w_*(\bm\mu+\bm u_n/\sqrt n)-\bm w_*(\bm\mu)}
			{\sqrt n\,a(n^{-1/2})}
			\le\frac{\psi(\bm u+h\bm v)-\psi(\bm u)}h.
		\end{aligned}
	\end{equation*}
	Letting $h\downarrow0$, the two endpoints converge to $\bm v^\top\nabla\psi(\bm u)$. Taking $\bm v$ to be each coordinate vector proves
	\begin{equation}
		\frac{\bm w_*(\bm\mu+\bm u_n/\sqrt n)-\bm w_*(\bm\mu)}
		{\sqrt n\,a(n^{-1/2})}
		\to\nabla\psi(\bm u).
		\label{eq:div_scaled_subgradient_limit}
	\end{equation}
	
	Now let $Q_n\in\mathcal Q_{n,M}^{\mathrm{KL}}(P)$, suppose $\bm t_{Q_n}\to\bm t$, and let $\hat{\bm\mu}_n$ be any sample-based statistic such that
	\begin{equation*}
		\sqrt n(\hat{\bm\mu}_n-\bm\mu_{Q_n})
		\xrightarrow{d}\bm Z\sim\mathcal N(\bm0,\Sigma).
	\end{equation*}
	With $\bm z=\sqrt n(\hat{\bm\mu}_n-\bm\mu_{Q_n})$, the exact support-function identity is
	\begin{equation*}
		\begin{aligned}
			\frac{R_{Q_n}(\bm w_*(\hat{\bm\mu}_n))}{a(n^{-1/2})}
			={}&\frac{r_{\bm\mu}(\bm t_{Q_n}/\sqrt n)
				-r_{\bm\mu}((\bm t_{Q_n}+\bm z)/\sqrt n)}
			{a(n^{-1/2})}\\
			&+\frac{(\bm w_*(\bm\mu+(\bm t_{Q_n}+\bm z)/\sqrt n)
				-\bm w_*(\bm\mu))^\top\bm z}
			{\sqrt n\,a(n^{-1/2})}.
		\end{aligned}
	\end{equation*}
	On a Skorokhod representation, the centered statistic $\bm z=\sqrt n(\hat{\bm\mu}_n-\bm\mu_{Q_n})$ converges almost surely to $\bm Z$. \eqref{eq:fast_general_local_expansion} and \eqref{eq:div_scaled_subgradient_limit} apply almost surely, because $\psi$ is differentiable at $\bm t+\bm Z$ almost surely. Thus
	\begin{equation}
		\frac{R_{Q_n}(\bm w_*(\hat{\bm\mu}_n))}{a(n^{-1/2})}
		\xrightarrow{d}\ell_\psi(\bm t,\bm Z).
		\label{eq:div_plugin_loss_distribution_limit}
	\end{equation}
	The same argument applies to any subsequence along which $\bm t_{Q_n}$ converges.
	
	Choose $L$ with $\bm t_Q\in B_2(L)$ for all sufficiently large $n$ and $\sup_{B_2(L)}\psi>0$. The deterministic subgradient argument leading to \eqref{eq:fast_general_plugin_remainder_bound} applies to any vectors $\bm t,\bm z$: expand the regret using the selected maximizer at $\bm\mu+(\bm t+\bm z)/\sqrt n$, and compare its support value at $\bm\mu+(\bm t+2\bm z)/\sqrt n$. It gives
	\begin{equation*}
		\begin{aligned}
			&\quad f_*(\bm\mu+\bm t/\sqrt n)
			-\bm w_*(\bm\mu+(\bm t+\bm z)/\sqrt n)^\top(\bm\mu+\bm t/\sqrt n)\\
			&\le r_{\bm\mu}(\bm t/\sqrt n)
			+r_{\bm\mu}((\bm t+2\bm z)/\sqrt n)
			-2r_{\bm\mu}((\bm t+\bm z)/\sqrt n)\\
			&\le b_L(n^{-1/2})+
			b_L\left(n^{-1/2}(1+2\|\bm z\|_2/L)\right),
			\qquad \bm t\in B_2(L).
		\end{aligned}
	\end{equation*}
	Since $b_L(n^{-1/2})/a(n^{-1/2})$ converges to $\sup_{B_2(L)}\psi\in(0,\infty)$, \eqref{eq:fast_radial_polynomial_growth} gives a constant $C$ such that, for $\bm t\in B_2(L)$ and $\bm z\in\mathbb R^d$, for sufficiently large $n$,
	\begin{equation}
		\begin{aligned}
			0&\le\frac{f_*(\bm\mu+\bm t/\sqrt n)-\bm w_*(\bm\mu+(\bm t+\bm z)/\sqrt n)^\top(\bm\mu+\bm t/\sqrt n)}{a(n^{-1/2})}\\
			&\le \frac{b_L(n^{-1/2})+b_L\left(n^{-1/2}(1+2\|\bm z\|_2/L)\right)}{a(n^{-1/2})}\\
			&\le \frac{b_L(n^{-1/2})}{a(n^{-1/2})}\left[1+C_b\left(1+\frac{2\|\bm z\|_2}{L}\right)^{c_b}\right]\\
			&\le C(1+\|\bm z\|_2)^{c_b}.
		\end{aligned}
		\label{eq:div_scaled_plugin_envelope}
	\end{equation}
	Here $C_b, c_b<\infty$ are defined in \eqref{eq:fast_radial_polynomial_growth} for $L$. Note that 
	\begin{equation*}
		\begin{aligned}
			&\quad\frac{f_*(\bm\mu+\bm t/\sqrt n)
				-\bm w_*(\bm\mu+(\bm t+\bm z)/\sqrt n)^\top
				(\bm\mu+\bm t/\sqrt n)}
			{a(n^{-1/2})}\\&=
			\frac{r_{\bm\mu}(\bm t/\sqrt n)
				-r_{\bm\mu}((\bm t+\bm z)/\sqrt n)}
			{a(n^{-1/2})}+
			\left[
			\frac{\bm w_*(\bm\mu+(\bm t+\bm z)/\sqrt n)-\bm w_*(\bm\mu)}
			{\sqrt n\,a(n^{-1/2})}
			\right]^\top\bm z\\
			&\to\ell_\psi(\bm t,\bm z)
		\end{aligned}
	\end{equation*}
	as $n\to\infty$. Passing to the limit in \eqref{eq:div_scaled_plugin_envelope} gives $\ell_\psi(\bm t,\bm z)\le C(1+\|\bm z\|_2)^{c_b}<\infty$ wherever $\psi$ is differentiable at $\bm t+\bm z$. Convexity of $\psi$ also gives $\ell_\psi\ge0$. Hence $\sup_{\bm t\in B_2(L)}\mathbb E\ell_\psi(\bm t,\bm Z)$ is finite.
	
Let $\bm t_j\to\bm t$ and enlarge the fixed radius $L$ to contain this sequence and its limit.
For almost every $\bm z$, $\psi$ is differentiable at $\bm t+\bm z$ and at every $\bm t_j+\bm z$.
The same subgradient inequalities used in \eqref{eq:div_scaled_subgradient_limit}, now applied to the fixed function $\psi$, give
$\nabla\psi(\bm t_j+\bm z)\to\nabla\psi(\bm t+\bm z)$.
Hence $\ell_\psi(\bm t_j,\bm Z)\to\ell_\psi(\bm t,\bm Z)$ almost surely.
The envelope \eqref{eq:div_scaled_plugin_envelope} on the larger ball is Gaussian-integrable, so dominated convergence yields
\begin{equation}
    \mathbb E\ell_\psi(\bm t_j,\bm Z)\to\mathbb E\ell_\psi(\bm t,\bm Z).
    \label{eq:div_expected_limiting_loss_continuity}
\end{equation}
	Thus the function $\bm{t}\mapsto \mathbb{E}\ell_\psi(\bm t,\bm Z)$ is continuous, and the supremum in \eqref{eq:div_gaussian_plugin_constant} is finite, attained, and equal to the supremum over $\|\Sigma^{-1/2}\bm t\|_2<M$.
	
	For $\hat{\bm\mu}_n=\widetilde{\bm\mu}_n$, Lemma~\ref{lem:div_kl_truncation} gives $
	\sqrt n(\widetilde{\bm\mu}_n-\bm\mu_{Q_n})=\bm Z_{n,Q_n}
	\xrightarrow{d}\bm Z$, i.e., the required convergence in distribution. Its moment bounds and \eqref{eq:div_scaled_plugin_envelope} give
	\begin{equation*}
		\sup_{n\ \mathrm{large}}\sup_Q
		\mathbb E_{Q^{\otimes n}}\left[
		\left(\frac{R_Q(\bm w_*(\widetilde{\bm\mu}_n))}{a(n^{-1/2})}\right)^2\right]
		\le C\sup_{n\ \mathrm{large}}\sup_Q
		\mathbb E_{Q^{\otimes n}}(1+\|\bm Z_{n,Q}\|_2)^{2c_b}<\infty.
	\end{equation*}
	Thus the normalized losses are uniformly integrable. More generally, the polynomial envelope requires only a uniformly bounded moment of order $q>c_b$ to obtain a bounded $(1+\rho)$th loss moment for some $\rho>0$. Lemma~\ref{lem:div_kl_truncation} supplies every fixed $q$, so this step remains valid even when $a(n^{-1/2})=o(n^{-1})$. Applying \eqref{eq:div_plugin_loss_distribution_limit} along any sequence $Q_n\in\mathcal Q_{n,M}^{\mathrm{KL}}(P)$ with $\bm t_{Q_n}\to\bm t$ gives
	\begin{equation}
		\frac{\mathbb E_{Q_n^{\otimes n}}R_{Q_n}(\bm w_*(\widetilde{\bm\mu}_n))}{a(n^{-1/2})}
		\to\mathbb E\ell_\psi(\bm t,\bm Z).
		\label{eq:div_truncated_sequential_risk_limit}
	\end{equation}
	Since all possible $\bm t_Q$, for all sufficiently large $n$ and all $Q\in\mathcal Q_{n,M}^{\mathrm{KL}}(P)$, lie in the same fixed compact ball $B_2(L)$, the subsequence argument below gives
	\begin{equation}
		\sup_Q\left|
		\frac{\mathbb E_{Q^{\otimes n}}R_Q(\bm w_*(\widetilde{\bm\mu}_n))}
		{a(n^{-1/2})}
		-\mathbb E\ell_\psi(\bm t_Q,\bm Z)\right|
		\to0.
		\label{eq:div_uniform_truncated_plugin_limit}
	\end{equation}
	Otherwise, if the supremum did not converge to zero, there would be a subsequence and distributions $Q_n$ along which the difference $\left|\mathbb E_{Q_n^{\otimes n}}R_{Q_n}(\bm w_*(\widetilde{\bm\mu}_n))/a(n^{-1/2})-\mathbb E\ell_\psi(\bm t_{Q_n},\bm Z)\right|$ stays bounded away from zero. The bound on $\bm t_{Q_n}$ gives a further subsequence with $\bm t_{Q_n}\to\bm t$. The distributional limit and uniform integrability yield
	\begin{equation*}
		\frac{\mathbb E_{Q_n^{\otimes n}}R_{Q_n}(\bm w_*(\widetilde{\bm\mu}_n))}
		{a(n^{-1/2})}\to\mathbb E\ell_\psi(\bm t,\bm Z),
		\qquad
		\mathbb E\ell_\psi(\bm t_{Q_n},\bm Z)\to\mathbb E\ell_\psi(\bm t,\bm Z),
	\end{equation*}
	a contradiction.
	
	For the ordinary mean $\bar{\bm{x}}$, \eqref{eq:div_truncation_lone_bound} gives
	\begin{equation}
		\sup_Q\mathbb E_{Q^{\otimes n}}\|\bm\Delta_n\|_2
		\le\sup_Q\mathbb E_Q\|\bm x_1-\bm x_1^{(n)}\|_2
		=O(n^{-1}).
		\label{eq:div_mean_truncation_lone}
	\end{equation}
	For every $h>0$, Markov's inequality gives
	\begin{equation*}
		\sup_Q Q^{\otimes n}\{\sqrt n\|\bm\Delta_n\|_2>h\}
		\le\frac{\sqrt n}{h}\sup_Q\mathbb E\|\bm\Delta_n\|_2
		=O(n^{-1/2})\to0.
	\end{equation*}
	Since
	\begin{equation*}
		\sqrt n(\bar{\bm x}-\bm\mu_Q)
		=\bm Z_{n,Q}+\sqrt n\bm\Delta_n,
	\end{equation*}
	Slutsky's theorem gives $\sqrt n(\bar{\bm x}-\bm\mu_{Q_n})
	\xrightarrow{d}\bm Z\sim\mathcal{N}(0,\Sigma)$. Equation \eqref{eq:div_plugin_loss_distribution_limit} therefore applies to $\hat{\bm\mu}_n=\bar{\bm x}$, i.e.,
	\begin{equation*}
		\frac{R_{Q_n}(\bm w_*(\bar{\bm x}))}{a(n^{-1/2})}\xrightarrow{d}\ell_\psi(\bm t,\bm Z).
	\end{equation*}
	For each fixed $K>0$, convergence in distribution implies convergence of the expectations after applying the bounded continuous function $x\mapsto\min\{\max\{x,0\},K\}$, i.e.,
	\begin{equation*}
		\lim_{n\to\infty}\mathbb E_{Q_n^{\otimes n}}\min\left\{\frac{R_{Q_n}(\bm w_*(\bar{\bm x}))}{a(n^{-1/2})},K\right\}=\mathbb E\min\{\ell_\psi(\bm t,\bm Z),K\}.
	\end{equation*}
	$R_{Q_n}(\bm w_*(\bar{\bm x}))$ and $\ell_\psi$ are nonnegative. Letting $K\to\infty$ and using MCT therefore gives
	\begin{equation}
		\liminf_{n\to\infty}
		\frac{\mathbb E_{Q_n^{\otimes n}}R_{Q_n}(\bm w_*(\bar{\bm x}))}
		{a(n^{-1/2})}
		\ge\mathbb E\ell_\psi(\bm t,\bm Z)
		\quad\text{if }\bm t_{Q_n}\to\bm t.
		\label{eq:div_mean_sequential_risk_lower}
	\end{equation}
	By \eqref{eq:div_mean_ellipsoid_bound}, every cluster point of $\bm t_Q$ lies in $\Sigma^{1/2}B_2(M)$. Conversely, for every fixed $\bm t$ with $\|\Sigma^{-1/2}\bm t\|_2<M$, the exponential perturbation $P_{\bm t/\sqrt n}$ has mean $\bm\mu+\bm t/\sqrt n$ and, by \eqref{eq:ng_local_kl_expansion}, satisfies
	\begin{equation*}
		D_{\mathrm{KL}}(P_{\bm t/\sqrt n}\|P)
		=\frac{\bm t^\top\Sigma^{-1}\bm t}{2n}
		+O(n^{-3/2})<\frac{M^2}{2n}
	\end{equation*}
	for large $n$. Fix any $\bm t$ with $\|\Sigma^{-1/2}\bm t\|_2<M$ and take $Q_n=P_{\bm t/\sqrt n}$. The mean parametrization gives $\bm\mu_{Q_n}=\bm\mu+\bm t/\sqrt n$, and $\bm t_{Q_n}=\sqrt n(\bm\mu_{Q_n}-\bm\mu)=\bm t$. Since $Q_n\in\mathcal Q_{n,M}^{\mathrm{KL}}(P)$ for all sufficiently large $n$, the definitions of the maximal regrets imply
	\begin{equation*}
		W_{n,M}^{\mathrm{tr}}(P)
		\ge\mathbb E_{Q_n^{\otimes n}}R_{Q_n}(\bm w_*(\widetilde{\bm\mu}_n)),\qquad
		W_{n,M}^{\mathrm{mean}}(P)
		\ge\mathbb E_{Q_n^{\otimes n}}R_{Q_n}(\bm w_*(\bar{\bm x})).
	\end{equation*}
	Dividing by $a(n^{-1/2})$ and applying \eqref{eq:div_truncated_sequential_risk_limit} and \eqref{eq:div_mean_sequential_risk_lower} yields
	\begin{equation*}
		\begin{aligned}
			&\liminf_{n\to\infty}\frac{W_{n,M}^{\mathrm{tr}}(P)}{a(n^{-1/2})}
			\ge\lim_{n\to\infty}
			\frac{\mathbb E_{Q_n^{\otimes n}}R_{Q_n}(\bm w_*(\widetilde{\bm\mu}_n))}{a(n^{-1/2})}
			=\mathbb E\ell_\psi(\bm t,\bm Z),\\
			&\liminf_{n\to\infty}\frac{W_{n,M}^{\mathrm{mean}}(P)}{a(n^{-1/2})}
			\ge\liminf_{n\to\infty}
			\frac{\mathbb E_{Q_n^{\otimes n}}R_{Q_n}(\bm w_*(\bar{\bm x}))}{a(n^{-1/2})}
			\ge\mathbb E\ell_\psi(\bm t,\bm Z).
		\end{aligned}
	\end{equation*}
	Both left sides are independent of the chosen $\bm t$, so each is at least $\sup_{\|\Sigma^{-1/2}\bm t\|_2<M}\mathbb E\ell_\psi(\bm t,\bm Z)$. 
	By \eqref{eq:div_expected_limiting_loss_continuity},
	\begin{equation*}
		\sup_{\|\Sigma^{-1/2}\bm t\|_2<M}\mathbb E\ell_\psi(\bm t,\bm Z)
		=\sup_{\|\Sigma^{-1/2}\bm t\|_2\le M}\mathbb E\ell_\psi(\bm t,\bm Z)
		=C_M^{\mathrm{plug}}(\psi,\Sigma),
	\end{equation*}
	and therefore
	\begin{equation*}
		\begin{aligned}
			&\liminf_{n\to\infty}\frac{W_{n,M}^{\mathrm{tr}}(P)}{a(n^{-1/2})}
			\ge C_M^{\mathrm{plug}}(\psi,\Sigma),\\
			&\liminf_{n\to\infty}\frac{W_{n,M}^{\mathrm{mean}}(P)}{a(n^{-1/2})}
			\ge C_M^{\mathrm{plug}}(\psi,\Sigma).
		\end{aligned}
	\end{equation*}
For the reverse bound, choose $Q_n$ within $a(n^{-1/2})/n$ of the supremum defining $W_{n,M}^{\mathrm{tr}}$.
The normalized risks are bounded by \eqref{eq:div_truncated_uniform_risk_upper}.
Along a subsequence attaining their limsup, compactness gives a further subsequence with $\bm t_{Q_n}\to\bm t$; \eqref{eq:div_mean_ellipsoid_bound} puts $\bm t\in\Sigma^{1/2}B_2(M)$.
Equations~\eqref{eq:div_uniform_truncated_plugin_limit} and
\eqref{eq:div_expected_limiting_loss_continuity} then give
\[
    \limsup_{n\to\infty}\frac{W_{n,M}^{\mathrm{tr}}(P)}{a(n^{-1/2})}
       =\mathbb E\ell_\psi(\bm t,\bm Z)
       \le C_M^{\mathrm{plug}}(\psi,\Sigma).
\]
	This proves \eqref{eq:div_asymptotic_plugin_comparison}. Finally, Theorem~\ref{thm:div_kl_minimax} and $\mathcal R_{n,M}^{\mathrm{KL}}(P)\le W_{n,M}^{\mathrm{tr}}(P)$ imply
	\begin{equation*}
		C_M^{\mathrm{plug}}(\psi,\Sigma)
		\ge V_{M,I_d}(\psi\circ\Sigma^{1/2})>0.
	\end{equation*}
	The two remaining comparisons follow because $W_{n,M}^{\mathrm{tr}}(P)/a(n^{-1/2})$ converges to this positive constant and $W_{n,M}^{\mathrm{mean}}(P)/a(n^{-1/2})$ has at least the same lower limit.
\end{proof}
\subsection{Proof of Proposition~\ref{prop:div_geometric_tail_criterion}}
\label{app:proof-tail-criterion}
\begin{proof}
	By definition,
	\begin{equation*}
		R_Q(\bm w_*(\bar{\bm x}))=r_{\bm\mu}(\bm\mu_Q-\bm\mu)-r_{\bm\mu}(\bar{\bm{x}}-\bm\mu)+(\bm{w}_*(\bar{\bm{x}})-\bm{w}_*(\bm{\mu}))^\top(\bar{\bm{x}}-\bm{\mu}_Q).
	\end{equation*}
	Combined with the subgradient inequality 
	\begin{equation*}
		r_{\bm\mu}(2\bar{\bm x}-\bm\mu_Q-\bm\mu)\ge r_{\bm\mu}(\bar{\bm{x}}-\bm\mu)+(\bm{w}_*(\bar{\bm{x}})-\bm{w}_*(\bm{\mu}))^\top(\bar{\bm{x}}-\bm{\mu}_Q),
	\end{equation*}
	we have
	\begin{equation*}
		R_Q(\bm w_*(\bar{\bm x}))\le r_{\bm\mu}(\bm\mu_Q-\bm\mu)+r_{\bm\mu}(2\bar{\bm x}-\bm\mu_Q-\bm\mu).
	\end{equation*}
	Since $2\bar{\bm x}-\bm\mu_Q-\bm\mu=(4\widetilde{\bm\mu}_n-2\bm\mu_Q-2\bm\mu)/2+(4\bm\Delta_n)/2$, convexity of $r_{\bm{\mu}}$ yields
	\begin{equation}
		0\le R_Q(\bm w_*(\bar{\bm x}))
		\le r_{\bm\mu}(\bm\mu_Q-\bm\mu)+\frac12 r_{\bm\mu}(4\widetilde{\bm\mu}_n-2\bm\mu_Q-2\bm\mu)
		+\frac12 r_{\bm\mu}(4\bm\Delta_n).
		\label{eq:div_raw_clipped_convex_decomposition}
	\end{equation}
	Write $\bm Z_{n,Q}=\sqrt n(\widetilde{\bm\mu}_n-\bm\mu_Q)$ as in \eqref{eq:div_truncated_mean}. Then
	\begin{equation*}
		\bm\mu_Q-\bm\mu=\frac{\bm t_Q}{\sqrt n},\qquad
		4\widetilde{\bm\mu}_n-2\bm\mu_Q-2\bm\mu
		=\frac{2\bm t_Q+4\bm Z_{n,Q}}{\sqrt n}.
	\end{equation*}
	Choose $L$ with $\bm t_Q\in B_2(L)$ and $\sup_{B_2(L)}\psi>0$. By \eqref{eq:fast_radial_scale_limits} and \eqref{eq:fast_radial_polynomial_growth}, the normalized sum
	\begin{equation*}
		A_{n,Q}:=\frac{r_{\bm\mu}(\bm\mu_Q-\bm\mu)
			+r_{\bm\mu}(4\widetilde{\bm\mu}_n-2\bm\mu_Q-2\bm\mu)/2}
		{a(n^{-1/2})}
	\end{equation*}
	satisfies, by $\|2\bm t_Q+4\bm Z_{n,Q}\|_2\le L(2+4\|\bm Z_{n,Q}\|_2/L)$,
	\begin{equation*}
		0\le A_{n,Q}
		\le\frac{b_L(n^{-1/2})}{a(n^{-1/2})}
		\left[1+\frac{C_b}{2}(2+4\|\bm Z_{n,Q}\|_2/L)^{c_b}\right]
		\le C(1+\|\bm Z_{n,Q}\|_2)^{c_b}.
	\end{equation*}
	for some $C$, $c_b$ which may depend on $L$. Lemma~\ref{lem:div_kl_truncation} therefore gives
	\begin{equation*}
		\sup_{n\ \mathrm{large}}\sup_Q\mathbb E A_{n,Q}^2<\infty.
	\end{equation*}
	Taking expectations in \eqref{eq:div_raw_clipped_convex_decomposition} proves the first inequality in \eqref{eq:div_geometric_tail_upper}. Moreover,
	\begin{equation*}
		0\le r_{\bm\mu}(4\bm\Delta_n)
		\le8R_\Omega\|\bm\Delta_n\|_2.
	\end{equation*}
	Equation \eqref{eq:div_mean_truncation_lone} proves the second inequality. Under \eqref{eq:div_geometric_tail_rate_condition}, the positive minimax limit in Theorem~\ref{thm:div_kl_minimax} supplies the matching lower bound.
	
	Suppose now that \eqref{eq:div_geometric_tail_constant_condition} holds. For nonnegative $x,A,B$ with $x\le A+B$, splitting according to $A>K/2$ gives
	\begin{equation*}
		x\mathbb I\{x>K\}\le2A\mathbb I\{A>K/2\}+2B.
	\end{equation*}
	Indeed, on $\{A\le K/2,x>K\}$ we have $B>K/2\ge A$, hence $x\le2B$; on $\{A>K/2\}$ use $x\le A+B$. Set
	\[
	L_{n,Q}:=\frac{R_Q(\bm w_*(\bar{\bm x}))}{a(n^{-1/2})}.
	\]
	Apply the inequality with $x=L_{n,Q}$, $A=A_{n,Q}$, and
	$B=r_{\bm\mu}(4\bm\Delta_n)/(2a(n^{-1/2}))$.
	Equation~\eqref{eq:div_raw_clipped_convex_decomposition} then implies
	\begin{equation*}
		\begin{aligned}
			\sup_Q\mathbb E[L_{n,Q}\mathbb I\{L_{n,Q}>K\}]
			&\le2\sup_Q\mathbb E[A_{n,Q}\mathbb I\{A_{n,Q}>K/2\}]\\
			&\quad+\frac{\sup_Q\mathbb E r_{\bm\mu}(4\bm\Delta_n)}{a(n^{-1/2})}\\
			&\le\frac4K\sup_Q\mathbb E A_{n,Q}^2
			+\frac{\sup_Q\mathbb E r_{\bm\mu}(4\bm\Delta_n)}{a(n^{-1/2})}.
		\end{aligned}
	\end{equation*}
	By \eqref{eq:div_geometric_tail_constant_condition} and the uniform second-moment bound for $A_{n,Q}$, for every $K>0$,
	\begin{equation*}
		\limsup_{n\to\infty}\sup_Q\mathbb E_{Q^{\otimes n}}\left[
		\frac{R_Q(\bm w_*(\bar{\bm x}))}{a(n^{-1/2})}
		\mathbb I\left\{\frac{R_Q(\bm w_*(\bar{\bm x}))}{a(n^{-1/2})}>K\right\}\right]
		\le\frac{C}{K}.
	\end{equation*}
	Consequently,
	\begin{equation}
		\lim_{K\to\infty}\limsup_{n\to\infty}\sup_Q\mathbb E_{Q^{\otimes n}}\left[\frac{R_Q(\bm w_*(\bar{\bm x}))}{a(n^{-1/2})}\mathbb I\left\{\frac{R_Q(\bm w_*(\bar{\bm x}))}{a(n^{-1/2})}>K\right\}\right]=0.
		\label{eq:asymp_uni_inte_of_mean_regret}
	\end{equation}
	
	Fix any sequence $Q_n\in\mathcal Q_{n,M}^{\mathrm{KL}}(P)$ such that $\bm t_{Q_n}\to\bm t$. Lemma~\ref{lem:div_kl_truncation} gives
	\begin{equation*}
		\sqrt n(\widetilde{\bm\mu}_n-\bm\mu_{Q_n})
		\xrightarrow{d}\bm Z,\qquad \bm Z\sim\mathcal N(\bm0,\Sigma).
	\end{equation*}
	Since $\bm\Delta_n=\bar{\bm x}-\widetilde{\bm\mu}_n$, Markov's inequality and \eqref{eq:div_mean_truncation_lone} imply, for every $\varepsilon>0$,
	\begin{equation*}
		\mathbb P_{Q_n^{\otimes n}}\{\sqrt n\|\bm\Delta_n\|_2>\varepsilon\}
		\le\frac{\sqrt n}{\varepsilon}\mathbb E_{Q_n^{\otimes n}}\|\bm\Delta_n\|_2
		=O(n^{-1/2})\to0.
	\end{equation*}
	Slutsky's theorem therefore gives
	\begin{equation*}
		\sqrt n(\bar{\bm x}-\bm\mu_{Q_n})
		=\sqrt n(\widetilde{\bm\mu}_n-\bm\mu_{Q_n})+\sqrt n\bm\Delta_n
		\xrightarrow{d}\bm Z.
	\end{equation*}
	Thus by \eqref{eq:div_plugin_loss_distribution_limit}, taking $\hat{\bm{\mu}}_n=\bar{\bm{x}}$ yields
	\begin{equation*}
		\frac{R_{Q_n}(\bm w_*(\bar{\bm x}))}{a(n^{-1/2})}
		\xrightarrow{d}\ell_\psi(\bm t,\bm Z).
	\end{equation*}
	For every fixed $K>0$, applying the bounded continuous function $x\mapsto\min\{\max\{x,0\},K\}$ gives
	\begin{equation*}
		\mathbb E_{Q_n^{\otimes n}}\min\left\{
		\frac{R_{Q_n}(\bm w_*(\bar{\bm x}))}{a(n^{-1/2})},K\right\}
		\to\mathbb E\min\{\ell_\psi(\bm t,\bm Z),K\}.
	\end{equation*}
	Since both losses are nonnegative, by \eqref{eq:asymp_uni_inte_of_mean_regret},
	\begin{equation}
		\begin{aligned}
			0&\le\frac{\mathbb E_{Q_n^{\otimes n}}R_{Q_n}(\bm w_*(\bar{\bm x}))}{a(n^{-1/2})}
			-\mathbb E_{Q_n^{\otimes n}}\min\left\{
			\frac{R_{Q_n}(\bm w_*(\bar{\bm x}))}{a(n^{-1/2})},K\right\}\\
			&\le\mathbb E_{Q_n^{\otimes n}}\left[
			\frac{R_{Q_n}(\bm w_*(\bar{\bm x}))}{a(n^{-1/2})}
			\mathbb I\left\{\frac{R_{Q_n}(\bm w_*(\bar{\bm x}))}{a(n^{-1/2})}>K\right\}\right],\\
			0&\le\mathbb E\ell_\psi(\bm t,\bm Z)-\mathbb E\min\{\ell_\psi(\bm t,\bm Z),K\}
			\le\mathbb E[\ell_\psi(\bm t,\bm Z)\mathbb I\{\ell_\psi(\bm t,\bm Z)>K\}].
		\end{aligned}
		\label{eq:both_expectation_can_be_unbounded}
	\end{equation}
	Asymptotic uniform integrability makes 
	\begin{equation*}
		\limsup_{n\to\infty}\mathbb E_{Q_n^{\otimes n}}\left[
		\frac{R_{Q_n}(\bm w_*(\bar{\bm x}))}{a(n^{-1/2})}
		\mathbb I\left\{\frac{R_{Q_n}(\bm w_*(\bar{\bm x}))}{a(n^{-1/2})}>K\right\}\right]\to0,
	\end{equation*}
	as $K\to\infty$. Passing to the limit in \eqref{eq:div_scaled_plugin_envelope}
	and using finiteness of every fixed Gaussian moment gives
	\begin{equation*}
		0\le\mathbb E\ell_\psi(\bm t,\bm Z)
		\le C\mathbb E(1+\|\bm Z\|_2)^{c_b}<\infty.
	\end{equation*}
	Hence $\mathbb E[\ell_\psi(\bm t,\bm Z)\mathbb I\{\ell_\psi(\bm t,\bm Z)>K\}]\to0$ as $K\to\infty$. Letting first $n\to\infty$ and then $K\to\infty$ in \eqref{eq:both_expectation_can_be_unbounded} therefore gives
	\begin{equation}
		\frac{\mathbb E_{Q_n^{\otimes n}}R_{Q_n}(\bm w_*(\bar{\bm x}))}{a(n^{-1/2})}
		\to\mathbb E\ell_\psi(\bm t,\bm Z).
		\label{eq:convergence_of_expectation_of_mean_regret}
	\end{equation}
	
	All possible $\bm t_Q$ belong to the same fixed compact ball $B_2(L)$ for sufficiently large $n$. Suppose there exists some $\varepsilon>0$, a subsequence $n_j$, and $Q_{n_j}\in\mathcal Q_{n_j,M}^{\mathrm{KL}}(P)$ such that
	\begin{equation*}
		\left|\frac{\mathbb E_{Q_{n_j}^{\otimes n_j}}R_{Q_{n_j}}(\bm w_*(\bar{\bm x}))}{a(n_j^{-1/2})}
		-\mathbb E\ell_\psi(\bm t_{Q_{n_j}},\bm Z)\right|\ge\varepsilon.
	\end{equation*}
	Passing to a further subsequence implies that there exists some $\bm{t}\in B_2(L)$, such that $\bm t_{Q_{n_j}}\to\bm t\in B_2(L)$. \eqref{eq:convergence_of_expectation_of_mean_regret} and \eqref{eq:div_expected_limiting_loss_continuity} imply
	\begin{equation*}
		\frac{\mathbb E_{Q_{n_j}^{\otimes n_j}}R_{Q_{n_j}}(\bm w_*(\bar{\bm x}))}{a(n_j^{-1/2})}
		\to\mathbb E\ell_\psi(\bm t,\bm Z),\qquad
		\mathbb E\ell_\psi(\bm t_{Q_{n_j}},\bm Z)\to\mathbb E\ell_\psi(\bm t,\bm Z),
	\end{equation*}
	a contradiction. Thus
	\begin{equation*}
		\sup_Q\left|
		\frac{\mathbb E_{Q^{\otimes n}}R_Q(\bm w_*(\bar{\bm x}))}{a(n^{-1/2})}
		-\mathbb E\ell_\psi(\bm t_Q,\bm Z)\right|\to0.
	\end{equation*}
	In particular,
	\begin{equation*}
		\left|\frac{W_{n,M}^{\mathrm{mean}}(P)}{a(n^{-1/2})}
		-\sup_Q\mathbb E\ell_\psi(\bm t_Q,\bm Z)\right|
		\le\sup_Q\left|
		\frac{\mathbb E_{Q^{\otimes n}}R_Q(\bm w_*(\bar{\bm x}))}{a(n^{-1/2})}
		-\mathbb E\ell_\psi(\bm t_Q,\bm Z)\right|\to0.
	\end{equation*}
	
	By \eqref{eq:div_mean_ellipsoid_bound},
	\begin{equation*}
		\sup_Q\|\Sigma^{-1/2}\bm t_Q\|_2^2\le M^2+o(1).
	\end{equation*}
	For each $n$, choose $Q_n\in\mathcal Q_{n,M}^{\mathrm{KL}}(P)$ such that
	\begin{equation*}
		0\le\sup_{Q\in\mathcal Q_{n,M}^{\mathrm{KL}}(P)}\mathbb E\ell_\psi(\bm t_Q,\bm Z)
		-\mathbb E\ell_\psi(\bm t_{Q_n},\bm Z)\le\frac1n.
	\end{equation*}
	All possible $\bm t_Q$ belong to the same fixed compact ball $B_2(L)$ for sufficiently large $n$. Continuity in \eqref{eq:div_expected_limiting_loss_continuity} therefore implies
	\begin{equation*}
		0\le\sup_{Q\in\mathcal Q_{n,M}^{\mathrm{KL}}(P)}\mathbb E\ell_\psi(\bm t_Q,\bm Z)
		\le\sup_{\bm t\in B_2(L)}\mathbb E\ell_\psi(\bm t,\bm Z)<\infty.
	\end{equation*}
	Choose a subsequence $n_j$ along which
	\begin{equation*}
		\sup_{Q\in\mathcal Q_{n_j,M}^{\mathrm{KL}}(P)}\mathbb E\ell_\psi(\bm t_Q,\bm Z)
		\to\limsup_{n\to\infty}
		\sup_{Q\in\mathcal Q_{n,M}^{\mathrm{KL}}(P)}\mathbb E\ell_\psi(\bm t_Q,\bm Z).
	\end{equation*}
	Since $\bm t_{Q_{n_j}}\in B_2(L)$, compactness gives a further subsequence, still denoted by $n_j$, and $\bm t\in B_2(L)$ such that $\bm t_{Q_{n_j}}\to\bm t$. By \eqref{eq:div_mean_ellipsoid_bound},
	\begin{equation*}
		\|\Sigma^{-1/2}\bm t\|_2^2
		=\lim_{j\to\infty}\|\Sigma^{-1/2}\bm t_{Q_{n_j}}\|_2^2
		\le M^2,
		\qquad \bm t\in\Sigma^{1/2}B_2(M).
	\end{equation*}
	Equation~\eqref{eq:div_expected_limiting_loss_continuity} gives
	\begin{equation*}
		\mathbb E\ell_\psi(\bm t_{Q_{n_j}},\bm Z)
		\to\mathbb E\ell_\psi(\bm t,\bm Z).
	\end{equation*}
	Combining these limits with the approximation error $1/n_j\to0$ yields
	\begin{equation*}
		\begin{aligned}
			\limsup_{n\to\infty}\sup_{Q\in\mathcal Q_{n,M}^{\mathrm{KL}}(P)}
			\mathbb E\ell_\psi(\bm t_Q,\bm Z)
			&=\lim_{j\to\infty}\mathbb E\ell_\psi(\bm t_{Q_{n_j}},\bm Z)=\mathbb E\ell_\psi(\bm t,\bm Z)\\
			&\le\sup_{\bm t\in\Sigma^{1/2}B_2(M)}\mathbb E\ell_\psi(\bm t,\bm Z)
			=C_M^{\mathrm{plug}}(\psi,\Sigma).
		\end{aligned}
	\end{equation*}
	Conversely, fix $\bm t$ with $\|\Sigma^{-1/2}\bm t\|_2<M$. Equation~\eqref{eq:ng_local_kl_expansion} gives
	\begin{equation*}
		nD_{\mathrm{KL}}(P_{\bm t/\sqrt n}\|P)
		=\frac12\bm t^\top\Sigma^{-1}\bm t+O(n^{-1/2})<\frac{M^2}{2}
	\end{equation*}
	for sufficiently large $n$. Thus $Q_n=P_{\bm t/\sqrt n}\in\mathcal Q_{n,M}^{\mathrm{KL}}(P)$ and $\bm t_{Q_n}=\bm t$, so
	\begin{equation*}
		\liminf_{n\to\infty}\sup_Q\mathbb E\ell_\psi(\bm t_Q,\bm Z)
		\ge\sup_{\|\Sigma^{-1/2}\bm t\|_2<M}\mathbb E\ell_\psi(\bm t,\bm Z)
		=C_M^{\mathrm{plug}}(\psi,\Sigma).
	\end{equation*}
	The last equality follows from \eqref{eq:div_expected_limiting_loss_continuity}, by approximating each boundary point $\bm t$ with $(1-j^{-1})\bm t$, $j\ge2$. Consequently,
	\begin{equation*}
		\frac{W_{n,M}^{\mathrm{mean}}(P)}{a(n^{-1/2})}
		\to C_M^{\mathrm{plug}}(\psi,\Sigma).
	\end{equation*}
	Together with \eqref{eq:div_asymptotic_plugin_comparison}, this proves \eqref{eq:div_equal_plugin_constants}.
	
	Finally, \eqref{eq:div_mean_truncation_lone} and $r_{\bm\mu}(4\bm\Delta_n)\le8R_\Omega\|\bm\Delta_n\|_2$ give
	\begin{equation*}
		0\le\frac{\sup_Q\mathbb E_{Q^{\otimes n}}r_{\bm\mu}(4\bm\Delta_n)}{a(n^{-1/2})}
		\le\frac{C}{n a(n^{-1/2})}.
	\end{equation*}
	If $\liminf_{n\to\infty}n a(n^{-1/2})>0$, the right side is bounded, proving \eqref{eq:div_geometric_tail_rate_condition} and hence rate-optimality. If $n a(n^{-1/2})\to\infty$, the right side converges to zero, proving \eqref{eq:div_geometric_tail_constant_condition} and hence \eqref{eq:div_equal_plugin_constants}.
\end{proof}
\subsection{Proof of Corollary~\ref{cor:div_local_order_tail_strength}}
\label{app:proof-power-tail}
\begin{proof}
	Choose $0<c_1<c$. Applying \eqref{eq:div_entropy_inequality} to
	$g(\bm{x})=c_1\|\bm x\|_2^\zeta\mathbb I\{\|\bm x\|_2>\tau_n\}$ gives
	\begin{equation*}
		\sup_Q\mathbb E_Q[
		\|\bm x\|_2^\zeta\mathbb I\{\|\bm x\|_2>\tau_n\}]
		\le\frac1{c_1}\left\{\frac{M^2}{2n}
		+e^{-(c-c_1)\tau_n^\zeta}\mathbb E_Pe^{c\|\bm x\|_2^\zeta}\right\}
		=O(n^{-1}).
	\end{equation*}
	Here $e^{-(c-c_1)\log^{2\zeta}(n+1)}=o(n^{-k})$ for every fixed $k>0$. Let $\bm D_{n,i}:=\bm x_i-\bm x_i^{(n)}$. Since
	\begin{equation*}
		\|\bm D_{n,i}\|_2
		\le\|\bm x_i\|_2\mathbb I\{\|\bm x_i\|_2>\tau_n\},
	\end{equation*}
	we have, uniformly in $Q$,
	\begin{equation*}
		\mathbb E_Q\|\bm D_{n,1}\|_2^\zeta\le\frac Cn,
		\qquad
		\|\mathbb E_Q\bm D_{n,1}\|_2
		\le\frac{C}{n\tau_n^{\zeta-1}}.
	\end{equation*}
	For $\zeta\ge2$ we also have
	\begin{equation*}
		\mathbb E_Q\|\bm D_{n,1}\|_2^2
		\le\mathbb E_Q[\|\bm x\|_2^2\mathbb I\{\|\bm x\|_2>\tau_n\}]
		\le\frac{C}{n\tau_n^{\zeta-2}}.
	\end{equation*}
	
	For independent centered real variables $U_1,\ldots,U_n$ with finite $\zeta$th absolute moments, the von Bahr--Esseen inequality \citep{vonBahrEsseen1965} gives, for $1\le\zeta\le2$,
	\begin{equation}
		\mathbb E\left|\sum_iU_i\right|^\zeta\le2\sum_i\mathbb E|U_i|^\zeta.
		\label{eq:von_bahr_esseen_ineq}  
	\end{equation}
	For $\zeta>2$, Rosenthal's inequality \citep{Rosenthal1970} gives
	\begin{equation}
		\mathbb E\left|\sum_iU_i\right|^\zeta
		\le C_{\zeta}\left\{\left(\sum_i\mathbb E U_i^2\right)^{\zeta/2}+\sum_i\mathbb E|U_i|^\zeta\right\}.
		\label{eq:rosenthal_ineq}
	\end{equation}
	Applying these scalar inequalities coordinatewise gives the corresponding vector inequalities, with constants depending also on the fixed $d$. Note that $\bm D_{n,1},\ldots,\bm D_{n,n}$ are i.i.d. under $Q^{\otimes n}$. Thus
	\begin{equation*}
		\bm\Delta_n
		=\mathbb E_Q\bm D_{n,1}
		+\frac1n\sum_{i=1}^n(\bm D_{n,i}-\mathbb E_Q\bm D_{n,1}),
	\end{equation*}
	where the items in the second term are independent and centered. For $\zeta\ge1$, the inequality $\|\bm a+\bm b\|_2^\zeta\le2^{\zeta-1}(\|\bm a\|_2^\zeta+\|\bm b\|_2^\zeta)$ gives
	\begin{equation*}
		\mathbb E_{Q^{\otimes n}}\|\bm\Delta_n\|_2^\zeta
		\le2^{\zeta-1}\|\mathbb E_Q\bm D_{n,1}\|_2^\zeta+\frac{2^{\zeta-1}}{n^\zeta}\mathbb E_{Q^{\otimes n}}
		\left\|\sum_{i=1}^n(\bm D_{n,i}-\mathbb E_Q\bm D_{n,1})\right\|_2^\zeta.
	\end{equation*}
	Jensen's inequality also gives
	\begin{equation*}
		\begin{aligned}
			\mathbb E_Q\|\bm D_{n,1}-\mathbb E_Q\bm D_{n,1}\|_2^\zeta\le2^{\zeta-1}\left(\mathbb E_Q\|\bm D_{n,1}\|_2^\zeta+\|\mathbb E_Q\bm D_{n,1}\|_2^\zeta\right)\le2^\zeta\mathbb E_Q\|\bm D_{n,1}\|_2^\zeta.
		\end{aligned}
	\end{equation*}
	For $\zeta>2$, the second moment in Rosenthal's inequality satisfies
	\begin{equation*}
		\mathbb E_Q\|\bm D_{n,1}-\mathbb E_Q\bm D_{n,1}\|_2^2
		=\mathbb E_Q\|\bm D_{n,1}\|_2^2-\|\mathbb E_Q\bm D_{n,1}\|_2^2
		\le\mathbb E_Q\|\bm D_{n,1}\|_2^2.
	\end{equation*}
	Applying \eqref{eq:von_bahr_esseen_ineq} and absorbing these factors into $C_{\zeta}$ yields, for $1\le\zeta\le2$,
	\begin{equation*}
		\mathbb E\|\bm\Delta_n\|_2^\zeta
		\le C_{\zeta}\left\{\|\mathbb E_Q\bm D_{n,1}\|_2^\zeta
		+n^{1-\zeta}\mathbb E_Q\|\bm D_{n,1}\|_2^\zeta\right\}
		=O(n^{-\zeta}).
	\end{equation*}
	For $\zeta>2$,
	\begin{equation*}
		\mathbb E\|\bm\Delta_n\|_2^\zeta
		\le C_{\zeta}\|\mathbb E_Q\bm D_{n,1}\|_2^\zeta
		+C_{\zeta} n^{1-\zeta}\mathbb E_Q\|\bm D_{n,1}\|_2^\zeta+C_{\zeta} n^{-\zeta/2}(\mathbb E_Q\|\bm D_{n,1}\|_2^2)^{\zeta/2}
		=O(n^{-\zeta}).
	\end{equation*}
	All these estimates are uniform in $Q$, proving $\sup_Q\mathbb E_{Q^{\otimes n}}\|\bm\Delta_n\|_2^\zeta=O(n^{-\zeta})$.
	
	The uniform local expansion on $B_2(1)$ gives, for some $r>0$,
	\begin{equation*}
		\sup_{\|\bm v\|_2\le1}\frac{r_{\bm\mu}(\varepsilon\bm v)}{\varepsilon^\alpha}
		\le1+\sup_{B_2(1)}\psi<\infty,
		\qquad 0<\varepsilon\le r.
	\end{equation*}
	For $0<\|\bm h\|_2\le r$, take $\varepsilon=\|\bm h\|_2$ and $\bm v=\bm h/\|\bm h\|_2$ to obtain $r_{\bm\mu}(\bm h)\le C\|\bm h\|_2^\alpha$. For $s\in[1,\alpha]$, use $\|\bm h\|_2^\alpha\le r^{\alpha-s}\|\bm h\|_2^s$ on $\|\bm h\|_2\le r$, and $\|\bm h\|_2\le r^{1-s}\|\bm h\|_2^s$ on $\|\bm h\|_2>r$. Together with the global bound $r_{\bm\mu}(\bm h)\le2R_\Omega\|\bm h\|_2$, this yields
	\begin{equation*}
		0\le r_{\bm\mu}(\bm h)\le C_s\|\bm h\|_2^s,
		\qquad \bm h\in\mathbb R^d.
	\end{equation*}
	Taking $s=\min\{\alpha,\zeta\}$ and applying $\sup_Q\mathbb E_{Q^{\otimes n}}\|\bm\Delta_n\|_2^\zeta=O(n^{-\zeta})$ gives
	\begin{equation}
		\sup_Q\mathbb E r_{\bm\mu}(4\bm\Delta_n)
		\le C\left(\sup_Q\mathbb E\|\bm\Delta_n\|_2^\zeta\right)^{\min\{\alpha,\zeta\}/\zeta}
		=O(n^{-\min\{\alpha,\zeta\}}).
		\label{eq:div_power_geometric_tail}
	\end{equation}
	By Proposition~\ref{prop:div_geometric_tail_criterion},
	\begin{equation*}
		W_{n,M}^{\mathrm{mean}}(P)
		\le C\{n^{-\alpha/2}+n^{-\min\{\alpha,\zeta\}}\}
		\le C'\{n^{-\alpha/2}+n^{-\zeta}\}.
	\end{equation*}
	If $\alpha\le2\zeta$, then $\min\{\alpha,\zeta\}\ge\alpha/2$, and \eqref{eq:div_geometric_tail_rate_condition} holds. If $\alpha<2\zeta$, then $\min\{\alpha,\zeta\}>\alpha/2$, and \eqref{eq:div_geometric_tail_constant_condition} holds. The strict rate separation for $\alpha>2\zeta$ is verified in Example~\ref{ex:div_matching_tail_geometry}.
\end{proof}
\subsection{Proof of Corollary~\ref{cor:div_raw_smooth_rate}}
\label{app:proof-raw-smooth}
\begin{proof}
	Choose $r>0$ so that $f_*$ is $C^2$ near $\bm\mu+B_2(r)$, and let $C_r$ bound the operator norm of its Hessian there. All suprema over $Q$ in this proof are over $\mathcal Q_{n,M}^{\mathrm{KL}}(P)$. \eqref{eq:div_mean_truncation_lone} gives
	\begin{equation*}
		\sup_Q\mathbb E_{Q^{\otimes n}}\|\bar{\bm x}-\widetilde{\bm\mu}_n\|_2=O(n^{-1}).
	\end{equation*}
	Recall $\bm Z_{n,Q}=\sqrt n(\widetilde{\bm\mu}_n-\bm\mu_Q)$ from \eqref{eq:div_truncated_mean}. Taking moment order $q=2$ in \eqref{eq:div_truncation_moment_bound} gives
	\begin{equation*}
		\sup_Q\mathbb E_{Q^{\otimes n}}\|\widetilde{\bm\mu}_n-\bm\mu_Q\|_2^2
		=\frac1n\sup_Q\mathbb E_{Q^{\otimes n}}\|\bm Z_{n,Q}\|_2^2
		=O(n^{-1}).
	\end{equation*}
	By \eqref{eq:div_mean_shift_bound}, $\sup_Q\|\bm\mu_Q-\bm\mu\|_2=O(n^{-1/2})$, so for sufficiently large $n$, $\sup_Q\|\bm\mu_Q-\bm\mu\|_2\le r/8$. Consequently, for every such $n$ and every $Q$ in the neighborhood,
	\begin{equation*}
		\{\|\widetilde{\bm\mu}_n-\bm\mu\|_2>r/4\}
		\subseteq\{\|\widetilde{\bm\mu}_n-\bm\mu_Q\|_2>r/8\}=\{\|\bm Z_{n,Q}\|_2>r\sqrt n/8\}.
	\end{equation*}
	Markov's inequality and \eqref{eq:div_truncation_moment_bound} with moment order $q=4$ therefore imply
	\begin{equation*}
		\begin{aligned}
			\sup_Q Q^{\otimes n}\{\|\widetilde{\bm\mu}_n-\bm\mu\|_2>r/4\}
			&\le\sup_Q Q^{\otimes n}\{\|\bm Z_{n,Q}\|_2>r\sqrt n/8\}\\
			&\le\frac{8^4}{r^4n^2}\sup_Q\mathbb E_{Q^{\otimes n}}\|\bm Z_{n,Q}\|_2^4
			=O(n^{-2}),
		\end{aligned}
	\end{equation*}
	where the implied constant may depend on the fixed $r$.
	
	On the event $E_n:=\{\|\widetilde{\bm\mu}_n-\bm\mu\|_2\le r/4,\ \|\bar{\bm x}-\widetilde{\bm\mu}_n\|_2\le r/4\}$, we have
	\begin{equation*}
		\|\bar{\bm x}-\bm\mu\|_2
		\le\|\bar{\bm x}-\widetilde{\bm\mu}_n\|_2
		+\|\widetilde{\bm\mu}_n-\bm\mu\|_2
		\le r/2<r,
		\qquad \|\bm\mu_Q-\bm\mu\|_2\le r/8<r.
	\end{equation*}
	Thus $\bar{\bm x},\bm\mu_Q\in \bm\mu+B_2(r)$, and convexity of the ball gives
	\begin{equation*}
		(1-s)\bar{\bm x}+s\bm\mu_Q\in \bm\mu+B_2(r),\qquad 0\le s\le1.
	\end{equation*}
	Taylor's formula implies
	\begin{equation*}
		\begin{aligned}
			R_Q(\bm w_*(\bar{\bm x}))&=f_*(\bm\mu_Q)-f_*(\bar{\bm x})
			-\nabla f_*(\bar{\bm x})^\top(\bm\mu_Q-\bar{\bm x})\\
			&=\int_0^1(1-s)(\bm\mu_Q-\bar{\bm x})^\top
			\nabla^2f_*\bigl(\bar{\bm x}+s(\bm\mu_Q-\bar{\bm x})\bigr)
			(\bm\mu_Q-\bar{\bm x})\,\dd s\\
			&\le\frac{C_r}{2}\|\bar{\bm x}-\bm\mu_Q\|_2^2\\
			&\le C_r\|\widetilde{\bm\mu}_n-\bm\mu_Q\|_2^2
			+C_r\|\bar{\bm x}-\widetilde{\bm\mu}_n\|_2^2\\
			&\le C_r\|\widetilde{\bm\mu}_n-\bm\mu_Q\|_2^2
			+\frac{C_rr}{4}\|\bar{\bm x}-\widetilde{\bm\mu}_n\|_2.
		\end{aligned}
	\end{equation*}
	On the complement, compactness gives
	$R_Q(\bm w_*(\bar{\bm x}))\le2R_\Omega\|\bm\mu_Q\|_2$, and Markov's inequality gives
	\begin{equation*}
		\begin{aligned}
			\sup_Q Q^{\otimes n}(E_n^c)
			&\le O(n^{-2})
			+\frac4r\sup_Q\mathbb E_{Q^{\otimes n}}
			\|\bar{\bm x}-\widetilde{\bm\mu}_n\|_2
			=O(n^{-1}).
		\end{aligned}
	\end{equation*}
	Taking expectations proves \eqref{eq:div_raw_smooth_rate}. If $H_{\bm\mu}\ne0$, then $\Sigma\succ0$ and $H_{\bm\mu}\succeq0$ imply
	$\bm t\mapsto\bm t^\top\Sigma^{1/2}H_{\bm\mu}\Sigma^{1/2}\bm t/2$ is nonaffine on every $B_2(M)$, $M>0$. Hence \eqref{eq:exact_limit_for_KL_neighborhood} has a positive finite fixed-$M$ constant, and
	\begin{equation*}
		\mathcal R_{n,M}^{\mathrm{KL}}(P)\asymp n^{-1}
		\asymp\sup_{Q\in\mathcal Q_{n,M}^{\mathrm{KL}}(P)}
		\mathbb E_{Q^{\otimes n}}R_Q(\bm w_*(\bar{\bm x})).
	\end{equation*}
\end{proof}

\subsection{Detailed construction for Example~\ref{ex:div_matching_tail_geometry}}
\label{app:tail-geometry-construction}
Fix $\zeta\ge1$ and $\alpha>1$. Let $x_1$ and $x_2$ be independent under $P$, with
\begin{equation}
	\mathbb{P}\{x_1=1/2\}=\mathbb{P}\{x_1=3/2\}=\frac12,
	\qquad
	p_{x_2}(z)=\frac{e^{-|z|^\zeta}}{\int_{\mathbb R}e^{-|u|^\zeta}\,\dd u}.
	\label{eq:div_tail_reference_law}
\end{equation}
Then $\bm\mu=\bm e_1$, $\Sigma\succ0$, and \eqref{eq:div_reference_tail_power} holds for a sufficiently small $c>0$. Take $\Omega=B_{\alpha/(\alpha-1)}^2$, so $f_*(\bm z)=\|\bm z\|_\alpha$. For bounded $\bm t$ and sufficiently small $\varepsilon>0$,
\begin{equation*}
	\begin{split}
		r_{\bm e_1}(\varepsilon\bm t)
		&=(1+\varepsilon t_1)
		\left\{\left(1+
		\left|\frac{\varepsilon t_2}{1+\varepsilon t_1}\right|^\alpha
		\right)^{1/\alpha}-1\right\}\\
		&=\frac{\varepsilon^\alpha|t_2|^\alpha}{\alpha}
		+o(\varepsilon^\alpha),
	\end{split}
\end{equation*}
uniformly on every fixed ball. Thus $\psi(\bm t)=|t_2|^\alpha/\alpha$, $a(\varepsilon)=\varepsilon^\alpha$, and Theorem~\ref{thm:div_kl_minimax} gives, for every $M>0$,
\begin{equation*}
	\mathcal R_{n,M}^{\mathrm{KL}}(P)
	\asymp W_{n,M}^{\mathrm{tr}}(P)\asymp n^{-\alpha/2}.
\end{equation*}
For this family, we claim that the bound in Corollary~\ref{cor:div_local_order_tail_strength} is sharp in rate:
\begin{equation}
	W_{n,M}^{\mathrm{mean}}(P)\asymp n^{-\alpha/2}+n^{-\zeta}.
	\label{eq:div_matching_tail_geometry_rates}
\end{equation}

To prove the additional $n^{-\zeta}$ lower bound, set
\begin{equation}
	B_n=\{n\le|x_2|\le n+1\},\qquad
	\varepsilon_n=\frac{M^2}{4n^{\zeta+1}},\qquad
	Q_n=(1-\varepsilon_n)P+\varepsilon_n P(\,\cdot\mid B_n).
	\label{eq:div_tail_mixture}
\end{equation}
For sufficiently large $n$, $0<\varepsilon_n<1$ and $0<P(B_n)<1$. The conditional law satisfies
\begin{equation*}
	\frac{\dd P(\,\cdot\mid B_n)}{\dd P}
	=\frac{\mathbb I_{B_n}}{P(B_n)},\qquad
	\frac{\dd Q_n}{\dd P}
	=1-\varepsilon_n+\frac{\varepsilon_n}{P(B_n)}\mathbb I_{B_n}.
\end{equation*}
Thus $Q_n\ll P$, and by convexity of $u\mapsto u\log u$ on $[0,\infty)$, with $0\log0=0$,
\begin{equation*}
	\left(1-\varepsilon_n+\frac{\varepsilon_n\mathbb I_{B_n}}{P(B_n)}\right)
	\log\left(1-\varepsilon_n+\frac{\varepsilon_n\mathbb I_{B_n}}{P(B_n)}\right)
	\le1\log 1+\varepsilon_n\frac{\mathbb I_{B_n}}{P(B_n)}
	\log\left(\frac{\mathbb I_{B_n}}{P(B_n)}\right).
\end{equation*}
Integrating with respect to $P$ yields
\begin{equation*}
	D_{\mathrm{KL}}(Q_n\|P)
	\le\varepsilon_n\int_{B_n}\frac{1}{P(B_n)}
	\log\frac1{P(B_n)}\,\dd P
	=\varepsilon_n\log\frac1{P(B_n)}.
\end{equation*}
Let
\[
c_\zeta:=\log\left(\frac12\int_{\mathbb R}e^{-|u|^\zeta}\,\dd u\right).
\]
Then
\begin{equation*}
		n^\zeta+c_\zeta\le\log\frac1{P(B_n)}\le(n+1)^\zeta+c_\zeta,
\end{equation*}
i.e.,
\begin{equation*}
    \log\frac1{P(B_n)}\sim n^\zeta.
\end{equation*}
Consequently,
\begin{equation}
	D_{\mathrm{KL}}(Q_n\|P)
	\le\varepsilon_n\log\frac1{P(B_n)}
	\sim\frac{M^2}{4n}<\frac{M^2}{2n}.
	\label{eq:div_tail_kl_budget}
\end{equation}
By symmetry in $x_2$ and independence of $x_1$, note that $\bm\mu_{Q_n}=\bm e_1$. Since $\mathbb E_P\|\bm x\|_2<\infty$ and $\|\bm x\|_2\le\sqrt{(3/2)^2+(n+1)^2}$ on $B_n$, we also have $\mathbb E_{Q_n}\|\bm x\|_2<\infty$. Hence $Q_n\in\mathcal Q_{n,M}^{\mathrm{KL}}(P)$ for all sufficiently large $n$. 

In the mixture representation, each observation independently selects the component $P$ with probability $1-\varepsilon_n$ or $P(\,\cdot\mid B_n)$ with probability $\varepsilon_n$, and is then drawn from the selected component. Let $N_n$ denote the number selecting $P(\,\cdot\mid B_n)$. Then
\begin{equation}
	N_n\sim\operatorname{Bin}(n,\varepsilon_n),\qquad
	\mathbb P\{N_n=1\}
	=n\varepsilon_n(1-\varepsilon_n)^{n-1}
	\sim\frac{M^2}{4n^\zeta},
	\label{eq:div_tail_single_draw_probability}
\end{equation}
since $n\varepsilon_n\to0$. Conditional on $N_n=1$, the exceptional index is the unique sample index whose selected component is $P(\,\cdot\mid B_n)$. It is uniformly distributed over $\{1,\ldots,n\}$. Since the sample mean is invariant under permutations, its conditional law is the same for every exceptional index. It therefore suffices to calculate with exceptional index $1$. By independence, centering, and Chebyshev's inequality, for every $h>0$,
\begin{equation*}
	\begin{aligned}
		&\mathbb E\left[\frac1n\sum_{i=2}^n x_{i,2}\,\middle|\,N_n=1,
		\text{the exceptional index is }1\right]=0,\\
		&\operatorname{Var}\left(\frac1n\sum_{i=2}^n x_{i,2}\,\middle|\,N_n=1,
		\text{the exceptional index is }1\right)
		=\frac{n-1}{n^2}\operatorname{Var}_P(x_2),\\
		&\mathbb P\left\{\left|\frac1n\sum_{i=2}^n x_{i,2}\right|>h\,\middle|\,N_n=1,
		\text{the exceptional index is }1\right\}
		\le\frac{(n-1)\operatorname{Var}_P(x_2)}{n^2h^2}\to0.
	\end{aligned}
\end{equation*}
On this conditional event, $|x_{1,2}|/n\in[1,1+1/n]$, and the triangle inequality gives
\begin{equation*}
	\begin{aligned}
		\bigl||\bar x_2|-1\bigr|\le\left||\bar x_2|-\frac{|x_{1,2}|}{n}\right|
		+\left|\frac{|x_{1,2}|}{n}-1\right|\le\left|\frac1n\sum_{i=2}^n x_{i,2}\right|+\frac1n.
	\end{aligned}
\end{equation*}
Thus $|\bar x_2|\to1$ in conditional probability, also when conditioning only on $N_n=1$ by the permutation invariance above. Independence of the two coordinates under $P$ and the definition of $B_n$ give $P\{x_1=1/2\mid B_n\}=P\{x_1=3/2\mid B_n\}=\frac12$.
Consequently, for every selection of mixture components, the first coordinates are independent with the same mean $1$ and variance $1/4$. Hence
\begin{equation*}
	\mathbb E[\bar x_1\mid N_n=1]=1,\quad
	\operatorname{Var}(\bar x_1\mid N_n=1)=\frac1{4n},\quad
	\mathbb P\{|\bar x_1-1|>h\mid N_n=1\}\le\frac1{4nh^2}\to0.
\end{equation*}
We have proved
\begin{equation}
	\bar x_1\to1,\qquad |\bar x_2|\to1
	\quad\text{in probability conditional on }N_n=1.
	\label{eq:div_tail_single_draw_mean_limit}
\end{equation}
For $\bm z\ne\bm0$, the maximizer over $\Omega=B_{\alpha/(\alpha-1)}^2$ is
\begin{equation*}
	w_{*,j}(\bm z)=\frac{\operatorname{sgn}(z_j)|z_j|^{\alpha-1}}
	{\|\bm z\|_\alpha^{\alpha-1}},\qquad j=1,2.
\end{equation*}
Since $\bm\mu_{Q_n}=\bm e_1$ and $f_*(\bm e_1)=1$, for $z_1>0$,
\begin{equation}
	\begin{aligned}
		R_{Q_n}(\bm w_*(\bm z))
		&=f_*(\bm e_1)-\bm w_*(\bm z)^\top\bm e_1\\
		&=1-\frac{z_1^{\alpha-1}}{(z_1^\alpha+|z_2|^\alpha)^{(\alpha-1)/\alpha}}
		=1-\left(1+|z_2/z_1|^\alpha\right)^{-(\alpha-1)/\alpha}.
	\end{aligned}
	\label{eq:div_tail_plugin_regret_formula}
\end{equation}
Since $\bar x_1\ge1/2$, this regret belongs to $[0,1]$. Continuity and the conditional convergence $|\bar x_2|/\bar x_1\to1$ imply convergence in conditional probability to $1-2^{-(\alpha-1)/\alpha}$. Boundedness then gives
\begin{equation}
	\mathbb E[R_{Q_n}(\bm w_*(\bar{\bm x}))\mid N_n=1]
	\to1-2^{-(\alpha-1)/\alpha}>0.
	\label{eq:div_tail_single_draw_regret_limit}
\end{equation}
By $Q_n\in\mathcal Q_{n,M}^{\mathrm{KL}}(P)$ and nonnegativity of regret,
\begin{equation*}
	\begin{aligned}
		W_{n,M}^{\mathrm{mean}}(P)
		&\ge\mathbb E_{Q_n^{\otimes n}}R_{Q_n}(\bm w_*(\bar{\bm x}))\\
		&\ge\mathbb P\{N_n=1\}\,
		\mathbb E[R_{Q_n}(\bm w_*(\bar{\bm x}))\mid N_n=1]\\
		&=\frac{M^2}{4n^\zeta}
		\{1-2^{-(\alpha-1)/\alpha}+o(1)\}.
	\end{aligned}
\end{equation*}
The minimax lower bound supplies $W_{n,M}^{\mathrm{mean}}(P)\gtrsim n^{-\alpha/2}$. Combining the two lower bounds with \eqref{eq:div_power_tail_regret_bound} proves \eqref{eq:div_matching_tail_geometry_rates}. In particular, if $\alpha>2\zeta$,
\begin{equation*}
	\frac{W_{n,M}^{\mathrm{mean}}(P)}{W_{n,M}^{\mathrm{tr}}(P)}
	\asymp n^{\alpha/2-\zeta}\to\infty.
\end{equation*}
Hence the truncated policy has strictly smaller maximal regret for every sufficiently large $n$ in this regime. In particular, for $(\zeta,\alpha)=(1,4)$, $x_2$ is standard Laplace, $\Omega=B_{4/3}^2$, and
\begin{equation*}
	W_{n,M}^{\mathrm{mean}}(P)\asymp n^{-1},
	\qquad W_{n,M}^{\mathrm{tr}}(P)\asymp\mathcal R_{n,M}^{\mathrm{KL}}(P)\asymp n^{-2}.
\end{equation*}

For the boundary $(\zeta,\alpha)=(1,2)$, $\Omega=B_2^2$, and $\bm\mu=\bm e_1$. We show that the same construction gives a strict constant gap between ordinary SAA and truncated SAA for every fixed $M>0$. In \eqref{eq:div_tail_reference_law}, $x_2$ now has density $e^{-|z|}/2$, so $\Sigma=\operatorname{diag}(1/4,2)$, $\psi(\bm t)=t_2^2/2$, and $H_{\bm e_1}=\operatorname{diag}(0,1)$. Consequently, Corollary~\ref{cor:div_kl_smooth} gives the truncated-SAA constant
\begin{equation*}
	\lim_{n\to\infty}nW_{n,M}^{\mathrm{tr}}(P)
	=\frac12\operatorname{tr}(\Sigma H_{\bm e_1})=1.
\end{equation*}

We first calculate the contribution from samples drawn entirely from $P$. Since $\bm\mu_P=\bm\mu_{Q_n}=\bm e_1$, \eqref{eq:div_tail_plugin_regret_formula} with $\alpha=2$ gives
\begin{equation*}
	R_P(\bm w_*(\bar{\bm x}))
	=1-\frac{\bar x_1}{\sqrt{\bar x_1^2+\bar x_2^2}}
	=1-\left(1+\frac{\bar x_2^2}{\bar x_1^2}\right)^{-1/2}.
\end{equation*}
Under $P^{\otimes n}$, $\bar x_1$ and $\bar x_2$ are independent. By \eqref{eq:div_tail_reference_law},
\begin{equation*}
	\frac12\le\bar x_1\le\frac32,\qquad
	\mathbb E_{P^{\otimes n}}\bar x_1=1,\qquad
	\operatorname{Var}_{P^{\otimes n}}(\bar x_1)=\frac1{4n},
\end{equation*}
so $\bar x_1\to1$ in probability. The Laplace moments $\mathbb E_Px_2^2=2$ and $\mathbb E_Px_2^4=24$ give
\begin{equation*}
		\mathbb E_{P^{\otimes n}}\bar x_2^2=\frac2n,\qquad
		\mathbb E_{P^{\otimes n}}\bar x_2^4
		=\frac{n\mathbb E_Px_2^4+3n(n-1)(\mathbb E_Px_2^2)^2}{n^4}
		=\frac{12}{n^2}+\frac{12}{n^3}.
\end{equation*}
Taylor's formula yields the uniform remainder bound
\begin{equation*}
	\left|1-(1+u)^{-1/2}-\frac u2\right|
	\le\frac38u^2,\qquad u\ge0.
\end{equation*}
Substituting $u=\bar x_2^2/\bar x_1^2$ and using $\bar x_1^{-4}\le16$ therefore gives
\begin{equation*}
		\left|n\mathbb E_{P^{\otimes n}}R_P(\bm w_*(\bar{\bm x}))
		-\frac n2\mathbb E_{P^{\otimes n}}\frac{\bar x_2^2}{\bar x_1^2}\right|\le6n\mathbb E_{P^{\otimes n}}\bar x_2^4=O(n^{-1}).
\end{equation*}
Independence of $\bar{x}_1$ and $\bar{x}_2$, and $\mathbb E_{P^{\otimes n}}\bar x_2^2=2/n$ imply
\begin{equation}
		n\mathbb E_{P^{\otimes n}}R_P(\bm w_*(\bar{\bm x}))
		=\frac n2\mathbb E_{P^{\otimes n}}[\bar x_1^{-2}]
		\mathbb E_{P^{\otimes n}}\bar x_2^2+O(n^{-1})=\mathbb E_{P^{\otimes n}}[\bar x_1^{-2}]+O(n^{-1})\to1.
	\label{eq:div_boundary_reference_risk}
\end{equation}
The last limit follows from $\bar x_1^{-2}\to1$ in probability and $0\le\bar x_1^{-2}\le4$.

Next, \eqref{eq:div_tail_mixture} gives $\varepsilon_n=M^2/(4n^2)$, and \eqref{eq:div_tail_single_draw_probability} implies
\begin{equation*}
	\mathbb P\{N_n=0\}=(1-\varepsilon_n)^n\to1,\qquad
	n\mathbb P\{N_n=1\}\to\frac{M^2}{4}.
\end{equation*}
Conditional on $N_n=0$, the sample law is $P^{\otimes n}$, so \eqref{eq:div_boundary_reference_risk} supplies the contribution $1$ after multiplying the risk by $n$. Conditional on $N_n=1$, \eqref{eq:div_tail_single_draw_mean_limit} gives $|\bar x_2|/\bar x_1\to1$ in probability. The regret formula \eqref{eq:div_tail_plugin_regret_formula} and its expectation limit \eqref{eq:div_tail_single_draw_regret_limit} then yield
\begin{equation*}
	n\mathbb P\{N_n=1\}\,
	\mathbb E[R_{Q_n}(\bm w_*(\bar{\bm x}))\mid N_n=1]
	\to\frac{M^2}{4}(1-2^{-1/2}).
\end{equation*}
The KL bound \eqref{eq:div_tail_kl_budget} ensures that $Q_n\in\mathcal Q_{n,M}^{\mathrm{KL}}(P)$ for all sufficiently large $n$. Since $\bm\mu_{Q_n}=\bm\mu_P$, retaining the nonnegative contributions from $N_n=0$ and $N_n=1$ gives
\begin{equation*}
	\begin{aligned}
		W_{n,M}^{\mathrm{mean}}(P)
		&\ge\mathbb E_{Q_n^{\otimes n}}R_{Q_n}(\bm w_*(\bar{\bm x}))\\
		&\ge(1-\varepsilon_n)^n\mathbb E_{P^{\otimes n}}R_P(\bm w_*(\bar{\bm x}))+\mathbb P\{N_n=1\}\,
		\mathbb E[R_{Q_n}(\bm w_*(\bar{\bm x}))\mid N_n=1].
	\end{aligned}
\end{equation*}
Multiplying by $n$ and using the preceding limits, we obtain
\begin{equation*}
		\liminf_{n\to\infty}nW_{n,M}^{\mathrm{mean}}(P)
		\ge1+\frac{M^2}{4}(1-2^{-1/2})>1=\lim_{n\to\infty}nW_{n,M}^{\mathrm{tr}}(P).
\end{equation*}
Thus the $N_n=0$ contribution matches the truncated-SAA constant, while the $N_n=1$ contribution makes the ordinary-SAA lower bound strictly larger for every $M>0$.

\begin{remark}[Scope of the optimality criteria]
	\label{rem:div_comparison_criteria_scope}
	The criteria \eqref{eq:div_geometric_tail_rate_condition} and \eqref{eq:div_geometric_tail_constant_condition} retain the direction of the truncation error through $r_{\bm\mu}(4\bm\Delta_n)$. Corollary~\ref{cor:div_local_order_tail_strength} gives sufficient conditions depending only on $\alpha$ and $\zeta$; Example~\ref{ex:div_matching_tail_geometry} proves their sharpness in rate over the stated class of reference distributions. For a particular $(P,\Omega)$, the directional criterion may apply even when these scalar conditions do not.
	
	If $P$ has bounded support, then $\bm x_i^{(n)}=\bm x_i$ $P$-almost surely for all sufficiently large $n$, and hence $Q$-almost surely for every $Q\ll P$. Thus $W_{n,M}^{\mathrm{mean}}(P)=W_{n,M}^{\mathrm{tr}}(P)$ eventually. For unbounded support, \eqref{eq:div_mean_truncation_lone} still gives the same sequential Gaussian fluctuation limit, while \eqref{eq:div_geometric_tail_constant_condition} ensures convergence of the expected normalized regrets. 
\end{remark}

\section{Auxiliary Local-Asymptotic Lemmas}
\label{app:auxiliary-kl}
\begin{lemma}[Duchi--Ruan tilts with bounded scores]
	\label{lem:ng_bounded_score_tilt}
	Let $P$ be a Borel probability measure on $\mathbb R^d$ with
	\begin{equation*}
		\mathbb E_P\|\bm x\|_2^2<\infty,\qquad
		\bm\mu=\mathbb E_P\bm x,\qquad
		\Sigma=\operatorname{Cov}_P(\bm x)\succ0.
	\end{equation*}
	Following \citet[Section~3.1]{DuchiRuan2021}, let $h_{\mathrm{DR}}:\mathbb R\to[-1,1]$ be $C^3$, with $h_{\mathrm{DR}}(s)=s$ for $|s|\le1/2$, $h_{\mathrm{DR}}'\ge0$, and bounded first three derivatives. For $L>0$, define
	\begin{equation*}
		\begin{aligned}
			\bm g_L(\bm x)
			&:=(\bm x-\bm\mu)\mathbb I\{\|\bm x-\bm\mu\|_2\le L\}
			-\mathbb E_P[(\bm x-\bm\mu)\mathbb I\{\|\bm x-\bm\mu\|_2\le L\}],\\
			\Gamma_L&:=\mathbb E_P[\bm g_L\bm g_L^\top].
		\end{aligned}
	\end{equation*}
	Then $\Gamma_L\xrightarrow{L\to\infty}\Sigma$ in operator norm, so $\Gamma_L\succ0$ for sufficiently large $L$. For such $L$, set
	\begin{equation*}
		J_L:=\mathbb E_P[(\bm x-\bm\mu)\bm g_L^\top]\Gamma_L^{-1/2}.
	\end{equation*}
	We have $J_L\xrightarrow{L\to\infty}\Sigma^{1/2}$ in operator norm, and $J_L$ is invertible for sufficiently large $L$.
	
	Fix such an $L$ and $M>0$. In the following assertions, $n\to\infty$ with $L,M$ fixed, and the thresholds for sufficiently large $n$ may depend on $P,L,M$. The suprema over $\bm u\in B_2(M)$ express uniformity in $\bm u$ for these fixed $L,M$. For all sufficiently large $n$ and $\bm u\in B_2(M)$, define
	\begin{equation}
		\frac{\dd Q_{n,\bm u}^{(L)}}{\dd P}(\bm x)
		:=\frac{1+h_{\mathrm{DR}}\bigl(n^{-1/2}\bm u^\top\Gamma_L^{-1/2}\bm g_L(\bm x)\bigr)}
		{\mathbb E_P\bigl[1+h_{\mathrm{DR}}\bigl(n^{-1/2}\bm u^\top\Gamma_L^{-1/2}\bm g_L(\bm x)\bigr)\bigr]}.
		\label{eq:ng_duchi_ruan_tilt}
	\end{equation}
	Uniformly over this ball, these measures are equivalent to $P$, and
	\begin{equation}
		\frac{\dd Q_{n,\bm u}^{(L)}}{\dd P}(\bm x)
		=1+\frac{\bm u^\top\Gamma_L^{-1/2}\bm g_L(\bm x)}{\sqrt n},
		\qquad
		\mathbb E_{Q_{n,\bm u}^{(L)}}\bm x
		=\bm\mu+\frac{J_L\bm u}{\sqrt n}.
		\label{eq:div_bounded_score_submodel}
	\end{equation}
	For every proper lower semicontinuous convex function $\phi:[0,\infty)\to[0,\infty]$ with $\phi(1)=\phi'(1)=0$, $\phi\in C^2$ near $1$, and $0<\phi''(1)<\infty$,
	\begin{equation}
		\sup_{\bm u\in B_2(M)}\left|
		n\mathbb E_P\phi\left(\frac{\dd Q_{n,\bm u}^{(L)}}{\dd P}\right)
		-\frac{\phi''(1)}2\|\bm u\|_2^2\right|\xrightarrow{n\to\infty}0.
		\label{eq:ng_bounded_score_divergence}
	\end{equation}
	In particular, for $\phi(s)=s\log s-s+1$, the expectation is $D_{\mathrm{KL}}(Q_{n,\bm u}^{(L)}\|P)$ and $\phi''(1)=1$.
	
	Under $P^{\otimes n}$, write
	\begin{equation*}
		L_n^{(L)}(\bm u)
		:=\frac{\dd (Q_{n,\bm u}^{(L)})^{\otimes n}}{\dd P^{\otimes n}}
		(\bm x_1,\ldots,\bm x_n).
	\end{equation*}
	These likelihoods have strictly positive continuous sample paths on $B_2(M)$, with
	\begin{equation*}
		\sup_{n\ \mathrm{large}}\mathbb E_{P^{\otimes n}}\|L_n^{(L)}\|_\infty^2<\infty,
	\end{equation*}
	and
	\begin{equation}
		L_n^{(L)}\xrightarrow[n\to\infty]{d}
		\bigl[\bm u\mapsto e^{\bm u^\top\bm Y-\|\bm u\|_2^2/2}\bigr]
		\quad\text{in }C(B_2(M)),\qquad \bm Y\sim\mathcal N(\bm0,I_d).
		\label{eq:ng_bounded_score_likelihood}
	\end{equation}
	The $L\to\infty$ limits concern the approximation of $\bm x-\bm\mu$ by the bounded scores. When combined with the sampling limits, one first lets $n\to\infty$ for each fixed $L$, and then lets $L\to\infty$.
\end{lemma}

\begin{proof}
	Centering and truncation give $\mathbb E_P\bm g_L=\bm0$, $\|\bm g_L\|_2\le2L$, and $\mathbb E_P\|\bm g_L\|_2^2\le\operatorname{tr}(\Sigma)$. Note that
	\begin{equation*}
		\bm g_L-(\bm x-\bm\mu)=-(\bm x-\bm\mu)\mathbb I\{\|\bm x-\bm\mu\|_2>L\}+\mathbb E_P\!\left[(\bm x-\bm\mu)\mathbb I\{\|\bm x-\bm\mu\|_2>L\}\right],
	\end{equation*}
	which, by DCT, yields
	\begin{equation*}
		\mathbb E_P\|\bm g_L-(\bm x-\bm\mu)\|_2^2
		\le\mathbb E_P[\|\bm x-\bm\mu\|_2^2\mathbb I\{\|\bm x-\bm\mu\|_2>L\}]\xrightarrow{L\to\infty}0.
	\end{equation*}
	Cauchy-Schwarz inequality therefore gives
	\begin{equation*}
		\begin{aligned}
			\left\|\mathbb E_P[(\bm x-\bm\mu)\bm g_L^\top]-\Sigma\right\|_{\mathrm{op}}&=\left\|\mathbb E_P\!\left[(\bm x-\bm\mu)(\bm g_L-(\bm x-\bm\mu))^\top\right]\right\|_{\mathrm{op}}\\
			&\le\mathbb E_P\!\left[\|\bm x-\bm\mu\|_2\|\bm g_L-(\bm x-\bm\mu)\|_2\right]\\
			&\le\sqrt{\operatorname{tr}(\Sigma)}\left(\mathbb E_P\|\bm g_L-(\bm x-\bm\mu)\|_2^2\right)^{1/2}\xrightarrow{L\to\infty}0,
		\end{aligned}
	\end{equation*}
	i.e., $\mathbb E_P[(\bm x-\bm\mu)\bm g_L^\top]\xrightarrow{L\to\infty}\Sigma$. Write
	\begin{equation*}
		\bm g_L\bm g_L^\top-(\bm x-\bm\mu)(\bm x-\bm\mu)^\top=(\bm g_L-(\bm x-\bm\mu))\bm g_L^\top+(\bm x-\bm\mu)(\bm g_L-(\bm x-\bm\mu))^\top,
	\end{equation*}
	which further gives
	\begin{equation*}
		\begin{aligned}
			\|\Gamma_L-\Sigma\|_{\mathrm{op}}&\le\left(\mathbb E_P\|\bm g_L-(\bm x-\bm\mu)\|_2^2\right)^{1/2}\left[\left(\mathbb E_P\|\bm g_L\|_2^2\right)^{1/2}+\left(\mathbb E_P\|\bm x-\bm\mu\|_2^2\right)^{1/2}\right]\\
			&\le2\sqrt{\operatorname{tr}(\Sigma)}\left(\mathbb E_P\|\bm g_L-(\bm x-\bm\mu)\|_2^2\right)^{1/2}\xrightarrow{L\to\infty}0.
		\end{aligned}
	\end{equation*}
	Since $\Sigma\succ0$, 
	\begin{equation*}
		J_L\xrightarrow{L\to\infty}\Sigma\Sigma^{-1/2}=\Sigma^{1/2},\quad
		\lambda_{\min}(\Gamma_L)\xrightarrow{L\to\infty}\lambda_{\min}(\Sigma)>0,
		\quad \det J_L\xrightarrow{L\to\infty}\det(\Sigma^{1/2})>0.
	\end{equation*}
	Fix $L$ with $\Gamma_L\succ0$ and $J_L$ invertible, and fix $M>0$. Throughout the remainder of the proof, all limits are as $n\to\infty$ with these $L,M$ fixed. For sufficiently large $n$,
	\begin{equation*}
		\sup_{\bm u\in B_2(M)}\sup_{\bm x\in\mathbb R^d}
		\left|\frac{\bm u^\top\Gamma_L^{-1/2}\bm g_L(\bm x)}{\sqrt n}\right|
		\le\frac{2ML\|\Gamma_L^{-1/2}\|_{\mathrm{op}}}{\sqrt n}\le\frac12.
	\end{equation*}
	Thus $h_{\mathrm{DR}}$ equals the identity at every argument in \eqref{eq:ng_duchi_ruan_tilt}, and its denominator equals
	\begin{equation*}
		1+\frac{\bm u^\top\Gamma_L^{-1/2}\mathbb E_P\bm g_L}{\sqrt n}=1.
	\end{equation*}
	Consequently,
	\begin{equation*}
		\frac12\le\frac{\dd Q_{n,\bm u}^{(L)}}{\dd P}\le\frac32,
		\qquad
		\mathbb E_{Q_{n,\bm u}^{(L)}}\|\bm x\|_2^2
		\le\frac32\mathbb E_P\|\bm x\|_2^2<\infty,
	\end{equation*}
	and
	\begin{equation*}
		\mathbb E_{Q_{n,\bm u}^{(L)}}\bm x-\bm\mu
		=\frac{\mathbb E_P[(\bm x-\bm\mu)\bm g_L^\top]\Gamma_L^{-1/2}\bm u}{\sqrt n}
		=\frac{J_L\bm u}{\sqrt n}.
	\end{equation*}
	This proves \eqref{eq:div_bounded_score_submodel}.
	
	Taylor's formula at $1$ gives, for sufficiently small $|z|$,
	\begin{equation*}
		\left|\phi(1+z)-\frac{\phi''(1)}2z^2\right|
		\le\frac{z^2}{2}\sup_{|v|\le|z|}|\phi''(1+v)-\phi''(1)|.
	\end{equation*}
	Since $\mathbb E_P(\bm u^\top\Gamma_L^{-1/2}\bm g_L)^2=\|\bm u\|_2^2$, when $n\to\infty$,
	\begin{equation*}
		\begin{aligned}
			&\sup_{\bm u\in B_2(M)}\left|
			n\mathbb E_P\phi\left(\frac{\dd Q_{n,\bm u}^{(L)}}{\dd P}\right)
			-\frac{\phi''(1)}2\|\bm u\|_2^2\right|\le\frac{M^2}{2}\sup_{|v|\le2ML\|\Gamma_L^{-1/2}\|_{\mathrm{op}}/\sqrt n}
			|\phi''(1+v)-\phi''(1)|
			\to0.
		\end{aligned}
	\end{equation*}
	This proves \eqref{eq:ng_bounded_score_divergence}.
	
	For $|z|\le1/2$, $|\log(1+z)-z+z^2/2|\le C|z|^3$. Applying this to each likelihood factor gives, uniformly over $\bm u\in B_2(M)$,
	\begin{equation}
		\begin{aligned}
			\log L_n^{(L)}(\bm u)
			&=\bm u^\top\frac1{\sqrt n}\sum_{i=1}^n\Gamma_L^{-1/2}\bm g_L(\bm x_i)\\
			&\quad-\frac12\bm u^\top\left(\frac1n\sum_{i=1}^n
			\Gamma_L^{-1/2}\bm g_L(\bm x_i)\bm g_L(\bm x_i)^\top\Gamma_L^{-1/2}\right)\bm u
			+O_{L,M}(n^{-1/2}).
		\end{aligned}
		\label{eq:duchi_ruan_likelihood_expansion}
	\end{equation}
	Note that $\Gamma_L^{-1/2}\bm g_L(\bm x)$ has mean zero and covariance $I_d$. CLT and LLN yield
	\begin{equation*}
		\begin{aligned}
			&\frac1{\sqrt n}\sum_{i=1}^n\Gamma_L^{-1/2}\bm g_L(\bm x_i)
			\xrightarrow[n\to\infty]{d}\mathcal N(\bm0,I_d),\\
			&\frac1n\sum_{i=1}^n\Gamma_L^{-1/2}\bm g_L(\bm x_i)\bm g_L(\bm x_i)^\top\Gamma_L^{-1/2}
			\xrightarrow[n\to\infty]{\mathbb P}I_d.
		\end{aligned}
	\end{equation*}
	The latter convergence is in operator norm, since $d$ is fixed. Equip $C(B_2(M))$ with the supremum norm. Since $|\bm u^\top A\bm u|\le M^2\|A\|_{\mathrm{op}}$ for $\bm u\in B_2(M)$, \eqref{eq:duchi_ruan_likelihood_expansion} and LLN give
	\begin{equation}
		\begin{aligned}
			&\quad\sup_{\bm u\in B_2(M)}
			\left|\log L_n^{(L)}(\bm u)
			-\left\{\bm u^\top\frac1{\sqrt n}
			\sum_{i=1}^n\Gamma_L^{-1/2}\bm g_L(\bm x_i)
			-\frac12\|\bm u\|_2^2\right\}\right|\\
			&\le\frac{M^2}{2}
			\left\|\frac1n\sum_{i=1}^n
			\Gamma_L^{-1/2}\bm g_L(\bm x_i)\bm g_L(\bm x_i)^\top
			\Gamma_L^{-1/2}-I_d\right\|_{\mathrm{op}}
			+C_{L,M}n^{-1/2}
			\xrightarrow[n\to\infty]{\mathbb P}0.
		\end{aligned}
		\label{eq:duchi_ruan_supremum_err_bound}
	\end{equation}
	The map
	\begin{equation*}
		\bm y\longmapsto
		\left[\bm u\mapsto\bm u^\top\bm y-\frac12\|\bm u\|_2^2\right],
		\qquad \mathbb R^d\longrightarrow C(B_2(M)),
	\end{equation*}
	is continuous, because
	\begin{equation*}
		\sup_{\bm u\in B_2(M)}|\bm u^\top(\bm y-\bm y')|
		\le M\|\bm y-\bm y'\|_2.
	\end{equation*}
	Thus CLT and the continuous mapping theorem imply
	\begin{equation*}
		\left[\bm u\mapsto\bm u^\top\frac1{\sqrt n}
		\sum_{i=1}^n\Gamma_L^{-1/2}\bm g_L(\bm x_i)
		-\frac12\|\bm u\|_2^2\right]
		\xrightarrow[n\to\infty]{d}
		\left[\bm u\mapsto\bm u^\top\bm Y-\frac12\|\bm u\|_2^2\right]
		\quad\text{in }C(B_2(M)),
	\end{equation*}
	where $\bm Y\sim\mathcal N(\bm0,I_d)$. Together with \eqref{eq:duchi_ruan_supremum_err_bound}, Slutsky's theorem in $C(B_2(M))$ yields
	\begin{equation*}
		\log L_n^{(L)}\xrightarrow[n\to\infty]{d}
		\left[\bm u\mapsto\bm u^\top\bm Y-\frac12\|\bm u\|_2^2\right]
		\quad\text{in }C(B_2(M)).
	\end{equation*}
	Finally, exponentiation is continuous on this space. Applying the continuous mapping theorem once more proves
	\begin{equation*}
		L_n^{(L)}\xrightarrow[n\to\infty]{d}
		\left[\bm u\mapsto e^{\bm u^\top\bm Y-\|\bm u\|_2^2/2}\right]
		\quad\text{in }C(B_2(M)),
	\end{equation*}
	which is \eqref{eq:ng_bounded_score_likelihood}. All these limits hold for fixed $L,M$ under $P^{\otimes n}$.
	
	To prove its moment bound, use $\log(1+z)\le z$ to obtain
	\begin{equation*}
		\|L_n^{(L)}\|_\infty^2
		\le\exp\left(2M\left\|\frac1{\sqrt n}\sum_{i=1}^n
		\Gamma_L^{-1/2}\bm g_L(\bm x_i)\right\|_2\right).
	\end{equation*}
	For every fixed $\bm v\in\mathbb R^d$, $e^z\le1+z+z^2e^{|z|}/2$ and centering imply
	\begin{equation*}
		\mathbb E_Pe^{\bm v^\top\Gamma_L^{-1/2}\bm g_L/\sqrt n}
		\le1+\frac{\|\bm v\|_2^2}{2n}
		e^{2L\|\Gamma_L^{-1/2}\bm v\|_2/\sqrt n}.
	\end{equation*}
	For sufficiently large $n$, independence therefore gives
	\begin{equation*}
		\mathbb E_{P^{\otimes n}}\exp\left(
		\bm v^\top\frac1{\sqrt n}\sum_{i=1}^n\Gamma_L^{-1/2}\bm g_L(\bm x_i)\right)
		\le e^{\|\bm v\|_2^2}.
	\end{equation*}
	Finally, $e^{c\|\bm z\|_2}\le\sum_{\bm v\in\{-1,1\}^d}e^{c\bm v^\top\bm z}$ yields, for all sufficiently large $n$,
	\begin{equation*}
		\mathbb E_{P^{\otimes n}}\|L_n^{(L)}\|_\infty^2
		\le\sum_{\bm v\in\{-1,1\}^d}
		\mathbb E_{P^{\otimes n}}\exp\left(
		2M\bm v^\top\frac1{\sqrt n}\sum_{i=1}^n\Gamma_L^{-1/2}\bm g_L(\bm x_i)\right)
		\le2^d e^{4M^2d}<\infty.
	\end{equation*}
\end{proof}

\begin{lemma}[Mean parametrization and the local KL expansion]
	\label{lem:ng_mean_parametrization}
	Under \eqref{eq:ng_reference_distribution}, there exist $\rho_P,\varepsilon_P,C_P>0$ and a $C^\infty$ map $\bm\theta(\bm h)$ on a neighborhood of $B_2(\varepsilon_P)$ such that $\bm\theta(\bm0)=\bm0$, $\nabla\Lambda_P(\bm\theta(\bm h))=\bm h$, and
	\begin{equation}
		\|\bm\theta(\bm h)\|_2\le\min\{\rho_P,C_P\|\bm h\|_2\},
		\qquad \|\bm h\|_2\le\varepsilon_P.
		\label{eq:ng_tilt_parameter_bound}
	\end{equation}
	Moreover,
	\begin{equation}
		\|\bm\theta(\bm h)-\Sigma^{-1}\bm h\|_2\le C_P\|\bm h\|_2^2,
		\qquad \|\bm h\|_2\le\varepsilon_P,
		\label{eq:ng_tilt_parameter_expansion}
	\end{equation}
	and, for $\|\bm\theta\|_2\le\rho_P$,
	\begin{equation}
		0\le\Lambda_P(\bm\theta)\le\frac{C_P}{2}\|\bm\theta\|_2^2,
		\qquad
		\left|\Lambda_P(\bm\theta)-\frac12\bm\theta^\top\Sigma\bm\theta\right|
		\le\frac{C_P}{6}\|\bm\theta\|_2^3.
		\label{eq:ng_cumulant_quadratic_bound}
	\end{equation}
	For $\|\bm h\|_2\le\varepsilon_P$, define
	\begin{equation}
		\frac{\dd P_{\bm h}}{\dd P}(\bm x)
		:=\exp\{\bm\theta(\bm h)^\top(\bm x-\bm\mu)
		-\Lambda_P(\bm\theta(\bm h))\}.
		\label{eq:ng_linear_tilt}
	\end{equation}
	Then $P_{\bm h}$ is a probability measure equivalent to $P$, with $P_{\bm0}=P$ and $\mathbb E_{P_{\bm h}}\bm x=\bm\mu+\bm h$. Uniformly as $\bm h\to\bm0$,
	\begin{equation}
		\Lambda_P(\bm\theta(\bm h))
		=\frac12\bm h^\top\Sigma^{-1}\bm h+O(\|\bm h\|_2^3),
		\label{eq:ng_mean_cumulant_expansion}
	\end{equation}
	and
	\begin{equation}
		D_{\mathrm{KL}}(P_{\bm h}\|P)
		=\bm\theta(\bm h)^\top\bm h-\Lambda_P(\bm\theta(\bm h))
		=\frac12\bm h^\top\Sigma^{-1}\bm h+O(\|\bm h\|_2^3).
		\label{eq:ng_local_kl_expansion}
	\end{equation}
\end{lemma}

\begin{proof}
	Fix $0<\rho<\eta/(4\sqrt d)$. Then
	\begin{equation*}
		\mathbb E_P e^{4\rho\|\bm x-\bm\mu\|_2}
		\le\sum_{\bm v\in\{-1,1\}^d}
		\mathbb E_P e^{4\rho\bm v^\top(\bm x-\bm\mu)}
		=\sum_{\bm v\in\{-1,1\}^d}e^{\Lambda_P(4\rho\bm v)}<\infty.
	\end{equation*}
	For each integer $k\ge0$,
	\begin{equation*}
		\sup_{\|\bm\theta\|_2\le7\rho/2}
		\|\bm x-\bm\mu\|_2^k e^{\bm\theta^\top(\bm x-\bm\mu)}
		\le k!(2/\rho)^k e^{4\rho\|\bm x-\bm\mu\|_2}.
	\end{equation*}
	Differentiation under the integral in \eqref{eq:ng_cumulant_definition} gives $\Lambda_P\in C^\infty$ on a neighborhood of $B_2(3\rho)$, with
	\begin{equation*}
		\nabla\Lambda_P(\bm\theta)
		=\frac{\mathbb E_P[(\bm x-\bm\mu)e^{\bm\theta^\top(\bm x-\bm\mu)}]}
		{\mathbb E_P e^{\bm\theta^\top(\bm x-\bm\mu)}}
	\end{equation*}
	and
	\begin{equation*}
		\nabla^2\Lambda_P(\bm\theta)
		=\frac{\mathbb E_P[(\bm x-\bm\mu)(\bm x-\bm\mu)^\top
			e^{\bm\theta^\top(\bm x-\bm\mu)}]}
		{\mathbb E_P e^{\bm\theta^\top(\bm x-\bm\mu)}}
		-\nabla\Lambda_P(\bm\theta)\nabla\Lambda_P(\bm\theta)^\top.
	\end{equation*}
	In particular, $\Lambda_P(\bm0)=0$, $\nabla\Lambda_P(\bm0)=\bm0$, and $\nabla^2\Lambda_P(\bm0)=\Sigma\succ0$.
	
	Choose $0<\rho_P\le\rho$ with
	\begin{equation*}
		\inf_{\|\bm\theta\|_2\le3\rho_P}
		\lambda_{\min}(\nabla^2\Lambda_P(\bm\theta))\ge\underline\lambda/2>0.
	\end{equation*}
	The inverse function theorem applied to $\nabla\Lambda_P$ at $\bm0$ gives a $C^\infty$ inverse $\bm\theta(\bm h)$ on a neighborhood of $\bm0$. Choose $\varepsilon_P>0$ so that this inverse is defined on a neighborhood of $B_2(\varepsilon_P)$ and $\bm\theta(B_2(\varepsilon_P))\subseteq B_2(\rho_P)$. The chain rule gives
	\begin{equation*}
		\nabla^2\Lambda_P(\bm\theta(\bm h))D\bm\theta(\bm h)=I_d,
	\end{equation*}
	so $D\bm\theta(\bm0)=\Sigma^{-1}$. Taylor's formula yields
	\begin{equation*}
		\bm\theta(\bm h)=\int_0^1D\bm\theta(s\bm h)[\bm h]\,\dd s
	\end{equation*}
	and
	\begin{equation*}
		\bm\theta(\bm h)-\Sigma^{-1}\bm h
		=\int_0^1(1-s)D^2\bm\theta(s\bm h)[\bm h,\bm h]\,\dd s.
	\end{equation*}
	Take $C_P\ge1$ bounding $\|D\bm\theta\|_{\mathrm{op}}$, $\|D^2\bm\theta\|_{\mathrm{op}}$, $\|\nabla^2\Lambda_P\|_{\mathrm{op}}$, and $\|D^3\Lambda_P\|_{\mathrm{op}}$ on the respective compact balls. Then \eqref{eq:ng_tilt_parameter_bound} and \eqref{eq:ng_tilt_parameter_expansion} follow.
	
	By the definition of $P_{\bm h}$,
	\begin{equation*}
		\int\dd P_{\bm h}
		=e^{-\Lambda_P(\bm\theta(\bm h))}
		\mathbb E_P e^{\bm\theta(\bm h)^\top(\bm x-\bm\mu)}=1.
	\end{equation*}
	The derivative formulas for $\Lambda_P$ give
	\begin{equation*}
		\mathbb E_{P_{\bm h}}(\bm x-\bm\mu)
		=\nabla\Lambda_P(\bm\theta(\bm h))=\bm h
	\end{equation*}
	and $P_{\bm0}=P$ since $\bm\theta(\bm0)=\bm0$ and $\Lambda_P(\bm0)=0$. Also, $0<\dd P_{\bm h}/\dd P<\infty$ $P$-almost surely, proving mutual absolute continuity.
	
	Finally, $\Lambda_P(\bm\theta)\ge\mathbb E_P[\bm\theta^\top(\bm x-\bm\mu)]=0$, and Taylor's formula gives
	\begin{equation*}
		\Lambda_P(\bm\theta)
		=\int_0^1(1-s)\bm\theta^\top\nabla^2\Lambda_P(s\bm\theta)\bm\theta\,\dd s
	\end{equation*}
	and
	\begin{equation*}
		\Lambda_P(\bm\theta)-\frac12\bm\theta^\top\Sigma\bm\theta
		=\frac12\int_0^1(1-s)^2D^3\Lambda_P(s\bm\theta)
		[\bm\theta,\bm\theta,\bm\theta]\,\dd s.
	\end{equation*}
	For $\|\bm\theta\|_2\le3\rho_P$, the bound on $D^3\Lambda_P$ gives
	\begin{equation*}
		\left|\Lambda_P(\bm\theta)-\frac12\bm\theta^\top\Sigma\bm\theta\right|
		\le\frac12\int_0^1(1-s)^2 C_P\|\bm\theta\|_2^3\,\dd s
		=\frac{C_P}{6}\|\bm\theta\|_2^3.
	\end{equation*}
	These identities and the bound on $\nabla^2\Lambda_P$ prove \eqref{eq:ng_cumulant_quadratic_bound}.
	
	Since $C_P\ge1$, $\|\Sigma\|_{\mathrm{op}}\le C_P$, and $\|\Sigma^{-1}\|_{\mathrm{op}}\le C_P$, they also give, for $\|\bm h\|_2\le\varepsilon_P$,
	\begin{equation*}
		\begin{split}
			\left|\Lambda_P(\bm\theta(\bm h))-\frac12\bm h^\top\Sigma^{-1}\bm h\right|
			&\le\left|\Lambda_P(\bm\theta(\bm h))-\frac12\bm\theta(\bm h)^\top\Sigma\bm\theta(\bm h)\right|+\frac{1}{2}\left|\bm h^\top\Sigma^{-1}\bm h-\bm\theta(\bm h)^\top\Sigma\bm\theta(\bm h)\right|\\
			&\le \frac{C_P}{6}\|\bm\theta(\bm h)\|_2^3+\frac{1}{2}\left|(\bm\theta(\bm h)-\Sigma^{-1}\bm h)^\top\Sigma(\bm\theta(\bm h)+\Sigma^{-1}\bm h)\right|\\
			&\le \frac{C_P}{6}\|\bm\theta(\bm h)\|_2^3+\frac{\|\Sigma\|_{\mathrm{op}}}{2}\|\bm\theta(\bm h)-\Sigma^{-1}\bm h\|_2\bigl(\|\bm\theta(\bm h)\|_2+\|\Sigma^{-1}\bm h\|_2\bigr).
		\end{split}
	\end{equation*}
	Consequently, \eqref{eq:ng_mean_cumulant_expansion} follows from
	\begin{equation*}
		\left|\Lambda_P(\bm\theta(\bm h))-\frac12\bm h^\top\Sigma^{-1}\bm h\right|
		\le(C_P^4/6+C_P^3)\|\bm h\|_2^3
		\le2C_P^4\|\bm h\|_2^3.
	\end{equation*}
	Finally, taking the expectation of the log density under $P_{\bm h}$ gives
	\begin{equation*}
		\begin{aligned}
			D_{\mathrm{KL}}(P_{\bm h}\|P)
			&=\mathbb E_{P_{\bm h}}\left[\bm\theta(\bm h)^\top(\bm x-\bm\mu)
			-\Lambda_P(\bm\theta(\bm h))\right]\\
			&=\bm\theta(\bm h)^\top\bm h-\Lambda_P(\bm\theta(\bm h)).
		\end{aligned}
	\end{equation*}
	By \eqref{eq:ng_tilt_parameter_expansion},
	\begin{equation*}
		|\bm\theta(\bm h)^\top\bm h-\bm h^\top\Sigma^{-1}\bm h|
		\le\|\bm\theta(\bm h)-\Sigma^{-1}\bm h\|_2\|\bm h\|_2
		\le C_P\|\bm h\|_2^3.
	\end{equation*}
	Together with \eqref{eq:ng_mean_cumulant_expansion}, this yields, uniformly as $\bm h\to\bm0$,
	\begin{equation*}
		D_{\mathrm{KL}}(P_{\bm h}\|P)
		=\frac12\bm h^\top\Sigma^{-1}\bm h+O(\|\bm h\|_2^3).
	\end{equation*}
\end{proof}

\begin{lemma}[A Bayes lower bound from a local likelihood limit]
	\label{lem:ng_likelihood_lower_bound}
	Let $P$ be a reference law with $\mathbb E_P\|\bm x\|_2<\infty$ and $\mathbb E_P\bm x=\bm\mu$. Suppose \eqref{eq:fast_general_local_expansion} holds, and fix $M>0$ and an invertible matrix $\Gamma\in\mathbb R^{d\times d}$. Consider a submodel $Q_{n,\bm u}\ll P$, $\bm u\in B_2(M)$, with
	\begin{equation*}
		\mathbb E_{Q_{n,\bm u}}\bm x=\bm\mu+\frac{\Gamma\bm u}{\sqrt n},
		\qquad
		L_n(\bm u):=\frac{\dd Q_{n,\bm u}^{\otimes n}}{\dd P^{\otimes n}}
		(\bm x_1,\ldots,\bm x_n).
	\end{equation*}
	Assume $L_n$ has strictly positive continuous sample paths on $B_2(M)$ for all sufficiently large $n$, and, under $P^{\otimes n}$,
	\begin{equation}
		L_n\xrightarrow{d} L\quad\text{in }C(B_2(M)),
		\qquad
		L(\bm u)=e^{\bm u^\top\bm Y-\|\bm u\|_2^2/2},
		\quad\bm Y\sim\mathcal N(\bm0,I_d),
		\label{eq:ng_likelihood_process_limit}
	\end{equation}
	with $\sup_{n\ \mathrm{large}}\mathbb E_{P^{\otimes n}}\|L_n\|_\infty^2<\infty$. Then
	\begin{equation}
		\liminf_{n\to\infty}
		\frac{\displaystyle\inf_{\hat{\bm w}_n}
			\sup_{\bm u\in B_2(M)}
			\mathbb E_{Q_{n,\bm u}^{\otimes n}}
			\ell_n(\Gamma\bm u,\hat{\bm w}_n)}{a(n^{-1/2})}
		\ge V_{M,I_d}(\psi\circ\Gamma).
		\label{eq:ng_submodel_lower_bound}
	\end{equation}
	The same lower bound applies to any distribution neighborhood containing this submodel for all sufficiently large $n$. The infimum is over measurable sample-based policies into $\Omega$.
\end{lemma}

\begin{proof}
	Since $C(B_2(M))$ is complete and separable, the Skorokhod representation theorem supplies copies of $L_n,L$ for which $\|L_n-L\|_\infty\to0$ almost surely. Their laws, including the stated moment bounds, are preserved. Also,
	\begin{equation*}
		\mathbb E\|L\|_\infty^2\le\mathbb E e^{2M\|\bm Y\|_2}<\infty.
	\end{equation*}
	For $K>0$,
	\begin{equation*}
		\mathbb E\|L_n-L\|_\infty
		\le\mathbb E\min\{\|L_n-L\|_\infty,K\}
		+K^{-1}\mathbb E\|L_n-L\|_\infty^2.
	\end{equation*}
	Letting first $n\to\infty$ and then $K\to\infty$ gives, on this coupling,
	\begin{equation}
		\mathbb E\|L_n-L\|_\infty\to0.
		\label{eq:ng_likelihood_coupling}
	\end{equation}
	
	For $\pi\in\mathcal P(B_2(M))$ and positive continuous $l$, define
	\begin{equation*}
		J_\pi(l):=\int\psi(\Gamma\bm u)l(\bm u)\,\pi(\dd\bm u)
		-\left(\int l\,\dd\pi\right)
		\psi\left(\Gamma\frac{\int\bm u l(\bm u)\,\pi(\dd\bm u)}{\int l\,\dd\pi}\right).
	\end{equation*}
	All integrals are over $B_2(M)$. Since $\int l\,\dd\pi>0$, the measure $l\,\dd\pi/\int l\,\dd\pi$ is a probability measure on $B_2(M)$, and
	\begin{equation*}
		J_\pi(l)
		=\left(\int l\,\dd\pi\right)
		\left\{\frac{\int\psi(\Gamma\bm u)l(\bm u)\,\pi(\dd\bm u)}{\int l\,\dd\pi}
		-\psi\left(\Gamma\frac{\int\bm u l\,\dd\pi}{\int l\,\dd\pi}\right)\right\}\ge0.
	\end{equation*}
	For $l=L_n$, Bayes' formula identifies $L_n\,\dd\pi/\int L_n\,\dd\pi$ as the posterior distribution, while $\int L_n\,\dd\pi$ is the density of the prior mixture of the sample laws relative to $P^{\otimes n}$. Thus $J_\pi(L_n)$ is the posterior Jensen gap multiplied by this mixture density.
	
	The finite convex function $\psi\circ\Gamma$ is bounded and Lipschitz on $B_2(M)$. Choose $C$ bounding its absolute value and a Lipschitz constant there. For positive continuous $l,k$, their weighted means belong to $B_2(M)$, since
	\begin{equation*}
		\left\|\frac{\int\bm u l\,\dd\pi}{\int l\,\dd\pi}\right\|_2
		\le\frac{\int\|\bm u\|_2l\,\dd\pi}{\int l\,\dd\pi}\le M,
		\qquad
		\left\|\frac{\int\bm u k\,\dd\pi}{\int k\,\dd\pi}\right\|_2\le M.
	\end{equation*}
	The identity
	\begin{equation*}
		\left(\int l\,\dd\pi\right)
		\left(\frac{\int\bm u l\,\dd\pi}{\int l\,\dd\pi}
		-\frac{\int\bm u k\,\dd\pi}{\int k\,\dd\pi}\right)=\int\bm u(l-k)\,\dd\pi
		-\left(\int(l-k)\,\dd\pi\right)
		\frac{\int\bm u k\,\dd\pi}{\int k\,\dd\pi}
	\end{equation*}
	therefore gives
	\begin{equation*}
		\begin{aligned}
			&\left(\int l\,\dd\pi\right)
			\left\|\frac{\int\bm u l\,\dd\pi}{\int l\,\dd\pi}
			-\frac{\int\bm u k\,\dd\pi}{\int k\,\dd\pi}\right\|_2\\
			&\quad\le\left\|\int\bm u(l-k)\,\dd\pi\right\|_2
			+M\left|\int(l-k)\,\dd\pi\right|\\
			&\quad\le2M\|l-k\|_\infty.
		\end{aligned}
	\end{equation*}
The first integral in $J_\pi(l)-J_\pi(k)$ is bounded by $C\|l-k\|_\infty$.
For the posterior term, add and subtract the expression with mass $\int l\,\dd\pi$ and the $k$-weighted mean. The preceding mean bound and the common Lipschitz constant give $(2CM+C)\|l-k\|_\infty$.
	Combining the two bounds gives
	\begin{equation}
		|J_\pi(l)-J_\pi(k)|\le2C(1+M)\|l-k\|_\infty,
		\qquad \pi\in\mathcal P(B_2(M)).
		\label{eq:ng_bayes_functional_stability}
	\end{equation}
	The constant is independent of $\pi$, and no uniform positive lower bound on $\int l\,\dd\pi$ or $\int k\,\dd\pi$ is needed. On the coupling used in \eqref{eq:ng_likelihood_coupling}, preservation of the marginal laws gives, for every $\pi$,
	\begin{equation*}
		\begin{aligned}
			\left|\mathbb E_{P^{\otimes n}}J_\pi(L_n)
			-\mathbb E_{\mathcal N(\bm0,I_d)}J_\pi(L)\right|
			&=\left|\mathbb E[J_\pi(L_n)-J_\pi(L)]\right|\le\mathbb E|J_\pi(L_n)-J_\pi(L)|\\
			&\le2C(1+M)\mathbb E\|L_n-L\|_\infty.
		\end{aligned}
	\end{equation*}
	The expectations on the right are taken on the common coupling space. Taking the supremum over $\pi$ and using \eqref{eq:ng_likelihood_coupling} yields
	\begin{equation}
		\sup_{\pi\in\mathcal P(B_2(M))}
		\left|\mathbb E_{P^{\otimes n}}J_\pi(L_n)
		-\mathbb E_{\mathcal N(\bm0,I_d)}J_\pi(L)\right|
		\le2C(1+M)\mathbb E\|L_n-L\|_\infty\to0.
		\label{eq:ng_uniform_bayes_convergence}
	\end{equation}
	
	Under the prior $\bm T\sim\pi$ and conditional sample law $Q_{n,\bm T}^{\otimes n}$, write $\mathbb E_\pi$ for the joint prior-sampling expectation. Since $\|\bm T\|_2\le M$ almost surely, conditional Jensen's inequality gives $\left\|\mathbb E_\pi[\bm T\mid\bm x_1,\ldots,\bm x_n]\right\|_2\le\mathbb E_\pi[\|\bm T\|_2\mid\bm x_1,\ldots,\bm x_n]\le M$. For a fixed action $\bm w\in\Omega$, the posterior expected loss is
	\begin{equation*}
		\begin{aligned}
			&\quad\mathbb E_\pi[\ell_n(\Gamma\bm T,\bm w)\mid\bm x_1,\ldots,\bm x_n]\\
			&=\mathbb E_\pi\left[f_*\left(\bm\mu+\frac{\Gamma\bm T}{\sqrt n}\right)
			\middle|\bm x_1,\ldots,\bm x_n\right]
			-\bm w^\top\left(\bm\mu+
			\frac{\Gamma\mathbb E_\pi[\bm T\mid\bm x_1,\ldots,\bm x_n]}{\sqrt n}\right).
		\end{aligned}
	\end{equation*}
	The first term is independent of $\bm w$. Maximizing the linear term over $\Omega$ therefore gives the Bayes rule
	\begin{equation*}
		\bm w_*\left(\bm\mu+
		\frac{\Gamma\mathbb E_\pi[\bm T\mid\bm x_1,\ldots,\bm x_n]}{\sqrt n}\right),
	\end{equation*}
	and its integrated regret satisfies
	\begin{equation*}
		\begin{aligned}
			\inf_{\hat{\bm w}_n}\mathbb E_\pi\ell_n(\Gamma\bm T,\hat{\bm w}_n)&=\mathbb E_\pi\left[
			f_*\left(\bm\mu+\frac{\Gamma\bm T}{\sqrt n}\right)
			-f_*\left(\bm\mu+
			\frac{\Gamma\mathbb E_\pi[\bm T\mid\bm x_1,\ldots,\bm x_n]}{\sqrt n}\right)\right]\\
			&=\mathbb E_\pi\left[
			r_{\bm\mu}(\Gamma\bm T/\sqrt n)
			-r_{\bm\mu}\left(
			\Gamma\mathbb E_\pi[\bm T\mid\bm x_1,\ldots,\bm x_n]/\sqrt n\right)\right].
		\end{aligned}
	\end{equation*}
	Bayes' formula gives
	\begin{equation*}
		\mathbb E_\pi[\bm T\mid\bm x_1,\ldots,\bm x_n]
		=\frac{\int\bm uL_n(\bm u)\,\pi(\dd\bm u)}{\int L_n\,\dd\pi}.
	\end{equation*}
	For the first term of $J_\pi(L_n)$, using $\mathbb E_{P^{\otimes n}}L_n(\bm u)=1$ yields
	\begin{equation*}
		\mathbb E_{P^{\otimes n}}\int\psi(\Gamma\bm u)L_n(\bm u)\,\pi(\dd\bm u)
		=\int\psi(\Gamma\bm u)\,\pi(\dd\bm u)
		=\mathbb E_\pi\psi(\Gamma\bm T).
	\end{equation*}
	For the second term, the density of the prior mixture of the sample laws relative to $P^{\otimes n}$ is $\int L_n\,\dd\pi$, so
	\begin{equation*}\mathbb E_{P^{\otimes n}}\left[
		\left(\int L_n\,\dd\pi\right)
		\psi\left(\Gamma\frac{\int\bm uL_n(\bm u)\,\pi(\dd\bm u)}{\int L_n\,\dd\pi}\right)\right]=\mathbb E_\pi\psi\left(
		\Gamma\mathbb E_\pi[\bm T\mid\bm x_1,\ldots,\bm x_n]\right).
	\end{equation*}
	Consequently,
	\begin{equation*}
		\mathbb E_{P^{\otimes n}}J_\pi(L_n)
		=\mathbb E_\pi\left[\psi(\Gamma\bm T)
		-\psi\left(\Gamma\mathbb E_\pi[\bm T\mid\bm x_1,\ldots,\bm x_n]\right)\right].
	\end{equation*}
	Both arguments belong to the same fixed ball:
	\begin{equation*}
		\|\Gamma\bm T\|_2\le\|\Gamma\|_{\mathrm{op}}M,
		\qquad
		\left\|\Gamma\mathbb E_\pi[\bm T\mid\bm x_1,\ldots,\bm x_n]\right\|_2
		\le\|\Gamma\|_{\mathrm{op}}M.
	\end{equation*}
	Write
	\[
	\bm m_n:=\mathbb E_\pi[\bm T\mid\bm x_1,\ldots,\bm x_n].
	\]
	The triangle inequality and the uniform local expansion therefore give
	\begin{equation*}
		\begin{aligned}
			&\quad\left|\frac{\inf_{\hat{\bm w}_n}\mathbb E_\pi
				\ell_n(\Gamma\bm T,\hat{\bm w}_n)}{a(n^{-1/2})}
			-\mathbb E_{P^{\otimes n}}J_\pi(L_n)\right|\\
			&\le\mathbb E_\pi\left|
			\frac{r_{\bm\mu}(\Gamma\bm T/\sqrt n)}{a(n^{-1/2})}
			-\psi(\Gamma\bm T)\right|+\mathbb E_\pi\left|
			\frac{r_{\bm\mu}(\Gamma\bm m_n/\sqrt n)}{a(n^{-1/2})}
			-\psi(\Gamma\bm m_n)\right|\\
			&\le2\sup_{\bm t\in B_2(\|\Gamma\|_{\mathrm{op}}M)}
			\left|\frac{r_{\bm\mu}(\bm t/\sqrt n)}{a(n^{-1/2})}-\psi(\bm t)\right|
			\to0,
		\end{aligned}
	\end{equation*}
	uniformly in $\pi$. In the Gaussian experiment, the same mixture calculation gives
	\begin{equation*}
		\mathbb E_{\mathcal N(\bm0,I_d)}J_\pi(L)
		=\mathbb E_{\pi,I_d}
		[\psi(\Gamma\bm T)-\psi(\Gamma\bm m_{\pi,I_d}(\bm Y))].
	\end{equation*}
	For every sample-based policy $\hat{\bm w}_n$ and every $\pi\in\mathcal P(B_2(M))$,
	\begin{equation*}
		\sup_{\bm u\in B_2(M)}
		\mathbb E_{Q_{n,\bm u}^{\otimes n}}\ell_n(\Gamma\bm u,\hat{\bm w}_n)
		\ge\int_{B_2(M)}
		\mathbb E_{Q_{n,\bm u}^{\otimes n}}\ell_n(\Gamma\bm u,\hat{\bm w}_n)\,\pi(\dd\bm u)=\mathbb E_\pi\ell_n(\Gamma\bm T,\hat{\bm w}_n).
	\end{equation*}
	Taking the infimum over the same policy class on both sides and then the supremum over $\pi$ gives
	\begin{equation*}
		\inf_{\hat{\bm w}_n}\sup_{\bm u\in B_2(M)}
		\mathbb E_{Q_{n,\bm u}^{\otimes n}}\ell_n(\Gamma\bm u,\hat{\bm w}_n)
		\ge\sup_{\pi\in\mathcal P(B_2(M))}\inf_{\hat{\bm w}_n}
		\mathbb E_\pi\ell_n(\Gamma\bm T,\hat{\bm w}_n).
	\end{equation*}
	The preceding uniform approximation and \eqref{eq:ng_uniform_bayes_convergence} imply
	\begin{equation*}
		\sup_{\pi\in\mathcal P(B_2(M))}
		\left|\frac{\inf_{\hat{\bm w}_n}\mathbb E_\pi
			\ell_n(\Gamma\bm T,\hat{\bm w}_n)}{a(n^{-1/2})}
		-\mathbb E_{\mathcal N(\bm0,I_d)}J_\pi(L)\right|
		\to0.
	\end{equation*}
	Using $|\sup u-\sup v|\le\sup|u-v|$ on $\mathcal P(B_2(M))$, we therefore obtain
	\begin{equation*}
		\sup_{\pi\in\mathcal P(B_2(M))}
		\frac{\inf_{\hat{\bm w}_n}\mathbb E_\pi
			\ell_n(\Gamma\bm T,\hat{\bm w}_n)}{a(n^{-1/2})}
		\to
		\sup_{\pi\in\mathcal P(B_2(M))}
		\mathbb E_{\mathcal N(\bm0,I_d)}J_\pi(L).
	\end{equation*}
	By the Gaussian identity above and the definition of $V_{M,I_d}$,
	\begin{equation*}
		\begin{aligned}
			\sup_{\pi\in\mathcal P(B_2(M))}
			\mathbb E_{\mathcal N(\bm0,I_d)}J_\pi(L)
			&=\sup_{\pi\in\mathcal P(B_2(M))}
			\mathbb E_{\pi,I_d}\left[
			(\psi\circ\Gamma)(\bm T)
			-(\psi\circ\Gamma)(\bm m_{\pi,I_d}(\bm Y))\right]\\
			&=V_{M,I_d}(\psi\circ\Gamma).
		\end{aligned}
	\end{equation*}
	Consequently,
	\begin{equation*}
		\begin{aligned}
			&\quad\liminf_{n\to\infty}
			\frac{\inf_{\hat{\bm w}_n}\sup_{\bm u\in B_2(M)}
				\mathbb E_{Q_{n,\bm u}^{\otimes n}}\ell_n(\Gamma\bm u,\hat{\bm w}_n)}
			{a(n^{-1/2})}\\
			&\ge\lim_{n\to\infty}\sup_{\pi\in\mathcal P(B_2(M))}
			\frac{\inf_{\hat{\bm w}_n}\mathbb E_\pi
				\ell_n(\Gamma\bm T,\hat{\bm w}_n)}{a(n^{-1/2})}
			=V_{M,I_d}(\psi\circ\Gamma).
		\end{aligned}
	\end{equation*}
	This proves \eqref{eq:ng_submodel_lower_bound} using the minimax--Bayes inequality.
\end{proof}

\begin{lemma}[Gaussian comparison on an ellipsoid]
	\label{lem:div_gaussian_comparison}
	Suppose \eqref{eq:fast_general_local_expansion} holds. For an invertible matrix $\Gamma\in\mathbb R^{d\times d}$, consider the observation $\bm Y\sim\mathcal N(\bm u,I_d)$, $\bm u\in B_2(M)$, and loss
	$\ell_n(\Gamma\bm u,\bm w)$. Then
	\begin{equation}
		\lim_{n\to\infty}
		\frac{\displaystyle\inf_\delta\sup_{\bm u\in B_2(M)}
			\mathbb E_{\mathcal N(\bm u,I_d)}
			\int_\Omega\ell_n(\Gamma\bm u,\bm w)\,
			\delta(\dd\bm w\mid\bm Y)}
		{a(n^{-1/2})}
		=V_{M,I_d}(\psi\circ\Gamma).
		\label{eq:div_gaussian_ellipsoid}
	\end{equation}
	Here $\delta$ ranges over Markov kernels from $\mathbb R^d$ to $\Omega$, with $\delta(\dd\bm w\mid\bm y)$ specifying the conditional distribution of $\bm w$ given $\bm Y=\bm y$. The infimum is unchanged by restricting the kernels $\delta$ to $\Omega_{n,L}$ for any $L\ge\|\Gamma\|_{\mathrm{op}}M$. This restriction means $\delta(\Omega_{n,L}\mid\bm y)=1$ for every $\bm y\in\mathbb R^d$. The map
	$(M,\Gamma)\mapsto V_{M,I_d}(\psi\circ\Gamma)$ is continuous for $M>0$ and invertible $\Gamma$.
\end{lemma}

\begin{proof}
	Fix $n$ and $L\ge\|\Gamma\|_{\mathrm{op}}M$. The parameter space $B_2(M)$ is compact, and the experiment $\bm Y\sim\mathcal N(\bm u,I_d)$ has a strictly positive Lebesgue density. For $\bm u,\bm v\in B_2(M)$, the Gaussian KL formula and Pinsker's inequality give
	\begin{equation*}
		D_{\mathrm{KL}}(\mathcal N(\bm u,I_d)\|\mathcal N(\bm v,I_d))
		=\frac12\|\bm u-\bm v\|_2^2,
		\qquad
		\|\mathcal N(\bm u,I_d)-\mathcal N(\bm v,I_d)\|_{\mathrm{TV}}
		\le\frac12\|\bm u-\bm v\|_2.
	\end{equation*}
	The loss $(\bm u,\bm w)\mapsto\ell_n(\Gamma\bm u,\bm w)$ is jointly continuous and affine in $\bm w$, and
	\begin{equation*}
		0\le\ell_n(\Gamma\bm u,\bm w)
		\le2R_\Omega(\|\bm\mu\|_2+L/\sqrt n),
		\qquad (\bm u,\bm w)\in B_2(M)\times\Omega.
	\end{equation*}
	For either decision space $\Omega$ or $\Omega_{n,L}$, compactness and convexity hold; for $\Omega_{n,L}$ these were proved in Lemma~\ref{lem:restricted_decision_set}. The class of all Markov kernels from $\mathbb R^d$ into the chosen decision space is nonempty, convex, and weakly closed in the sense used in Lemma~\ref{lem:cov_compact_minimax}. The same statistical minimax argument therefore applies to this experiment and loss.
	
	For $\pi\in\mathcal P(B_2(M))$, let $\bm T\sim\pi$ and $\bm Y=\bm T+\bm G$, where $\bm G\sim\mathcal N(\bm0,I_d)$ is independent of $\bm T$. Then $\bm m_{\pi,I_d}(\bm Y)=\mathbb E_\pi[\bm T\mid\bm Y]$, and $\|\bm m_{\pi,I_d}(\bm Y)\|_2\le\mathbb E_\pi[\|\bm T\|_2\mid\bm Y]\le M$. For $\bm w\in\Omega$, the posterior expected loss is
	\begin{equation*}
		\mathbb E_\pi[\ell_n(\Gamma\bm T,\bm w)\mid\bm Y]
		=\mathbb E_\pi\left[f_*\left(\bm\mu+\frac{\Gamma\bm T}{\sqrt n}\right)\middle|\bm Y\right]
		-\bm w^\top\left(\bm\mu+\frac{\Gamma\bm m_{\pi,I_d}(\bm Y)}{\sqrt n}\right).
	\end{equation*}
	The first term does not depend on $\bm w$. Hence a Bayes rule is
	\begin{equation*}
		\bm w_*\left(\bm\mu+\frac{\Gamma\bm m_{\pi,I_d}(\bm Y)}{\sqrt n}\right)
		\in\arg\max_{\bm w\in\Omega}
		\left(\bm\mu+\frac{\Gamma\bm m_{\pi,I_d}(\bm Y)}{\sqrt n}\right)^\top\bm w
		\subseteq\Omega_{n,L},
	\end{equation*}
	since $\|\Gamma\bm m_{\pi,I_d}(\bm Y)\|_2\le\|\Gamma\|_{\mathrm{op}}M\le L$. Thus every prior admits a Bayes rule in $\Omega_{n,L}$, and its Bayes risk is the same for the restricted (i.e.,$\Omega_{n,L}$) and unrestricted (i.e., $\Omega$) decision spaces. Applying the minimax identity of the form of \eqref{eq:minimax_regret_equals_sup_bayes} to both spaces proves the assertion about restricting the kernels.
	
	The unrestricted Bayes risk equals
	\begin{equation*}
		\begin{aligned}
			\inf_\delta\mathbb E_\pi\int_\Omega\ell_n(\Gamma\bm T,\bm w)\,\delta(\dd\bm w\mid\bm Y)&=\mathbb E_\pi\left[f_*\left(\bm\mu+\frac{\Gamma\bm T}{\sqrt n}\right)
			-f_*\left(\bm\mu+\frac{\Gamma\bm m_{\pi,I_d}(\bm Y)}{\sqrt n}\right)\right]\\
			&=\mathbb E_\pi\left[r_{\bm\mu}\left(\frac{\Gamma\bm T}{\sqrt n}\right)
			-r_{\bm\mu}\left(\frac{\Gamma\bm m_{\pi,I_d}(\bm Y)}{\sqrt n}\right)\right].
		\end{aligned}
	\end{equation*}
	Consequently, the minimax identity of the form of \eqref{eq:minimax_regret_equals_sup_bayes} gives the exact equality
	\begin{equation*}
		\begin{aligned}
			&\inf_\delta\sup_{\bm u\in B_2(M)}
			\mathbb E_{\mathcal N(\bm u,I_d)}
			\int_\Omega \ell_n(\Gamma\bm u,\bm w)\,
			\delta(\dd\bm w\mid\bm Y)\\
			&\quad=\sup_{\pi\in\mathcal P(B_2(M))}
			\mathbb E_\pi\left[
			r_{\bm\mu}\left(\frac{\Gamma\bm T}{\sqrt n}\right)
			-r_{\bm\mu}\left(
			\frac{\Gamma\bm m_{\pi,I_d}(\bm Y)}{\sqrt n}\right)
			\right].
		\end{aligned}
	\end{equation*}
	Both $\Gamma\bm T$ and $\Gamma\bm m_{\pi,I_d}(\bm Y)$ belong to $B_2(L)$. The uniform expansion \eqref{eq:fast_general_local_expansion} therefore implies
	\begin{equation*}
		\begin{aligned}
			&\quad\sup_{\pi\in\mathcal P(B_2(M))}\left|
			\frac{\mathbb E_\pi[r_{\bm\mu}(\Gamma\bm T/\sqrt n)
				-r_{\bm\mu}(\Gamma\bm m_{\pi,I_d}(\bm Y)/\sqrt n)]}{a(n^{-1/2})}
			-\mathbb E_\pi[\psi(\Gamma\bm T)-\psi(\Gamma\bm m_{\pi,I_d}(\bm Y))]\right|\\
			&\le2\sup_{\bm t\in B_2(L)}
			\left|\frac{r_{\bm\mu}(\bm t/\sqrt n)}{a(n^{-1/2})}-\psi(\bm t)\right|\to0.
		\end{aligned}
	\end{equation*}
	Since
	\begin{equation*}
		\sup_{\pi\in\mathcal P(B_2(M))}
		\mathbb E_\pi[\psi(\Gamma\bm T)-\psi(\Gamma\bm m_{\pi,I_d}(\bm Y))]
		=V_{M,I_d}(\psi\circ\Gamma),
	\end{equation*}
	taking suprema over the same priors proves \eqref{eq:div_gaussian_ellipsoid}.
	
	To prove continuity, let $(M_k,\Gamma_k)\to(M,\Gamma)$ with $M_k,M>0$ and $\Gamma_k,\Gamma$ invertible. Parametrize the prior by $\bm T=M\bm U$, where $\bm U\sim\pi\in\mathcal P(B_2(1))$. The correspondence between priors on $B_2(M)$ and $B_2(1)$ is bijective. Given $\bm U=\bm u$, $\bm{Y}\sim\mathcal N(M\bm u,I_d)$, and
	\begin{equation*}
		l_M(\bm u,\bm y)=e^{M\bm u^\top\bm y-M^2\|\bm u\|_2^2/2},
		\qquad
		\mathbb E_\pi[\bm U\mid\bm Y=\bm y]
		=\frac{\int\bm u l_M(\bm u,\bm y)\,\pi(\dd\bm u)}{\int l_M(\bm u,\bm y)\,\pi(\dd\bm u)}.
	\end{equation*}
	All integrals below are over $B_2(1)$. For $l=l_{M_k}$ or $l=l_M$, define
	\begin{equation*}
		\begin{aligned}
			\mathcal J_\pi(\Gamma,M;l)
			:=\mathbb E_{\mathcal N(\bm0,I_d)}\Bigg[
			&\int\psi(\Gamma M\bm u)l(\bm u,\bm Y)\,\pi(\dd\bm u)\\
			&-\left(\int l(\bm u,\bm Y)\,\pi(\dd\bm u)\right)
			\psi\left(\Gamma M\frac{\int\bm u l(\bm u,\bm Y)\,\pi(\dd\bm u)}
			{\int l(\bm u,\bm Y)\,\pi(\dd\bm u)}\right)\Bigg].
		\end{aligned}
	\end{equation*}
	The arguments $(\Gamma,M)$ determine the function $\bm u\mapsto\psi(\Gamma M\bm u)$, whereas $l$ determines the mixture density and posterior mean. The mixture-density calculation used for $J_\pi$ gives
	\begin{equation*}
		V_{M,I_d}(\psi\circ\Gamma)
		=\sup_{\pi\in\mathcal P(B_2(1))}\mathcal J_\pi(\Gamma,M;l_M),
		\qquad
		V_{M_k,I_d}(\psi\circ\Gamma_k)
		=\sup_{\pi\in\mathcal P(B_2(1))}\mathcal J_\pi(\Gamma_k,M_k;l_{M_k}).
	\end{equation*}
	
	Choose a fixed $L$ with $M_k\|\Gamma_k\|_{\mathrm{op}}\le L$ and $M\|\Gamma\|_{\mathrm{op}}\le L$ for all sufficiently large $k$. The finite convex function $\psi$ is bounded and Lipschitz on $B_2(L)$. Since $M_k\Gamma_k\to M\Gamma$,
	\begin{equation*}
		\sup_{\bm u\in B_2(1)}|\psi(\Gamma_kM_k\bm u)-\psi(\Gamma M\bm u)|
		\le C\|M_k\Gamma_k-M\Gamma\|_{\mathrm{op}}\to0.
	\end{equation*}
	Also, $\bm u\mapsto\psi(\Gamma_kM_k\bm u)$ and $\bm u\mapsto\psi(\Gamma M\bm u)$ have a common bound and Lipschitz constant on $B_2(1)$, after increasing $C$ if necessary. For each fixed $\bm y$,
	\begin{equation*}
		\sup_{\bm u\in B_2(1)}|\log l_{M_k}(\bm u,\bm y)-\log l_M(\bm u,\bm y)|
		\le |M_k-M|\|\bm y\|_2+\frac12|M_k^2-M^2|\to0.
	\end{equation*}
	Thus $\sup_{\bm u\in B_2(1)}|l_{M_k}(\bm u,\bm y)-l_M(\bm u,\bm y)|\to0$. Since $M_k$ is bounded, for a constant $C$ and all sufficiently large $k$,
	\begin{equation*}
		\sup_{\bm u\in B_2(1)}|l_{M_k}(\bm u,\bm y)-l_M(\bm u,\bm y)|
		\le2e^{C\|\bm y\|_2},
		\qquad
		\mathbb E_{\mathcal N(\bm0,I_d)}e^{C\|\bm Y\|_2}<\infty.
	\end{equation*}
	DCT gives
	\begin{equation*}
		\mathbb E_{\mathcal N(\bm0,I_d)}
		\sup_{\bm u\in B_2(1)}|l_{M_k}(\bm u,\bm Y)-l_M(\bm u,\bm Y)|\to0.
	\end{equation*}
	
	Fix $\pi\in\mathcal P(B_2(1))$. Adding and subtracting $\mathcal J_\pi(\Gamma,M;l_{M_k})$ gives
	\begin{equation*}
		\begin{aligned}
			&\quad\mathcal J_\pi(\Gamma_k,M_k;l_{M_k})-\mathcal J_\pi(\Gamma,M;l_M)\\
			&=\bigl[\mathcal J_\pi(\Gamma_k,M_k;l_{M_k})-\mathcal J_\pi(\Gamma,M;l_{M_k})\bigr]
			+\bigl[\mathcal J_\pi(\Gamma,M;l_{M_k})-\mathcal J_\pi(\Gamma,M;l_M)\bigr].
		\end{aligned}
	\end{equation*}
	For the first difference, only $\Gamma_kM_k$ is replaced by $\Gamma M$. The posterior mean is unchanged, and
	\begin{equation*}
		\left\|\frac{\int\bm u l_{M_k}(\bm u,\bm Y)\,\pi(\dd\bm u)}
		{\int l_{M_k}(\bm u,\bm Y)\,\pi(\dd\bm u)}\right\|_2
		\le\frac{\int\|\bm u\|_2 l_{M_k}(\bm u,\bm Y)\,\pi(\dd\bm u)}
		{\int l_{M_k}(\bm u,\bm Y)\,\pi(\dd\bm u)}\le1.
	\end{equation*}
	Thus both $\bm u$ and this posterior mean belong to $B_2(1)$. Each of the two terms in the functional changes in absolute value by at most
	\begin{equation*}
		\sup_{\bm u\in B_2(1)}|\psi(\Gamma_kM_k\bm u)-\psi(\Gamma M\bm u)|
		\int l_{M_k}(\bm u,\bm Y)\,\pi(\dd\bm u).
	\end{equation*}
	Moreover, Tonelli's theorem and the likelihood-ratio identity give
	\begin{equation*}
		\mathbb E_{\mathcal N(\bm0,I_d)}\int l_{M_k}(\bm u,\bm Y)\,\pi(\dd\bm u)
		=\int\mathbb E_{\mathcal N(\bm0,I_d)}l_{M_k}(\bm u,\bm Y)\,\pi(\dd\bm u)
		=\int1\,\pi(\dd\bm u)=1.
	\end{equation*}
	Consequently,
	\begin{equation*}
		|\mathcal J_\pi(\Gamma_k,M_k;l_{M_k})-\mathcal J_\pi(\Gamma,M;l_{M_k})|
		\le2\sup_{\bm u\in B_2(1)}|\psi(\Gamma_kM_k\bm u)-\psi(\Gamma M\bm u)|.
	\end{equation*}
	
	For the second difference, $(\Gamma,M)$ is fixed and only $l_{M_k}$ is replaced by $l_M$. For each fixed $\bm y$, the functions $l_{M_k}(\cdot,\bm y)$ and $l_M(\cdot,\bm y)$ are positive and continuous on $B_2(1)$. The function $\bm u\mapsto\psi(\Gamma M\bm u)$ is bounded and Lipschitz there, with both bounds at most $C$. Thus \eqref{eq:ng_bayes_functional_stability} applies with the radius in that estimate equal to $1$, its matrix equal to $\Gamma M$, and its two weight functions equal to $l_{M_k}(\cdot,\bm y)$ and $l_M(\cdot,\bm y)$. The two integrands defining $\mathcal J_\pi(\Gamma,M;l_{M_k})$ and $\mathcal J_\pi(\Gamma,M;l_M)$ are precisely the corresponding $J_\pi$ functionals. Taking expectation of the pointwise bound and using $|\mathbb E A-\mathbb E B|\le\mathbb E|A-B|$ gives
	\begin{equation*}
		|\mathcal J_\pi(\Gamma,M;l_{M_k})-\mathcal J_\pi(\Gamma,M;l_M)|
		\le4C\,\mathbb E_{\mathcal N(\bm0,I_d)}
		\sup_{\bm u\in B_2(1)}|l_{M_k}(\bm u,\bm Y)-l_M(\bm u,\bm Y)|.
	\end{equation*}
	Both bounds are independent of $\pi$. Taking suprema over the same priors gives
	\begin{equation*}
		\begin{aligned}
			|V_{M_k,I_d}(\psi\circ\Gamma_k)-V_{M,I_d}(\psi\circ\Gamma)|
			&\le\sup_{\pi\in\mathcal P(B_2(1))}|\mathcal J_\pi(\Gamma_k,M_k;l_{M_k})-\mathcal J_\pi(\Gamma,M;l_M)|\\
			&\le2\sup_{\bm u\in B_2(1)}|\psi(\Gamma_kM_k\bm u)-\psi(\Gamma M\bm u)|\\
			&\quad+4C\,\mathbb E_{\mathcal N(\bm0,I_d)}
			\sup_{\bm u\in B_2(1)}|l_{M_k}(\bm u,\bm Y)-l_M(\bm u,\bm Y)|\to0.
		\end{aligned}
	\end{equation*}
	This proves the asserted continuity.
\end{proof}

\end{document}